\documentclass[a4paper,11pt]{amsart}

\usepackage[T1]{fontenc}
\usepackage{amsmath,amssymb,amsthm}
\usepackage{mathrsfs}
\usepackage{enumerate}
\usepackage[all]{xy}
\usepackage[margin=25mm]{geometry}
\usepackage{hyperref}
\hypersetup{hidelinks}
\theoremstyle{plain}
\newtheorem{theorem}{Theorem}[section]
\newtheorem{maintheorem}[theorem]{Main Theorem}
\newtheorem{proposition}[theorem]{Proposition}
\newtheorem{lemma}[theorem]{Lemma}
\newtheorem{corollary}[theorem]{Corollary}

\theoremstyle{definition}
\newtheorem{definition}[theorem]{Definition}
\newtheorem{example}[theorem]{Example}

\theoremstyle{remark}
\newtheorem{remark}[theorem]{Remark}

\numberwithin{equation}{section}

\title{Derived algebraic cobordism and Chow motives over prisms}
\date{\today}
\author{Yuki Kato}
\thanks{The author was supported by Grants-in-Aid for Scientific Research No.~23K03080, Japan Society for the Promotion of Science.}
\address{National Institute of Technology, Kurume College,
1-1-1 Komorino, Kurume, Fukuoka 830-8555, Japan}
\email{kato\_051@kurume-nct.ac.jp}

\newcommand{\Vect}{\mathrm{Vect}}
\newcommand{\Spec}{\operatorname{Spec}}
\newcommand{\Span}{\mathrm{Span}}
\newcommand{\Corr}{\mathrm{Corr}}
\newcommand{\fib}{\operatorname{fib}}
\newcommand{\dMGL}{\mathrm{dMGL}}
\newcommand{\Sp}{\mathrm{Sp}}

\newcommand{\A}{\mathbb{A}}
\newcommand{\Gm}{\mathbb{G}_m}
\newcommand{\Ga}{\mathbb{G}_a}
\newcommand{\tm}{\widetilde{\mathfrak{m}}}
\newcommand{\NAlg}{\mathrm{NAlg}}
\newcommand{\CAlg}{\mathrm{CAlg}}
\newcommand{\Shv}{\operatorname{Shv}}
\newcommand{\Map}{\operatorname{Map}}
\newcommand{\Exc}{\operatorname{Exc}}
\newcommand{\Ext}{\operatorname{Ext}}
\newcommand{\Mod}{\operatorname{Mod}}
\newcommand{\Ainf}{{\mathbb{A}_{\inf}}}
\newcommand{\Prism}{\mathrm{Prism}}
\newcommand{\conn}{\operatorname{conn}}

\newcommand{\cofib}{\operatorname{cofib}}

\begin{document}
\begin{abstract}
We construct non-$\mathbb{A}^1$-invariant derived algebraic cobordism for $\infty$-categories equipped with a universal line datum and apply this construction to derived formal schemes over prisms. Over a perfect prism $(\mathbb{A}_{\inf R}, I_R)$, cotangent and relative tangential excisive approximations induce symmetric monoidal almost tilting equivalences between the arithmetic, prismatic, and tilted cobordism categories, as well as derived base-change specializations along $\theta$ and $\rho_\flat$. On iterated loop spaces, the resulting \v{C}ech--Amitsur totalization towers satisfy derived $p$-complete hyperdescent and stabilize at an explicit finite stage in each homotopy degree. We also define universal $\mathbb{G}_a$- and $\mathbb{G}_m$-oriented derived span motives; in characteristic zero, $\mathbb{G}_a$-oriented correspondences recover classical rational Grothendieck Chow motives, and over a perfectoid base both oriented theories satisfy the same tilting and base-change equivalences.
\end{abstract}
\maketitle

\section{Introduction}
\label{sec:introduction}
In motivic homotopy theory, Voevodsky's algebraic cobordism spectrum $\mathrm{MGL}$ is the motivic analogue of the complex cobordism spectrum $\mathrm{MU}$, representing the universal oriented cohomology theory in the stable motivic homotopy category $\mathrm{SH}(S)$ over a base scheme $S$~\cite{V,Panin2003Oriented, GeSn}. Garkusha and Panin~\cite{GP} refined the approach via spans into framed correspondences. By using motivic stable homotopy theory and framed correspondences, Elmanto, Hoyois, Khan, Sosnilo, and Yakerson~\cite{ModMGL} developed algebraic cobordism as the monoidal unit of the stable $\infty$-category with finite syntomic correspondences.

Lowrey--Sch\"{u}rg~\cite{LS16} and Annala~\cite{Ann21} introduced bivariant algebraic cobordism theories for derived schemes without $\A^1$-invariance. In this paper, we first define non-$\A^1$-invariant derived algebraic cobordism $\dMGL_{\mathcal{C}}$ on any $\infty$-category $\mathcal{C}$ equipped with a universal line datum and an admissible class of morphisms $\mathbf{S}$, as an $\infty$-categorical refinement of the algebraic cobordism via spans~\cite{KatoMGL} introduced by the author (Section~\ref{sec:span_formalism}, Definition~\ref{def:relation-dmgl}, and Propositions~\ref{thm:localized-span-universal-property} and~\ref{thm:full-orientation-from-line-thom}). Motivated by the characterization of periodic algebraic $K$-theory via multiplicative orientations without $\A^1$-invariance~\cite[Proposition~8.10 and Theorem~8.11]{AHIConnerFloyd}, we also formulate universal additive and multiplicative laws in the spectral span category; these laws, corresponding to the formal group laws in stable homotopy theory, define universal $\Ga$-oriented derived Chow motives and universal $\Gm$-oriented derived span motives. Over a field of characteristic zero with rational coefficients, our $\Ga$-oriented span relations annihilate derived double-point cobordism relations and agree with the additive formal group law specialization $\mathbb L_* \otimes \mathbb Q \twoheadrightarrow \mathbb Q$ of Levine--Morel's Lazard ring (Proposition~\ref{prop:lazard-ga-compatibility}), recovering classical rational Grothendieck Chow motives (Theorem~\ref{thm:p1-chow-comparison}).

In this paper, we fix a prime number $p>0$. Let $R^\circ$ be a $p$-torsion-free integral perfectoid ring of mixed characteristic $(0,p)$ with tilt $R^{\flat\circ}$, and set $\mathbb{A}_{\inf R}=\Ainf(R^\circ)=W(R^{\flat\circ})$ and $I_R=\ker(\theta\colon\mathbb{A}_{\inf R}\to R^\circ)$ for Fontaine's map $\theta$. By Bhatt and Scholze's equivalence between perfect prisms and integral perfectoid rings~\cite[Theorem~3.10]{BhattScholze2017Prisms}, the perfect prism $(\mathbb{A}_{\inf R},I_R)$ specializes along $\theta$ to $R^\circ$ and along the reduction map $\rho_\flat\colon\mathbb{A}_{\inf R}\to\mathbb{A}_{\inf R}/(p)$ to $R^{\flat\circ}$. The formal scheme $\mathfrak F_R=\operatorname{Spf}(\mathbb{A}_{\inf R})$, with its Cartier divisors $V(I_R)$ and $V(p)$, gives their common base. Choosing a distinguished generator $\xi$ of $I_R$ gives $\delta(\xi)\in\mathbb{A}_{\inf R}^\times$ for the $\delta$-ring structure on $\mathbb{A}_{\inf R}$. Setting $\bar\xi=\rho_\flat(\xi)\in R^{\flat\circ}$, $\omega^\flat=\bar\xi^{1/p}\in R^{\flat\circ}$, $\widetilde{\omega}=[\omega^\flat]\in\mathbb{A}_{\inf R}$, and $\omega=\theta(\widetilde{\omega})\in R^\circ$ defines the adic towers $R_n=R^\circ/(\omega^n)$ and $R_n^\flat=R^{\flat\circ}/((\omega^\flat)^n)$ for $n\ge 1$.

For any $1\le k\le p$, the elements $(\xi,\widetilde{\omega}^k)$ and $(p,\widetilde{\omega}^k)$ form regular sequences of length $2$ in $\mathbb{A}_{\inf R}$ generating the same ideal, fitting into the commutative specialization diagram
\[
\xymatrix@C=3.0em@R=1.7em{
& \mathbb{A}_{\inf R} \ar@{->>}[dl]_-{\theta} \ar@{->>}[dr]^-{\rho_\flat} & \\
R^\circ \ar@{->>}[dr] & & R^{\flat\circ} \ar@{->>}[dl] \\
& R_k \cong R_k^\flat. &
}
\]

Applying this $\infty$-categorical construction of $\dMGL$ and its oriented refinements to derived formal schemes over the perfect prism $(\mathbb{A}_{\inf R}, I_R)$ and to their arithmetic and tilted almost formal specializations over $R^\circ$ and $R^{\flat\circ}$ (Definitions~\ref{def:nonunital-comparison-sources} and~\ref{def:derived-gm-span-motives}, and Subsection~\ref{subsec:tilting-oriented-motives}), we obtain our main almost tilting equivalences, proved in Theorems~\ref{thm:tilting-main}(1) and~\ref{thm:tilting-oriented-motives}(1). Throughout, write $(-)^a$ for almostification and $\operatorname{Spf}^a$ for almost formal spectra; all tensor products $\otimes$ and completions $(-)^\wedge$ are derived (Section~\ref{sec:fiber-specialization-and-deformation}).

\begin{maintheorem}
\label{thm:intro-tilting}
Let $A^\circ$ be an integral perfectoid $R^\circ$-algebra with tilt $A^{\flat\circ}$ over $R^{\flat\circ}$ and perfect prism $\widetilde{A}=\Ainf(A^\circ)=W(A^{\flat\circ})$. For any split augmented derived $(p,I_R)$-complete formal scheme $\operatorname{Spf}(\widetilde{A}\oplus I)\to X_0=\operatorname{Spf}(\widetilde{A})$ with connective derived $(p,I_R)$-complete augmentation ideal $I$, boundary specializations $I_\theta=(I\otimes_{\mathbb{A}_{\inf R},\theta}R^\circ)^a$ and $I_\flat=(I\otimes_{\mathbb{A}_{\inf R},\rho_\flat}R^{\flat\circ})^a$, and any quasi-smooth formal schemes $X,Y\in\mathscr D_{\mathrm{cmp}}(X_0)$ over $X_0$ with boundary specializations $i_\theta^*X,i_\theta^*Y$ and $i_\flat^*X,i_\flat^*Y$ (Subsection~\ref{subsec:fiber-specializations}), for the transported cobordism categories $\mathcal M_{\mathcal D,\mathrm{cmp}}^{a,\mathrm{nu}}$ along the adic tower and their intrinsic mod-$\omega_{\mathcal D}$ associated graded versions $\mathcal M_{\mathcal D,\mathrm{gr}}^{a,\mathrm{nu}}$ (Definition~\ref{def:nonunital-comparison-sources}), there are symmetric monoidal almost tilting equivalences of $\infty$-categories and of derived $\bar\omega_{\mathcal D}$-complete module $\infty$-categories $\Mod_E^\wedge$ over $E=\mathbf{dMGL}_{\mathcal D,\mathrm{cmp}}^{a,\mathrm{nu}}=\operatorname{End}(\mathbf 1)$ ($\mathcal D\in\{\mathrm{Prism},\theta,\flat\}$):
\[
\begin{aligned}
\mathcal M_{\theta,\mathrm{cmp}}^{a,\mathrm{nu}}(A^\circ)
&\;\longleftarrow\;
\mathcal M_{\Prism,\mathrm{cmp}}^{a,\mathrm{nu}}(\widetilde{A})
\;\longrightarrow\;
\mathcal M_{\flat,\mathrm{cmp}}^{a,\mathrm{nu}}(A^{\flat\circ}), \\
\mathcal M_{\theta,\mathrm{gr}}^{a,\mathrm{nu}}(A^\circ)
&\;\longleftarrow\;
\mathcal M_{\Prism,\mathrm{gr}}^{a,\mathrm{nu}}(\widetilde{A})
\;\longrightarrow\;
\mathcal M_{\flat,\mathrm{gr}}^{a,\mathrm{nu}}(A^{\flat\circ}), \\
\Mod_{\mathbf{dMGL}_{\theta,\mathrm{cmp}}^{a,\mathrm{nu}}(A^\circ)}^\wedge
&\;\longleftarrow\;
\Mod_{\mathbf{dMGL}_{\Prism,\mathrm{cmp}}^{a,\mathrm{nu}}(\widetilde{A})}^\wedge
\;\longrightarrow\;
\Mod_{\mathbf{dMGL}_{\flat,\mathrm{cmp}}^{a,\mathrm{nu}}(A^{\flat\circ})}^\wedge,
\end{aligned}
\]
compatible with the symmetric monoidal associated graded realization functors $\mathcal M_{\mathcal D,\mathrm{cmp}}^{a,\mathrm{nu}}\to\mathcal M_{\mathcal D,\mathrm{gr}}^{a,\mathrm{nu}}$ ($\mathcal D\in\{\mathrm{Prism},\theta,\flat\}$), together with almost tilting equivalences of cobordism spectra along the adic tower across all $n\ge 1$:
\[
\dMGL_{\theta,\mathrm{cmp}}^{a,\mathrm{nu}}(i_\theta^*X\to i_\theta^*Y;I_\theta)
\;\longleftarrow\;
\dMGL_{\Prism,\mathrm{cmp}}^{a,\mathrm{nu}}(X\to Y;I)
\;\longrightarrow\;
\dMGL_{\flat,\mathrm{cmp}}^{a,\mathrm{nu}}(i_\flat^*X\to i_\flat^*Y;I_\flat),
\]
compatible with refinements of hypercovers.
\end{maintheorem}

\begin{maintheorem}
\label{thm:intro-tilting-oriented}
Under the hypotheses and notation of Main Theorem~\ref{thm:intro-tilting}, for the transported categories of universal $\mathbb G_m$-oriented derived span motives and $\mathbb G_a$-oriented derived Chow motives along the adic tower, together with their intrinsic associated graded versions (Subsection~\ref{subsec:tilting-oriented-motives}), there are symmetric monoidal almost tilting equivalences:
\[
\begin{aligned}
\mathcal{DM}_{\mathbb G_m,\theta}^{a,\mathrm{nu}}(A^\circ)
&\;\longleftarrow\;
\mathcal{DM}_{\mathbb G_m,\Prism}^{a,\mathrm{nu}}(\widetilde{A})
\;\longrightarrow\;
\mathcal{DM}_{\mathbb G_m,\flat}^{a,\mathrm{nu}}(A^{\flat\circ}), \\
\mathcal{DM}_{\mathrm{Chow},\theta}^{a,\mathrm{nu}}(A^\circ)
&\;\longleftarrow\;
\mathcal{DM}_{\mathrm{Chow},\Prism}^{a,\mathrm{nu}}(\widetilde{A})
\;\longrightarrow\;
\mathcal{DM}_{\mathrm{Chow},\flat}^{a,\mathrm{nu}}(A^{\flat\circ}),
\end{aligned}
\]
compatible with the associated graded equivalences $\mathcal{DM}_{\mathbb G_m,\theta,\mathrm{gr}}^{a,\mathrm{nu}}(A^\circ)\leftarrow\mathcal{DM}_{\mathbb G_m,\Prism,\mathrm{gr}}^{a,\mathrm{nu}}(\widetilde A)\rightarrow\mathcal{DM}_{\mathbb G_m,\flat,\mathrm{gr}}^{a,\mathrm{nu}}(A^{\flat\circ})$ (resp.\ $\mathcal{DM}_{\mathrm{Chow},\theta,\mathrm{gr}}^{a,\mathrm{nu}}(A^\circ)\leftarrow\mathcal{DM}_{\mathrm{Chow},\Prism,\mathrm{gr}}^{a,\mathrm{nu}}(\widetilde A)\rightarrow\mathcal{DM}_{\mathrm{Chow},\flat,\mathrm{gr}}^{a,\mathrm{nu}}(A^{\flat\circ})$), identifying the excisive approximations of the oriented tangent motives of $(I,I_\theta,I_\flat)$ along the adic tower across all $n\ge 1$.
\end{maintheorem}
 
We outline the main ingredients of the proof (Sections~\ref{sec:span_formalism}--\ref{sec:oriented-span-motives}). For an $\infty$-category $\mathcal{C}$ equipped with a universal line datum and a class of admissible morphisms closed under homotopy pullbacks, derived l.c.i.\ spans compose by homotopy pullbacks, and derived self-intersections of zero sections define a preorientation (Proposition~\ref{prop:zero-section-class}). Localizing the stable motivic envelope at admissible pushout excision, projective bundle decompositions, presentation hyperdescent, and normalized projective Thom comparison maps defines $\dMGL_{\mathcal{C}}$ (Propositions~\ref{thm:localized-span-universal-property} and~\ref{thm:full-orientation-from-line-thom}). Applying this construction to $\mathcal{T}_R$ and to the arithmetic and tilted almost formal categories $\mathcal{C}_{R^\circ}^a$ and $\mathcal{C}_{R^{\flat\circ}}^a$ via Gabber--Ramero almost unitalization and flat firm modules~\cite{GR,kato2023mathematics,KatoCalcNUAlg}, with quasiregular semiperfectoid (QRSP) rings~\cite{BMS19,BhattScholze2017Prisms} providing local splitting models, relates the prismatic and boundary theories under $i_\theta^*$ and $i_\flat^*$.

Along the adic towers $\{R_n\}_{n\ge 1}$ and $\{R_n^\flat\}_{n\ge 1}$, each square-zero transition kernel $I_n=(\omega^n)/(\omega^{n+1})\cong R_1$ reduces the deformation tangent data across all $n\ge 1$ to the shared residue ring $R_1\simeq R_1^\flat$ (Lemma~\ref{lem:reduction} and Corollary~\ref{cor:reduction-controls}). On unobstructed affine formal schemes, the adic lift spaces at finite levels therefore form nonempty torsors under derivation groups identified via $R_1\simeq R_1^\flat$ (Corollary~\ref{cor:torsor-comparison}), represented by simplicial affine Bar hypercovers over derived affine vector group schemes (Proposition~\ref{prop:affine-bar-models}). Adjoining local section sieves of these lift torsors across all $n\ge 1$ to the quasisyntomic--Zariski topology $\tau_{\mathrm{geom}}$ (Definition~\ref{def:master-prismatic-theory}) defines the deformation topology $\tau_{\mathrm{pld}}$ and, via universal completed perfections, the comparison topology $\tau_{\mathrm{cmp}}$ with hypercomplete spectral sheaf category $\mathcal{E}_R=\Shv_{\mathrm{hyp}}(\mathcal{T}_R,\tau_{\mathrm{cmp}};\Sp)$ (Proposition~\ref{prop:common-coefficient-category}). Over perfect prisms, the boundary specializations identify the coefficient sites and preserve all motivic relation and orientation generators (Propositions~\ref{prop:perfect-boundary-coefficient-sites} and~\ref{prop:boundary-specialization-functoriality}).

Combining Goodwillie's excisive approximation functor~\cite{Goodwillie3,Heuts2018,KatoCalcNUAlg} with the twisted-arrow end description of pointed diagram lifts (Proposition~\ref{prop:coherent-dual-lift-complex}) gives the proof of Main Theorems~\ref{thm:intro-tilting} and~\ref{thm:intro-tilting-oriented} together with their integral $\mathcal{O}$-linear versions. Applying cotangent excisive approximation $P^1_{\mathrm{cot}}$ to the split augmentation $B\oplus I\to B$ and relative tangential excisive approximation $P^1_{\mathrm{rel}}$ across the square-zero adic transitions $0\to I_n\to R_{n+1}\to R_n\to 0$ ($n\ge 1$) endows any split lift torsor with a basepoint given by the zero derivation coherent across hypercovers and passes the coefficients to the completed associated graded almost rings $\widehat G_\theta\xleftarrow{\sim}\widehat G_{\mathrm{Prism}}\xrightarrow{\sim}\widehat G_\flat\simeq(R_1[\![\bar\omega_{\mathcal D}]\!])^a$ (Subsection~\ref{subsec:reduction-and-p1}). Since Lemma~\ref{lem:reduction} and Corollary~\ref{cor:reduction-controls} identify all acting derivation groups via $R_1\simeq R_1^\flat$, the levelwise lift spectra combine across all $n\ge 1$ into symmetric monoidal almost tilting equivalences for $\mathcal M_{\mathcal D,\mathrm{cmp}}^{a,\mathrm{nu}}$, $\mathcal{DM}_{\mathbb G_m,\mathcal D}^{a,\mathrm{nu}}$, and $\mathcal{DM}_{\mathrm{Chow},\mathcal D}^{a,\mathrm{nu}}$ together with their associated graded versions $\mathcal M_{\mathcal D,\mathrm{gr}}^{a,\mathrm{nu}}$, $\mathcal{DM}_{\mathbb G_m,\mathcal D,\mathrm{gr}}^{a,\mathrm{nu}}$, and $\mathcal{DM}_{\mathrm{Chow},\mathcal D,\mathrm{gr}}^{a,\mathrm{nu}}$. Retaining the integral bases $\mathbb A_{\inf R}$, $(R^\circ)^a$, and $(R^{\flat\circ})^a$ under $P^1_{\mathrm{cot}}$ defines the integral $\mathcal{O}$-linear categories $\mathcal M_{\mathcal D,\mathrm{cmp}}^{\mathrm{lin}}$, $\mathcal{DM}_{\mathbb G_m,\mathcal D}^{\mathrm{lin}}$, and $\mathcal{DM}_{\mathrm{Chow},\mathcal D}^{\mathrm{lin}}$ (Definition~\ref{def:integral-linear-comparison-sources}); dualizability of $R^\circ$, $R^{\flat\circ}$, and $R_k$ ($1\le k\le p$) over $\mathbb{A}_{\inf R}$ (Lemma~\ref{lem:tilting-frame-koszul}), together with Lemma~\ref{lem:p1-fiber-comparison}, gives derived base-change equivalences along $\theta$ and $\rho_\flat$ (Theorems~\ref{thm:tilting-main}(2) and~\ref{thm:tilting-oriented-motives}(2)), recovers the split spectra along the adic tower from the truncated transition layers, and gives mod-$(\omega^k)$ tilting equivalences up to $R_p^a=(R^\circ/(p))^a\simeq (R^{\flat\circ}/(\bar\xi))^a=(R_p^\flat)^a$ (Propositions~\ref{prop:koszul-tilting-and-tower-recovery} and~\ref{prop:oriented-koszul-tilting}). The boundary categories $\mathcal M_{\nu,\mathrm{cmp}}^{a,\mathrm{nu}}$, $\mathcal M_{\nu,\mathrm{cmp}}^{\mathrm{lin}}$, $\mathcal{DM}_{\mathbb G_m,\nu}^{a,\mathrm{nu}}$, $\mathcal{DM}_{\mathrm{Chow},\nu}^{a,\mathrm{nu}}$, $\mathcal{DM}_{\mathbb G_m,\nu}^{\mathrm{lin}}$, and $\mathcal{DM}_{\mathrm{Chow},\nu}^{\mathrm{lin}}$ are built on the presentation span categories transported along the coefficient site and admit symmetric monoidal realization functors to the intrinsic boundary categories, whereas the mod-$\omega_{\mathcal D}$ associated graded categories and the evaluated cobordism and tangent spectra are intrinsic on both boundaries.

We also obtain three consequences. First, Main Theorem~\ref{thm:intro-tilting} implies almost conservativity on derived completions (Corollary~\ref{cor:mixed-to-equal-conservativity}): a prismatic morphism induces an almost equivalence on derived $\omega$-completed cobordism over $(R^\circ)^a$ if and only if its tilt induces an almost equivalence on derived $\omega^\flat$-completed cobordism over $(R^{\flat\circ})^a$. Second, along $\tau_{\mathrm{cmp}}$-hypercovers pulled back to geometric iterated loop spaces $\Omega_X^n U$ ($n\ge 0$), cosimplicial totalizations compute $\tau_{\mathrm{cmp}}$-sheafified invariants and satisfy derived $p$-complete hyperdescent in all homotopy degrees (Theorem~\ref{thm:derived-p-complete-descent}). For a pointed derived affine scheme $U=\Spec(A\oplus J)$ over $\Spec A$ with connective non-unital augmentation ideal $J\in\NAlg_A$ and cotangent fiber $Q_A(J)$, passing to the $k$-fold loop space $\Omega_A^k U=\Spec(A\oplus J^{(k)})$ ($k\ge 1$) suspends $J$ in $\NAlg_A$, shifting the cotangent fiber to $Q_A(J^{(k)})\simeq Q_A(J)[k]$ and raising the connectivity of the augmentation ideal to $\conn(J^{(k)})=\conn(Q_A(J))+k$ (unlike nonzero unital rings, which remain $0$-connective). On the \v{C}ech--Amitsur nerve of $\Spec A\to\Omega_A^{k-1}U$, quasi-smoothness of a diagram $\sigma$ and excisive approximation on modules linearize the lift functor to the exact tensor product $K_\sigma^\vee\otimes_A(-)$ with bounded-below dual lift complex $K_\sigma^\vee$, so cosimplicial Dold--Kan extracts the $q$-fold tensor powers $(J^{(k)})^{\otimes_A q}$ and yields positive totalization slope $a+k-1$ for $a=\conn(Q_A(J))$; consequently, whenever $a+k>1$---in particular at $k=2$ for any connective ideal ($a\ge 0$) and at $k=1$ for $1$-connective ideals ($a\ge 1$)---the linearized totalization tower stabilizes at an explicit finite stage in every homotopy degree $i\in\mathbb Z$ (Theorem~\ref{thm:finite-stage-detection}). Third, over the neutral locus where lift torsors split quasisyntomically, $\tau_{\mathrm{geom}}$ and $\tau_{\mathrm{pld}}$ generate the same hypercomplete $\infty$-topoi, so localization induces an equivalence of hypercomplete additive span theories (Proposition~\ref{prop:neutral-locus-equivalence}).

\medskip
\noindent\textbf{Organization of the paper.}
Section~\ref{sec:span_formalism} develops the abstract span formalism and oriented derived algebraic cobordism. Sections~\ref{sec:prismatic-master-theory} and~\ref{sec:fiber-specialization-and-deformation} construct the prismatic cobordism theory over prisms, establish the Koszul relations for $\mathbb A_{\inf R}$, and introduce the boundary almost formal categories. Section~\ref{sec:deformation-perfectoid} develops square-zero lifting torsors, affine Bar models, coherent diagram lifts, cotangent and relative tangential excisive approximations $P^1_{\mathrm{cot}}$ and $P^1_{\mathrm{rel}}$, and the comparison topologies. Section~\ref{sec:linearized-comparison} constructs the excisive approximations of the boundary realizations and proves Main Theorem~\ref{thm:intro-tilting} and Corollary~\ref{cor:mixed-to-equal-conservativity}. Section~\ref{sec:hyperdescent-and-neutral-locus} establishes hyperdescent on iterated loop spaces, stabilization at a finite stage, and the comparison over the neutral locus. Section~\ref{sec:oriented-span-motives} introduces universal $\mathbb{G}_a$- and $\mathbb{G}_m$-oriented derived span motives, establishes compatibility of the $\mathbb G_a$-relations with the additive specialization of Levine--Morel's Lazard ring to recover classical rational Grothendieck Chow motives in characteristic zero (Proposition~\ref{prop:lazard-ga-compatibility} and Theorem~\ref{thm:p1-chow-comparison}), and proves integral prismatic base change and almost tilting equivalences along the adic tower after excisive approximation.

\medskip
\noindent\textbf{Notation and conventions.}
Throughout the paper, we fix Grothendieck universes.
\begin{itemize}
\item We let $\mathcal S$ denote the $\infty$-category of spaces, and $\Sp$ denote the stable $\infty$-category of spectra with sphere spectrum $\mathbb S$, smash product $\otimes_{\mathbb S}$ or $\wedge$, and connective subcategory $\Sp_{\ge 0}$. More generally, $(-)_{\ge 0}$ denotes connective objects, $(-)^{\simeq}$ the maximal sub-$\infty$-groupoid, or core, and $\mathcal P\mathrm{r}_{\mathrm{st}}^{\mathrm{L}}$ the $\infty$-category of presentable stable $\infty$-categories and functors preserving small colimits. We write $\Omega^\infty\colon\Sp\to\mathcal S$ for the infinite loop space functor, $\widetilde{\Sigma}^\infty\colon\mathcal S_*\to\Sp$ for the reduced suspension spectrum functor on pointed spaces, $X_+ = X\amalg *$ for adding a disjoint basepoint, and $\mathbb S/p^m = \cofib(p^m\colon\mathbb S\to\mathbb S)$ for the mod-$p^m$ Moore spectrum.
\item For an $\infty$-category $\mathcal D$ and object $Y\in\mathcal D$, $\mathcal D_{/Y}$ denotes the overcategory; when $\mathcal D$ has a terminal object $*$, $\mathcal D_*$ denotes the $\infty$-category of pointed objects. We write $\Map_{\mathcal D}(X,Y)\in\mathcal S$ for mapping spaces, $\operatorname{map}_{\mathcal E}(X,Y)\in\Sp$ for mapping spectra in a stable $\infty$-category $\mathcal E$, $\underline{\operatorname{map}}$ or $\underline{\mathrm{Hom}}$ for internal mapping objects in closed symmetric monoidal or enriched settings, and $\underline{\Map}$ exclusively for mapping stacks.
\item For a spectrum or module object $M$ and $n\in\mathbb Z$, $M[n]=\Sigma^n M$ denotes the $n$-fold suspension, and $\tau_{\ge 0}$ denotes connective truncation. The connectivity of $M$ is $\conn(M)=\sup\{c\in\mathbb Z:\pi_j(M)=0\text{ for all }j<c\}\in\mathbb Z\cup\{-\infty,+\infty\}$ with $\conn(0)=+\infty$. Tor-amplitude is defined in cohomological degrees, so a perfect complex of Tor-amplitude $[-1,0]$ is locally represented by a two-term homological complex $P_1\to P_0$ of finite projective modules in degrees $1$ and $0$. For a ring or almost ring $A$, $D(A)\simeq\Mod_A$ denotes its unbounded derived $\infty$-category, with connective subcategory $\Mod_A^{\ge 0}$, perfect dualizable subcategory $\operatorname{Perf}(A)$, and relative tensor product $\otimes_A$. All tensor products and completions are derived.
\item We write $\mathrm{CMon}(-)$ for $\mathbb E_\infty$-monoids, $(-)^{\mathrm{gp}}$ for group completion, $\operatorname{Grp}(\mathcal S)$ for group objects in $\mathcal S$, and $\mathrm{Alg}_{\mathbb E_1}(-)$ for $\mathbb E_1$-algebras. We write $\CAlg(\mathcal C)$ for commutative algebra objects in a symmetric monoidal $\infty$-category $\mathcal C$, $\CAlg_A$ for connective animated commutative $A$-algebras, $\delta\text{-}\CAlg_{\mathbb A_{\inf R}}$ for commutative $\delta$-$\mathbb A_{\inf R}$-algebras, and $\NAlg_A$ for connective non-unital animated $A$-algebras.
\item All fiber products of derived schemes, formal schemes, and presheaves are homotopy fiber products. We write $\operatorname{Tot}_{[n]\in\Delta} = \varprojlim_{[n]\in\Delta}$ for cosimplicial totalization, $|-|=\operatorname*{colim}_{[n]\in\Delta^{\mathrm{op}}}$ for geometric realization, and $[X/G]$ for the homotopy quotient, which is defined by the bar construction of the action groupoid of a group object $G$ acting on $X$ in an $\infty$-topos.
\item Throughout the arithmetic sections, $W(-)$ denotes Witt vectors, $[-]$ the Teichm\"{u}ller lift, and a $\delta$-ring carries its Frobenius lift $\varphi(x)=x^p+p\delta(x)$. We write $(-)^\wedge_I$ for derived $I$-completion (Subsection~\ref{subsec:almost-formal-base}) and $E^\wedge = \varprojlim_{m\ge 1}(E\wedge\mathbb S/p^m)$ for derived $p$-completion on $\Sp$.
\end{itemize}

\section{Derived algebraic cobordism of \texorpdfstring{$\infty$}{infinity}-categories with a universal line datum}
\label{sec:span_formalism}

Throughout this section, let $\mathcal{C}$ be a small $\infty$-category that admits finite limits, universal finite coproducts, and the additively admissible pushouts of Subsection~\ref{subsec:admissible-spans}, with terminal object $* \in \mathcal{C}$. Let $\mathrm{Pre}(\mathcal C)=\operatorname{Fun}(\mathcal C^{\mathrm{op}},\mathcal S)$ denote the $\infty$-topos of space-valued presheaves equipped with the Yoneda embedding $j\colon\mathcal C\to\mathrm{Pre}(\mathcal C)$, and let $\mathcal{P}_{\mathrm{Sp}}(\mathcal{C}) = \operatorname{Fun}(\mathcal{C}^{\mathrm{op}}, \Sp)$ denote the stable $\infty$-category of spectral presheaves endowed with the objectwise smash product $\otimes_{\mathbb{S}}$ and the connective subcategory $\mathcal{P}_{\mathrm{Sp}}(\mathcal{C})_{\ge 0}$. For any Grothendieck topology $\tau$ on $\mathcal{C}$, let $\mathcal{H}_\tau(\mathcal{C})$ denote the $\infty$-topos of $\tau$-hypercomplete sheaves with the left exact hypercomplete sheafification functor $a_\tau^{\mathrm{hyp}} \colon \mathrm{Pre}(\mathcal{C}) \to \mathcal{H}_\tau(\mathcal{C})$.

\subsection{Universal line data, vector bundles, and projective bundles}
\label{subsec:universal-line-datum}

For $\mathbf{BG}_m \in \mathrm{CMon}(\mathrm{Pre}(\mathcal{C})) = \CAlg(\mathrm{Pre}(\mathcal{C})^\times)$ with multiplication $\otimes$ and unit $1 \colon * \to \mathbf{BG}_m$, the slice $\infty$-category $\mathrm{Pre}(\mathcal{C})_{/\mathbf{BG}_m}$ carries the symmetric monoidal structure
\[
(A \to \mathbf{BG}_m) \boxtimes (B \to \mathbf{BG}_m)
= \bigl(A \times B \longrightarrow \mathbf{BG}_m \times \mathbf{BG}_m \xrightarrow{\;\otimes\;} \mathbf{BG}_m\bigr)
\]
with unit $1 \colon * \to \mathbf{BG}_m$, satisfying $\CAlg(\mathrm{Pre}(\mathcal{C})_{/\mathbf{BG}_m}, \boxtimes) \simeq \mathrm{CMon}(\mathrm{Pre}(\mathcal{C}))_{/\mathbf{BG}_m}$. Since colimits in $\mathrm{Pre}(\mathcal{C})_{/\mathbf{BG}_m}$ are created in the $\infty$-topos $\mathrm{Pre}(\mathcal{C})$, $\boxtimes$ preserves small colimits in each variable. Tensoring with $(\Sp_{\ge 0}, \otimes_{\mathbb{S}})$ in $\CAlg(\mathcal{P}\mathrm{r}^{\mathrm{L}})$~\cite[Section~4.8.1]{HA} gives the presentably symmetric monoidal $\infty$-category $\Sp_{\ge 0}(\mathrm{Pre}(\mathcal{C})_{/\mathbf{BG}_m})$, together with the strong symmetric monoidal left adjoint $\Sigma^\infty_+$ and its lax symmetric monoidal right adjoint $\Omega^\infty \colon \Sp_{\ge 0}(\mathrm{Pre}(\mathcal{C})_{/\mathbf{BG}_m}) \to \mathrm{Pre}(\mathcal{C})_{/\mathbf{BG}_m}$~\cite[Corollary~7.3.2.7]{HA}. In particular, applying $\Omega^\infty$ to any $\mathcal{L}_{\mathrm{Sp}} \in \CAlg(\Sp_{\ge 0}(\mathrm{Pre}(\mathcal{C})_{/\mathbf{BG}_m}))$ equips its underlying morphism $\pi \colon \mathcal{L} \to \mathbf{BG}_m$ with a commutative monoid structure over $\mathbf{BG}_m$ with unit $1_{\mathcal{L}} \colon * \to \mathcal{L}$ lifting $1$.

\begin{definition}
\label{def:universal-line-datum}
A \emph{universal line datum} on $\mathcal{C}$ consists of:
\begin{enumerate}[(1)]
\item a commutative monoid $\mathbf{BG}_m \in \mathrm{CMon}(\mathrm{Pre}(\mathcal{C}))$ with unit $1 \colon * \to \mathbf{BG}_m$, together with a morphism $\pi \colon \mathcal{L} \to \mathbf{BG}_m$ in $\mathrm{Pre}(\mathcal{C})$;
\item a commutative algebra $\mathcal{L}_{\mathrm{Sp}} \in \CAlg\bigl(\Sp_{\ge 0}(\mathrm{Pre}(\mathcal{C})_{/\mathbf{BG}_m})\bigr)$ equipped with an equivalence $\Omega^\infty(\mathcal{L}_{\mathrm{Sp}}) \to (\pi \colon \mathcal{L} \to \mathbf{BG}_m)$ in $\mathrm{Pre}(\mathcal{C})_{/\mathbf{BG}_m}$;
\item a subpresheaf $\mathcal{L}^\circ \subset \mathcal{L}$ closed under the induced multiplication and containing $1_{\mathcal{L}}$, such that $1_{\mathcal{L}} \colon * \to \mathcal{L}^\circ$ is an equivalence in $\mathrm{Pre}(\mathcal{C})$ and $* \times_{\mathbf{BG}_m} \mathcal{L}^\circ \subset * \times_{\mathbf{BG}_m} \mathcal{L}$ coincides with the subpresheaf of multiplicative units $(* \times_{\mathbf{BG}_m} \mathcal{L})^\times$.
\end{enumerate}
\end{definition}

Fix a universal line datum $(\mathbf{BG}_m, \mathcal{L}, \mathcal{L}_{\mathrm{Sp}}, \mathcal{L}^\circ)$ on $\mathcal{C}$. Write $1_! \colon \mathrm{Pre}(\mathcal{C}) \rightleftarrows \mathrm{Pre}(\mathcal{C})_{/\mathbf{BG}_m} \colon 1^*$ for the adjunction induced by $1 \colon * \to \mathbf{BG}_m$, and $1_!^{\mathrm{Sp}} \colon \mathcal{P}_{\mathrm{Sp}}(\mathcal{C})_{\ge 0} \rightleftarrows \Sp_{\ge 0}(\mathrm{Pre}(\mathcal{C})_{/\mathbf{BG}_m}) \colon 1^*_{\mathrm{Sp}}$ for its tensor product with $\Sp_{\ge 0}$.

\begin{proposition}
\label{prop:ring-from-line}
\begin{enumerate}[(1)]
\item The right adjoint $1^*_{\mathrm{Sp}}$ is lax symmetric monoidal, giving a connective $\mathbb{E}_\infty$-ring presheaf $\mathbb{A}^1_{\mathcal{C},\mathrm{Sp}} = 1^*_{\mathrm{Sp}}(\mathcal{L}_{\mathrm{Sp}}) \in \CAlg(\mathcal{P}_{\mathrm{Sp}}(\mathcal{C})_{\ge 0})$ whose underlying multiplicative $\mathbb{E}_\infty$-monoid presheaf is $\mathbb{A}^1_{\mathcal{C}} = \Omega^\infty \mathbb{A}^1_{\mathcal{C},\mathrm{Sp}} \simeq * \times_{\mathbf{BG}_m} \mathcal{L}$.
\item The presheaf of units $\mathbb{G}_{m,\mathcal{C}}^{\mathrm{pre}} = \mathrm{GL}_1(\mathbb{A}^1_{\mathcal{C},\mathrm{Sp}})$ admits an equivalence $\Omega_1(\mathbf{BG}_m) \to \mathbb{G}_{m,\mathcal{C}}^{\mathrm{pre}}$ in $\mathrm{CMon}(\mathrm{Pre}(\mathcal{C}))^{\mathrm{gp}}$.
\end{enumerate}
\end{proposition}
\begin{proof}
Since $1 \colon * \to \mathbf{BG}_m$ lies in $\mathrm{CMon}(\mathrm{Pre}(\mathcal{C}))$, the left adjoint $1_!$ from $(\mathrm{Pre}(\mathcal{C}), \times)$ to $(\mathrm{Pre}(\mathcal{C})_{/\mathbf{BG}_m}, \boxtimes)$ is strong symmetric monoidal and preserves small colimits. Under $(\mathrm{Pre}(\mathcal{C}), \times) \otimes (\Sp_{\ge 0}, \otimes_{\mathbb{S}}) \simeq (\mathcal{P}_{\mathrm{Sp}}(\mathcal{C})_{\ge 0}, \otimes_{\mathbb{S}})$ in $\CAlg(\mathcal{P}\mathrm{r}^{\mathrm{L}})$, tensoring $1_!$ with $\Sp_{\ge 0}$ makes $1_!^{\mathrm{Sp}}$ strong symmetric monoidal, so $1^*_{\mathrm{Sp}}$ is lax symmetric monoidal by~\cite[Corollary~7.3.2.7]{HA}. Passing to right adjoints in $1_!^{\mathrm{Sp}} \circ \Sigma^\infty_+ \simeq \Sigma^\infty_+ \circ 1_!$ gives $\Omega^\infty \circ 1^*_{\mathrm{Sp}} \simeq 1^* \circ \Omega^\infty$ as lax symmetric monoidal functors, so $\mathbb{A}^1_{\mathcal{C}} \simeq 1^*(\mathcal{L}) = * \times_{\mathbf{BG}_m} \mathcal{L}$ in $\mathrm{CMon}(\mathrm{Pre}(\mathcal{C}))$.

For~(2), Definition~\ref{def:universal-line-datum}(3) identifies $\mathrm{GL}_1(\mathbb{A}^1_{\mathcal{C},\mathrm{Sp}}) = (\mathbb{A}^1_{\mathcal{C}})^\times$ with $* \times_{\mathbf{BG}_m} \mathcal{L}^\circ$ in $\mathrm{CMon}(\mathrm{Pre}(\mathcal{C}))$. Pulling back the equivalence $1_{\mathcal{L}} \colon * \to \mathcal{L}^\circ$ in $\mathrm{CMon}(\mathrm{Pre}(\mathcal{C}))_{/\mathbf{BG}_m}$ along $1 \colon * \to \mathbf{BG}_m$ gives the equivalence $\Omega_1(\mathbf{BG}_m) = * \times_{\mathbf{BG}_m} * \to * \times_{\mathbf{BG}_m} \mathcal{L}^\circ = \mathbb{G}_{m,\mathcal{C}}^{\mathrm{pre}}$ in $\mathrm{CMon}(\mathrm{Pre}(\mathcal{C}))^{\mathrm{gp}}$.
\end{proof}

Write $A_U = \mathbb{A}^1_{\mathcal{C},\mathrm{Sp}}(U)$ for $U \in \mathcal{C}$, and let $\operatorname{Perf}(A_U) \subset \Mod_{A_U}$ denote the $\infty$-category of perfect $A_U$-modules.

\begin{definition}
\label{def:presheaf-linear-algebra}
For each integer $r \ge 0$, define the following presheaves on $\mathcal{C}$, with functoriality given by base change along $\mathbb{A}^1_{\mathcal{C},\mathrm{Sp}}$~\cite[Theorem~4.5.3.1, Theorem~4.8.5.16, and Proposition~7.2.4.2]{HA}:
\begin{enumerate}[(1)]
\item the commutative monoid presheaf $\operatorname{Perf}_{\mathcal{C}}^{\simeq} \in \mathrm{CMon}(\mathrm{Pre}(\mathcal{C}))$, given by $\operatorname{Perf}_{\mathcal{C}}^{\simeq}(U) = \operatorname{Perf}(A_U)^{\simeq}$ under $\oplus$;
\item the full subpresheaf $\mathrm{Free}_r \subset \operatorname{Perf}_{\mathcal{C}}^{\simeq}$ whose value at $U$ is the connected component of $x_r(U) = A_U^{\oplus r}$ with direct sum $\oplus \colon \mathrm{Free}_r \times \mathrm{Free}_s \to \mathrm{Free}_{r+s}$, together with $\mathbf{GL}_{r,\mathcal{C}}^{\mathrm{pre}} = \Omega_{x_r}\mathrm{Free}_r \in \operatorname{Grp}(\mathrm{Pre}(\mathcal{C}))$;
\item the projection $\pi^{\mathrm{pre}} \colon \mathcal{V}^{\mathrm{pre}} \to \operatorname{Perf}_{\mathcal{C}}^{\simeq}$, given by $\mathcal{V}^{\mathrm{pre}}(U) = (\operatorname{Perf}(A_U)_{A_U/})^{\simeq} \to \operatorname{Perf}(A_U)^{\simeq}$ with zero section $s^{\mathrm{pre}}(P) = (0 \colon A_U \to P)$, and its restriction $\pi_r^{\mathrm{pre}} \colon \mathcal{V}_r^{\mathrm{pre}} \to \mathrm{Free}_r$ with zero section $s_r^{\mathrm{pre}}$;
\item the subpresheaf $\mathrm{Epi}_r^{\mathrm{pre}} \to \mathrm{Free}_r \times \mathrm{Free}_1$ of
\[
U \longmapsto \operatorname{Fun}(\Delta^1, \operatorname{Perf}(A_U))^{\simeq} \times_{\operatorname{Perf}(A_U)^{\simeq} \times \operatorname{Perf}(A_U)^{\simeq}} \bigl(\mathrm{Free}_r(U) \times \mathrm{Free}_1(U)\bigr)
\]
spanned at any $U$ by the connected components of morphisms $q \colon P \to L$ admitting a section in $\operatorname{Perf}(A_U)$.
\end{enumerate}
\end{definition}

\begin{lemma}
\label{lem:presheaf-duality-epi}
For each integer $r \ge 0$:
\begin{enumerate}[(1)]
\item There are equivalences $[*/\mathbf{GL}_{r,\mathcal{C}}^{\mathrm{pre}}] \to \mathrm{Free}_r$ in $\mathrm{Pre}(\mathcal{C})$ and $\mathbf{GL}_{1,\mathcal{C}}^{\mathrm{pre}} \to \mathbb{G}_{m,\mathcal{C}}^{\mathrm{pre}}$ in $\operatorname{Grp}(\mathrm{Pre}(\mathcal{C}))$.
\item Dualization $P \mapsto P^\vee = \underline{\mathrm{Hom}}_{A_U}(P, A_U)$ induces an involution $(-)^\vee \colon \operatorname{Perf}_{\mathcal{C}}^{\simeq} \to \operatorname{Perf}_{\mathcal{C}}^{\simeq}$ in $\mathrm{CMon}(\mathrm{Pre}(\mathcal{C}))$ restricting to $(-)^\vee \colon \mathrm{Free}_r \to \mathrm{Free}_r$.
\item Writing $\mathrm{Epi}^{\mathrm{pre}}(\mathcal{U}_r^\vee)$ and $\mathrm{Epi}^{\mathrm{pre}}(\mathcal{U}_r^\vee \oplus \mathcal{O})$ for the pullbacks of $\mathrm{Epi}_r^{\mathrm{pre}} \to \mathrm{Free}_r$ along $(-)^\vee \colon \mathrm{Free}_r \to \mathrm{Free}_r$ and of $\mathrm{Epi}_{r+1}^{\mathrm{pre}} \to \mathrm{Free}_{r+1}$ along $(-)^\vee \oplus x_1$, there are natural morphisms over $\mathrm{Free}_r$,
\[
i_{\infty,r}^{\mathrm{pre}} \colon \mathrm{Epi}^{\mathrm{pre}}(\mathcal{U}_r^\vee) \longrightarrow \mathrm{Epi}^{\mathrm{pre}}(\mathcal{U}_r^\vee \oplus \mathcal{O}), \qquad
\iota_r^{\mathrm{pre}} \colon \mathcal{V}_r^{\mathrm{pre}} \longrightarrow \mathrm{Epi}^{\mathrm{pre}}(\mathcal{U}_r^\vee \oplus \mathcal{O}),
\]
sending $(P^\vee \xrightarrow{q} L)$ to $(P^\vee \oplus A_U \xrightarrow{(q,0)} L)$ and $(A_U \xrightarrow{x} P)$ to $(P^\vee \oplus A_U \xrightarrow{(x^\vee,\mathrm{id}_{A_U})} A_U)$, respectively, together with $z_r^{\mathrm{pre}} = \iota_r^{\mathrm{pre}} \circ s_r^{\mathrm{pre}}$. The morphism $\iota_r^{\mathrm{pre}}$ is a monomorphism identifying $\mathcal{V}_r^{\mathrm{pre}}$ with the components of $(P^\vee \oplus A_U \xrightarrow{(q_1,q_0)} L)$ where $q_0 \colon A_U \to L$ is an equivalence, and $\mathcal{V}_r^{\mathrm{pre}} \times_{\mathrm{Epi}^{\mathrm{pre}}(\mathcal{U}_r^\vee \oplus \mathcal{O})} \mathrm{Epi}^{\mathrm{pre}}(\mathcal{U}_r^\vee)$ is empty at any $U \in \mathcal{C}$ with $A_U \not\simeq 0$.
\end{enumerate}
\end{lemma}
\begin{proof}
For~(1), pointwise connectedness of $\mathrm{Free}_r$ gives $[*/\mathbf{GL}_{r,\mathcal{C}}^{\mathrm{pre}}] \xrightarrow{\sim} \mathrm{Free}_r$ via the bar construction~\cite[Theorem~5.2.6.15]{HA}, and $\mathbf{GL}_{1,\mathcal{C}}^{\mathrm{pre}}(U) \simeq \mathrm{Aut}_{A_U}(A_U) \simeq \mathrm{GL}_1(A_U) = \mathbb{G}_{m,\mathcal{C}}^{\mathrm{pre}}(U)$. For~(2), dualization on $\operatorname{Perf}(A_U)$~\cite[Section~4.6.1 and Proposition~7.2.4.4]{HA} commutes with base change and preserves finite biproducts and $A_U^{\oplus r}$; passing to cores gives the duality involution $(-)^\vee$ on $\operatorname{Perf}_{\mathcal{C}}^{\simeq}$ and $\mathrm{Free}_r$.

For~(3), naturality of $i_{\infty,r}^{\mathrm{pre}}$ and $\iota_r^{\mathrm{pre}}$ follows from the preservation of duals, biproducts, and split epimorphisms under base change. Fix $U \in \mathcal{C}$ and $P \in \mathrm{Free}_r(U)$. Since the space of pairs $(L, q_0)$ with $L \in \mathrm{Free}_1(U)$ and $q_0 \colon A_U \to L$ an equivalence is contractible, the components of $\mathrm{Epi}^{\mathrm{pre}}(\mathcal{U}_r^\vee \oplus \mathcal{O})(U)$ over $P$ where $q_0$ is an equivalence map equivalently to their fiber over $(A_U, \mathrm{id}_{A_U})$, namely $\Map_{A_U}(P^\vee, A_U) \simeq \Map_{A_U}(A_U, P)$. Thus $\iota_r^{\mathrm{pre}}$ is a pointwise inclusion of connected components, hence a monomorphism. Any point in $\mathcal{V}_r^{\mathrm{pre}}(U) \times_{\mathrm{Epi}^{\mathrm{pre}}(\mathcal{U}_r^\vee \oplus \mathcal{O})(U)} \mathrm{Epi}^{\mathrm{pre}}(\mathcal{U}_r^\vee)(U)$ equips $q_0 \colon A_U \to L$ simultaneously with an inverse and a nullhomotopy $q_0 \simeq 0$, forcing $A_U \simeq 0$.
\end{proof}

Fix a Grothendieck topology $\tau$ on $\mathcal{C}$.

\begin{definition}
\label{def:vb-pb-hypersheaf}
For each integer $r \ge 0$:
\begin{enumerate}[(1)]
\item Set $\mathbf{GL}_{r,\mathcal{C}} = a_\tau^{\mathrm{hyp}}(\mathbf{GL}_{r,\mathcal{C}}^{\mathrm{pre}})$, $\mathbb{G}_{m,\mathcal{C}} = \mathbf{GL}_{1,\mathcal{C}}$, and $\mathbf{BGL}_{r,\mathcal{C}} = a_\tau^{\mathrm{hyp}}(\mathrm{Free}_r)$, equipped with direct sum $\oplus \colon \mathbf{BGL}_{r,\mathcal{C}} \times \mathbf{BGL}_{s,\mathcal{C}} \to \mathbf{BGL}_{r+s,\mathcal{C}}$ and the dual involution $(-)^\vee \colon \mathbf{BGL}_{r,\mathcal{C}} \to \mathbf{BGL}_{r,\mathcal{C}}$ induced by $a_\tau^{\mathrm{hyp}}$. Applying $a_\tau^{\mathrm{hyp}}$ to Definition~\ref{def:presheaf-linear-algebra} and Lemma~\ref{lem:presheaf-duality-epi} defines over $\mathbf{BGL}_{r,\mathcal{C}}$ the universal total space $\pi_{\mathcal{U}_r} \colon \mathbb{V}(\mathcal{U}_r) = a_\tau^{\mathrm{hyp}}(\mathcal{V}_r^{\mathrm{pre}}) \to \mathbf{BGL}_{r,\mathcal{C}}$ with zero section $s_{\mathcal{U}_r} = a_\tau^{\mathrm{hyp}}(s_r^{\mathrm{pre}})$, the universal projective bundle $p_{\mathcal{U}_r} \colon \mathbb{P}(\mathcal{U}_r) = a_\tau^{\mathrm{hyp}}(\mathrm{Epi}_r^{\mathrm{pre}}) \to \mathbf{BGL}_{r,\mathcal{C}}$ with tautological line bundle $\mathcal{O}_{\mathbb{P}(\mathcal{U}_r)}(1)$ classified by the projection $\mathbb{P}(\mathcal{U}_r) \to \mathbf{BGL}_{1,\mathcal{C}}$, and the universal morphisms $i_{\infty,\mathcal{U}_r} = a_\tau^{\mathrm{hyp}}(i_{\infty,r}^{\mathrm{pre}})$, $\iota_{\mathcal{U}_r} = a_\tau^{\mathrm{hyp}}(\iota_r^{\mathrm{pre}})$, and $z_{\mathcal{U}_r} = a_\tau^{\mathrm{hyp}}(z_r^{\mathrm{pre}})$.
\item For $X \in \mathcal{H}_\tau(\mathcal{C})$, set $\Vect_r(X) = \Map_{\mathcal{H}_\tau(\mathcal{C})}(X, \mathbf{BGL}_{r,\mathcal{C}})$, writing $\mathcal{U}_r \in \Vect_r(\mathbf{BGL}_{r,\mathcal{C}})$ for $\mathrm{id}_{\mathbf{BGL}_{r,\mathcal{C}}}$ and $\mathcal{O}_X^{\oplus r}$ for the bundle classified by $X \to * \xrightarrow{a_\tau^{\mathrm{hyp}}(x_r)} \mathbf{BGL}_{r,\mathcal{C}}$, with $\mathcal{O}_X = \mathcal{O}_X^{\oplus 1}$. Direct sums $\mathcal{E} \oplus \mathcal{F}$ and duals $\mathcal{E}^\vee$ are defined via $\oplus$ and $(-)^\vee$. For $\mathcal{E} \in \Vect_r(X)$ classified by $c_{\mathcal{E}} \colon X \to \mathbf{BGL}_{r,\mathcal{C}}$, pullback along $c_{\mathcal{E}}$ defines $\pi_{\mathcal{E}} \colon \mathbb{V}_X(\mathcal{E}) \to X$, $p_{\mathcal{E}} \colon \mathbb{P}_X(\mathcal{E}) \to X$, $\mathcal{O}_{\mathbb{P}_X(\mathcal{E})}(1)$, and the morphisms
\[
\begin{aligned}
s_{\mathcal{E}} &\colon X \longrightarrow \mathbb{V}_X(\mathcal{E}), &
z_{\mathcal{E}} &\colon X \longrightarrow \mathbb{P}_X(\mathcal{E}^\vee \oplus \mathcal{O}_X), \\
i_\infty &\colon \mathbb{P}_X(\mathcal{E}^\vee) \longrightarrow \mathbb{P}_X(\mathcal{E}^\vee \oplus \mathcal{O}_X), &
\iota_{\mathcal{E}} &\colon \mathbb{V}_X(\mathcal{E}) \longrightarrow \mathbb{P}_X(\mathcal{E}^\vee \oplus \mathcal{O}_X).
\end{aligned}
\]
\item The \emph{unimodular locus} is $\mathbb{V}_X(\mathcal{E}^\vee)^\times = \mathbb{P}_X(\mathcal{E}) \times_{\mathbf{BGL}_{1,\mathcal{C}}} *$, formed along the classifying map of $\mathcal{O}_{\mathbb{P}_X(\mathcal{E})}(1)$ and $a_\tau^{\mathrm{hyp}}(x_1) \colon * \to \mathbf{BGL}_{1,\mathcal{C}}$.
\end{enumerate}
\end{definition}

\begin{proposition}
\label{prop:vb-pb-properties}
Let $X \in \mathcal{H}_\tau(\mathcal{C})$ and $\mathcal{E} \in \Vect_r(X)$.
\begin{enumerate}[(1)]
\item Pullback along any morphism $g \colon X' \to X$ in $\mathcal{H}_\tau(\mathcal{C})$ commutes with $\mathbb{V}_X(\mathcal{E})$, $\mathbb{P}_X(\mathcal{E})$, $\mathbb{V}_X(\mathcal{E}^\vee)^\times$, $\mathcal{O}_{\mathbb{P}_X(\mathcal{E})}(1)$, and $s_{\mathcal{E}}, z_{\mathcal{E}}, i_\infty, \iota_{\mathcal{E}}$.
\item In $\mathcal{H}_\tau(\mathcal{C})$, there are equivalences $[*/\mathbf{GL}_{r,\mathcal{C}}] \to \mathbf{BGL}_{r,\mathcal{C}}$ and $[\mathbb{V}_X(\mathcal{E}^\vee)^\times / \mathbb{G}_{m,\mathcal{C}}] \to \mathbb{P}_X(\mathcal{E})$ over $X$.
\item The morphism $\iota_{\mathcal{E}} \colon \mathbb{V}_X(\mathcal{E}) \to \mathbb{P}_X(\mathcal{E}^\vee \oplus \mathcal{O}_X)$ is a monomorphism in $\mathcal{H}_\tau(\mathcal{C})$. For $r = 0$, $z_0 \colon X \to \mathbb{P}_X(\mathcal{O}_X)$ is an equivalence; for $r = 1$, $\mathbb{P}(\mathcal{U}_1) \simeq \mathbf{BGL}_{1,\mathcal{C}}$ with $\mathcal{O}_{\mathbb{P}(\mathcal{U}_1)}(1) \simeq \mathcal{U}_1$.
\item If any $U \in \mathcal{C}$ with $A_U \simeq 0$ is covered by the empty $\tau$-sieve, then $\mathbb{P}_X(0) \simeq \varnothing$ and $X \times_{\mathbb{P}_X(\mathcal{E}^\vee \oplus \mathcal{O}_X)} \mathbb{P}_X(\mathcal{E}^\vee) \simeq \varnothing$ along $z_{\mathcal{E}}$ and $i_\infty$.
\item If $1 \colon * \to a_\tau^{\mathrm{hyp}}(\mathbf{BG}_m)$ is an effective epimorphism in $\mathcal{H}_\tau(\mathcal{C})$, then $a_\tau^{\mathrm{hyp}}(\mathbf{BG}_m) \simeq \mathbf{BGL}_{1,\mathcal{C}}$.
\end{enumerate}
\end{proposition}
\begin{proof}
Claim~(1) follows from Definition~\ref{def:vb-pb-hypersheaf}(2) and the pasting law for pullbacks. For~(2), since $a_\tau^{\mathrm{hyp}}$ preserves colimits and finite limits, Lemma~\ref{lem:presheaf-duality-epi}(1) gives $\mathbf{BGL}_{r,\mathcal{C}} \simeq [*/\mathbf{GL}_{r,\mathcal{C}}]$, and descent along $a_\tau^{\mathrm{hyp}}(x_1) \colon * \to \mathbf{BGL}_{1,\mathcal{C}} \simeq [*/\mathbb{G}_{m,\mathcal{C}}]$ gives $\mathbb{P}_X(\mathcal{E}) \simeq [\mathbb{V}_X(\mathcal{E}^\vee)^\times / \mathbb{G}_{m,\mathcal{C}}]$.

For~(3), left exactness of $a_\tau^{\mathrm{hyp}}$, pullback along $c_{\mathcal{E}}$, and Lemma~\ref{lem:presheaf-duality-epi}(3) show that $\iota_{\mathcal{E}}$ is a monomorphism. Since $A_U$ is commutative, any split epimorphism $A_U \to L$ ($L \in \mathrm{Free}_1(U)$) is an equivalence; hence $\iota_0^{\mathrm{pre}} \circ s_0^{\mathrm{pre}}$ and both projections $\mathrm{Epi}_1^{\mathrm{pre}} \to \mathrm{Free}_1$ are equivalences, giving $z_0 \colon X \xrightarrow{\sim} \mathbb{P}_X(\mathcal{O}_X)$ and $\mathbb{P}(\mathcal{U}_1) \simeq \mathbf{BGL}_{1,\mathcal{C}}$ with $\mathcal{O}_{\mathbb{P}(\mathcal{U}_1)}(1) \simeq \mathcal{U}_1$.

For~(4), since $\mathrm{Epi}_0^{\mathrm{pre}}(U)$ and the presheaf fiber product of Lemma~\ref{lem:presheaf-duality-epi}(3) are empty whenever $A_U \not\simeq 0$, the empty-sieve condition implies that both hypersheafify to $\varnothing$. Pulling back along $c_{\mathcal{E}}$ and $s_{\mathcal{E}}$ via $z_{\mathcal{E}} = \iota_{\mathcal{E}} \circ s_{\mathcal{E}}$ gives $\mathbb{P}_X(0) \simeq \varnothing$ and $X \times_{\mathbb{P}_X(\mathcal{E}^\vee \oplus \mathcal{O}_X)} \mathbb{P}_X(\mathcal{E}^\vee) \simeq \varnothing$. Claim~(5) follows from descent along $1 \colon * \to a_\tau^{\mathrm{hyp}}(\mathbf{BG}_m)$ and Proposition~\ref{prop:ring-from-line}(2).
\end{proof}

\begin{definition}
\label{def:geometric-topology}
A Grothendieck topology $\tau$ on $\mathcal{C}$ is \emph{locally trivializing} for the universal line datum if:
\begin{enumerate}[(1)]
\item the functor $a_\tau^{\mathrm{hyp}} \circ j \colon \mathcal{C}
      \hookrightarrow \mathcal{H}_\tau(\mathcal{C})$ is fully faithful
      and preserves finite coproducts;
\item the morphism $1 \colon * \to a_\tau^{\mathrm{hyp}}(\mathbf{BG}_m)$ is an effective epimorphism in $\mathcal{H}_\tau(\mathcal{C})$; and
\item an object $U \in \mathcal{C}$ satisfies $U \simeq \varnothing$ in $\mathcal{C}$ if and only if $A_U = \mathbb{A}^1_{\mathcal{C},\mathrm{Sp}}(U) \simeq 0$.
\end{enumerate}
\end{definition}

\subsection{Projective Thom objects and non-\texorpdfstring{$\mathbb{A}^1$}{A1}-localized geometry}
\label{subsec:projective-thom-objects}

Throughout the remainder of Section~\ref{sec:span_formalism}, let $\tau$ be a Grothendieck topology on $\mathcal{C}$ that is locally trivializing for the universal line datum (Definition~\ref{def:geometric-topology}), identifying any $U \in \mathcal{C}$ with $(a_\tau^{\mathrm{hyp}} \circ j)(U) \in \mathcal{H}_\tau(\mathcal{C})$. For $X \in \mathcal{C}$, a vector bundle $\mathcal{E} \in \Vect_r(X)$ is \emph{represented} if its associated objects $\mathbb{V}_X(\mathcal{E})$, $\mathbb{V}_X(\mathcal{E}^\vee)$, $\mathbb{V}_X(\mathcal{E}^\vee)^\times$, $\mathbb{P}_X(\mathcal{E})$, $\mathbb{P}_X(\mathcal{E}^\vee)$, and $\mathbb{P}_X(\mathcal{E}^\vee \oplus \mathcal{O}_X)$, along with their structural morphisms, lie in the essential image of $\mathcal{C} \hookrightarrow \mathcal{H}_\tau(\mathcal{C})$.

\begin{definition}
\label{def:projective-thom-objects} Let $\tau$ be locally trivializing
for the universal line datum, and let $X \in
\mathcal{H}_\tau(\mathcal{C})$, $\mathcal{E} \in \Vect_r(X)$, and
$i_\infty \colon \mathbb{P}_X(\mathcal{E}^\vee) \to
\mathbb{P}_X(\mathcal{E}^\vee \oplus \mathcal{O}_X)$ be the hyperplane
morphism of Definition~\ref{def:vb-pb-hypersheaf}(2).
\begin{enumerate}[(1)]
\item The \emph{projective Thom object} of $\mathcal{E}$ is the pointed cofiber
\[
\operatorname{Th}_X(\mathcal{E})
= \operatorname{cofib}\Bigl(\mathbb{P}_X(\mathcal{E}^\vee)_+ \xrightarrow{\ i_\infty\ } \mathbb{P}_X(\mathcal{E}^\vee \oplus \mathcal{O}_X)_+\Bigr)
\;\in\; \mathcal{H}_\tau(\mathcal{C})_*.
\]
\item For any spectrum-valued additive bivariant theory $E$ on $(\mathcal{C}, \mathbf{S})$ (Subsection~\ref{subsec:admissible-spans}) that admits a $\tau$-hypersheaf extension $\underline{E}(- \to Y) \colon \mathcal{H}_\tau(\mathcal{C})^{\mathrm{op}} \to \Sp$ in its contravariant variable, its reduced extension to $\mathcal{H}_\tau(\mathcal{C})_*^{\mathrm{op}} \to \Sp$ sends pointed cofibers to fibers; hence, we set
\[
\underline{E}\bigl(\operatorname{Th}_X(\mathcal{E}) \to Y\bigr)
= \operatorname{fib}\Bigl(\underline{E}\bigl(\mathbb{P}_X(\mathcal{E}^\vee \oplus \mathcal{O}_X) \to Y\bigr) \xrightarrow{\ i_\infty^*\ } \underline{E}\bigl(\mathbb{P}_X(\mathcal{E}^\vee) \to Y\bigr)\Bigr).
\]
\end{enumerate}
\end{definition}

\subsection{Admissible spans and correspondence spectra}
\label{subsec:admissible-spans}

\begin{definition}
\label{def:admissible-class}
Let $\tau$ be a locally trivializing Grothendieck topology on $\mathcal{C}$ for the universal line datum. A class $\mathbf{S} \subset \operatorname{Mor}(\mathcal{C})$ of morphisms in $\mathcal{C}$ is called \emph{admissible} if:
\begin{enumerate}[(1)]
\item the class $\mathbf{S}$ contains all equivalences, is closed under composition and finite coproducts, and is stable under homotopy pullbacks in $\mathcal{C}$;
\item for any $X \in \mathcal{C}$, the initial morphism $\varnothing \to X$ and the codiagonal $\nabla \colon X \amalg X \to X$ lie in $\mathbf{S}$; and
\item membership in $\mathbf{S}$ is $\tau$-local on the target, and for any $U \in \mathcal{C}$, any $r \ge 0$, any $\mathcal{E} \in \Vect_r(U)$, and any morphism $V \to \mathbb{V}_U(\mathcal{E})$ (resp.\ $V \to \mathbb{P}_U(\mathcal{E}^\vee \oplus \mathcal{O}_U)$) in $\mathcal{H}_\tau(\mathcal{C})$ with $V \in \mathcal{C}$, the homotopy pullback $U \times_{\mathbb{V}_U(\mathcal{E})} V$ (resp.\ $U \times_{\mathbb{P}_U(\mathcal{E}^\vee \oplus \mathcal{O}_U)} V$) in $\mathcal{H}_\tau(\mathcal{C})$ is representable in $\mathcal{C}$ and its projection to $V$ lies in $\mathbf{S}$.
\end{enumerate}
\end{definition}

When $\mathcal{E} \in \Vect_r(U)$ is represented in $\mathcal{C}$, Condition~(3) of Definition~\ref{def:admissible-class} holds automatically whenever $s_{\mathcal{E}}$ and $z_{\mathcal{E}}$ lie in $\mathbf{S}$. More generally, if pullbacks of the trivial zero sections $s_{\mathcal{O}_U^{\oplus r}}$ and $z_{\mathcal{O}_U^{\oplus r}}$ satisfy effective $\tau$-descent into $\mathcal{C}$, Condition~(3) for arbitrary $\mathcal{E} \in \Vect_r(U)$ follows from the case $\mathcal{E} = \mathcal{O}_U^{\oplus r}$ by $\tau$-local trivialization of $\mathbf{BGL}_{r,\mathcal{C}} \simeq [*/\mathbf{GL}_{r,\mathcal{C}}]$.

\begin{definition}
\label{def:admissible-spans-stacks}
Let $\tau$ be a locally trivializing topology on $\mathcal{C}$ and let $\mathbf{S} \subset \operatorname{Mor}(\mathcal{C})$ be an admissible class (Definition~\ref{def:admissible-class}).
\begin{enumerate}[\rm (1)]
\item A morphism $f \colon Y \to X$ in $\mathcal{H}_\tau(\mathcal{C})$ is \emph{$\mathbf{S}$-representable}, written $f \in \widetilde{\mathbf{S}}$, if for any $U \in \mathcal{C}$ and any morphism $U \to X$ in $\mathcal{H}_\tau(\mathcal{C})$, the homotopy fiber product $Y \times_X U$ in $\mathcal{H}_\tau(\mathcal{C})$ is representable by an object of $\mathcal{C}$ and the projection $Y \times_X U \to U$ belongs to $\mathbf{S}$.
\item For any full subcategory $\mathcal{X}_0 \subset \mathcal{H}_\tau(\mathcal{C})$ closed under finite limits and finite coproducts in $\mathcal{H}_\tau(\mathcal{C})$, set $\widetilde{\mathbf{S}}_0 = \widetilde{\mathbf{S}} \cap \operatorname{Mor}(\mathcal{X}_0)$. We write $\Span_{\widetilde{\mathbf{S}}}(\mathcal{X}_0)$ for the span $(\infty,2)$-category~\cite[Section~5]{Hoyois2017SixOps} whose objects are those of $\mathcal{X}_0$ and whose $1$-morphisms from $X$ to $Y$ are spans $[X \xleftarrow{s} Z \xrightarrow{t} Y]$ in $\mathcal{X}_0$ with backward leg $s \in \widetilde{\mathbf{S}}_0$ and arbitrary forward leg $t$, composed via homotopy pullbacks in $\mathcal{X}_0$.
\end{enumerate}
\end{definition}

\begin{lemma}
\label{lem:admissible-stacks-properties}
Let $\tau$ be a locally trivializing topology on $\mathcal{C}$ and let $\mathbf{S}$ be an admissible class on $(\mathcal{C}, \tau)$.
\begin{enumerate}[\rm (1)]
\item The class $\widetilde{\mathbf{S}}$ contains all equivalences, is closed under composition, finite coproducts, and arbitrary homotopy pullbacks in $\mathcal{H}_\tau(\mathcal{C})$, restricts on $\mathcal{C} \hookrightarrow \mathcal{H}_\tau(\mathcal{C})$ to $\mathbf{S}$, and contains the initial morphism $\varnothing \to X$ and the codiagonal $\nabla \colon X \amalg X \to X$ for any $X \in \mathcal{H}_\tau(\mathcal{C})$.
\item For any stack $X \in \mathcal{H}_\tau(\mathcal{C})$ and any vector bundle $\mathcal{E} \in \Vect_r(X)$ ($r \ge 0$), the zero section $s_{\mathcal{E}} \colon X \to \mathbb{V}_X(\mathcal{E})$ and the projective zero section $z_{\mathcal{E}} \colon X \to \mathbb{P}_X(\mathcal{E}^\vee \oplus \mathcal{O}_X)$ belong to $\widetilde{\mathbf{S}}$.
\item For any full subcategory $\mathcal{X}_0 \subset \mathcal{H}_\tau(\mathcal{C})$ closed under finite limits and finite coproducts in $\mathcal{H}_\tau(\mathcal{C})$, the span $(\infty,2)$-category $\Span_{\widetilde{\mathbf{S}}}(\mathcal{X}_0)$ is well-defined and its underlying $(\infty,1)$-category $\Span_{\widetilde{\mathbf{S}}}(\mathcal{X}_0)_1$ is semiadditive.
\end{enumerate}
\end{lemma}
\begin{proof}
For~(1), stability of $\widetilde{\mathbf{S}}$ under equivalences, composition, and homotopy pullbacks follows from Definition~\ref{def:admissible-class}(1) and the pasting law, and left exactness of $a_\tau^{\mathrm{hyp}} \circ j$ identifies $\widetilde{\mathbf{S}} \cap \operatorname{Mor}(\mathcal{C}) = \mathbf{S}$ via pullback along $\mathrm{id}_X$. Since finite coproducts in $\mathcal{H}_\tau(\mathcal{C})$ are universal and preserved by $a_\tau^{\mathrm{hyp}} \circ j$ (Definition~\ref{def:geometric-topology}(1)), testing along morphisms $U \to X$ ($U \in \mathcal{C}$) reduces closure under finite coproducts and membership of $\varnothing \to X$ and $\nabla \colon X \amalg X \to X$ to Definition~\ref{def:admissible-class}(1)--(2).

For~(2), any morphism $f \colon V \to \mathbb{V}_X(\mathcal{E})$ ($V \in \mathcal{C}$) factors through $\mathbb{V}_V(g^*\mathcal{E}) \to \mathbb{V}_X(\mathcal{E})$ for $g = \pi_{\mathcal{E}} \circ f$. By Proposition~\ref{prop:vb-pb-properties}(1) and the pasting law, $X \times_{\mathbb{V}_X(\mathcal{E})} V \simeq V \times_{\mathbb{V}_V(g^*\mathcal{E})} V$, so Definition~\ref{def:admissible-class}(3) gives $s_{\mathcal{E}} \in \widetilde{\mathbf{S}}$, and analogously $z_{\mathcal{E}} \in \widetilde{\mathbf{S}}$. Claim~(3) follows from~(1), since finite coproducts in $\mathcal{X}_0$ are universal and hence form biproducts in $\Span_{\widetilde{\mathbf{S}}}(\mathcal{X}_0)_1$.
\end{proof}
\begin{definition}\label{def:span-category-correspondences}
Let $(\mathcal{C}, \mathbf{S})$ be an $\infty$-category equipped with an admissible class.
\begin{enumerate}[(1)]
\item Let $\Span_{\mathbf S}(\mathcal{C})$ be the $(\infty, 2)$-category of spans~\cite[Section~5]{Hoyois2017SixOps} whose objects are those of $\mathcal{C}$ and whose $1$-morphisms from $X$ to $Y$ are spans $[X \xleftarrow{s} Z \xrightarrow{t} Y]$ with $s \in \mathbf S$ and arbitrary forward leg $t$, composed via homotopy pullbacks in $\mathcal{C}$.
\item By Definition~\ref{def:admissible-class}(1)--(2) and universality of finite coproducts in $\mathcal{C}$, the underlying $(\infty,1)$-category $\Span_{\mathbf S}(\mathcal{C})_1$ is semiadditive: finite coproducts in $\mathcal{C}$ form categorical biproducts in $\Span_{\mathbf S}(\mathcal{C})_1$. Hence, $\Span_{\mathbf S}(\mathcal{C})_1$ is enriched over the symmetric monoidal $\infty$-category $(\mathrm{CMon}(\mathcal{S}), \otimes)$ of $\mathbb{E}_\infty$-spaces~\cite{MR3345192}, where for $X, Y \in \mathcal{C}$ the mapping $\mathbb{E}_\infty$-space $\Corr_{\mathbf S}(X, Y) = \Span_{\mathbf S}(\mathcal{C})(X, Y)^\simeq$ carries the addition induced by disjoint union of middle terms with unit given by the empty span.
\item Applying enrichment base change along the symmetric monoidal group completion functor $(-)^{\mathrm{gp}} \colon (\mathrm{CMon}(\mathcal{S}), \otimes) \to (\Sp_{\ge 0}, \wedge)$ defines the $\Sp_{\ge 0}$-enriched \emph{additive spectral span category} $\mathcal{A}_{\mathcal{C}}$ with the same objects as $\mathcal{C}$, mapping spectra $\mathcal{A}_{\mathcal{C}}(X,Y) = \Corr_{\mathbf S}(X,Y)^{\mathrm{gp}} \in \Sp_{\ge 0}$, and coherently associative composition pairings $\Corr_{\mathbf S}(X, Y)^{\mathrm{gp}} \wedge \Corr_{\mathbf S}(Y, Z)^{\mathrm{gp}} \to \Corr_{\mathbf S}(X, Z)^{\mathrm{gp}}$.
\end{enumerate}
\end{definition}

A \emph{spectrum-valued additive bivariant theory} $E$ on $(\mathcal{C}, \mathbf{S})$ is a spectrally enriched category $\mathcal{D}_E$ receiving an $\Sp$-enriched functor $\Phi_E \colon \mathcal{A}_{\mathcal{C}} \to \mathcal{D}_E$ on the same objects. Its bivariant spectra are $E(X \to Y) = \operatorname{map}_{\mathcal{D}_E}(X, Y) \in \Sp$, where $X \to Y$ denotes the ordered pair of source and target objects rather than a morphism in $\mathcal{C}$; for $X = Y$, we write $E(X \xrightarrow{\mathrm{id}} X)$ or $E(X \to X)$. Composition in $\mathcal{D}_E$ gives associative pairings $\circ \colon E(Y \to Z) \wedge E(X \to Y) \to E(X \to Z)$ and units $\mathbf{1}_X = 1 \in \pi_0 E(X \xrightarrow{\mathrm{id}} X)$. Any morphism $f \colon P \to X$ in $\mathcal{C}$ defines a contravariant graph span $f^* = [P \xleftarrow{\mathrm{id}_P} P \xrightarrow{f} X] \in \mathcal{A}_{\mathcal{C}}(P, X)$; precomposition with $\Phi_E(f^*)$ induces the pullback $f^* \colon E(X \to Y) \to E(P \to Y)$, $\alpha \mapsto \alpha \circ \Phi_E(f^*)$. Similarly, any morphism $g \colon X \to Y$ in $\mathbf{S}$ defines a covariant span $g_* = [Y \xleftarrow{g} X \xrightarrow{\mathrm{id}_X} X] \in \mathcal{A}_{\mathcal{C}}(Y, X)$, whose image induces the pushforward $g_* \colon E(X \to Z) \to E(Y \to Z)$. When $E$ admits a $\tau$-hypersheaf extension $\underline{E}(- \to Y)$ in its contravariant variable, this extension evaluates on $\mathcal{H}_\tau(\mathcal{C})$ as in Definition~\ref{def:projective-thom-objects}(2).

Consider a homotopy pushout square in $\mathcal{C}$ formed by $i_1 \colon X_0 \to X_1$ and $i_2 \colon X_0 \to X_2$ with pushout $X = X_1 \amalg_{X_0} X_2$. This square is \emph{additively $\mathbf S$-admissible} if, for any $Y \in \mathcal{C}$, base change along $X_1 \to X$, $X_2 \to X$, and $X_0 \to X$ induces a Cartesian square of group completions of correspondence spaces $\Corr_{\mathbf S}(-, Y)^{\mathrm{gp}}$ in $\mathrm{CMon}(\mathcal{S})^{\mathrm{gp}}$ such that $\pi_0\Corr_{\mathbf S}(X_1, Y)^{\mathrm{gp}} \oplus \pi_0\Corr_{\mathbf S}(X_2, Y)^{\mathrm{gp}} \to \pi_0\Corr_{\mathbf S}(X_0, Y)^{\mathrm{gp}}$ is surjective.

\begin{proposition}
\label{prop:admissible-pushout-excision}
For any additively $\mathbf S$-admissible homotopy pushout square $X = X_1 \amalg_{X_0} X_2$ in $\mathcal{C}$ and any $Y \in \mathcal{C}$, restriction induces a Cartesian square in $\Sp$:
\[
\xymatrix{
\Corr_{\mathbf S}(X, Y)^{\mathrm{gp}} \ar[r] \ar[d] & \Corr_{\mathbf S}(X_1, Y)^{\mathrm{gp}} \ar[d] \\
\Corr_{\mathbf S}(X_2, Y)^{\mathrm{gp}} \ar[r] & \Corr_{\mathbf S}(X_0, Y)^{\mathrm{gp}}.
}
\]
\end{proposition}
\begin{proof}
By pullback stability of $\mathbf S$ (Definition~\ref{def:admissible-class}(1)), base change along the structural maps from $X_1$, $X_2$, and $X_0$ to $X$ defines a Cartesian square in $\mathrm{CMon}(\mathcal S)^{\mathrm{gp}} \simeq \Sp_{\ge 0}$. By the Mayer--Vietoris exact sequence, surjectivity on $\pi_0$ implies that the homotopy pullback in $\Sp$ has vanishing $\pi_{-1}$ and hence is connective, giving a Cartesian square in $\Sp$.
\end{proof}

For spans $\alpha=[X \xleftarrow{a} A \to Z_A]$ and $\beta=[Y \xleftarrow{b} B \to Z_B]$ in $\Span_{\mathbf S}(\mathcal C)$ and a morphism $f \colon X \to Y$ in $\mathbf S$, let $f^*\beta$ denote the pullback span along $f$, let $f_*\alpha$ be the span obtained by composing $a$ with $f$, and let $\cdot_X$ and $\cdot_Y$ denote the relative products formed by homotopy fiber products of middle terms over $X$ and $Y$.

\begin{proposition}
\label{prop:general-projection-formula}
Let $f \colon X \to Y$ belong to $\mathbf S$, and let $\alpha=[X \xleftarrow{a} A \to Z_A]$ and $\beta=[Y \xleftarrow{b} B \to Z_B]$ be spans in $\Span_{\mathbf S}(\mathcal C)$. There is an equivalence
\[
f_*\bigl(\alpha \mathbin{\cdot_X} f^*\beta\bigr) \longrightarrow f_*\alpha \mathbin{\cdot_Y} \beta
\qquad\text{in }\Corr_{\mathbf S}(Y, Z_A \times Z_B),
\]
inducing corresponding identities in additive span spectra and on homotopy groups.
\end{proposition}
\begin{proof}
Since $f^*\beta$ has middle term $X \times_Y B$, the middle term of $f_*(\alpha \mathbin{\cdot_X} f^*\beta)$ is the homotopy pullback $A \times_X (X \times_Y B) \simeq A \times_Y B$, which coincides with the middle term of $f_*\alpha \mathbin{\cdot_Y} \beta$ compatibly with the projections to $Y$ and $Z_A \times Z_B$~\cite[Section~2]{KatoMGL}. Passing to group completions and homotopy groups implies the remaining assertions.
\end{proof}

\begin{proposition}
\label{prop:zero-section-class}
Let $\mathcal{X}_0 \subset \mathcal{H}_\tau(\mathcal{C})$ be a full subcategory closed under finite limits and finite coproducts in $\mathcal{H}_\tau(\mathcal{C})$, and let $\mathcal{E} \in \Vect_r(X)$ be a vector bundle on $X \in \mathcal{X}_0$ such that $\mathbb{V}_X(\mathcal{E}) \in \mathcal{X}_0$. The homotopy self-intersection $Z_{\mathcal{E}} = X \times_{\mathbb{V}_X(\mathcal{E})} X$ of the zero section $s_{\mathcal{E}} \colon X \to \mathbb{V}_X(\mathcal{E})$ defines a span $[X \xleftarrow{p_1} Z_{\mathcal{E}} \xrightarrow{p_2} X]$ in $\Span_{\widetilde{\mathbf{S}}}(\mathcal{X}_0)$ with a homotopy $p_1 \simeq p_2$ over $X$, whose class $\widetilde{e}(\mathcal{E})\in\pi_0\Corr_{\widetilde{\mathbf{S}}_0}(X,X)^{\mathrm{gp}}$ is functorial under pullbacks along arbitrary morphisms $f \colon Y \to X$ in $\mathcal{X}_0$ with $\mathbb{V}_Y(f^*\mathcal{E}) \in \mathcal{X}_0$.
\end{proposition}
\begin{proof}
Since $s_{\mathcal{E}} \in \widetilde{\mathbf{S}}$ by Lemma~\ref{lem:admissible-stacks-properties}(2) and $\mathcal{X}_0$ is closed under finite limits in $\mathcal{H}_\tau(\mathcal{C})$, pullback stability implies that $Z_{\mathcal{E}} = X \times_{\mathbb{V}_X(\mathcal{E})} X \in \mathcal{X}_0$ and $p_1 \colon Z_{\mathcal{E}} \to X$ lies in $\widetilde{\mathbf{S}}_0$. The bundle projection $\pi_{\mathcal{E}} \colon \mathbb{V}_X(\mathcal{E}) \to X$ satisfies $\pi_{\mathcal{E}} \circ s_{\mathcal{E}} \simeq \mathrm{id}_X$, which gives $p_1 \simeq \pi_{\mathcal{E}} \circ s_{\mathcal{E}} \circ p_1 \simeq \pi_{\mathcal{E}} \circ s_{\mathcal{E}} \circ p_2 \simeq p_2$. Thus $(Z_{\mathcal{E}}, p_1, p_2)$ lies in the image of the diagonal inclusion $(\mathcal{X}_{0/X}^{\widetilde{\mathbf{S}}_0})^{\simeq} \to \Corr_{\widetilde{\mathbf{S}}_0}(X, X)$ and defines $\widetilde{e}(\mathcal{E}) \in \pi_0 \Corr_{\widetilde{\mathbf{S}}_0}(X, X)^{\mathrm{gp}}$, acting on $\Corr_{\widetilde{\mathbf{S}}_0}(X, -)$ by relative fiber product over $X$.

For any morphism $f \colon Y \to X$ in $\mathcal{X}_0$ with $\mathbb{V}_Y(f^*\mathcal{E}) \in \mathcal{X}_0$, the zero section $s_{f^*\mathcal{E}} \colon Y \to \mathbb{V}_Y(f^*\mathcal{E})$ is the homotopy pullback of $s_{\mathcal{E}}$ along $\mathbb{V}_Y(f^*\mathcal{E}) \to \mathbb{V}_X(\mathcal{E})$ (Proposition~\ref{prop:vb-pb-properties}(1)), so the pasting law for homotopy pullbacks gives $Y \times_X (X \times_{\mathbb{V}_X(\mathcal{E})} X) \simeq Y \times_{\mathbb{V}_Y(f^*\mathcal{E})} Y$ over $Y$ and hence $f^*\widetilde{e}(\mathcal{E}) = \widetilde{e}(f^*\mathcal{E})$.
\end{proof}

\begin{remark}
\label{rem:first-chern-class-notation}
For a line bundle $\mathcal{L} \in \Vect_1(X)$, we define its first Chern class by $\widetilde{c}_1(\mathcal{L})=\widetilde{e}(\mathcal{L})\in\pi_0\Corr_{\widetilde{\mathbf{S}}_0}(X,X)^{\mathrm{gp}}$.
\end{remark}

\subsection{Orientations and oriented cohomology theories}
\label{subsec:orientable-cohomology}

\begin{definition}
\label{def:zero-section-preorientation}
Let $\mathcal{C}$ be equipped with a universal line datum (Definition~\ref{def:universal-line-datum}), the projective bundle construction, and a locally trivializing topology $\tau$ (Definition~\ref{def:geometric-topology}), and let $E$ be a spectrum-valued additive bivariant theory on $(\mathcal{C}, \mathbf{S})$ (Subsection~\ref{subsec:admissible-spans}).
\begin{enumerate}[(1)]
\item A \emph{preorientation} on $E$ is an assignment of a class $\widetilde{c}_1(\mathcal{L}) \in \pi_0 E(X \xrightarrow{\mathrm{id}} X)$ to any line bundle $\mathcal{L}$ on $X$, compatible with pullbacks.
\item Suppose that $E$ is equipped with a preorientation and admits a $\tau$-hypersheaf extension $\underline{E}(- \to Y)$. A \emph{vector bundle twist} on $E$ is an assignment $\mathcal{E} \mapsto \Sigma^{\mathcal{E}}$ of spectrum autoequivalences compatible with pullbacks and satisfying $\Sigma^0 \simeq \mathrm{id}$ and $\Sigma^{\mathcal{E} \oplus \mathcal{F}} \simeq \Sigma^{\mathcal{E}} \Sigma^{\mathcal{F}}$, with inverse $\Sigma^{-\mathcal{E}}$. A \emph{vector bundle orientation} on $E$ consists of a vector bundle twist together with natural Thom equivalences
\[
\operatorname{th}_{\mathcal{E}} \colon \Sigma^{-\mathcal{E}} \underline{E}(X \to Y) \longrightarrow \underline{E}\bigl(\operatorname{Th}_X(\mathcal{E}) \to Y\bigr)
\]
induced by the projective zero sections, which are normalized at the zero bundle, compatible with pullbacks, multiplicative under direct sums, and satisfy $z_{\mathcal{L}}^*(\operatorname{th}_{\mathcal{L}}(1)) = \widetilde{c}_1(\mathcal{L})$ for the unit $1 = \mathbf{1}_X \in \pi_0 E(X \xrightarrow{\mathrm{id}} X)$ along the projective zero section $z_{\mathcal{L}} \colon X \hookrightarrow \mathbb{P}_X(\mathcal{L}^\vee \oplus \mathcal{O}_X)$ of any line bundle $\mathcal{L}$.
\end{enumerate}
\end{definition}

\begin{definition}
\label{def:oriented-cohomology-theory}
Under the hypotheses of Definition~\ref{def:zero-section-preorientation}, a spectrum-valued additive bivariant theory $E$ on $(\mathcal{C}, \mathbf{S})$ admitting a $\tau$-hypersheaf extension $\underline{E}(- \to Y)$ is called an \emph{oriented cohomology theory} if it satisfies:
\begin{enumerate}[(1)]
\item the theory $E$ transforms additively $\mathbf{S}$-admissible pushout squares into Cartesian squares of spectra (Proposition~\ref{prop:admissible-pushout-excision});
\item the theory $E$ carries a preorientation and a compatible vector bundle orientation (Definition~\ref{def:zero-section-preorientation}); and
\item for each integer $r \ge 0$, any projective bundle $p_{\mathcal{E}} \colon \mathbb{P}_X(\mathcal{E}) \to X$ associated with $\mathcal{E} \in \Vect_{r+1}(X)$ ($X \in \mathcal{H}_\tau(\mathcal{C})$) with tautological line bundle $\mathcal{O}(1) = \mathcal{O}_{\mathbb{P}_X(\mathcal{E})}(1) \in \Vect_1(\mathbb{P}_X(\mathcal{E}))$, and any $Y \in \mathcal{C}$, the morphism
\[
\theta_{\mathcal{E}} \colon \bigoplus_{i=0}^r \underline{E}(X \to Y) \longrightarrow \underline{E}\bigl(\mathbb{P}_X(\mathcal{E}) \to Y\bigr), \qquad (\alpha_i)_{0 \le i \le r} \longmapsto \sum_{i=0}^r p_{\mathcal{E}}^*(\alpha_i) \circ \widetilde{c}_1(\mathcal{O}(1))^i,
\]
is an equivalence of spectra.
\end{enumerate}

When $E$ is defined on $\mathcal{A}_{\mathcal{X}_0}$ and satisfies the projective bundle property for the universal rank-$r$ bundle $\mathcal{U}_r$ on $\mathbf{BGL}_{r,\mathcal{C}}$ ($r \ge 1$), let $t = \widetilde{c}_1(\mathcal{O}(1))$ for $\mathcal{O}(1) = \mathcal{O}_{\mathbb{P}(\mathcal{U}_r)}(1)$. Expressing $t^r$ in the basis $\{1, t, \dots, t^{r-1}\}$ via the relation $t^r - p_{\mathcal{U}_r}^*c_1^{\mathrm{univ}} t^{r-1} + \cdots + (-1)^r p_{\mathcal{U}_r}^*c_r^{\mathrm{univ}} = 0$ uniquely determines the \emph{universal Chern classes} $c_i^{\mathrm{univ}} \in \pi_0 E(\mathbf{BGL}_{r,\mathcal{C}} \xrightarrow{\mathrm{id}} \mathbf{BGL}_{r,\mathcal{C}})$ for $1 \le i \le r$, with $c_1^{\mathrm{univ}} = \widetilde{c}_1(\mathcal{U}_1)$ since $\mathbb{P}(\mathcal{U}_1) \simeq \mathbf{BGL}_{1,\mathcal{C}}$ and $\mathcal{O}_{\mathbb{P}(\mathcal{U}_1)}(1) \simeq \mathcal{U}_1$. Pulling back along the classifying map $c_{\mathcal{E}} \colon X \to \mathbf{BGL}_{r,\mathcal{C}}$ defines the Chern classes $c_i(\mathcal{E}) = c_{\mathcal{E}}^* c_i^{\mathrm{univ}}$ for any $\mathcal{E} \in \Vect_r(X)$, which satisfy the Whitney sum formula whenever $E$ satisfies the splitting principle and multiplicativity under direct sums~\cite[Definition~4.10 and Theorem~4.15]{KatoMGL}.
\end{definition}

\subsection{Additive and multiplicative Chern laws}
\label{subsec:chern-law-refinements}

\begin{definition}
\label{def:additive-multiplicative-chern-laws}
Let $E$ be an oriented cohomology theory on $(\mathcal{C}, \mathbf{S})$ (Definition~\ref{def:oriented-cohomology-theory}), with first Chern class $\widetilde{c}_1(\mathcal{L}) \in \pi_0 E(X \xrightarrow{\mathrm{id}} X)$ for a line bundle $\mathcal{L}$ on $X$ and multiplication given by composition.
\begin{enumerate}[(1)]
\item The theory $E$ satisfies the \emph{additive Chern law}, or is \emph{$\mathbb{G}_a$-oriented}, if $\widetilde{c}_1(\mathcal{O}_X) = 0$ and $\widetilde{c}_1(\mathcal{L} \otimes \mathcal{L}') = \widetilde{c}_1(\mathcal{L}) + \widetilde{c}_1(\mathcal{L}')$ for all line bundles $\mathcal{L}, \mathcal{L}'$ on $X$.
\item Given a central coefficient class $\beta$ compatible with pullbacks, $E$ satisfies the \emph{multiplicative Chern law with parameter $\beta$}, or is \emph{$\mathbb{G}_m$-oriented}, if $\widetilde{c}_1(\mathcal{O}_X) = 0$ and $\widetilde{c}_1(\mathcal{L} \otimes \mathcal{L}') = \widetilde{c}_1(\mathcal{L}) + \widetilde{c}_1(\mathcal{L}') - \beta_X \widetilde{c}_1(\mathcal{L}) \widetilde{c}_1(\mathcal{L}')$ for all line bundles $\mathcal{L}, \mathcal{L}'$ on $X$, where $\beta_X$ is the pullback of $\beta$ to $X$. When $\beta$ is invertible, the classes $u(\mathcal{L}) = \mathbf{1}_X - \beta_X \widetilde{c}_1(\mathcal{L})$ satisfy $u(\mathcal{L} \otimes \mathcal{L}') = u(\mathcal{L}) \circ u(\mathcal{L}')$.
\end{enumerate}
In a graded refinement, these identities are interpreted homogeneously in the endomorphism ring equipped with the total grading, and the same conditions apply to any theory equipped with a preorientation.
\end{definition}

\begin{definition}
\label{def:chern-law-span-quotient}
Let $\mathcal{A}$ be an additive spectral span category, with stable presheaf envelope $\mathcal{M}_0(\mathcal{A}) = \operatorname{Fun}_{\Sp}(\mathcal{A}^{\mathrm{op}}, \Sp)$ and Yoneda embedding $m_0 \colon \mathcal{A} \to \mathcal{M}_0(\mathcal{A})$. For a line bundle $\mathcal{L}$ on $X$, let $\widetilde{c}_1(\mathcal{L}) \in \pi_0 \mathcal{A}(X,X) \simeq \pi_0 \operatorname{map}_{\mathcal{M}_0(\mathcal{A})}(m_0(X), m_0(X))$ be the class of its zero section, with products given by composition, and write $\operatorname{Pic}(X) = \pi_0 \Vect_1(X)$ for the group of isomorphism classes of line bundles on $X$. For a fixed central class $\beta$ compatible with pullbacks, consider the families of relation endomorphisms $\delta \colon \mathbb{S} \to \mathcal{A}(X,X)$ as $X$ ranges over the objects of $\mathcal{A}$ and $\mathcal{L},\mathcal{L}'$ range over $\operatorname{Pic}(X)$ given by:
\begin{itemize}
\item the endomorphisms $\widetilde{c}_1(\mathcal{O}_X)$ and $\widetilde{c}_1(\mathcal{L}\otimes\mathcal{L}')-\widetilde{c}_1(\mathcal{L})-\widetilde{c}_1(\mathcal{L}')$ for the $\mathbb{G}_a$-relations $\mathcal{I}_{\mathbb{G}_a}$; and
\item the endomorphisms $\widetilde{c}_1(\mathcal{O}_X)$ and $\widetilde{c}_1(\mathcal{L}\otimes\mathcal{L}')-\widetilde{c}_1(\mathcal{L})-\widetilde{c}_1(\mathcal{L}')+\beta_X\widetilde{c}_1(\mathcal{L})\widetilde{c}_1(\mathcal{L}')$ for the $\mathbb{G}_m$-relations $\mathcal{I}_{\mathbb{G}_m,\beta}$.
\end{itemize}
In the presentable $\infty$-category of symmetric monoidal $\Sp$-enriched $\infty$-categories, let $q_{\mathcal{I}} \colon \mathcal{A} \to \mathcal{A} / \mathcal{I}$ ($\mathcal{I} \in \{\mathcal{I}_{\mathbb{G}_a}, \mathcal{I}_{\mathbb{G}_m, \beta}\}$) denote the universal symmetric monoidal spectral quotient category obtained by coning off the two-sided tensor ideal generated by the morphisms $\delta$~\cite{MR3345192}. Enriched left Kan extension along $q_{\mathcal{I}}$ defines the symmetric monoidal base-change functor $(q_{\mathcal{I}})_! \colon \mathcal{M}_0(\mathcal{A}) \to \operatorname{Fun}_{\Sp}((\mathcal{A}/\mathcal{I})^{\mathrm{op}}, \Sp)$ preserving small colimits~\cite[Proposition~4.8.1.10 and Corollary~4.8.1.12]{HA}, which universally equips the images of the generating morphisms $\delta$ with coherent nullhomotopies compatible with composition and external products.
\end{definition}

\subsection{Derived algebraic cobordism}
\label{subsec:derived-algebraic-cobordism}

Choose an essentially small full subcategory $\mathcal{C}_0 \subset \mathcal{C}$ generating $\mathcal{H}_\tau(\mathcal{C})$ under $\tau$-colimits and closed under fiber products along $\mathbf{S}$, finite coproducts, and the selected additively $\mathbf{S}$-admissible pushouts. Let $\mathcal{X}_0 \subset \mathcal{H}_\tau(\mathcal{C})$ be the smallest essentially small full subcategory containing $\mathcal{C}_0$ and $\{\mathbf{BGL}_{r,\mathcal{C}}, \mathbb{V}(\mathcal{U}_r), \mathbb{P}(\mathcal{U}_r)\}_{r \ge 0}$ and closed under finite limits and finite coproducts in $\mathcal{H}_\tau(\mathcal{C})$; by Proposition~\ref{prop:vb-pb-properties}(1), $\mathcal{X}_0$ is closed under $\mathbb{V}_X(\mathcal{E})$, $\mathbb{P}_X(\mathcal{E})$, and $Z_{\mathcal{E}}$ for all $X \in \mathcal{X}_0$ and $\mathcal{E} \in \Vect_r(X)$. Let $\mathcal{A}_{\mathcal{X}_0}$ denote the symmetric monoidal additive spectral span $\infty$-category on $\mathcal{X}_0$ with backward legs in $\widetilde{\mathbf{S}}_0$ and arbitrary forward legs, and set $\mathcal{M}_0(\mathcal{C}) = \operatorname{Fun}_{\mathrm{Sp}}(\mathcal{A}_{\mathcal{X}_0}^{\mathrm{op}}, \mathrm{Sp})$ with Day convolution symmetric monoidal structure and spectral Yoneda functor $m_0 \colon \mathcal{A}_{\mathcal{X}_0} \to \mathcal{M}_0(\mathcal{C})$.

\begin{definition}
\label{def:relation-dmgl}
In the ungraded setting, take the trivial vector bundle twist $\Sigma^{\mathcal{E}} = \mathrm{id}$. Let $\mathcal{J}_{\dMGL} \subset \mathcal{M}_0(\mathcal{C})$ be the smallest localizing tensor-ideal subcategory containing the cofibers of the following four small families of morphisms, and let $W_{\dMGL}$ be the strongly saturated tensor-ideal class of morphisms whose cofibers lie in $\mathcal{J}_{\dMGL}$:
\begin{enumerate}[\rm (1)]
\item the assembly morphisms $m_0(X_1) \amalg_{m_0(X_{12})} m_0(X_2) \to m_0(X)$ induced via contravariant graph spans by the selected additively $\mathbf{S}$-admissible pushout squares $X = X_1 \amalg_{X_{12}} X_2$ in $\mathcal{C}_0$;
\item for any $X \in \mathcal{X}_0$ and any rank-$(r+1)$ vector bundle $\mathcal{E}$ on $X$ ($r \ge 0$), the projective bundle morphism
\[
\vartheta_{\mathcal{E}} = \bigl(m_0(p_{\mathcal{E}}^*) \circ m_0(\widetilde{c}_1(\mathcal{O}(1)))^j\bigr)_{0 \le j \le r} \colon m_0\bigl(\mathbb{P}_X(\mathcal{E})\bigr) \longrightarrow \bigoplus_{j=0}^{r} m_0(X);
\]
\item the assembly morphisms $\varinjlim_{[n]\in\Delta^{\mathrm{op}}} \bigoplus_{\alpha \in I_n} m_0(U_{n,\alpha}) \to m_0(X)$ for a small generating set of $\tau$-hypercovers $\bigl\{\coprod_{\alpha \in I_n} U_{n,\alpha} \to X\bigr\}_{[n]\in\Delta^{\mathrm{op}}}$ in $\mathcal{X}_0$ containing a $\mathcal{C}_0$-hypercover of any $X \in \mathcal{X}_0$;
\item for each $r \ge 0$ and any pullback $\mathcal{E} = c_{\mathcal{E}}^*\mathcal{U}_r$ along $c_{\mathcal{E}} \colon X \to \mathbf{BGL}_{r,\mathcal{C}}$ with $X \in \mathcal{X}_0$, setting $m_0\bigl(\operatorname{Th}_X(\mathcal{E})\bigr) = \operatorname{cofib}\bigl(m_0(\mathbb{P}_X(\mathcal{E}^\vee)) \xrightarrow{m_0(i_\infty^*)} m_0(\mathbb{P}_X(\mathcal{E}^\vee \oplus \mathcal{O}_X))\bigr)$, the morphism
\[
\varrho_{\mathcal{E}} \colon m_0\bigl(\operatorname{Th}_X(\mathcal{E})\bigr) \longrightarrow m_0(X)
\]
induced by $(z_{\mathcal{E}})_* \in \mathcal{A}_{\mathcal{X}_0}\bigl(\mathbb{P}_X(\mathcal{E}^\vee \oplus \mathcal{O}_X), X\bigr)$ via the disjointness identity $(z_{\mathcal{E}})_* \circ i_\infty^* \simeq 0$ (Proposition~\ref{prop:vb-pb-properties}(4)).
\end{enumerate}
We set $\mathcal{M}(\mathcal{C}) = \mathcal{M}_0(\mathcal{C})/\mathcal{J}_{\dMGL} \simeq W_{\dMGL}^{-1}\mathcal{M}_0(\mathcal{C})$, write $\ell \colon \mathcal{M}_0(\mathcal{C}) \to \mathcal{M}(\mathcal{C})$ for the localization functor and $m = \ell \circ m_0 \colon \mathcal{A}_{\mathcal{X}_0} \to \mathcal{M}(\mathcal{C})$, and set $\dMGL_{\mathcal{C}}(X \to Y) = \operatorname{map}_{\mathcal{M}(\mathcal{C})}(m(X), m(Y))$ for $X,Y \in \mathcal{X}_0$, with canonical map $q_{\mathcal{C}} \colon \Corr_{\widetilde{\mathbf{S}}_0}(X,Y)^{\mathrm{gp}} \to \dMGL_{\mathcal{C}}(X \to Y)$. By (3), for any $Y \in \mathcal{X}_0$, the $\tau$-hypersheaf $U \mapsto \dMGL_{\mathcal{C}}(U \to Y)$ on $\mathcal{C}_0$ extends uniquely by $\tau$-hypersheaf right Kan extension to $\underline{\dMGL}_{\mathcal{C}}(- \to Y) \colon \mathcal{H}_\tau(\mathcal{C})^{\mathrm{op}} \to \mathrm{Sp}$, agreeing with $\dMGL_{\mathcal{C}}(X \to Y)$ for all $X \in \mathcal{X}_0$.
\end{definition}

\begin{proposition}
\label{thm:localized-span-universal-property}
The derived cobordism quotient category $\mathcal{M}(\mathcal{C})=\mathcal{M}_0(\mathcal{C})/\mathcal{J}_{\dMGL}$ is a presentable stable $\infty$-category and inherits a unique presentably symmetric monoidal structure for which the localization functor $\ell\colon\mathcal{M}_0(\mathcal{C})\to\mathcal{M}(\mathcal{C})$ is exact and symmetric monoidal. The mapping spectra $\dMGL_{\mathcal{C}}(X\to Y)=\operatorname{map}_{\mathcal{M}(\mathcal{C})}(m(X),m(Y))$ inherit associative composition pairings and external products. For any presentable stable $\infty$-category $\mathcal{D}$, precomposition with $\ell$ restricts to an equivalence between:
\begin{enumerate}[(a)]
\item functors $\mathcal{M}(\mathcal{C}) \to \mathcal{D}$ preserving small colimits; and
\item functors $\mathcal{M}_0(\mathcal{C}) \to \mathcal{D}$ preserving small colimits that invert the admissible pushout assembly morphisms, the projective bundle morphisms $\vartheta_{\mathcal{E}}$, the $\tau$-presentation hypercover assembly morphisms, and the projective Thom comparison morphisms $\varrho_{\mathcal{E}}$.
\end{enumerate}
\end{proposition}
\begin{proof}
By~\cite[Sections~1.4.4 and~4.8.1]{HA}, $\mathcal{M}_0(\mathcal{C})$ is presentable, stable, and forms the free stable cocompletion of $\mathcal{A}_{\mathcal{X}_0}$ generated by the representables $m_0(X)$. Within the chosen universe, the generating morphisms in (1)--(4) of Definition~\ref{def:relation-dmgl} and their Day convolution translates against representables $\{m_0(Z)\}_{Z \in \mathcal{X}_0}$ form a small set, so $W_{\dMGL}$ is strongly saturated and of small generation. Thus, the accessible exact localization $\ell\colon \mathcal{M}_0(\mathcal{C}) \to W_{\dMGL}^{-1}\mathcal{M}_0(\mathcal{C}) \simeq \mathcal{M}_0(\mathcal{C})/\mathcal{J}_{\dMGL}$ exists by~\cite[Propositions~5.5.4.15 and~5.5.4.20]{HT}. Since $W_{\dMGL}$ is closed under Day convolution against all objects,~\cite[Proposition~2.2.1.9]{HA} endows $\mathcal{M}(\mathcal{C})$ with a unique presentably symmetric monoidal structure making $\ell$ symmetric monoidal, preserving composition and external products on mapping spectra, and satisfying the universal property.
\end{proof}

For any $X \in \mathcal{X}_0$ and line bundle $\mathcal{L}$ on $X$, we also write $\widetilde{c}_1(\mathcal{L})\in\pi_0\underline{\dMGL}_{\mathcal{C}}(X\to X)$ for the image under $q_{\mathcal{C}}\colon\Corr_{\widetilde{\mathbf{S}}_0}(X,X)^{\mathrm{gp}}\to\underline{\dMGL}_{\mathcal{C}}(X\to X)$ of the first Chern class of Proposition~\ref{prop:zero-section-class}.

\begin{proposition}
\label{thm:full-orientation-from-line-thom}
The derived algebraic cobordism theory $\underline{\dMGL}_{\mathcal{C}}$ satisfies the following properties for any $Y \in \mathcal{X}_0$:
\begin{enumerate}[(1)]
\item For any $X \in \mathcal{X}_0$ and any line bundle $\mathcal{L} \in \Vect_1(X)$ with zero section $s_{\mathcal{L}}\colon X \to \mathbb{V}_X(\mathcal{L})$, we have $\widetilde{c}_1(\mathcal{L}) = q_{\mathcal{C}}((s_{\mathcal{L}})_*) \circ q_{\mathcal{C}}(s_{\mathcal{L}}^*)$ in $\underline{\dMGL}_{\mathcal{C}}(X\to X)$, and precomposition with $\widetilde{c}_1(\mathcal{L})$ acts on $\underline{\dMGL}_{\mathcal{C}}(X\to Y)$ as the derived self-intersection operator $s_{\mathcal{L}}^*(s_{\mathcal{L}})_*$.
\item For any $X \in \mathcal{H}_\tau(\mathcal{C})$ and any rank-$(r+1)$ vector bundle $\mathcal{E}$ on $X$ ($r \ge 0$) with projective bundle $p_{\mathcal{E}}\colon \mathbb{P}_X(\mathcal{E}) \to X$, the projective bundle morphism
\[
\theta_{\mathcal{E}} \colon \bigoplus_{j=0}^r \underline{\dMGL}_{\mathcal{C}}(X\to Y) \longrightarrow \underline{\dMGL}_{\mathcal{C}}\bigl(\mathbb{P}_X(\mathcal{E})\to Y\bigr),
\qquad (\alpha_j)_{0\le j\le r} \longmapsto \sum_{j=0}^r p_{\mathcal{E}}^*(\alpha_j) \circ \widetilde{c}_1(\mathcal{O}(1))^j,
\]
is an equivalence.
\item With the trivial vector bundle twist $\Sigma^{\mathcal{E}} = \mathrm{id}$ of Definition~\ref{def:relation-dmgl}, for any $X \in \mathcal{H}_\tau(\mathcal{C})$ and any vector bundle $\mathcal{E}$ on $X$, $\underline{\dMGL}_{\mathcal{C}}$ satisfies the normalized projective Thom equivalence
\[
\operatorname{th}_{\mathcal{E}} \colon \underline{\dMGL}_{\mathcal{C}}(X \to Y) \longrightarrow \underline{\dMGL}_{\mathcal{C}}\bigl(\operatorname{Th}_X(\mathcal{E}) \to Y\bigr).
\]
\end{enumerate}
In particular, $\underline{\dMGL}_{\mathcal{C}}$ restricts to an oriented cohomology theory on $(\mathcal{C}_0,\mathbf{S})$ and is initial among all $\tau$-hypercomplete oriented realizations of $\mathcal{A}_{\mathcal{X}_0}$ satisfying admissible pushout excision, projective bundle decompositions, $\tau$-presentation hyperdescent, and normalized projective Thom equivalences.
\end{proposition}
\begin{proof}
For~(1), in $\mathcal{A}_{\mathcal{X}_0}$ the spans $(s_{\mathcal{L}})_*$ and $s_{\mathcal{L}}^*$ compose via homotopy pullback to $[X \leftarrow X\times_{\mathbb{V}_X(\mathcal{L})}X \to X]$, which represents $\widetilde{c}_1(\mathcal{L})$ (Proposition~\ref{prop:zero-section-class}). Functoriality of $q_{\mathcal{C}}$ then gives $\widetilde{c}_1(\mathcal{L}) = q_{\mathcal{C}}((s_{\mathcal{L}})_*)\circ q_{\mathcal{C}}(s_{\mathcal{L}}^*)$, whose precomposition action on $\underline{\dMGL}_{\mathcal{C}}(X\to Y)$ is $s_{\mathcal{L}}^*(s_{\mathcal{L}})_*$.

For~(2) and~(3), when $X \in \mathcal{X}_0$, Definition~\ref{def:relation-dmgl}(2) and~(4) invert $\vartheta_{\mathcal{E}}$ and $\varrho_{\mathcal{E}}$ in $\mathcal{M}(\mathcal{C})$. Applying $\operatorname{map}_{\mathcal{M}(\mathcal{C})}(-, m(Y))$ gives the equivalences $\theta_{\mathcal{E}}$ and $\operatorname{th}_{\mathcal{E}}$, where $\operatorname{Th}_X(\mathcal{E}) \simeq c_{\mathcal{E}}^*\operatorname{Th}_{\mathbf{BGL}_{r,\mathcal{C}}}(\mathcal{U}_r)$ along the classifying map $c_{\mathcal{E}}\colon X\to\mathbf{BGL}_{r,\mathcal{C}}$ by Proposition~\ref{prop:vb-pb-properties}(1). For general $X \in \mathcal{H}_\tau(\mathcal{C})$, choose a $\tau$-hypercover $U_\bullet \to X$ with each $U_n$ a coproduct of objects in $\mathcal{C}_0 \subset \mathcal{X}_0$. By universality of colimits in $\mathcal{H}_\tau(\mathcal{C})$ and Proposition~\ref{prop:vb-pb-properties}(1), both $\mathbb{P}_X(\mathcal{E})$ and $\operatorname{Th}_X(\mathcal{E})$ are the geometric realizations of their pullbacks to $U_\bullet$. Since $\underline{\dMGL}_{\mathcal{C}}(- \to Y)$ sends colimits to limits and finite biproducts commute with limits in $\Sp$, both $\theta_{\mathcal{E}}$ and $\operatorname{th}_{\mathcal{E}}$ are equivalences for all $X \in \mathcal{H}_\tau(\mathcal{C})$.

At $r = 0$, $\mathbb{P}_X(0) = \varnothing$ and $z_0 \simeq \mathrm{id}_X$ give $\operatorname{th}_0 = \mathrm{id}$, and Proposition~\ref{prop:vb-pb-properties}(1) gives pullback compatibility of $\operatorname{th}_{\mathcal{E}}$. For $X, X' \in \mathcal{X}_0$ with $\mathcal{E} \in \Vect_r(X)$ and $\mathcal{F} \in \Vect_s(X')$, since $\ell$ is symmetric monoidal (Proposition~\ref{thm:localized-span-universal-property}), the zig-zag of equivalences
\[
m(\operatorname{Th}_{X\times X'}(\mathcal{E}\boxplus\mathcal{F}))
\xrightarrow{\ell(\varrho_{\mathcal{E}\boxplus\mathcal{F}})}
m(X\times X')
\leftarrow
m(X)\otimes m(X')
\xleftarrow{\ell(\varrho_{\mathcal{E}})\otimes\ell(\varrho_{\mathcal{F}})}
m(\operatorname{Th}_X(\mathcal{E}))\otimes m(\operatorname{Th}_{X'}(\mathcal{F}))
\]
identifies $\operatorname{th}_{\mathcal{E}\boxplus\mathcal{F}}$ with $\operatorname{th}_{\mathcal{E}}\otimes\operatorname{th}_{\mathcal{F}}$ compatibly with $(s_{\mathcal{E}\boxplus\mathcal{F}})_* = (s_{\mathcal{E}})_*\otimes(s_{\mathcal{F}})_*$ on $\mathbb{V}_X(\mathcal{E})\times\mathbb{V}_{X'}(\mathcal{F})\simeq\mathbb{V}_{X\times X'}(\mathcal{E}\boxplus\mathcal{F})$. By $\tau$-hypersheaf descent and diagonal pullback, $\operatorname{th}$ is multiplicative under direct sums on all of $\mathcal{H}_\tau(\mathcal{C})$.

Finally, for any line bundle $\mathcal{L}$ on $X \in \mathcal{C}_0$, since $\iota_{\mathcal{L}} \colon \mathbb{V}_X(\mathcal{L}) \hookrightarrow \mathbb{P}_X(\mathcal{L}^\vee \oplus \mathcal{O}_X)$ is a monomorphism with $z_{\mathcal{L}} = \iota_{\mathcal{L}} \circ s_{\mathcal{L}}$ (Proposition~\ref{prop:vb-pb-properties}(3)), we have $X \times_{\mathbb{P}_X(\mathcal{L}^\vee \oplus \mathcal{O}_X)} X \simeq X \times_{\mathbb{V}_X(\mathcal{L})} X$ and hence
\[
z_{\mathcal{L}}^*\bigl(\operatorname{th}_{\mathcal{L}}(1)\bigr)
= q_{\mathcal{C}}((z_{\mathcal{L}})_*) \circ q_{\mathcal{C}}(z_{\mathcal{L}}^*)
= q_{\mathcal{C}}((s_{\mathcal{L}})_*) \circ q_{\mathcal{C}}(s_{\mathcal{L}}^*)
= \widetilde{c}_1(\mathcal{L}).
\]
Together with~(1), (2), and admissible pushout excision (Definition~\ref{def:relation-dmgl}(1)), $\underline{\dMGL}_{\mathcal{C}}$ is an oriented cohomology theory on $(\mathcal{C}_0,\mathbf{S})$, and initiality follows from Proposition~\ref{thm:localized-span-universal-property}.
\end{proof}
\section{Almost algebraic cobordisms over prisms}
\label{sec:prismatic-master-theory}

We apply Section~\ref{sec:span_formalism} to construct a derived algebraic cobordism theory over prisms. The base affine formal schemes $U=\operatorname{Spf}(B)$ of prisms over $(\mathbb{A}_{\inf R}, I_R)$ carry the $\delta$-ring structure, whereas the fiber objects $X\to U$ forming admissible spans are quasi-compact, quasi-separated derived $(p,I_R)$-complete formal schemes over the underlying ring $B$.

\subsection{Animated prisms and boundary divisors}
\label{subsec:tilting-frame-and-prisms}
Fix a prime number $p$. We recall animated commutative rings, $\delta$-rings, and animated prisms.

\begin{definition}
\label{def:animated-rings-and-prisms}
\begin{enumerate}[\rm (1)]
\item Let $\mathrm{Poly}_{\mathbb{Z}}$ denote the category of finitely generated polynomial rings $\mathbb{Z}[x_1,\dots,x_n]$ ($n \ge 0$). Following~\cite[Section~5.1]{CS24} and~\cite[Section~5.5.8]{HT}, the $\infty$-category of \emph{animated commutative rings}, or \emph{animated rings}, is the sifted-colimit completion
\[
\mathrm{Ani}(\mathrm{Ring}) = \mathcal{P}_{\Sigma}(\mathrm{Poly}_{\mathbb{Z}}) \simeq \operatorname{Fun}^{\times}\bigl(\mathrm{Poly}_{\mathbb{Z}}^{\mathrm{op}}, \mathcal{S}\bigr),
\]
where the right-hand side is the full subcategory spanned by functors preserving finite products, equivalently obtained from the category of simplicial commutative rings by inverting weak homotopy equivalences. For $A \in \mathrm{Ani}(\mathrm{Ring})$, we write $\CAlg_A = \mathrm{Ani}(\mathrm{Ring})_{A/}$ for connective animated $A$-algebras, $\Mod_A^{\ge 0}$ for connective $A$-modules, and $L_{B/A} \in \Mod_B^{\ge 0}$ for the Andr\'e--Quillen relative cotangent complex constructed by animating K\"ahler differentials on polynomial $A$-algebras.

\item Following~\cite[Definitions~2.1 and~3.2]{BhattScholze2017Prisms}, a \emph{$\delta$-ring} is a pair $(A,\delta)$ consisting of a commutative ring $A$ and a set-theoretic map $\delta \colon A \to A$ satisfying $\delta(0)=\delta(1)=0$,
\[
\delta(xy) = x^p\delta(y) + y^p\delta(x) + p\delta(x)\delta(y), \qquad
\delta(x+y) = \delta(x) + \delta(y) - \sum_{i=1}^{p-1} \frac{1}{p}\binom{p}{i} x^i y^{p-i}
\]
for all $x,y \in A$; its associated ring endomorphism $\varphi_A(x) = x^p + p\delta(x)$ lifts the Frobenius endomorphism of $A/(p)$. A \emph{prism} is a pair $(A,I)$ consisting of a $\delta$-ring $A$ and an invertible ideal $I \subset A$ such that $A$ is derived $(p,I)$-complete and $p \in I + \varphi_A(I)A$. A prism $(A,I)$ is \emph{bounded} if $A/I$ has bounded $p^\infty$-torsion, \emph{orientable} if $I$ is principal, and \emph{perfect} if $\varphi_A \colon A \to A$ is an isomorphism.

\item Let $\mathrm{Poly}_{\mathbb{Z}}^{\delta}$ denote the category of finitely generated free $\delta$-rings $\mathbb{Z}\{x_1,\dots,x_n\}_{\delta}$. The $\infty$-category of \emph{animated $\delta$-rings} is $\mathrm{Ani}(\mathrm{Ring}_{\delta}) = \mathcal{P}_{\Sigma}(\mathrm{Poly}_{\mathbb{Z}}^{\delta})$~\cite[Definition~A.11]{BhattLurie2022Prismatization}, equipped with the conservative forgetful functor $\mathrm{Ani}(\mathrm{Ring}_{\delta}) \to \mathrm{Ani}(\mathrm{Ring})$ preserving sifted colimits and animated Frobenius endomorphism $\varphi_B \colon B \to B$. Following~\cite[Construction~2.1 and Definition~2.4]{BhattLurie2022Prismatization}, an \emph{animated prism} $(B, I \xrightarrow{\alpha} B)$ consists of an animated $\delta$-ring $B$ and a generalized Cartier divisor $\alpha \colon I \to B$, with invertible $I \in \Mod_B$ and derived Hodge--Tate quotient $B/I = \operatorname{cofib}(\alpha) \in \CAlg_B$, such that $B$ is derived $(p,I)$-complete and $p \in \operatorname{im}(\pi_0(I) \to \pi_0(B)) + \varphi_{\pi_0(B)}\bigl(\operatorname{im}(\pi_0(I) \to \pi_0(B))\bigr)\pi_0(B)$. An animated prism is \emph{perfect} if $\varphi_B \colon B \to B$ is an equivalence; any perfect animated prism is discrete and hence a classical perfect prism~\cite[Corollary~2.14 and Lemma~2.15]{BhattLurie2022Prismatization}.
\end{enumerate}
\end{definition}

\begin{remark}
\label{rem:animated-vs-Einfty-and-rigidity}
Unlike connective $\mathbb{E}_\infty$-rings in $\CAlg^{\mathrm{cn}}(\mathrm{Sp})$, animated rings retain strict polynomial operations in mixed and positive characteristic~\cite[Section~5.1]{CS24}; hence polynomial algebras $\mathbb{Z}[x_1,\dots,x_n]$ remain discrete and flat, Frobenius lifts and $\delta$-structures extend by animation, and quasi-smooth morphisms are characterized by perfect Andr\'e--Quillen cotangent complexes of Tor-amplitude in $[-1,0]$. Over any fixed prism $(A_0,I_0)$, rigidity of generalized Cartier divisors~\cite[Corollary~2.10]{BhattLurie2022Prismatization} identifies the $\infty$-category of animated prisms over $(A_0,I_0)$ with the $\infty$-category of derived $(p,I_0)$-complete animated $\delta$-$A_0$-algebras $B$, with induced prism structure $\bigl(B, I_0 \widehat{\otimes}_{A_0} B \to B\bigr)$ and Hodge--Tate specialization $B \otimes_{A_0} (A_0/I_0)$.
\end{remark}

\begin{definition}
\label{def:perfectoid-tilting-frame}
Consider the formal scheme $\mathfrak{F}_R=\operatorname{Spf}(\mathbb A_{\inf R})$ with its $(p,I_R)=(p,\xi)$-adic topology and its two Cartier divisors:
\begin{enumerate}[(1)]
\item the \emph{arithmetic locus} $V(I_R)=\operatorname{Spf}(R^\circ)$, carrying the induced $(p)$-adic topology; and
\item the \emph{tilted locus} $V(p)=\operatorname{Spf}(R^{\flat\circ})$, carrying the induced $(\bar\xi)$-adic topology.
\end{enumerate}
\end{definition}

\begin{definition}
\label{def:category-animated-prisms}
We consider three categories of affine formal schemes of prisms over $\mathfrak{F}_R$:
\begin{enumerate}[(1)]
\item Let $\mathcal T_R$ be the category of prisms $(B,I_R B)$ over $(\mathbb A_{\inf R},I_R)$~\cite[Definition~3.2]{BhattScholze2017Prisms}, regarded contravariantly via their affine formal schemes $U=\operatorname{Spf}(B)$. Here $B$ is a derived $(p,\xi)$-complete $\delta$-$\mathbb A_{\inf R}$-algebra such that $\xi$ is a non-zero-divisor in $B$. Morphisms $\operatorname{Spf}(B')\to\operatorname{Spf}(B)$ in $\mathcal T_R$ are $\delta$-$\mathbb A_{\inf R}$-algebra homomorphisms $B\to B'$.
\item Let $\mathcal T_R^{\mathrm{perf}}\subset\mathcal T_R$ be the full subcategory spanned by affine formal schemes $U=\operatorname{Spf}(B)$ of perfect prisms over $(\mathbb A_{\inf R},I_R)$, for which the Frobenius lift $\varphi_B\colon B\to B$ is an isomorphism. For any $U=\operatorname{Spf}(B)\in\mathcal T_R$, we set $B_\infty = (\varinjlim_{\varphi_B} B)^\wedge_{(p,I_R)}$ and $U^{\mathrm{perf}} = \operatorname{Spf}(B_\infty)$, equipped with the canonical morphism $\epsilon_U\colon U^{\mathrm{perf}}\to U$ induced by $B\to B_\infty$.
\item Let $\mathcal T_R^{\mathrm{an}}$ be the $\infty$-category of strict bounded animated prisms $(B,J_B)$ over $(\mathbb A_{\inf R},I_R)$~\cite[Definition~2.4, Remark~2.9, and Corollary~2.10]{BhattLurie2022Prismatization}, where $B$ is a connective, derived $(p,\xi)$-complete animated $\delta$-$\mathbb A_{\inf R}$-algebra, $(J_B\to B)\simeq (I_R\otimes_{\mathbb A_{\inf R}}B\to B)$ with $B/J_B\simeq B\otimes_{\mathbb A_{\inf R},\theta}R^\circ$, and $\pi_0(B/J_B)$ has bounded $p^\infty$-torsion.
\end{enumerate}
In $\mathcal T_R^{\mathrm{an}}$, a bounded homotopy fiber product of affine formal schemes over $\operatorname{Spf}(C)$ is given by $\operatorname{Spf}(B_1)\times_{\operatorname{Spf}(C)}^{h}\operatorname{Spf}(B_2) \simeq \operatorname{Spf}\bigl((B_1\otimes_CB_2)^\wedge_{(p,I_R)}\bigr)$ with the induced animated $\delta$-structure~\cite[Proposition~A.20]{BhattLurie2022Prismatization}. Formal gluing of affine formal schemes is performed in the associated hypersheaf categories.
\end{definition}

\begin{proposition}
\label{prop:universal-prism-perfection}
For any $U=\operatorname{Spf}(B)\in\mathcal T_R$, the ring $B_\infty$ is discrete and $U^{\mathrm{perf}}$ belongs to $\mathcal T_R^{\mathrm{perf}}$. The assignment $U\mapsto U^{\mathrm{perf}}$ defines a right adjoint $P\colon\mathcal T_R\to\mathcal T_R^{\mathrm{perf}}$ to the full inclusion $\mathcal T_R^{\mathrm{perf}}\to\mathcal T_R$: for any $V\in\mathcal T_R^{\mathrm{perf}}$, composition with $\epsilon_U$ induces an equivalence
\[
\Map_{\mathcal T_R}(V,U^{\mathrm{perf}})\longrightarrow\Map_{\mathcal T_R}(V,U).
\]
In particular, $\epsilon_V\colon V^{\mathrm{perf}}\to V$ is an equivalence for any $V\in\mathcal T_R^{\mathrm{perf}}$.
\end{proposition}
\begin{proof}
For $U=\operatorname{Spf}(B)\in\mathcal T_R$,~\cite[Lemma~3.9]{BhattScholze2017Prisms} shows that the derived $(p,I_R)$-completion $B_\infty=(\varinjlim_{\varphi_B}B)^\wedge_{(p,I_R)}$ is concentrated in degree $0$ and carries a unique $\delta$-$\mathbb A_{\inf R}$-algebra structure making $(B_\infty,I_R B_\infty)$ a perfect prism over $(\mathbb A_{\inf R},I_R)$ whose reduction $B_\infty/I_R B_\infty$ is an integral perfectoid $R^\circ$-algebra~\cite[Theorem~3.10]{BhattScholze2017Prisms}. Thus $U^{\mathrm{perf}}\in\mathcal T_R^{\mathrm{perf}}$, and $\iota_B\colon B\to B_\infty$ defines $\epsilon_U=\operatorname{Spf}(\iota_B)\colon U^{\mathrm{perf}}\to U$.

For any $V=\operatorname{Spf}(C)\in\mathcal T_R^{\mathrm{perf}}$, $C$ is derived $(p,I_R)$-complete and the Frobenius lift $\varphi_C\colon C\to C$ is an isomorphism of $\delta$-$\mathbb A_{\inf R}$-algebras. By the universal properties of derived $(p,I_R)$-completion and of $\varinjlim_{\varphi_B}B$, we obtain equivalences
\[
\Map_{\mathcal T_R}(V,U^{\mathrm{perf}})
\longrightarrow\Map_{\delta\text{-}\CAlg_{\mathbb A_{\inf R}}}(B_\infty,C)
\longrightarrow\varprojlim_{\varphi_B^*}\Map_{\delta\text{-}\CAlg_{\mathbb A_{\inf R}}}(B,C)
\longrightarrow\Map_{\mathcal T_R}(V,U),
\]
where the transition maps $\varphi_B^*(f)=\varphi_C\circ f$ are equivalences with inverse $\varphi_C^{-1}\circ (-)$. Thus $P$ is right adjoint to the full inclusion $\mathcal T_R^{\mathrm{perf}}\to\mathcal T_R$, implying that $\epsilon_V$ is an equivalence for any $V\in\mathcal T_R^{\mathrm{perf}}$.
\end{proof}

\subsection{Formal geometry parametrized by prisms}
\label{subsec:prismatic-geometry-and-master-localizations}

\begin{definition}
\label{def:category-prismatic-formal-schemes}
Let $U=\operatorname{Spf}(B)\in\mathcal T_R$ be the affine formal scheme of a prism over $(\mathbb{A}_{\inf R}, I_R)$, and let $\tau$ be a subcanonical Grothendieck topology on $\mathrm{dForm}(U)$ containing the Zariski topology and coarser than the completely flat topology. Given any connective derived $(p,I_R)$-complete animated $B$-algebra $A$ and variables $T_1,\dots,T_r$, we write $A\langle T_1,\dots,T_r\rangle = (A[T_1,\dots,T_r])^\wedge_{(p,I_R)}$.
\begin{enumerate}[(1)]
\item Let $\mathrm{dForm}(U)$ be the $\infty$-category of quasi-compact and quasi-separated derived $(p,I_R)$-complete formal schemes over $U$, locally modeled by $\operatorname{Spf}(A)$ for connective, derived $(p,I_R)$-complete animated $B$-algebras $A$ over the underlying ring $B$~\cite{SAG}, closed under homotopy pullbacks, universal finite coproducts, and additively admissible pushouts. The derived $(p,I_R)$-completed affine line $\widehat{\mathbb{A}}^1_U = \operatorname{Spf}(B\langle t\rangle)$ and multiplicative group $\widehat{\mathbb{G}}_{m,U} = \operatorname{Spf}(B\langle t, t^{-1}\rangle)$ endow $(\mathrm{dForm}(U), \tau)$ with a universal line datum $\mathbf{BG}_m = a_\tau^{\mathrm{hyp}}([* / \widehat{\mathbb{G}}_{m,U}])$, $\mathcal{L} = a_\tau^{\mathrm{hyp}}([\widehat{\mathbb{A}}^1_U / \widehat{\mathbb{G}}_{m,U}])$, and $\mathcal{L}^\circ = a_\tau^{\mathrm{hyp}}([\widehat{\mathbb{G}}_{m,U} / \widehat{\mathbb{G}}_{m,U}]) \simeq *$.
\item Any morphism $u\colon U'=\operatorname{Spf}(B')\to U=\operatorname{Spf}(B)$ in $\mathcal T_R$ induces the completed derived base-change functor $u^*\colon\mathrm{dForm}(U)\to\mathrm{dForm}(U')$ sending $\operatorname{Spf}(A)$ to $\operatorname{Spf}\bigl((A\otimes_BB')^\wedge_{(p,I_R)}\bigr)$. We denote by $\mathcal C_{\mathcal T_R}\to\mathcal T_R$ the resulting Cartesian fibration of pairs $(U,X)$ with $U\in\mathcal T_R$ and $X\in\mathrm{dForm}(U)$.
\item For any $U\in\mathcal T_R$, let $\mathbf S_{\mathrm{Prism},U}$ be the class of morphisms $f\colon X\to Y$ in $\mathrm{dForm}(U)$ that are quasi-compact, quasi-separated, topologically locally of finite presentation, and derived quasi-smooth, meaning that their completed relative cotangent complex $\widehat L_{X/Y}$ is perfect of Tor-amplitude $[-1,0]$ in cohomological degrees. Let $\mathbf S_{\mathrm{Prism}}$ be the class of vertical morphisms in $\mathcal C_{\mathcal T_R}$ whose fiber map belongs to $\mathbf S_{\mathrm{Prism},U}$.
\end{enumerate}
\end{definition}

\begin{proposition}
\label{prop:prismatic-axioms-verification}
For any $U = \operatorname{Spf}(B) \in \mathcal{T}_R$ and any subcanonical Grothendieck topology $\tau$ on $\mathrm{dForm}(U)$ between the Zariski topology and the completely flat topology, $\mathrm{dForm}(U)$ carries a universal line datum $(\mathbf{BG}_m,\mathcal{L},\mathcal{L}_{\mathrm{Sp}},\mathcal{L}^\circ)$ for which $\tau$ is locally trivializing, and $\mathbf{S}_{\mathrm{Prism},U}$ is an admissible class. For any finite locally free module $\mathcal{E}$ of rank $r$ on $X \in \mathrm{dForm}(U)$, the spaces $\mathbb{V}_X(\mathcal{E})$, $\mathbb{V}_X(\mathcal{E}^\vee)$, $\mathbb{V}_X(\mathcal{E}^\vee)^\times$, and $\mathbb{P}_X(\mathcal{E})$ are represented in $\mathrm{dForm}(U)$, and both $s_{\mathcal{E}} \colon X \hookrightarrow \mathbb{V}_X(\mathcal{E})$ and $z_{\mathcal{E}} \colon X \hookrightarrow \mathbb{P}_X(\mathcal{E}^\vee \oplus \mathcal{O}_X)$ belong to $\mathbf{S}_{\mathrm{Prism},U}$. For any morphism $u \colon U' \to U$ in $\mathcal{T}_R$, the base-change functor $u^* \colon \mathrm{dForm}(U) \to \mathrm{dForm}(U')$ preserves homotopy pullbacks, finite coproducts, $\mathbf{S}_{\mathrm{Prism},U}$, $\mathbb{V}_X(\mathcal{E})$, and $\mathbb{P}_X(\mathcal{E})$.
\end{proposition}
\begin{proof}
Write $\mathcal{H} = \mathcal{H}_\tau(\mathrm{dForm}(U))$, and let $\widehat{\mathbb{A}}^1_{U,\mathrm{Sp}}$ be the connective $\mathbb{E}_\infty$-ring $\tau$-hypersheaf $X \mapsto \tau_{\ge 0} R\Gamma(X, \mathcal{O}_X)$, with underlying space sheaf $\widehat{\mathbb{A}}^1_U$ and unit group $\widehat{\mathbb{G}}_{m,U} \simeq \mathrm{GL}_1(\widehat{\mathbb{A}}^1_{U,\mathrm{Sp}})$. Descent along $1 \colon * \twoheadrightarrow \mathbf{BG}_m = a_\tau^{\mathrm{hyp}}([*/\widehat{\mathbb{G}}_{m,U}])$ gives a symmetric monoidal equivalence $\Sp_{\ge 0}(\mathcal{H}_{/\mathbf{BG}_m}) \xrightarrow{\sim} \Mod_{\Sigma^\infty_+\widehat{\mathbb{G}}_{m,U}}(\Sp_{\ge 0}(\mathcal{H}))$. Under this equivalence, the $\mathrm{GL}_1$-counit $\Sigma^\infty_+\widehat{\mathbb{G}}_{m,U} \to \widehat{\mathbb{A}}^1_{U,\mathrm{Sp}}$ defines $\mathcal{L}_{\mathrm{Sp}} \in \CAlg(\Sp_{\ge 0}(\mathrm{Pre}(\mathrm{dForm}(U))_{/\mathbf{BG}_m}))$ lifting $\mathcal{L}$ with $1^*_{\mathrm{Sp}}(\mathcal{L}_{\mathrm{Sp}}) \simeq \widehat{\mathbb{A}}^1_{U,\mathrm{Sp}}$; together with $\mathcal{L}^\circ \simeq *$, this verifies Definition~\ref{def:universal-line-datum}.

Since $\tau$ is subcanonical and refines the Zariski topology, $\tau$ is locally trivializing (Definition~\ref{def:geometric-topology}), and derived $(p,I_R)$-complete flat descent~\cite[Section~4]{BMS19} identifies $\Vect_r(X)$ with rank-$r$ locally free $\mathcal{O}_X$-modules. Over $X = \operatorname{Spf}(A)$, derived $(p,I_R)$-completeness implies $(p,I_R)\pi_0(A) \subseteq \operatorname{Jac}(\pi_0(A))$~\cite[Tags~091V and~0G1U]{stacks-project}, so rank-$r$ projective modules lift Zariski locally from $\pi_0(A)/(p,I_R)$ to $A$, and a row in $\pi_0(S)^r$ is unimodular if and only if its entries generate the unit ideal. Standard Zariski gluing of completed symmetric algebras and homogeneous localizations $A\langle T_1/T_i, \dots, T_r/T_i\rangle$ ($1 \le i \le r$) therefore represents $\mathbb{V}_X(\mathcal{E})$, $\mathbb{V}_X(\mathcal{E}^\vee)$, $\mathbb{V}_X(\mathcal{E}^\vee)^\times$, and $\mathbb{P}_X(\mathcal{E})$ in $\mathrm{dForm}(U)$. Since derived $(p,I_R)$-completion is symmetric monoidal on connective modules~\cite[Tag~091V]{stacks-project}, $u^*$ preserves homotopy pullbacks, finite coproducts, $\mathbb{V}_X(\mathcal{E})$, and $\mathbb{P}_X(\mathcal{E})$.

Finally, equivalences, $\varnothing \to X$, and $\nabla \colon X \amalg X \to X$ have vanishing completed relative cotangent complexes, while perfection and Tor-amplitude $[-1,0]$ of $\widehat{L}_{X/Y}$ are stable under completed base change and composition by cotangent transitivity. For any rank-$r$ vector bundle $\mathcal{E}$ on $X$, both $s_{\mathcal{E}}$ and $z_{\mathcal{E}} = \iota_{\mathcal{E}} \circ s_{\mathcal{E}}$ have completed relative cotangent complex $\mathcal{E}^\vee[1]$ of Tor-amplitude $[-1,-1] \subset [-1,0]$ and hence lie in $\mathbf{S}_{\mathrm{Prism},U}$. Derived $(p,I_R)$-complete flat descent~\cite[Section~4]{BMS19} makes $\mathbf{S}_{\mathrm{Prism},U}$ $\tau$-local on the target, completing the verification of Definition~\ref{def:admissible-class}.
\end{proof}

\begin{definition}
\label{def:master-prismatic-theory}
Applying Definition~\ref{def:span-category-correspondences} fiberwise to $(\mathrm{dForm}(U),\mathbf S_{\mathrm{Prism},U})$ for any $U\in\mathcal T_R$ defines the additive correspondence spectra $\Corr_{\mathbf S_{\mathrm{Prism},U}}(X,Y)^{\mathrm{gp}}\in\Sp_{\geq0}$ for $X,Y\in\mathrm{dForm}(U)$, functorial in $U \in \mathcal{T}_R$ by Proposition~\ref{prop:prismatic-axioms-verification}. Let $\tau_{\mathrm{geom}}$ be the \emph{quasisyntomic--Zariski topology} on $\mathcal T_R$, generated by completely quasisyntomic covers, quasiregular semiperfectoid (QRSP) refinements~\cite[Section~4]{BMS19}, and Zariski open coverings. Enriched localization over the comparison topology $\tau_{\mathrm{cmp}}$ in Definitions~\ref{def:perfect-presentation-topology} and~\ref{def:hypersheafified-dmgla} defines the prismatic derived algebraic cobordism theory $\underline{\dMGL}_{\mathrm{Prism}}^a$, which satisfies the orientation properties of Proposition~\ref{thm:full-orientation-from-line-thom}.
\end{definition}

\section{Boundary fiber specializations and deformation geometry}
\label{sec:fiber-specialization-and-deformation}

We specialize the prismatic theory over $\mathcal{T}_R$ along $V(I_R) = \operatorname{Spf}(R^\circ)$ and $V(p) = \operatorname{Spf}(R^{\flat\circ})$. Since an integral perfectoid ring contains a compatible system of $p$-power roots of a pseudo-uniformizer $\omega$ (Subsection~\ref{subsec:almost-formal-base}), the ideal $\mathfrak{m} = \bigcup_{r\ge 0} \omega^{1/p^r} R^\circ$ satisfies $\mathfrak{m}^2 = \mathfrak{m}$, so integral perfectoid tilting and deformation theory are formulated in almost mathematics~\cite{Scholze2012Perfectoid,GR,BMS18}.

\subsection{Almost formal base categories and almost cobordism}
\label{subsec:almost-formal-base}

Since $R^{\flat\circ}$ is a perfect $\mathbb F_p$-algebra, $\bar\xi=\rho_\flat(\xi)\in R^{\flat\circ}$ admits unique compatible $p$-power roots $(\omega^\flat)^{1/p^r}=\bar\xi^{1/p^{r+1}}\in R^{\flat\circ}$ ($r\ge 0$) with $(\omega^\flat)^p=\bar\xi$. Using the Teichm\"{u}ller lift $[-]\colon R^{\flat\circ}\to W(R^{\flat\circ})=\mathbb A_{\inf R}$, set $\widetilde\omega=[\omega^\flat]\in \mathbb A_{\inf R}$ and $\omega^{1/p^r}=\theta([(\omega^\flat)^{1/p^r}])\in R^\circ$ ($r\ge 0$), with $\omega=\omega^1=\theta(\widetilde\omega)$.

\begin{lemma}
\label{lem:tilting-frame-koszul}
For $\omega$, $\omega^\flat$, and $\widetilde\omega$ as above:
\begin{enumerate}[\rm (1)]
\item There is a unique element $y_0\in\mathbb A_{\inf R}^\times$ satisfying $\xi=\widetilde\omega^p+p y_0$, with $\omega^p=-\theta(y_0)p$ in $R^\circ$ and $(\omega^\flat)^p=\bar\xi$ in $R^{\flat\circ}$. The elements $\omega\in R^\circ$, $\omega^\flat\in R^{\flat\circ}$, and $\widetilde\omega\in\mathbb A_{\inf R}$ are non-zero-divisors. The $(\omega)$-adic and $(p)$-adic topologies on $R^\circ$ coincide, and the $(\omega^\flat)$-adic and $(\bar\xi)$-adic topologies on $R^{\flat\circ}$ coincide.
\item For each integer $1\le k\le p$, the ideals $(\xi,\widetilde\omega^k)$ and $(p,\widetilde\omega^k)$ in $\mathbb A_{\inf R}$ are equal and are generated by regular sequences $(\xi,\widetilde\omega^k)$ and $(p,\widetilde\omega^k)$ of length $2$, inducing ring isomorphisms
\[
R_k = R^\circ/(\omega^k) \;\longleftarrow\; = \mathbb A_{\inf R}/(p,\widetilde\omega^k) \;\longrightarrow\; R^{\flat\circ}/((\omega^\flat)^k) = R_k^\flat
\]
whose composite defines $\iota_k\colon R_k\to R_k^\flat$ compatible with $\theta$ and $\rho_\flat$; in particular, $k=1$ gives $\iota_1\colon R_1\to R_1^\flat$ and $k=p$ gives the mod-$p$ isomorphism $\iota_p\colon R_p=R^\circ/(p)\to R^{\flat\circ}/(\bar\xi)=R_p^\flat$.
\item For each $1\le k\le p$, the quotients $R^\circ$, $R^{\flat\circ}$, and $R_k$ belong to $\operatorname{Perf}(\mathbb A_{\inf R})$ and are dualizable in $D(\mathbb A_{\inf R})$, while $R_k\in\operatorname{Perf}(R^\circ)$ and $R_k^\flat\in\operatorname{Perf}(R^{\flat\circ})$ are dualizable over $R^\circ$ and $R^{\flat\circ}$. For any $M\in D(\mathbb A_{\inf R})$, there are equivalences
\[
\bigl(M\otimes_{\mathbb A_{\inf R},\theta}R^\circ\bigr)\otimes_{R^\circ}R_k
\;\longleftarrow\;
M\otimes_{\mathbb A_{\inf R}}R_k
\;\longrightarrow\;
\bigl(M\otimes_{\mathbb A_{\inf R},\rho_\flat}R^{\flat\circ}\bigr)\otimes_{R^{\flat\circ}}R_k^\flat.
\]
\end{enumerate}
\end{lemma}
\begin{proof}
Since $\rho_\flat(\xi)=\bar\xi=\rho_\flat(\widetilde\omega^p)$ and $p$ is a non-zero-divisor in $W(R^{\flat\circ})$, there is a unique $y_0\in\mathbb A_{\inf R}$ with $\xi=\widetilde\omega^p+p y_0$. Using $\delta(\widetilde\omega^p)=0$, we obtain $\delta(\xi)\equiv (1-p^{p-1})y_0^p\pmod{(p,\widetilde\omega)}$. Since $\delta(\xi)\in\mathbb A_{\inf R}^\times$, $1-p^{p-1}\in\mathbb Z_p^\times$, and $(p,\widetilde\omega)^p\subseteq (p,\xi)\subseteq\operatorname{Jac}(\mathbb A_{\inf R})$, it follows that $y_0\in\mathbb A_{\inf R}^\times$.

Applying $\theta$ gives $\omega^p=-\theta(y_0)p$ with $-\theta(y_0)\in (R^\circ)^\times$, so $\omega$ is a non-zero-divisor in $R^\circ$. If $\widetilde\omega w=0$ in $\mathbb A_{\inf R}$, then $\xi w=p y_0 w$ forces $\theta(w)=0$, so $w=\xi w_1$ with $\widetilde\omega w_1=0$ as $\xi$ is a non-zero-divisor~\cite[Lemma~3.8]{BhattScholze2017Prisms}; $(p,\xi)$-adic separatedness gives $w\in\bigcap_{n\ge 1}\xi^n\mathbb A_{\inf R}=0$. Injectivity of $[-]\colon R^{\flat\circ}\to\mathbb A_{\inf R}$ then implies that $\omega^\flat$ is a non-zero-divisor in $R^{\flat\circ}$, proving~(1).

For~(2), since $\xi=\widetilde\omega^{p-k}\widetilde\omega^k+p y_0$ with $y_0\in\mathbb A_{\inf R}^\times$, we have $(\xi,\widetilde\omega^k)=(p,\widetilde\omega^k)$ for $1\le k\le p$, and $(\omega^p)=(p)$, $((\omega^\flat)^p)=(\bar\xi)$ at $k=p$. Since $\xi$ (resp.\ $p$) is a non-zero-divisor in $\mathbb A_{\inf R}$ and $\omega^k$ (resp.\ $(\omega^\flat)^k$) is a non-zero-divisor in $R^\circ$ (resp.\ $R^{\flat\circ}$), both $(\xi,\widetilde\omega^k)$ and $(p,\widetilde\omega^k)$ are regular sequences of length $2$ with common quotient $\iota_k\colon R_k\xrightarrow{\sim}R_k^\flat$. Part~(3) follows from the resulting two-term free resolutions of $R^\circ$, $R^{\flat\circ}$, $R_k$, and $R_k^\flat$, together with transitivity of derived tensor products.
\end{proof}

For each $n \geq 1$, the adic quotients $R_n = R^\circ/(\omega^n)$ and $R_n^\flat = R^{\flat\circ}/((\omega^\flat)^n)$ form adic towers with inverse limits $R^\circ$ and $R^{\flat\circ}$ and identified stages $R_k \simeq R_k^\flat$ for $1\le k\le p$.

The ideals $\mathfrak{m} = \varinjlim_{r\geq 0} \omega^{1/p^r} R^\circ \subset R^\circ$ and $\mathfrak{m}^\flat = \varinjlim_{r\geq 0} (\omega^\flat)^{1/p^r} R^{\flat\circ} \subset R^{\flat\circ}$ are filtered colimits of free rank-one modules and hence flat over $R^\circ$ and $R^{\flat\circ}$~\cite[Tag~05UT]{stacks-project}. The relations $(\omega^{1/p^{r+1}})^p = \omega^{1/p^r}$ imply $\mathfrak{m}^2 = \mathfrak{m}$ and $(\mathfrak{m}^\flat)^2 = \mathfrak{m}^\flat$. The tensor squares $\tm = \mathfrak{m} \otimes_{R^\circ} \mathfrak{m}$ and $\tm^\flat = \mathfrak{m}^\flat \otimes_{R^{\flat\circ}} \mathfrak{m}^\flat$ are likewise flat. An $R^\circ$-module $M$ is \emph{almost zero} if $\mathfrak{m}M = 0$, and \emph{almostification} $M \mapsto M^a$ is the exact symmetric monoidal localization quotienting by almost zero modules~\cite[Section~2.2]{GR}.

At each level $n \ge 1$, we distinguish between the flat modules $\mathfrak{m}_n = \mathfrak{m} \otimes_{R^\circ} R_n$ and $\tm_n = \tm \otimes_{R^\circ} R_n$ obtained by base change, and the idempotent image ideal $\mathfrak{m}R_n = \operatorname{im}(\mathfrak{m}_n \to R_n)$, with analogous notation $\mathfrak{m}_n^\flat, \tm_n^\flat$, and $\mathfrak{m}^\flat R_n^\flat$ over $R_n^\flat$. Under each isomorphism $\iota_k\colon R_k\to R_k^\flat$ ($1\le k\le p$) of Lemma~\ref{lem:tilting-frame-koszul}(2), the residues of $\omega^{1/p^r}=\theta([(\omega^\flat)^{1/p^r}])$ and $(\omega^\flat)^{1/p^r}=\rho_\flat([(\omega^\flat)^{1/p^r}])$ coincide, identifying $\mathfrak{m}R_k\xrightarrow{\sim}\mathfrak{m}^\flat R_k^\flat$, $\mathfrak{m}_k\xrightarrow{\sim}\mathfrak{m}_k^\flat$, and $\tm_k\xrightarrow{\sim}\tm_k^\flat$, and hence inducing a symmetric monoidal equivalence $\Mod_{R_k^a}\to\Mod_{(R_k^\flat)^a}$. The surjection $\mathfrak{m}_n \to \mathfrak{m}R_n$ is an almost isomorphism with kernel annihilated by $\mathfrak{m}$, while $\mathfrak{m}_n$ is flat over $R_n$~\cite[Remark~2.1.4(ii)]{GR}. Almostification over $R_n$ and $R_n^\flat$ is taken with respect to $\mathfrak{m}R_n$ and $\mathfrak{m}^\flat R_n^\flat$.

The localization $(-)^a \colon \Mod_{R^\circ} \to \Mod_{(R^\circ)^a}$ has right adjoint $M_* = \operatorname{Hom}_{(R^\circ)^a}((R^\circ)^a, M)$ and exact left adjoint $M_! = \tm \otimes_{R^\circ} M_*$~\cite[Proposition~2.2.21]{GR}, inducing on derived $\infty$-categories an adjoint triple
\[
(-)_! \;\dashv\; (-)^a \;\dashv\; R(-)_*.
\]
Since $(-)^a \colon D(R^\circ) \to D((R^\circ)^a)$ is both a left and a right adjoint, it preserves all small homotopy limits and all small colimits, and $R_n^a \in \operatorname{Perf}((R^\circ)^a)$ and $(R_n^\flat)^a \in \operatorname{Perf}((R^{\flat\circ})^a)$ remain dualizable.

For a non-unital animated $R^\circ$-algebra $I$, write $R^\circ \oplus I$ for its direct sum unitalization. For any almost module $M$ over an almost $(R^\circ)^a$-algebra $A$, the counit $(M_!)^a \to M$ is an equivalence. Let $A$ be a discrete or animated $(R^\circ)^a$-algebra. Both $\tau_{\ge 0}R(A)_*$ and Gabber--Ramero's almost unitalization
\[
A_{!!} = R^\circ \amalg_{R^\circ \oplus \tm} (R^\circ \oplus A_!)
\]
are unital animated $R^\circ$-algebras satisfying $(\tau_{\ge 0}R(A)_*)^a \xrightarrow{\sim} A$ and $(A_{!!})^a \xrightarrow{\sim} A$~\cite[Propositions~2.2.14 and~2.2.25]{GR},~\cite[Sections~2--3]{kato2023mathematics}. Since $\tm$ and $\tm_n$ are flat, $M \mapsto M_!$ is exact and embeds the derived almost category fully faithfully into ordinary derived modules, computing all derived tensor products, completed cotangent complexes, and obstruction groups via firm representatives.

\begin{definition}
\label{def:almost-formal-base-categories}
Let $\mathcal{C}_{R^\circ}^a$ be the $\infty$-category of quasi-compact, quasi-separated derived almost formal schemes over $\operatorname{Spf}^a((R^\circ)^a)$, where $\operatorname{Spf}^a(-)$ (resp.\ $\Spec^a(-)$) denotes formal (resp.\ affine) spectra in the opposite $\infty$-category of derived $\omega$-complete (resp.\ general) connective animated almost algebras:
\begin{enumerate}[(1)]
\item Objects of $\mathcal{C}_{R^\circ}^a$ are formed by finite Zariski gluing of affine spectra $\operatorname{Spf}^a(A)$ for connective animated $(R^\circ)^a$-algebras $A$ that are derived $\omega$-complete, meaning $A \xrightarrow{\sim} \varprojlim_{n\geq 1} \bigl(A \otimes_{(R^\circ)^a} R_n^a\bigr)$, with level-$n$ reduction $\Spec^a\bigl(A \otimes_{(R^\circ)^a} R_n^a\bigr)$, fiber products given by derived $\omega$-completed tensor products, universal finite coproducts, and additively admissible pushouts.

Writing $R^\circ\langle t_1,\dots,t_r\rangle = (R^\circ[t_1,\dots,t_r])^\wedge_\omega$ and $\tau_{\mathrm{Zar}}$ for the Zariski topology, the universal line datum on $(\mathcal{C}_{R^\circ}^a, \tau_{\mathrm{Zar}})$ is induced as in Proposition~\ref{prop:prismatic-axioms-verification} by $\widehat{\mathbb{A}}^{1,a}_{R^\circ} = \operatorname{Spf}^a(R^\circ\langle t\rangle^a)$ and $\widehat{\mathbb{G}}_{m,R^\circ}^a = \operatorname{Spf}^a(R^\circ\langle t, t^{-1}\rangle^a)$. For any $\operatorname{Spf}^a(B) \in \mathcal{C}_{R^\circ}^a$, the equivalence $\widehat{\mathbb{A}}^{1,a}_{R^\circ}(\operatorname{Spf}^a(B)) \to \Omega^\infty R(B)_*$ identifies unimodular rows in $B$ with unimodular rows in $\pi_0 R(B)_*$, and $\tau_{\ge 0}R(B)_* \simeq 0$ forces $B \simeq (\tau_{\ge 0}R(B)_*)^a \simeq 0$, so $\operatorname{Spf}^a(B) = \varnothing$ is covered by the empty Zariski sieve. Therefore, $\tau_{\mathrm{Zar}}$ is locally trivializing, and for any finite locally free module $\mathcal{E}$ on $X \in \mathcal{C}_{R^\circ}^a$, derived $\omega$-completed almost symmetric algebras and standard Zariski open affine coverings represent $\mathbb{V}_X(\mathcal{E})$, $\mathbb{V}_X(\mathcal{E}^\vee)^\times$, and $\mathbb{P}_X(\mathcal{E})$ in $\mathcal{C}_{R^\circ}^a$.
\item A morphism $f \colon X \to Y$ in $\mathcal{C}_{R^\circ}^a$ is \emph{quasi-smooth} if it is topologically locally of almost finite presentation and its completed relative almost cotangent complex $\widehat{L}_{X/Y}^a$ is perfect of Tor-amplitude $[-1,0]$ in cohomological degrees. The class $\mathbf{S}_\theta$ of all quasi-compact, quasi-separated quasi-smooth morphisms in $\mathcal{C}_{R^\circ}^a$ is an admissible class, verified identically to Proposition~\ref{prop:prismatic-axioms-verification}.
\end{enumerate}
The tilted category $\mathcal{C}_{R^{\flat\circ}}^a$ and its admissible class $\mathbf{S}_\flat$ are defined identically over $\operatorname{Spf}^a((R^{\flat\circ})^a)$ via $\omega^\flat$ and $R_n^\flat$. Applying Section~\ref{sec:span_formalism} to $(\mathcal{C}_{R^\circ}^a, \tau_{\mathrm{Zar}}, \mathbf{S}_\theta)$ and $(\mathcal{C}_{R^{\flat\circ}}^a, \tau_{\mathrm{Zar}}, \mathbf{S}_\flat)$ defines the arithmetic and tilted almost derived algebraic cobordism theories $\underline{\dMGL}_{\theta}^a$ and $\underline{\dMGL}_{\flat}^a$.
\end{definition}

For a ring $B$ and a finitely generated ideal $I=(f_1,\dots,f_r)\subset B$, a complex $K\in D(B)$ is derived $I$-complete if $\mathrm{RHom}_B(B_f,K)=0$ for all $f\in I$~\cite[Tag~091S]{stacks-project}, defining the full subcategory $D_{\mathrm{comp}}(B,I)\subseteq D(B)$, and its derived $I$-completion is $K^\wedge_I \simeq \varprojlim_{n\ge 1}\bigl(K\otimes_BK_\bullet(B;f_1^n,\dots,f_r^n)\bigr)$~\cite[Tags~091V and~0920]{stacks-project}. Here $K_\bullet(B;f_1,\dots,f_r)$ denotes the Koszul complex on $f_1,\dots,f_r$; when $f_1,\dots,f_r$ form a regular sequence in $B$, $K_\bullet(B;f_1^n,\dots,f_r^n)\simeq B/(f_1^n,\dots,f_r^n)$, which applies to the regular sequences $\omega\in R^\circ$, $\omega^\flat\in R^{\flat\circ}$, and $(p,\xi)$ in $\mathbb A_{\inf R}$~\cite[Lemma~3.8]{BhattScholze2017Prisms}. Throughout, an \emph{integral perfectoid $R^\circ$-algebra} means a discrete, adically complete perfectoid ring equipped with a continuous $R^\circ$-algebra structure, corresponding under Bhatt--Scholze's equivalence to a perfect prism over $(\mathbb{A}_{\inf R}, I_R)$ in $\mathcal{T}_R^{\mathrm{perf}}$~\cite[Theorem~3.10]{BhattScholze2017Prisms}.

\subsection{Arithmetic and tilted fiber specializations}
\label{subsec:fiber-specializations}

For an affine formal scheme $X = \operatorname{Spf}(A)$ over $\operatorname{Spf}(B)\in\mathcal T_R$, restricting along $\theta \colon \mathbb{A}_{\inf R} \twoheadrightarrow R^\circ$ and $\rho_\flat \colon \mathbb{A}_{\inf R} \twoheadrightarrow R^{\flat\circ}$ defines the \emph{arithmetic} and \emph{tilted fiber specializations} $i_\theta^*\operatorname{Spf}(A) = \operatorname{Spf}^a(A_\theta)\in\mathcal{C}_{R^\circ}^a$ and $i_\flat^*\operatorname{Spf}(A) = \operatorname{Spf}^a(A_\flat)\in\mathcal{C}_{R^{\flat\circ}}^a$, where
\[
A_\theta = \Bigl(\bigl(A \otimes_{\mathbb{A}_{\inf R},\theta} R^\circ\bigr)^a\Bigr)^\wedge_{\omega},
\qquad
A_\flat = \Bigl(\bigl(A \otimes_{\mathbb{A}_{\inf R},\rho_\flat} R^{\flat\circ}\bigr)^a\Bigr)^\wedge_{\omega^\flat},
\]
with level-$n$ reductions $\Spec^a\bigl((A \otimes_{\mathbb{A}_{\inf R}} R_n)^a\bigr)$ and $\Spec^a\bigl((A \otimes_{\mathbb{A}_{\inf R}} R_n^\flat)^a\bigr)$. Zariski gluing extends $i_\theta^*$ and $i_\flat^*$ to $\mathcal{C}_{\mathcal{T}_R}$ and to module pullbacks.

\begin{proposition}
\label{prop:prismatic-fiber-completion}
Let $A$ be a derived $(p,\xi)$-complete animated $\mathbb{A}_{\inf R}$-algebra. Then $A \otimes_{\mathbb{A}_{\inf R},\theta} R^\circ$ is derived $p$- and $\omega$-complete, and $A \otimes_{\mathbb{A}_{\inf R},\rho_\flat} R^{\flat\circ}$ is derived $\bar\xi$- and $\omega^\flat$-complete. The canonical completion maps are equivalences:
\[
\bigl(A \otimes_{\mathbb{A}_{\inf R},\theta} R^\circ\bigr)^a \longrightarrow A_\theta,
\qquad
\bigl(A \otimes_{\mathbb{A}_{\inf R},\rho_\flat} R^{\flat\circ}\bigr)^a \longrightarrow A_\flat.
\]
\end{proposition}
\begin{proof}
By Lemma~\ref{lem:tilting-frame-koszul}(3), $A \otimes_{\mathbb{A}_{\inf R},\theta} R^\circ \simeq \cofib(A\xrightarrow{\xi}A)$ and $A \otimes_{\mathbb{A}_{\inf R},\rho_\flat} R^{\flat\circ} \simeq \cofib(A\xrightarrow{p}A)$. Since $A$ is derived $(p,\xi)$-complete, these cofibers are derived $p$- and $\xi$-complete over $\mathbb A_{\inf R}$, and hence derived $p$-complete in $D(R^\circ)$ and derived $\bar\xi$-complete in $D(R^{\flat\circ})$~\cite[Tags~091S, 091U, and~0924]{stacks-project}. Lemma~\ref{lem:tilting-frame-koszul}(1) identifies the $(p)$-adic and $(\omega)$-adic topologies on $R^\circ$, and the $(\bar\xi)$-adic and $(\omega^\flat)$-adic topologies on $R^{\flat\circ}$, giving derived $\omega$- and $\omega^\flat$-completeness. Finally, since $(-)^a$ is symmetric monoidal and right adjoint to the exact functor $(-)_!$~\cite[Proposition~2.2.21 and Corollary~2.2.22]{GR}, it commutes with $-\otimes_{R^\circ}R_n$, $-\otimes_{R^{\flat\circ}}R_n^\flat$, and $\varprojlim_{n\ge 1}$.
\end{proof}

\begin{proposition}
\label{prop:prismatic-fiber-geometric-compatibility}
For any $U \in \mathcal{T}_R$, the functors $i_\theta^* \colon \mathrm{dForm}(U) \to \mathcal{C}_{R^\circ}^a$ and $i_\flat^* \colon \mathrm{dForm}(U) \to \mathcal{C}_{R^{\flat\circ}}^a$ satisfy:
\begin{enumerate}[\rm (1)]
\item Both $i_\theta^*$ and $i_\flat^*$ preserve finite homotopy limits and finite coproducts, and carry $\mathbf{S}_{\mathrm{Prism},U}$ into $\mathbf{S}_\theta$ and $\mathbf{S}_\flat$.
\item For any finite locally free module $\mathcal{E}$ on $X \in \mathrm{dForm}(U)$ and each $\nu \in \{\theta, \flat\}$, the functor $i_\nu^*$ commutes with the formation of relative total spaces, projective bundles, zero sections, projective zero sections, and derived zero loci:
\[
i_\nu^*\mathbb{V}_X(\mathcal{E}) \simeq \mathbb{V}_{i_\nu^*X}(i_\nu^*\mathcal{E}), \qquad
i_\nu^*\mathbb{P}_X(\mathcal{E}) \simeq \mathbb{P}_{i_\nu^*X}(i_\nu^*\mathcal{E}).
\]
\end{enumerate}
\end{proposition}
\begin{proof}
For (1), Proposition~\ref{prop:prismatic-fiber-completion} implies that $i_\theta^*$ and $i_\flat^*$ are given on affine formal schemes by $-\otimes_{\mathbb{A}_{\inf R}} R^\circ$ and $-\otimes_{\mathbb{A}_{\inf R}} R^{\flat\circ}$ followed by exact almostification $(-)^a$, preserving completed derived tensor products, homotopy pullbacks, terminal objects, and finite disjoint unions. For any $f \colon X \to Y$ in $\mathbf{S}_{\mathrm{Prism},U}$, base change gives $\widehat{L}_{i_\theta^*X / i_\theta^*Y}^a \simeq i_\theta^*\widehat{L}_{X/Y}$ and $\widehat{L}_{i_\flat^*X / i_\flat^*Y}^a \simeq i_\flat^*\widehat{L}_{X/Y}$, preserving perfection, Tor-amplitude $[-1,0]$, topological local almost finite presentation, quasi-compactness, and quasi-separatedness.

For (2), $\mathbb{V}_X(\mathcal{E})$ and $\mathbb{P}_X(\mathcal{E})$ are formed locally from completed symmetric algebras and standard Zariski open affine coverings, which commute with $i_\theta^*$ and $i_\flat^*$ by Definition~\ref{def:almost-formal-base-categories}(1); preservation of zero sections, projective zero sections, and derived zero loci follows from (1).
\end{proof}

\begin{definition}
\label{def:prismatically-liftable-pair}
An \emph{augmented common formal scheme} in $\mathcal{C}_{\mathcal{T}_R}$ is a pair $\bigl(X_0 = \operatorname{Spf}(B),\, X_I = \operatorname{Spf}(B \oplus I)\bigr)$ consisting of $X_0 \in \mathcal{T}_R$ and a split augmented derived $(p,\xi)$-complete animated $B$-algebra $B \oplus I \to B$ over the underlying ring $B$. Fix an integral perfectoid $R^\circ$-algebra $A^\circ$ with tilt $A^{\flat\circ}$ and perfect prism $\widetilde{A} = \Ainf(A^\circ) = W(A^{\flat\circ})$, and consider a pair $\bigl(\operatorname{Spf}^a((A^\circ)^a \oplus J_\theta),\, \operatorname{Spf}^a((A^{\flat\circ})^a \oplus J_\flat)\bigr)$ in $\mathcal{C}_{R^\circ}^a \times \mathcal{C}_{R^{\flat\circ}}^a$, where $J_\theta$ and $J_\flat$ are connective derived $\omega$- and $\omega^\flat$-complete non-unital animated almost augmentation ideals (Subsection~\ref{subsec:canonical-split-lifting}). This pair admits a \emph{common augmented prismatic lift} if there exists an augmented common formal scheme $\bigl(\operatorname{Spf}(\widetilde{A}),\, \operatorname{Spf}(\widetilde{A} \oplus I)\bigr)$ equipped with augmentation-preserving equivalences $i_\theta^*\operatorname{Spf}(\widetilde{A} \oplus I) \to \operatorname{Spf}^a((A^\circ)^a \oplus J_\theta)$ and $i_\flat^*\operatorname{Spf}(\widetilde{A} \oplus I) \to \operatorname{Spf}^a((A^{\flat\circ})^a \oplus J_\flat)$. It is \emph{prismatically liftable} if such lifts exist locally in the comparison topology with coherent transition data.
\end{definition}

\section{Deformation theory and comparison topologies}
\label{sec:deformation-perfectoid}

Throughout this section, a \emph{finite admissible diagram} is a finite diagram of quasi-compact, quasi-separated quasi-smooth formal schemes over the base formal scheme, whose vertices and admissible legs have perfect completed relative cotangent complexes of Tor-amplitude $[-1,0]$. It presents an admissible span $[X\xleftarrow{s}Z\to Y]$, an admissible pushout square, or a relation diagram. This family is closed under derived base change and the finite operations of Section~\ref{sec:span_formalism}. A lift of such a diagram consists of compatible lifts of its vertices, arrows, and coherence homotopies.

\subsection{Adic towers and square-zero transitions}
\label{subsec:adic-towers}

For $n\geq 1$, set $I_n = \ker(R_{n+1}\twoheadrightarrow R_n) = (\omega^n)/(\omega^{n+1})$ and $I_n^\flat = \ker(R_{n+1}^\flat\twoheadrightarrow R_n^\flat) = ((\omega^\flat)^n)/((\omega^\flat)^{n+1})$. Since $2n\ge n+1$, $I_n^2=0$ in $R_{n+1}$ and $(I_n^\flat)^2=0$ in $R_{n+1}^\flat$. Since $\omega\in R^\circ$ and $\omega^\flat\in R^{\flat\circ}$ are non-zero-divisors, multiplication by $\omega^n$ and $(\omega^\flat)^n$ induces module isomorphisms
\begin{equation}
\label{eq:square-zero-kernel-identification}
\cdot\,\omega^n\colon R_1 = R^\circ/(\omega) \longrightarrow I_n,
\qquad
\cdot\,(\omega^\flat)^n\colon R_1^\flat = R^{\flat\circ}/(\omega^\flat) \longrightarrow I_n^\flat.
\end{equation}

Let $\mathcal C_{R_n}^{qs,\mathrm{aff}}$ and $\mathcal C_{R_n^\flat}^{qs,\mathrm{aff}}$ denote the full $\infty$-subcategories of affine quasi-smooth derived almost $R_n$- and $R_n^\flat$-schemes. For $Z_n\in\mathcal C_{R_n}^{qs,\mathrm{aff}}$, write $\mathcal O_{Z_n}$ for a chosen coordinate animated $R_n$-algebra, $\mathcal O_{Z_n}^a$ for its almostification, and $L_{Z_n/R_n}^a$ for its relative almost cotangent complex. Let $\mathrm{Lift}(Z_n)$ be the $\infty$-groupoid of $Z_{n+1}\in\mathcal C_{R_{n+1}}^{qs,\mathrm{aff}}$ equipped with a reduction equivalence $Z_{n+1}\times_{R_{n+1}}R_n \to Z_n$, and define $\mathrm{Lift}(W_n)$ for $W_n\in\mathcal C_{R_n^\flat}^{qs,\mathrm{aff}}$ analogously. Set $I_{Z_n}^a = (I_n\otimes_{R_n}\mathcal O_{Z_n})^a$ and $I_{W_n}^{\flat,a} = (I_n^\flat\otimes_{R_n^\flat}\mathcal O_{W_n})^a$, realized via the flat firm modules $\tm_n$ and $\tm_n^\flat$ of Subsection~\ref{subsec:almost-formal-base} by the firm representatives $\tm_n\otimes_{R_n}(I_n\otimes_{R_n}\mathcal O_{Z_n})$ and $\tm_n^\flat\otimes_{R_n^\flat}(I_n^\flat\otimes_{R_n^\flat}\mathcal O_{W_n})$. Mapping spaces $\Map_{\mathcal O_{Z_n}^a}(-,-)$ and obstruction groups $\Ext^i_{\mathcal O_{Z_n}^a}(-,-)$ in the derived almost module category are computed via these firm models, distinguishing internal derived Hom objects $\underline{\mathrm{RHom}}_{\mathcal O_{Z_n}^a}$ from mapping spaces $\Map_{\mathcal O_{Z_n}^a} = \Omega^\infty R(\underline{\mathrm{RHom}}_{\mathcal O_{Z_n}^a})_*$.

\subsection{Obstruction theory and torsor comparison}
\label{subsec:obstruction-torsor-comparison}

\begin{proposition}
\label{prop:deformation-control}
For any $Z_n\in\mathcal C_{R_n}^{qs,\mathrm{aff}}$, there is a canonical obstruction class $[o_{Z_n}]\in\Ext^2_{\mathcal O_{Z_n}^a}(L_{Z_n/R_n}^a,I_{Z_n}^a)$ to lifting $Z_n$ along $R_{n+1}\to R_n$. The space $\mathrm{Lift}(Z_n)$ is nonempty if and only if $[o_{Z_n}]=0$, in which case $\mathrm{Lift}(Z_n)$ is a torsor under the derivation group space $\Map_{\mathcal{O}_{Z_n}^a}(L_{Z_n/R_n}^a, I_{Z_n}^a[1])$, trivialized upon choosing a base lift $Z_{n+1}\in\mathrm{Lift}(Z_n)$.
\end{proposition}
\begin{proof}
Since $I_n^2=0$, $R_{n+1}\to R_n$ is classified by a derivation $d_{R_n}\colon R_n\to R_n\oplus I_n[1]$, presenting $R_{n+1}$ as $R_n\times_{R_n\oplus I_n[1]}R_n$ along the zero section $s_0\colon R_n\to R_n\oplus I_n[1]$~\cite[Theorem~7.4.1.26]{HA}. By square-zero deformation theory for derived almost algebras~\cite[Section~3.2]{GR}, composing the cotangent transitivity map $L_{Z_n/R_n}^a\to (L_{R_n}^a\otimes_{R_n^a}\mathcal O_{Z_n}^a)[1]$ with $d_{R_n}$ defines $o_{Z_n}\in\Map_{\mathcal O_{Z_n}^a}(L_{Z_n/R_n}^a,I_{Z_n}^a[2])$ fitting into the homotopy fiber sequence
\[
\mathrm{Lift}(Z_n)\longrightarrow * \xrightarrow{\ o_{Z_n}\ } \Map_{\mathcal O_{Z_n}^a}\bigl(L_{Z_n/R_n}^a,I_{Z_n}^a[2]\bigr)
\]
over $0$. Thus $\mathrm{Lift}(Z_n)$ is nonempty if and only if $[o_{Z_n}]=0$, in which case $\mathrm{Lift}(Z_n)\simeq\operatorname{Path}(o_{Z_n},0)$ is a nonempty torsor under the based loop space $\Omega\Map_{\mathcal O_{Z_n}^a}(L_{Z_n/R_n}^a,I_{Z_n}^a[2]) \simeq \Map_{\mathcal O_{Z_n}^a}(L_{Z_n/R_n}^a,I_{Z_n}^a[1])$.
\end{proof}

\begin{lemma}
\label{lem:reduction}
For any $Z_n\in\mathcal C_{R_n}^{qs,\mathrm{aff}}$ with first adic reduction $Z_1=Z_n\times_{R_n}R_1\in\mathcal C_{R_1}^{qs,\mathrm{aff}}$ and any integer $k\ge 0$, there is an equivalence
\[
\Map_{\mathcal O_{Z_n}^a}\bigl(L_{Z_n/R_n}^a,I_{Z_n}^a[k]\bigr)
\longrightarrow
\Map_{\mathcal O_{Z_1}^a}\bigl(L_{Z_1/R_1}^a,\mathcal O_{Z_1}^a[k]\bigr),
\]
inducing isomorphisms $\Ext^i_{\mathcal O_{Z_n}^a}(L_{Z_n/R_n}^a,I_{Z_n}^a)\to\Ext^i_{\mathcal O_{Z_1}^a}(L_{Z_1/R_1}^a,\mathcal O_{Z_1}^a)$ for all $i\in\mathbb Z$.
\end{lemma}
\begin{proof}
Let $j\colon Z_1\to Z_n$ be the reduction morphism, with derived pullback $j^*$ and restriction of scalars $j_*$. Via $\cdot\,\omega^n\colon R_1\xrightarrow{\sim}I_n$ from~\eqref{eq:square-zero-kernel-identification}, we have $I_n\otimes_{R_n}\mathcal O_{Z_n} \simeq R_1\otimes_{R_n}\mathcal O_{Z_n} \simeq j_*\mathcal O_{Z_1}$, and almostification identifies $I_{Z_n}^a\simeq j_*\mathcal O_{Z_1}^a$. Since $j^*L_{Z_n/R_n}^a \simeq L_{Z_1/R_1}^a$, the $(j^*,j_*)$-adjunction induces an equivalence of mapping spectra $\operatorname{map}_{\Mod_{\mathcal O_{Z_n}^a}}(L_{Z_n/R_n}^a,I_{Z_n}^a) \to \operatorname{map}_{\Mod_{\mathcal O_{Z_1}^a}}(L_{Z_1/R_1}^a,\mathcal O_{Z_1}^a)$; applying $\Omega^\infty(-[k])$ gives the equivalence of mapping spaces for each $k\ge 0$, and taking $H^i = \pi_{-i}$ gives the $\Ext^i$ isomorphisms for all $i\in\mathbb Z$.
\end{proof}

\begin{corollary}
\label{cor:reduction-controls}
For any $Z_n\in\mathcal C_{R_n}^{qs,\mathrm{aff}}$ with first adic reduction $Z_1=Z_n\times_{R_n}R_1\in\mathcal C_{R_1}^{qs,\mathrm{aff}}$, there is an equivalence
\[
\Map_{\mathcal O_{Z_n}^a}\bigl(L_{Z_n/R_n}^a,I_{Z_n}^a[1]\bigr) \;\longrightarrow\; \Map_{\mathcal O_{Z_1}^a}\bigl(L_{Z_1/R_1}^a,\mathcal O_{Z_1}^a[1]\bigr).
\]
\end{corollary}

\begin{remark}
\label{rem:reduction-controls-functoriality}
By Corollary~\ref{cor:reduction-controls}, the acting derivation group space at level $n$ depends functorially only on the first adic reduction $Z_1$.
\end{remark}

\begin{lemma}
\label{lem:cotangent-dual-vanishing}
For any $Z_n\in\mathcal C_{R_n}^{qs,\mathrm{aff}}$ and any $i\geq 2$, $H^i\bigl((L_{Z_n/R_n}^a)^\vee\otimes_{\mathcal O_{Z_n}^a}I_{Z_n}^a\bigr) = 0$.
\end{lemma}
\begin{proof}
Since $L_{Z_n/R_n}^a$ is perfect of Tor-amplitude $[-1,0]$, its dual $(L_{Z_n/R_n}^a)^\vee$ has Tor-amplitude in cohomological degrees $[0,1]$. Since $I_n$ and $\mathcal O_{Z_n}$ are connective and $\tm_n$ is flat over $R_n$, $I_{Z_n}^a$ is concentrated in cohomological degrees $(-\infty,0]$. Thus $(L_{Z_n/R_n}^a)^\vee\otimes_{\mathcal O_{Z_n}^a}I_{Z_n}^a$ is concentrated in degrees $(-\infty,1]$, implying the vanishing for $i\ge 2$.
\end{proof}
\begin{corollary}
\label{cor:honest-chart-lift-existence}
If $Z_n$ is the almost realization of a specified ordinary affine quasi-smooth derived $R_n$-scheme $\widetilde{Z}_n = \operatorname{Spec}(A_n)$ compatible with $(-)_!$ and derived base change along the adic quotients, then $Z_n$ admits an affine quasi-smooth lift to $R_{n+1}$.
\end{corollary}
\begin{proof}
Over the affine ring $A_n$, the relative cotangent complex $L_{A_n/R_n}$ is a perfect $A_n$-module of Tor-amplitude $[-1,0]$ and $I_n \otimes_{R_n} A_n$ is connective. Dualizability of $L_{A_n/R_n}$ over $A_n$ and Lemma~\ref{lem:cotangent-dual-vanishing} imply $\operatorname{Ext}^2_{A_n}(L_{A_n/R_n}, I_n \otimes_{R_n} A_n) \cong H^2\bigl(L_{A_n/R_n}^\vee \otimes_{A_n} (I_n \otimes_{R_n} A_n)\bigr) = 0$. Thus $\widetilde{Z}_n$ admits an ordinary affine quasi-smooth lift $\widetilde{Z}_{n+1} = \operatorname{Spec}(A_{n+1})$ over $R_{n+1}$~\cite[Tag~08SP]{stacks-project}, whose almost realization defines $Z_{n+1} \in \operatorname{Lift}(Z_n)$.
\end{proof}
\begin{corollary}
\label{cor:torsor-comparison}
Let $Z_n\in\mathcal C_{R_n}^{qs,\mathrm{aff}}$ and $W_n\in\mathcal C_{R_n^\flat}^{qs,\mathrm{aff}}$ satisfy:
\begin{enumerate}[\rm (a)]
\item their first reductions $Z_1=Z_n\times_{R_n}R_1$ and $W_1=W_n\times_{R_n^\flat}R_1^\flat$ are identified via $\iota_1\colon R_1\xrightarrow{\sim} R_1^\flat$; and
\item both $\mathrm{Lift}(Z_n)$ and $\mathrm{Lift}(W_n)$ are nonempty.
\end{enumerate}
Then:
\begin{enumerate}[\rm (1)]
\item Setting $G_{Z_1} = \Map_{\mathcal O_{Z_1}^a}(L_{Z_1/R_1}^a,\mathcal O_{Z_1}^a[1])$, there are equivalences:
\[
\Map_{\mathcal O_{Z_n}^a}\bigl(L_{Z_n/R_n}^a,I_{Z_n}^a[1]\bigr)
\;\longrightarrow\;
G_{Z_1}
\;\longleftarrow\;
\Map_{\mathcal O_{W_n}^a}\bigl(L_{W_n/R_n^\flat}^a,I_{W_n}^{\flat,a}[1]\bigr).
\]
\item Both $\mathrm{Lift}(Z_n)$ and $\mathrm{Lift}(W_n)$ are nonempty $G_{Z_1}$-torsors, and any choice of base lifts $Z_{n+1}\in\mathrm{Lift}(Z_n)$ and $W_{n+1}\in\mathrm{Lift}(W_n)$ determines a $G_{Z_1}$-equivariant equivalence $\mathrm{Lift}(Z_n)\to\mathrm{Lift}(W_n)$ via $g\cdot Z_{n+1}\mapsto g\cdot W_{n+1}$.
\end{enumerate}
\end{corollary}
\begin{proof}
Part (1) follows from Corollary~\ref{cor:reduction-controls} and the identification $Z_1\simeq W_1$ over $\iota_1\colon R_1\xrightarrow{\sim} R_1^\flat$. Part (2) follows from Proposition~\ref{prop:deformation-control} and the orbit equivalences $G_{Z_1}\to\mathrm{Lift}(Z_n)$ and $G_{Z_1}\to\mathrm{Lift}(W_n)$ defined by $Z_{n+1}$ and $W_{n+1}$.
\end{proof}

To compare correspondence diagrams across these torsors, we present the lift torsors by simplicial affine Bar hypercovers in Subsection~\ref{subsec:affine-bar-models}, establish pointed comparisons via cotangent and relative tangential excisive approximation in Subsection~\ref{subsec:canonical-split-lifting}, and adjoin section sieves to form the deformation and comparison topologies in Subsection~\ref{subsec:lift-torsor-localization-sites}.

\subsection{Affine Bar models for lift torsors}
\label{subsec:affine-bar-models}

Let $A'\to A$ be a square-zero extension of connective animated rings with discrete kernel $M$, let $B$ be a connective quasi-smooth animated $A$-algebra, and set $M_B=B\otimes_AM$. As in Proposition~\ref{prop:deformation-control}, letting $o_B\in K_B=\Map_{\Mod_B}(L_{B/A},M_B[2])$ be the obstruction to lifting $B$ along $A'\to A$, the global lift $\infty$-groupoid $\operatorname{Lift}_{A'}(B) = \{(\widetilde B,\alpha)\mid \widetilde B\in\CAlg_{A'},\, \alpha\colon\widetilde B\otimes_{A'}A\xrightarrow{\sim} B\}^{\simeq}$ is identified with $\operatorname{Path}_{K_B}(o_B,0)$, forming a torsor under $G_B = \Omega K_B \simeq \Map_{\Mod_B}(L_{B/A},M_B[1])$ when $[o_B]=0$. Suppose that $M$ factors as a finite projective module over a quotient ring $A_1$ of $A$ and $L_{B/A}$ is presented by a cofiber sequence $P_1\xrightarrow{d}P_0\to L_{B/A}$ of finite projective $B$-modules concentrated in degree $0$, and set $B_1=B\otimes_AA_1$, $P_{1,i}=P_i\otimes_BB_1$, and $V_i=\underline{\Map}_{B_1}(P_{1,i},M_B)$ for $i\in\{0,1\}$, where $\underline{\Map}_{B_1}(-,-)$ denotes the mapping stack on derived affine schemes over $\Spec B_1$ and $\mathbf B(-)=[*/-]$ denotes the classifying stack functor. Since $M_B\simeq B_1\otimes_{A_1}M$ and $L_{B_1/A_1}\simeq L_{B/A}\otimes_BB_1$, adjunction along $B\to B_1$ identifies $K_B\simeq\Map_{\Mod_{B_1}}(L_{B_1/A_1},M_B[2])$ and $G_B\simeq\Map_{\Mod_{B_1}}(L_{B_1/A_1},M_B[1])$; passing to \'etale sheaves on derived affine schemes over $\Spec B_1$ defines the acting derivation group stack $\underline{G}_B = \underline{\Map}_{B_1}(L_{B_1/A_1},M_B[1])$ and the $\underline{G}_B$-torsor stack $\underline{\operatorname{Lift}}_{A'}(B)$ with global sections $\operatorname{Lift}_{A'}(B)$.

\begin{proposition}
\label{prop:affine-bar-models}
Let $A'\to A$ be a square-zero extension whose discrete kernel $M$ factors as a finite projective module over a quotient ring $A_1$ of $A$, and let $B$ be a connective quasi-smooth animated $A$-algebra with $[o_B]=0$ whose relative cotangent complex $L_{B/A}$ is presented by a cofiber sequence $P_1\xrightarrow{d}P_0\to L_{B/A}$ of finite projective $B$-modules concentrated in degree $0$.
\begin{enumerate}[\rm (1)]
\item Over $\Spec B_1$, the acting derivation group stack $\underline{G}_B$ fits into a homotopy fiber sequence of connective $\mathbb E_\infty$-group stacks
\[
V_0 \xrightarrow{\ d^*\ } V_1 \longrightarrow \underline{G}_B \longrightarrow \mathbf B V_0,
\]
where $V_0$ and $V_1$ are derived affine vector group schemes over $\Spec B_1$, presenting $\underline{G}_B \simeq [V_1 / V_0]$ whose almostification recovers the almost derivation group stack over $\Spec^a B_1^a$.
\item The lift stack $\underline{\operatorname{Lift}}_{A'}(B)$ is presented over $\Spec B_1$ by the simplicial affine Bar hypercover
\[
\underline{\operatorname{Lift}}_{A'}(B) \;\simeq\; \bigl[\,T_B \,/\, V_0\,\bigr]
\;\simeq\;
\operatorname*{colim}_{[q]\in\Delta^{\mathrm{op}}}\bigl(T_B\times_{\Spec B_1} V_0^{\times_{\Spec B_1} q}\bigr),
\]
where $T_B$ is an affine $V_1$-torsor over $\Spec B_1$ parametrizing termwise lifts, and each term $T_B\times_{\Spec B_1} V_0^{\times_{\Spec B_1} q}$ ($q\ge 0$) is a derived affine scheme over $\Spec B_1$.
\item This simplicial affine presentation commutes with derived base change and, after almostification, with the boundary specializations $i_\theta^*$ and $i_\flat^*$.
\end{enumerate}
\end{proposition}
\begin{proof}
For (1), applying $\underline{\Map}_{B_1}(-,M_B[1])$ to the base change of the cofiber sequence $P_{1,1}\to P_{1,0}\to L_{B_1/A_1}\to P_{1,1}[1]$ and rotating the resulting fiber sequence $V_1\to\underline{G}_B\to \mathbf B V_0\xrightarrow{\mathbf B(d^*)}\mathbf B V_1$ produces $V_0\xrightarrow{d^*}V_1\to\underline{G}_B\to \mathbf B V_0$, where $V_i\simeq P_{1,i}^\vee\otimes_{B_1}M_B$ is represented over $\Spec B_1$ by the derived affine scheme $\Spec\bigl(\operatorname{Sym}_{B_1}(P_{1,i}\otimes_{A_1}M^\vee)\bigr)$, exhibiting $\underline{G}_B\simeq [V_1/V_0]$.

For (2), in the delooping fiber sequence $\mathbf B V_1 \to \mathbf B\underline{G}_B \to \mathbf B^2 V_0 \simeq \underline{\Map}_{B_1}(P_{1,0},M_B[2])$, since $P_{1,0}$ is finite projective over $B_1$ in degree $0$ and $M_B[2]$ is $2$-connective, $P_{1,0}^\vee\otimes_{B_1}M_B[2]$ is $2$-connective and the global obstruction in $\pi_0(\mathbf B^2 V_0(B_1))\cong\Ext^2_{B_1}(P_{1,0},M_B)=0$ vanishes. Thus $\underline{\operatorname{Lift}}_{A'}(B)$ is induced from an affine $V_1$-torsor $T_B$ over $\Spec B_1$ via the associated bundle quotient $\underline{\operatorname{Lift}}_{A'}(B)\simeq T_B\times^{V_1}\underline{G}_B = [(T_B\times_{\Spec B_1}\underline{G}_B)/V_1]\simeq [T_B/V_0]$, whose Bar nerve is $[q]\mapsto T_B\times_{\Spec B_1} V_0^{\times_{\Spec B_1} q}$.

Part (3) follows from Proposition~\ref{prop:prismatic-fiber-geometric-compatibility}, as derived base change and exact almostification preserve finite projective modules, duals, and finite fiber products of affine schemes.
\end{proof}

\begin{example}
\label{ex:adic-transition-bar-model}
For each $n\ge 1$, the square-zero adic transition $R_{n+1}\to R_n$ has kernel $M=I_n\xrightarrow{\sim}R_1$, which is free of rank $1$ over the quotient ring $A_1=R_1$ of $A=R_n$, and likewise for $R_{n+1}^\flat\to R_n^\flat$ with kernel $I_n^\flat\xrightarrow{\sim}R_1^\flat$; Proposition~\ref{prop:affine-bar-models} therefore applies across all stages of the arithmetic and tilted adic towers.
\end{example}
\subsection{Non-unital algebras, cotangent fibers, and coherent diagram lifts}
\label{subsec:canonical-split-lifting}

Let $X_0=\operatorname{Spf}(B)\in\mathcal T_R$ be the affine formal scheme of a prism over $(\mathbb{A}_{\inf R}, I_R)$ with boundary specializations $i_\theta^*(X_0)=\operatorname{Spf}^a(B_\theta)\in\mathcal C_{R^\circ}^a$ and $i_\flat^*(X_0)=\operatorname{Spf}^a(B_\flat)\in\mathcal C_{R^{\flat\circ}}^a$, where $B_\theta=(B\otimes_{\mathbb A_{\inf R},\theta}R^\circ)^a$ and $B_\flat=(B\otimes_{\mathbb A_{\inf R},\rho_\flat}R^{\flat\circ})^a$.

Let $\NAlg_{\mathrm{Prism}}(B)$, abbreviated $\NAlg_{\mathrm{Prism}}$, be the pointed $\infty$-category of connective, derived $(p,I_R)$-complete non-unital animated $B$-algebras over the underlying ring $B$. By operadic Smith ideal theory~\cite{WY2024},~\cite[Subsection~4.1]{KatoCalcNUAlg}, direct sum unitalization $I\mapsto B\oplus I$ and the augmentation ideal functor $\fib(-\to B)$ define mutually inverse equivalences $\NAlg_{\mathrm{Prism}}(B) \rightleftarrows \CAlg_{B//B}^{\wedge}$ with split augmented, connective, derived $(p,I_R)$-complete animated $B$-algebras. Contravariantly, any $I\in\NAlg_{\mathrm{Prism}}(B)$ corresponds to the pointed formal scheme $X_I=\operatorname{Spf}(B\oplus I)$, whose structure projection $X_I\to X_0$ is induced by the zero section $B\to B\oplus I$ and whose basepoint section $X_0\to X_I$ is induced by the augmentation $\epsilon_I\colon B\oplus I\to B$. On $\Mod_B^{\ge 0,\wedge}$, the underlying module functor $\operatorname{ev}\colon\NAlg_{\mathrm{Prism}}(B)\to\Mod_B^{\ge 0,\wedge}$ admits the completed non-unital free animated algebra left adjoint $\operatorname{Free}_B$. The trivial square-zero functor $\operatorname{sqz}_B\colon\Mod_B^{\ge 0,\wedge}\to\NAlg_{\mathrm{Prism}}(B)$ admits the left adjoint completed relative cotangent fiber functor $Q_B\colon\NAlg_{\mathrm{Prism}}(B)\to\Mod_B^{\ge 0,\wedge}$ given by $Q_B(I) = \widehat{L}_{(B\oplus I)/B}\widehat{\otimes}_{B\oplus I} B$~\cite[Section~7.3.4]{HA},~\cite[Definition~4.3]{KatoCalcNUAlg}, with $\operatorname{ev}\circ\operatorname{sqz}_B\simeq\mathrm{id}$.

For $\nu\in\{\theta,\flat\}$, set $\omega_\theta=\omega$, $\omega_\flat=\omega^\flat$, $R_\theta=R^\circ$, $R_\flat=R^{\flat\circ}$, $\tm_\theta=\tm$, and $\tm_\flat=\tm^\flat$. Let $\NAlg_\nu^a(B_\nu)$, abbreviated $\NAlg_\nu^a$, be the pointed $\infty$-category of connective, derived $\omega_\nu$-complete non-unital animated almost $B_\nu$-algebras, realized via firm representatives $J_!=\tm_\nu\otimes_{R_\nu}J_*$ and almost unitalization $(B_\nu\oplus J)_{!!}$~\cite[Subsection~4.4]{KatoCalcNUAlg}. Let $\Mod_{B_\nu}^{\ge 0,\wedge}$ denote the $\infty$-category of connective, derived $\omega_\nu$-complete almost $B_\nu$-modules. Unitalization identifies $\NAlg_\nu^a(B_\nu)\simeq\CAlg_{B_\nu//B_\nu}^{a,\wedge}$, with adjunctions $\operatorname{Free}_{B_\nu}\colon\Mod_{B_\nu}^{\ge 0,\wedge}\rightleftarrows\NAlg_\nu^a(B_\nu)\colon\!\operatorname{ev}$ and $Q_{B_\nu}^a\colon\NAlg_\nu^a(B_\nu)\rightleftarrows\Mod_{B_\nu}^{\ge 0,\wedge}\colon\!\operatorname{sqz}_{B_\nu}$ satisfying $Q_{B_\nu}^a(J) = \widehat{L}_{(B_\nu\oplus J)/B_\nu}^a\widehat{\otimes}_{B_\nu\oplus J} B_\nu$ and $\operatorname{ev}\circ\operatorname{sqz}_{B_\nu}\simeq\mathrm{id}$. By Proposition~\ref{prop:prismatic-fiber-completion}, derived base change and almostification define the pointed specialization functors $i_\theta^*\colon\NAlg_{\mathrm{Prism}}(B)\to\NAlg_\theta^a(B_\theta)$, $I\mapsto (I\otimes_{\mathbb A_{\inf R},\theta}R^\circ)^a$, and $i_\flat^*\colon\NAlg_{\mathrm{Prism}}(B)\to\NAlg_\flat^a(B_\flat)$, $I\mapsto (I\otimes_{\mathbb A_{\inf R},\rho_\flat}R^{\flat\circ})^a$.

For each $\nu\in\{\theta,\flat\}$, these cotangent fiber adjunctions and boundary specializations fit into the square of $\infty$-categories, commuting in both directions by Proposition~\ref{prop:excisive-reflection-means}(3):
\[
\xymatrix@C=4.5em@R=2.2em{
\NAlg_{\mathrm{Prism}}(B) \ar@<0.6ex>[r]^-{Q_B} \ar[d]_-{i_\nu^*}
& \Mod_B^{\ge 0,\wedge} \ar@<0.6ex>[l]^-{\operatorname{sqz}_B} \ar[d]^-{i_\nu^*} \\
\NAlg_\nu^a(B_\nu) \ar@<0.6ex>[r]^-{Q_{B_\nu}^a}
& \Mod_{B_\nu}^{\ge 0,\wedge}. \ar@<0.6ex>[l]^-{\operatorname{sqz}_{B_\nu}}
}
\]

For $\mathcal C\in\{\NAlg_{\mathrm{Prism}}(B),\NAlg_\nu^a(B_\nu)\}$ with base $A\in\{B,B_\theta,B_\flat\}$ and cotangent fiber adjunction unit $\eta_I\colon I\to\operatorname{sqz}_A(Q_A(I))$, write $\Sigma_{\mathcal C} I = 0\amalg_I^{\wedge}0$ for the categorical suspension functor in $\mathcal C$.

\begin{proposition}
\label{prop:excisive-reflection-means}
Let $X_0=\operatorname{Spf}(B)\in\mathcal T_R$ have boundary specializations $i_\nu^*(X_0)=\operatorname{Spf}^a(B_\nu)$ ($\nu\in\{\theta,\flat\}$). Let $\mathcal C$ denote either $\NAlg_{\mathrm{Prism}}(B)$ or $\NAlg_\nu^a(B_\nu)$, with base $A\in\{B,B_\theta,B_\flat\}$, module category $\Mod_A^{\ge 0,\wedge}$, and adjunctions $(\operatorname{Free}_A\dashv\operatorname{ev})$ and $(Q_A\dashv\operatorname{sqz}_A)$.
\begin{enumerate}[\rm (1)]
\item The $\infty$-category $\mathcal C$ admits all finite homotopy limits and all small colimits, and under unitalization $I\mapsto A\oplus I$ there is an equivalence
\begin{equation}
\label{eq:suspension-loop-nonunital}
A\oplus\Sigma_{\mathcal C} I \;\longrightarrow\; A\widehat{\otimes}_{A\oplus I}A.
\end{equation}

\item Both left adjoints $Q_A\colon\mathcal C\to\Mod_A^{\ge 0,\wedge}$ and $\operatorname{Free}_A\colon\Mod_A^{\ge 0,\wedge}\to\mathcal C$ preserve $0$ and all small colimits, with $Q_A\circ\operatorname{Free}_A\xrightarrow{\sim}\mathrm{id}_{\Mod_A^{\ge 0,\wedge}}$, and for any $I\in\mathcal C$ and $m\ge 0$ there is an equivalence
\begin{equation}
\label{eq:explicit-excisive-reflection}
Q_A(\Sigma_{\mathcal C}^m I) \;\longrightarrow\; Q_A(I)[m].
\end{equation}

\item For each $\nu\in\{\theta,\flat\}$ and any $I\in\NAlg_{\mathrm{Prism}}(B)$, and $M\in\Mod_B^{\ge 0,\wedge}$, the functor $i_\nu^*\colon\NAlg_{\mathrm{Prism}}(B)\to\NAlg_\nu^a(B_\nu)$ preserves $0$ and finite completed colimits, and induces equivalences
\begin{equation}
\label{eq:boundary-reflection-commutation}
i_\nu^*\bigl(Q_B(I)\bigr) \;\longrightarrow\; Q_{B_\nu}^a\bigl(i_\nu^* I\bigr),
\qquad
i_\nu^*\bigl(\operatorname{sqz}_B(M)\bigr) \;\longrightarrow\; \operatorname{sqz}_{B_\nu}\bigl(i_\nu^* M\bigr).
\end{equation}
\end{enumerate}
\end{proposition}
\begin{proof}
For (1), the forgetful functor $\CAlg_{B//B}\to\Mod_B^{\ge 0}$ is conservative and preserves small homotopy limits, and $\Mod_B^{\ge 0,\wedge}\subset\Mod_B^{\ge 0}$ is closed under homotopy limits~\cite[Tag~091S]{stacks-project}, so $\NAlg_{\mathrm{Prism}}(B)$ admits finite homotopy limits computed via $\operatorname{ev}$. Dually, $\CAlg_{B//B}$ admits small colimits, and symmetric monoidal left adjoint derived completion $(-)^\wedge_{(p,I_R)}$~\cite[Tag~091V]{stacks-project} computes all small colimits in $\NAlg_{\mathrm{Prism}}(B)$. The same holds in $\NAlg_\nu^a(B_\nu)$ via exact firm representatives $M_!=\tm_\nu\otimes_{R_\nu}M_*$ and almost unitalization $C_{!!}$~\cite[Subsection~4.4]{KatoCalcNUAlg}. By the universal property of the completed algebraic cotangent complex $\widehat{L}_{(A\oplus I)/A}$~\cite[Section~7.3.4]{HA},~\cite[Definition~4.3]{KatoCalcNUAlg} and letting $\epsilon_{I*}\colon\Mod_A^{\wedge}\to\Mod_{A\oplus I}^{\wedge}$ denote restriction of scalars along the augmentation $\epsilon_I\colon A\oplus I\to A$, for any $I\in\mathcal C$ and $M\in\Mod_A^{\ge 0,\wedge}$ we have
\[
\Map_{\mathcal C}\bigl(I,\operatorname{sqz}_A(M)\bigr)
\;\simeq\;
\Map_{\Mod_{A\oplus I}^{\wedge}}\bigl(\widehat{L}_{(A\oplus I)/A},\epsilon_{I*}M\bigr)
\;\simeq\;
\Map_{\Mod_A^{\ge 0,\wedge}}\bigl(Q_A(I),M\bigr),
\]
establishing $Q_A\dashv\operatorname{sqz}_A$. Since $I\mapsto A\oplus I$ preserves pushouts and sends $0$ to $A$, we obtain~\eqref{eq:suspension-loop-nonunital}.

For (2), left adjoints preserve initial objects and all small colimits, while taking left adjoints of $\operatorname{ev}\circ\operatorname{sqz}_A\simeq\mathrm{id}_{\Mod_A^{\ge 0,\wedge}}$ gives $Q_A\circ\operatorname{Free}_A\simeq\mathrm{id}_{\Mod_A^{\ge 0,\wedge}}$. Since derived completion on $\Mod_A$ is exact and preserves connectivity of shifted modules $M[1]$ for $M\in\Mod_A^{\ge 0,\wedge}$, categorical suspension in $\Mod_A^{\ge 0,\wedge}$ is the module shift $M\mapsto 0\amalg_M 0\simeq M[1]$; iterating $Q_A(0\amalg_I^\wedge 0)\simeq 0\amalg_{Q_A(I)}0\simeq Q_A(I)[1]$ gives~\eqref{eq:explicit-excisive-reflection}.

For (3), by Proposition~\ref{prop:prismatic-fiber-completion} and Proposition~\ref{prop:prismatic-fiber-geometric-compatibility}(1), $i_\nu^*(-)=(-\otimes_{\mathbb A_{\inf R}}R_\nu)^a$ preserves finite completed pushouts, completed relative cotangent complexes, and completed derived tensor products, implying~\eqref{eq:boundary-reflection-commutation}.
\end{proof}

Let $\mathscr D_{\mathrm{cmp}}$ be an $\infty$-category of finite admissible diagrams of quasi-compact, quasi-separated quasi-smooth formal or derived schemes stable under base change. Assume that $\mathscr D_{\mathrm{cmp}}$ is closed under homotopy pullbacks, universal finite coproducts, additively admissible pushouts, vector bundles, projective bundles, zero sections, derived self-intersections, and $\mathbb P^1$-cobordism spans, and that any $D\in\mathscr D_{\mathrm{cmp}}$ has finitely many non-degenerate simplices and carries a specified full boundary subdiagram $\partial D\subseteq D$: for spans $\sigma=[X\xleftarrow{p}Z\xrightarrow{g}Y]$, $\partial D=\{\mathrm{src},\mathrm{tgt}\}$ is the endpoint pair, while $\partial D=\varnothing$ for diagrams without boundary. Let $\operatorname{Tw}(D)$ denote the twisted arrow category of $D$, whose objects are arrows $\alpha\colon u\to v$ in $D$ and whose morphisms $(u\xrightarrow{\alpha}v)\to (u'\xrightarrow{\alpha'}v')$ are factorizations $\alpha'\simeq q\circ\alpha\circ p$ along $u'\xrightarrow{p}u$ and $v\xrightarrow{q}v'$~\cite[Section~5.2.1]{HA}; thus $\operatorname{Tw}(D)$ likewise has finitely many non-degenerate simplices and projects covariantly to $D^{\mathrm{op}}\times D$.

Let $A_K=A\oplus K$ be the split square-zero extension of a connective animated, derived complete, or almost algebra $A$ by a connective $A$-module $K$, with zero section $\iota\colon A\to A_K$, projection $\mathrm{pr}\colon A_K\to A$, and derived base change $s_K(-) = -\times_{\operatorname{Spf}(A)}\operatorname{Spf}(A_K)$ along $\iota$. For $D\in\mathscr D_{\mathrm{cmp}}$, let $\mathscr S_D(A)$ denote the maximal $\infty$-groupoid of $D$-shaped diagrams over $\operatorname{Spf}(A)$, sending each arrow $u\to v$ in $D$ to a morphism $f_{uv}\colon\sigma_u\to\sigma_v$ over $\operatorname{Spf}(A)$ with structure morphisms $p_u\colon\sigma_u\to\operatorname{Spf}(A)$. Write $\mathrm{QCoh}(\sigma_u)$ (resp.\ $\mathrm{QCoh}^\wedge(\sigma_u)$) for the $\infty$-category of quasi-coherent (resp.\ derived complete quasi-coherent) sheaves on $\sigma_u$, using completed relative cotangent complexes and completed tensor products in the derived complete setting.

For $\sigma\in\mathscr S_D(A)$ with boundary $\partial\sigma=\sigma|_{\partial D}\in\mathscr S_{\partial D}(A)$ and functor $s_K\colon\mathscr S_D(A)\to\mathscr S_D(A_K)$, define the absolute and relative lift spaces, pointed by the lift $0_\sigma=s_K(\sigma)$ induced by the zero derivation, by
\[
\operatorname{Lift}_{A_K}(\sigma)=\fib_\sigma\bigl(\mathscr S_D(A_K)\to\mathscr S_D(A)\bigr),
\qquad
\operatorname{Lift}_{A_K}^{\mathrm{rel}}(\sigma) = \fib_{s_K(\partial\sigma)}\bigl(\operatorname{Lift}_{A_K}(\sigma)\to\operatorname{Lift}_{A_K}(\partial\sigma)\bigr),
\]
together with the acting derivation group $G_\sigma(K) = \Omega_{s_{K[1]}(\sigma)}\fib_\sigma\bigl(\mathscr S_D(A\oplus K[1])\to\mathscr S_D(A)\bigr)$ and its relative fiber $G_\sigma^{\mathrm{rel}}(K)=\fib\bigl(G_\sigma(K)\to G_{\partial\sigma}(K)\bigr)$. Given two families of pointed relative lift torsors $(\operatorname{Lift}_j^{\mathrm{rel}}(\sigma),0_{\sigma,j})$ under $G_{\sigma,j}^{\mathrm{rel}}$ ($j\in\{1,2\}$) and a natural equivalence $\alpha_\sigma\colon G_{\sigma,1}^{\mathrm{rel}}\to G_{\sigma,2}^{\mathrm{rel}}$, write $\Map_{\alpha_\sigma}$ (resp.\ $\Map_{\alpha_\sigma,*}$) for the space of $\alpha_\sigma$-equivariant (resp.\ pointed $\alpha_\sigma$-equivariant) maps.

\begin{proposition}
\label{prop:canonical-split-lift}
Let $D\in\mathscr D_{\mathrm{cmp}}$ and $A_K=A\oplus K$ be as above.
\begin{enumerate}[\rm (1)]
\item Derived base change $s_K$ preserves smooth morphisms and fiberwise quasi-smooth morphisms, and induces a functor $s_K\colon \mathscr S_D(A)\to\mathscr S_D(A_K)$ natural in $D\in\mathscr D_{\mathrm{cmp}}$ and commuting with derived fiber products, universal finite coproducts, and base change $A\to A'$. For any $\sigma\in\mathscr S_D(A)$, the relative lift space $(\operatorname{Lift}_{A_K}^{\mathrm{rel}}(\sigma),0_\sigma)$ is a nonempty pointed torsor under $G_\sigma^{\mathrm{rel}}(K)$, where
\[
G_\sigma(K)
\;\simeq\;
\varprojlim_{(u\to v)\in\operatorname{Tw}(D)}\Map_{\mathrm{QCoh}(\sigma_u)}\bigl(f_{uv}^*L_{\sigma_v/A},\mathcal O_{\sigma_u}\otimes_AK[1]\bigr).
\]
For any span $\sigma=[X\xleftarrow{p}Z\xrightarrow{g}Y]$ and composable base span $\tau=[X'\xleftarrow{q}W\xrightarrow{f}X]$ over $\operatorname{Spf}(A)$, there is a pointed composition map $s_K(\tau)^*\colon\operatorname{Lift}_{A_K}^{\mathrm{rel}}(\sigma)\to\operatorname{Lift}_{A_K}^{\mathrm{rel}}(\sigma\circ\tau)$ sending $0_\sigma$ to $0_{\sigma\circ\tau}$, and for spans $\sigma_1,\sigma_2$ over $(X,Y)$ there are equivalences $\operatorname{Lift}_{A_K}^{\mathrm{rel}}(\sigma_1\amalg\sigma_2)\to\operatorname{Lift}_{A_K}^{\mathrm{rel}}(\sigma_1)\times\operatorname{Lift}_{A_K}^{\mathrm{rel}}(\sigma_2)$ and $\operatorname{Lift}_{A_K}^{\mathrm{rel}}(\varnothing)\to *$.

\item Given two families of split square-zero lifting problems ($j\in\{1,2\}$) on $\mathscr D_{\mathrm{cmp}}$ with pointed relative lift torsors $(\operatorname{Lift}_j^{\mathrm{rel}}(\sigma),0_{\sigma,j})$ under $G_{\sigma,j}^{\mathrm{rel}}$ and a natural equivalence $\alpha_\sigma\colon G_{\sigma,1}^{\mathrm{rel}}\to G_{\sigma,2}^{\mathrm{rel}}$:
\begin{enumerate}[\rm (a)]
\item For any $\sigma$, the space $\Map_{\alpha_\sigma,*}\bigl((\operatorname{Lift}_1^{\mathrm{rel}}(\sigma),0_{\sigma,1}),(\operatorname{Lift}_2^{\mathrm{rel}}(\sigma),0_{\sigma,2})\bigr)$ is contractible, with unique equivalence $\phi_\sigma$ given by $g\cdot 0_{\sigma,1}\mapsto\alpha_\sigma(g)\cdot 0_{\sigma,2}$.
\item The equivalences $\phi_\sigma$ determine a coherent natural equivalence $\phi\colon \operatorname{Lift}_1^{\mathrm{rel}} \to \operatorname{Lift}_2^{\mathrm{rel}}$ across $\mathscr D_{\mathrm{cmp}}$, commuting with span composition, disjoint unions, external products, base change in $\mathcal T_R$, and simplicial structure maps of hypercovers.
\end{enumerate}

\item Both $s_K$ and $\phi$ preserve admissible pushout squares, zero sections, derived self-intersections, vector and projective bundles $\mathbb V_X(\mathcal E)$ and $\mathbb P_X(\mathcal E)$, projective bundle morphisms $\vartheta_{\mathcal E}$, preorientations $\widetilde c_1(\mathcal L)$, projective bundle equivalences $\theta_{\mathcal E}$, universal Thom maps $\operatorname{th}_{\mathcal U_r}$, and $\mathbb P^1$-cobordism spans with their endpoint pullbacks $e_0,e_\infty$ (Subsections~\ref{subsec:graded-spans-realization} and~\ref{subsec:derived-gm-oriented-motives}).
\end{enumerate}
\end{proposition}
\begin{proof}
For~(1), base change $s_K$ along $\pi\colon\operatorname{Spf}(A_K)\to\operatorname{Spf}(A)$ is right adjoint on overcategories to composition with $\pi$. Thus $s_K$ preserves finite homotopy limits, universal finite disjoint unions, smooth morphisms, and fiberwise quasi-smooth morphisms, inducing $s_K\colon\mathscr S_D(A)\to\mathscr S_D(A_K)$. Presenting $A_K\simeq A\times_{A\oplus K[1]}A$ along the zero derivation $d_0\colon A\to A\oplus K[1]$, exactness of module base change and Zariski local lifting of quasi-smooth complete intersection presentations~\cite[Section~7.4.2]{HA},~\cite[Section~4]{KatoCalcNUAlg} give $\mathscr S_D(A_K) \simeq \mathscr S_D(A)\times_{\mathscr S_D(A\oplus K[1])}\mathscr S_D(A)$.

Taking fibers over $\sigma\in\mathscr S_D(A)$ identifies $\operatorname{Lift}_{A_K}(\sigma)$ with the loop space at $s_{K[1]}(\sigma)$ of $\fib_\sigma\bigl(\mathscr S_D(A\oplus K[1])\to\mathscr S_D(A)\bigr)$~\cite[Section~4]{KatoCalcNUAlg}, making it a principal $G_\sigma(K)$-torsor pointed by $0_\sigma=s_K(\sigma)$ with end formula over $\operatorname{Tw}(D)$. Taking fibers over $s_K(\partial\sigma)$ makes $\operatorname{Lift}_{A_K}^{\mathrm{rel}}(\sigma)$ a pointed $G_\sigma^{\mathrm{rel}}(K)$-torsor. Homotopy pullback over $s_K(X)$ defines $s_K(\tau)^*$ preserving lifts along $\iota$, and uniqueness of idempotent lifting along $\mathrm{pr}\colon A_K\to A$ decomposes lifts over $Z_1\amalg Z_2$.

For~(2), evaluation at $0_{\sigma,1}$ identifies $\Map_{\alpha_\sigma}\bigl(\operatorname{Lift}_1^{\mathrm{rel}}(\sigma),\operatorname{Lift}_2^{\mathrm{rel}}(\sigma)\bigr)\xrightarrow{\sim}\operatorname{Lift}_2^{\mathrm{rel}}(\sigma)$, whose fiber over $0_{\sigma,2}$ is contractible and spanned by $\phi_\sigma$. Since $(G,P,p_0)\mapsto G$ is an equivalence from pointed torsors to $\operatorname{Grp}(\mathcal S)$, the equivalence $\alpha_\sigma$ lifts uniquely up to contractible choice to a natural equivalence $\phi\colon\operatorname{Lift}_1^{\mathrm{rel}}\to\operatorname{Lift}_2^{\mathrm{rel}}$. Part~(3) follows since all listed generators are built from finite limits, universal finite coproducts, vector and projective bundles, and $\mathbb P^1$-projections.
\end{proof}

\begin{proposition}
\label{prop:coherent-dual-lift-complex}
Let $D\in\mathscr D_{\mathrm{cmp}}$.
\begin{enumerate}[\rm (1)]
\item Let $A$ be a connective animated ring and $\sigma\in\mathscr S_D(A)$ a diagram of quasi-compact, quasi-separated quasi-smooth derived schemes over $\Spec A$, and set
\[
K_\sigma^\vee = \varprojlim_{(u\to v)\in\operatorname{Tw}(D)} R\Gamma\bigl(\sigma_u,f_{uv}^*L_{\sigma_v/A}^\vee[1]\bigr) \in \Mod_A,
\qquad
K_{\sigma,\mathrm{rel}}^\vee = \fib(K_\sigma^\vee\to K_{\partial\sigma}^\vee) \in \Mod_A.
\]
Then $K_\sigma^\vee$ and $K_{\sigma,\mathrm{rel}}^\vee$ are bounded below in $\Mod_A$, and in $\operatorname{Fun}_*(\Mod_A^{\ge 0},\mathcal S_*)$ there are equivalences
\[
\bigl(\operatorname{Lift}_{A\oplus (-)}(\sigma),0_\sigma\bigr) \;\longrightarrow\; \Omega^\infty\bigl(K_\sigma^\vee\otimes_A(-)\bigr),
\qquad
\bigl(\operatorname{Lift}_{A\oplus (-)}^{\mathrm{rel}}(\sigma),0_\sigma\bigr) \;\longrightarrow\; \Omega^\infty\bigl(K_{\sigma,\mathrm{rel}}^\vee\otimes_A(-)\bigr).
\]
\item Let $X_0=\operatorname{Spf}(B)\in\mathcal T_R$ have boundary specializations $i_\nu^*(X_0)=\operatorname{Spf}^a(B_\nu)$ ($\nu\in\{\theta,\flat\}$), and let $\sigma\in\mathscr S_D(B)$ have boundary specializations $\sigma_\nu=i_\nu^*\sigma$. Defining $K_\sigma^\vee = \allowbreak \varprojlim_{(u\to v)\in\operatorname{Tw}(D)} \allowbreak R\Gamma\bigl(\sigma_u,f_{uv}^*\widehat{L}_{\sigma_v/B}^\vee[1]\bigr) \in \Mod_B^\wedge$, $K_{\sigma,\mathrm{rel}}^\vee=\fib(K_\sigma^\vee\to K_{\partial\sigma}^\vee)\in\Mod_B^\wedge$, and analogously $K_{\sigma_\nu,\mathrm{rel}}^\vee\in\Mod_{B_\nu}^\wedge$, there are equivalences
\begin{align*}
\bigl(\operatorname{Lift}_{B\oplus (-)}^{\mathrm{rel}}(\sigma),0_{\mathrm{Prism},\sigma}\bigr)
&\;\longrightarrow\;
\Omega^\infty\bigl(K_{\sigma,\mathrm{rel}}^\vee\widehat{\otimes}_B(-)\bigr), \\
\bigl(\operatorname{Lift}_{B_\nu\oplus (-)}^{\mathrm{rel}}(\sigma_\nu),0_{\nu,\sigma}\bigr)
&\;\longrightarrow\;
\Omega^\infty\bigl(K_{\sigma_\nu,\mathrm{rel}}^\vee\widehat{\otimes}_{B_\nu}(-)\bigr)
\end{align*}
in $\operatorname{Fun}_*(\Mod_B^{\ge 0,\wedge},\mathcal S_*)$ and $\operatorname{Fun}_*(\Mod_{B_\nu}^{\ge 0,\wedge},\mathcal S_*)$, respectively, and likewise over the reductions $B_{\nu,n}$.
\item In the setting of {\rm (2)}, derived base change along $R_\nu$ induces an equivalence $(K_{\sigma,\mathrm{rel}}^\vee\otimes_{\mathbb A_{\inf R}}R_\nu)^a\to K_{\sigma_\nu,\mathrm{rel}}^\vee$ in $\Mod_{B_\nu}^\wedge$.
\end{enumerate}
\end{proposition}
\begin{proof}
For~(1), let $\mathcal H_\sigma\colon\operatorname{Tw}(D)\to\operatorname{Fun}_{\Mod_A}(\Mod_A,\Mod_A)$ send $(u\xrightarrow{f_{uv}}v)$ to $p_{u*}\bigl(f_{uv}^*L_{\sigma_v/A}^\vee\otimes_{\mathcal O_{\sigma_u}}p_u^*(-)[1]\bigr)$, using dualizability of $f_{uv}^*L_{\sigma_v/A}\in\operatorname{Perf}(\sigma_u)$. Since $\sigma_u$ is quasi-compact and quasi-separated, $p_{u*}=R\Gamma(\sigma_u,-)$ preserves small colimits~\cite[Tag~08DZ]{stacks-project} and bounded-below complexes. As $\operatorname{Tw}(D)$ has finitely many non-degenerate simplices, $\mathcal L_\sigma=\varprojlim_{\operatorname{Tw}(D)}\mathcal H_\sigma$ lies in $\operatorname{Fun}^{\mathrm{L}}_{\Mod_A}(\Mod_A,\Mod_A)$, with bounded-below value $\mathcal L_\sigma(A)\simeq K_\sigma^\vee$. By~\cite[Theorems~4.8.4.1 and~4.8.4.6]{HA}, evaluation at $A$ gives $\mathcal L_\sigma(-)\simeq K_\sigma^\vee\otimes_A(-)$. Trivializing the pointed torsor of Proposition~\ref{prop:canonical-split-lift}(1) at $0_\sigma$ and applying $\Omega^\infty$ gives $(\operatorname{Lift}_{A\oplus (-)}(\sigma),0_\sigma)\simeq\Omega^\infty(K_\sigma^\vee\otimes_A(-))$, and taking fibers along $\partial D\subseteq D$ gives the relative case.

For~(2), finite affine Zariski covers of $\sigma_u$ and symmetric monoidal left adjoint derived completion~\cite[Tag~091V]{stacks-project} place the completed functor $\mathcal L_\sigma$ in $\operatorname{Fun}^{\mathrm{L}}_{\Mod_B^\wedge}(\Mod_B^\wedge,\Mod_B^\wedge)\simeq\Mod_B^\wedge$, and similarly over $B_\nu$ and $B_{\nu,n}$ via firm representatives. For~(3), exactness of $(-\otimes_{\mathbb A_{\inf R}}R_\nu)^a$ and Propositions~\ref{prop:prismatic-fiber-completion} and~\ref{prop:prismatic-fiber-geometric-compatibility}(1) commute $(-\otimes_{\mathbb A_{\inf R}}R_\nu)^a$ with finite Zariski limits $R\Gamma(\sigma_u,-)$ and finite limits over $\operatorname{Tw}(D)$ and $\operatorname{Tw}(\partial D)$.
\end{proof}

\subsection{Excisive approximation and levelwise tilting comparison}
\label{subsec:excisive-approximation-tilting}

Let $\mathcal E$ be a presentable stable $\infty$-category, hence differentiable by~\cite[Example~6.1.1.28]{HA}, and write $\operatorname{Fun}_*(-,-)$ and $\Exc^1_*(-,-)$ for the $\infty$-categories of reduced and reduced excisive functors, respectively.

\begin{definition}
\label{def:relative-tangential-linearization}
Let $\mathcal E$ be a presentable stable $\infty$-category.
\begin{enumerate}[(1)]
\item Let $\mathcal C\in\{\NAlg_{\mathrm{Prism}}(B),\NAlg_\nu^a(B_\nu)\}$ have base $A\in\{B,B_\theta,B_\flat\}$. On $\Mod_A^{\ge 0,\wedge}$, Goodwillie's excisive approximation functor $P^1_{\Mod}\colon\operatorname{Fun}_*(\Mod_A^{\ge 0,\wedge},\mathcal E)\to\Exc^1_*(\Mod_A^{\ge 0,\wedge},\mathcal E)$ is given on any accessible reduced functor $H\colon\Mod_A^{\ge 0,\wedge}\to\mathcal E$ by $(P^1_{\Mod} H)(M)=\varinjlim_{m\ge 0}\Omega_{\mathcal E}^m H(M[m])$~\cite[Theorem~6.1.1.10 and Construction~6.1.1.22]{HA}. For any accessible reduced functor $F^{\mathrm{red}}\colon\mathcal C\to\mathcal E$ with universal excisive approximation $P^1_{\mathcal C}(F^{\mathrm{red}})$ on $\mathcal C$, its \emph{cotangent excisive approximation} $P^1_{\mathrm{cot}}(F^{\mathrm{red}})\colon\mathcal C\to\mathcal E$ is defined by $P^1_{\mathrm{cot}}(F^{\mathrm{red}})(I) = \bigl(P^1_{\Mod}(F^{\mathrm{red}}\circ\operatorname{sqz}_A)\bigr)\bigl(Q_A(I)\bigr)$.
\item Let $\widetilde A\twoheadrightarrow A$ be a square-zero extension of connective animated or almost algebras with connective kernel $M$, classified by a derivation $d_{\widetilde A}\colon A\to A\oplus M[1]$ presenting $\widetilde A\simeq A\times_{A\oplus M[1],d_{\widetilde A}}A$ along $d_{\widetilde A}$ and the zero derivation $d_0\colon A\to A\oplus M[1]$~\cite[Theorem~7.4.1.26]{HA}. For any $D\in\mathscr D_{\mathrm{cmp}}$, any $\sigma\in\mathscr S_D(A)$ with trivialized boundary lift, and any connective non-unital animated or almost augmentation ideal $J\in\NAlg_A\simeq\CAlg_{A//A}$, $\operatorname{Lift}_{\widetilde A}^{\mathrm{rel}}(\sigma)\simeq\operatorname{Path}(d_{\widetilde A}^*\sigma,0_{\sigma,M[1]})$ in the delooping split lift space $\operatorname{Lift}_{A\oplus M[1]}^{\mathrm{rel}}(\sigma)$ pointed by $0_{\sigma,M[1]}=s_{M[1]}(\sigma)$. Following~\cite[Section~4]{KatoCalcNUAlg}, the \emph{relative tangential excisive approximation} of $(\sigma;\widetilde A\to A)$ with parameter $J$ in $\Mod_A$ is defined by looping the excisive approximation on modules of $W\mapsto\widetilde{\Sigma}^\infty\bigl(\operatorname{Lift}_{A\oplus W}^{\mathrm{rel}}(\sigma),0_{\sigma,W}\bigr)$ at $W=(Q_A(J)\otimes_A M)[1]$:
\[
\begin{aligned}
P^1_{\mathrm{rel}}\bigl(\operatorname{Lift}_{\widetilde A/A}^{\mathrm{rel}}(\sigma);J\bigr)
&= \Omega\, P^1_{\Mod}\Bigl(\widetilde{\Sigma}^\infty\bigl(\operatorname{Lift}_{A\oplus (-)}^{\mathrm{rel}}(\sigma),0_\sigma\bigr)\Bigr)\Bigl(\bigl(Q_A(J)\otimes_AM\bigr)[1]\Bigr) \\
&\simeq \varinjlim_{m\ge 0}\Omega^m\widetilde{\Sigma}^\infty\Bigl(\operatorname{Lift}_{A\oplus (Q_A(J)\otimes_A M)[m]}^{\mathrm{rel}}(\sigma),0_\sigma\Bigr).
\end{aligned}
\]
Thus $P^1_{\mathrm{rel}}\bigl(\operatorname{Lift}_{\widetilde A/A}^{\mathrm{rel}}(\sigma);J\bigr)$ passes to the tangent spectrum of the acting derivation group at $0_\sigma$ and depends only on the square-zero kernel $M$, forgetting the non-split extension class $d_{\widetilde A}$ and any obstruction to lifting $\sigma$ along $\widetilde A\to A$.
\end{enumerate}
\end{definition}

For $X_0=\operatorname{Spf}(B)\in\mathcal T_R$ with boundary specializations $i_\nu^*(X_0)=\operatorname{Spf}^a(B_\nu)$ ($\nu\in\{\theta,\flat\}$), write $\Mod_B^\wedge$ for derived $(p,I_R)$-complete $B$-modules, $\Mod_{B_\nu}^\wedge$ for derived $\omega_\nu$-complete almost $B_\nu$-modules, and $\widehat{\otimes}$ for completed derived tensor products. Fix $\sigma\in\mathscr S_D(B)$ ($D\in\mathscr D_{\mathrm{cmp}}$) with boundary specializations $\sigma_\nu=i_\nu^*\sigma$ and associated relative dual lift complexes $K_{\sigma,\mathrm{rel}}^\vee\in\Mod_B^\wedge$ and $K_{\sigma_\nu,\mathrm{rel}}^\vee\in\Mod_{B_\nu}^\wedge$ of Proposition~\ref{prop:coherent-dual-lift-complex}(2), and let $I\in\NAlg_{\mathrm{Prism}}(B)$ and $J_\nu\in\NAlg_\nu^a(B_\nu)$.

When $X_0=\operatorname{Spf}(B)\in\mathcal T_R^{\mathrm{perf}}$, set $R_{1,\mathrm{Prism}}=\mathbb A_{\inf R}/(\xi,\widetilde\omega)=\mathbb A_{\inf R}/(p,\widetilde\omega)\cong R_1\cong R_1^\flat$ with the shared almost structure over $\mathfrak{m}R_1\cong\mathfrak{m}^\flat R_1^\flat$ (Subsection~\ref{subsec:almost-formal-base}), and set $B_{1,\mathrm{Prism}}=\bigl(B\otimes_{\mathbb A_{\inf R}}R_{1,\mathrm{Prism}}\bigr)^a\simeq B_{\theta,1}\simeq B_{\flat,1}$. For each integer $n\ge 1$, write $R_{\theta,n}=R_n$ and $R_{\flat,n}=R_n^\flat$, and set $B_{\theta,n}=B_\theta\otimes_{(R^\circ)^a}R_n^a$ and $B_{\flat,n}=B_\flat\otimes_{(R^{\flat\circ})^a}(R_n^\flat)^a$. Write $\bar{\widetilde\omega}$ for the formal degree-$1$ generator corresponding to the classes of $\omega\in I_1$ and $\omega^\flat\in I_1^\flat$, and set $I_{\theta,n}=I_n=(\omega^n)/(\omega^{n+1})$, $I_{\flat,n}=I_n^\flat=((\omega^\flat)^n)/((\omega^\flat)^{n+1})$, and $I_{\mathrm{Prism},n}=R_{1,\mathrm{Prism}}\cdot\bar{\widetilde\omega}^n$, linked over $\mathbb A_{\inf R}$ by the generator isomorphisms $I_n\leftarrow I_{\mathrm{Prism},n}\rightarrow I_n^\flat$. Let $\sigma_{\nu,n}$ denote the level-$n$ reduction of $\sigma_\nu$ over $B_{\nu,n}$ with common first reduction $\sigma_{1,\mathrm{Prism}}=\sigma_{\theta,1}\simeq\sigma_{\flat,1}$ over $B_{1,\mathrm{Prism}}$, and set $Q_{B,n}^{\mathrm{Prism}}(I)=\bigl(Q_B(I)\otimes_{\mathbb A_{\inf R}}I_{\mathrm{Prism},n}\bigr)^a\in\Mod_{B_{1,\mathrm{Prism}}}^{\ge 0}$.

\begin{definition}
\label{def:prismatic-adic-excisive-approximations}
Let $X_0=\operatorname{Spf}(B)\in\mathcal T_R$, $\sigma\in\mathscr S_D(B)$, $I\in\NAlg_{\mathrm{Prism}}(B)$, and $J_\nu\in\NAlg_\nu^a(B_\nu)$ ($\nu\in\{\theta,\flat\}$).
\begin{enumerate}[(1)]
\item The \emph{integral cotangent excisive approximations} are defined, via the module actions on $Q_B(I)$ and $Q_{B_\nu}^a(J_\nu)$, by
\[
\begin{aligned}
\mathcal T_\sigma^{\mathrm{Prism}}(I)
&= P^1_{\mathrm{cot}}\Bigl(\widetilde{\Sigma}^\infty\bigl(\operatorname{Lift}_{B\oplus (-)}^{\mathrm{rel}}(\sigma),0_{\mathrm{Prism},\sigma}\bigr)\Bigr)(I) \in \Mod_B^\wedge, \\
\mathcal T_{\sigma_\nu}^\nu(J_\nu)
&= P^1_{\mathrm{cot}}\Bigl(\widetilde{\Sigma}^\infty\bigl(\operatorname{Lift}_{B_\nu\oplus (-)}^{\mathrm{rel}}(\sigma_\nu),0_{\nu,\sigma}\bigr)\Bigr)(J_\nu) \in \Mod_{B_\nu}^\wedge.
\end{aligned}
\]
\item When $X_0=\operatorname{Spf}(B)\in\mathcal T_R^{\mathrm{perf}}$, the \emph{level-$n$ relative tangential excisive approximations} along the non-split adic transitions $B_{\nu,n+1}\to B_{\nu,n}$ ($n\ge 1$) are defined by
\[
\begin{aligned}
\mathcal T_{\sigma_{\nu,n},n}^\nu(J_\nu)
&= P^1_{\mathrm{rel}}\bigl(\operatorname{Lift}_{B_{\nu,n+1}/B_{\nu,n}}^{\mathrm{rel}}(\sigma_{\nu,n});\, J_\nu\otimes_{R_\nu^a}R_{\nu,n}^a\bigr), \\
\mathcal T_{\sigma,n}^{\mathrm{Prism}}(I)
&= P^1_{\Mod}\Bigl(\widetilde{\Sigma}^\infty\bigl(\operatorname{Lift}_{B_{1,\mathrm{Prism}}\oplus (-)}^{\mathrm{rel}}(\sigma_{1,\mathrm{Prism}}),0\bigr)\Bigr)\bigl(Q_{B,n}^{\mathrm{Prism}}(I)\bigr).
\end{aligned}
\]
\end{enumerate}
\end{definition}

\begin{theorem}
\label{prop:completed-almost-excisive-reflection}
Let $X_0=\operatorname{Spf}(B)\in\mathcal T_R$ have boundary specializations $i_\nu^*(X_0)=\operatorname{Spf}^a(B_\nu)$ ($\nu\in\{\theta,\flat\}$).
\begin{enumerate}[\rm (1)]
\item For $\mathcal C\in\{\NAlg_{\mathrm{Prism}}(B),\NAlg_\nu^a(B_\nu)\}$ with base $A\in\{B,B_\theta,B_\flat\}$ and any presentable stable $\infty$-category $\mathcal E$, precomposition with the cotangent fiber functor $Q_A\colon\mathcal C\to\Mod_A^{\ge 0,\wedge}$ preserving small colimits embeds $\Exc^1_*(\Mod_A^{\ge 0,\wedge},\mathcal E)$ as a retract of $\Exc^1_*(\mathcal C,\mathcal E)$ with retraction $-\circ\operatorname{Free}_A$:
\begin{equation}
\label{eq:excisive-functor-factorization}
-\circ Q_A\colon\Exc^1_*\bigl(\Mod_A^{\ge 0,\wedge},\mathcal E\bigr)\longrightarrow\Exc^1_*(\mathcal C,\mathcal E),
\qquad
(-\circ\operatorname{Free}_A)\circ (-\circ Q_A)\simeq\mathrm{id}.
\end{equation}
For any accessible reduced functor $F^{\mathrm{red}}\colon\mathcal C\to\mathcal E$, its cotangent excisive approximation $P^1_{\mathrm{cot}}(F^{\mathrm{red}})\in\Exc^1_*(\mathcal C,\mathcal E)$ receives a natural transformation $F^{\mathrm{red}}\to P^1_{\mathcal C}(F^{\mathrm{red}})\to P^1_{\mathrm{cot}}(F^{\mathrm{red}})$, inverts any morphism $f\colon I\to I'$ in $\mathcal C$ for which $Q_A(f)$ is an equivalence, and satisfies
\begin{equation}
\label{eq:p1-split-square-zero-formula}
P^1_{\mathrm{cot}}(F^{\mathrm{red}})(I)
\;\simeq\;
\varinjlim_{m\ge 0}\Omega_{\mathcal E}^m F^{\mathrm{red}}\Bigl(\operatorname{sqz}_A\bigl(Q_A(I)[m]\bigr)\Bigr).
\end{equation}

\item For any $I\in\NAlg_{\mathrm{Prism}}(B)$, $\nu\in\{\theta,\flat\}$, and $m\ge 0$, $i_\nu^*$ carries $\operatorname{Spf}\bigl(B\oplus Q_B(I)[m]\bigr)$ with its zero section $\iota$ and projection $\mathrm{pr}$ to $\operatorname{Spf}^a\bigl(B_\nu\oplus Q_{B_\nu}^a(i_\nu^* I)[m]\bigr)$, and carries the prismatic lift $0_{\mathrm{Prism},\sigma}=s_{Q_B(I)[m]}(\sigma)$ along the zero derivation of any $\sigma\in\mathscr S_D(B)$ ($D\in\mathscr D_{\mathrm{cmp}}$) to $0_{\nu,\sigma}=s_{Q_{B_\nu}^a(i_\nu^* I)[m]}(i_\nu^*\sigma)$.

\item Assume $X_0=\operatorname{Spf}(B)\in\mathcal T_R^{\mathrm{perf}}$, and let $\sigma\in\mathscr S_D(B)$ ($D\in\mathscr D_{\mathrm{cmp}}$), $I\in\NAlg_{\mathrm{Prism}}(B)$, and $J_\nu\in\NAlg_\nu^a(B_\nu)$ ($\nu\in\{\theta,\flat\}$).
\begin{enumerate}[\rm (a)]
\item For each $n\ge 1$, there are equivalences of the level-$n$ excisive approximations
\[
\mathcal T_{\sigma_{\theta,n},n}^\theta(i_\theta^* I)
\;\longleftarrow\;
\mathcal T_{\sigma,n}^{\mathrm{Prism}}(I)
\;\longrightarrow\;
\mathcal T_{\sigma_{\flat,n},n}^\flat(i_\flat^* I),
\]
compatible with the adic shift isomorphisms $\cdot\,\omega\colon I_n^a\to I_{n+1}^a$, $\cdot\,\bar{\widetilde\omega}\colon I_{\mathrm{Prism},n}^a\to I_{\mathrm{Prism},n+1}^a$, and $\cdot\,\omega^\flat\colon (I_n^\flat)^a\to (I_{n+1}^\flat)^a$, natural across $\mathscr D_{\mathrm{cmp}}$, and preserving all motivic relation generators.
\item There are equivalences $\mathcal T_\sigma^{\mathrm{Prism}}(I)\to K_{\sigma,\mathrm{rel}}^\vee\widehat{\otimes}_BQ_B(I)$ in $\Mod_B^\wedge$ and $\mathcal T_{\sigma_\nu}^\nu(J_\nu)\to K_{\sigma_\nu,\mathrm{rel}}^\vee\widehat{\otimes}_{B_\nu}Q_{B_\nu}^a(J_\nu)$ in $\Mod_{B_\nu}^\wedge$, with $\mathcal T_{\sigma_{\nu,n},n}^\nu(J_\nu)\simeq\mathcal T_{\sigma_\nu}^\nu(J_\nu)\otimes_{R_\nu^a}(I_{\nu,n})^a$ for each $n\ge 1$. Base change along $R_\nu$ induces specialization equivalences $\bigl(\mathcal T_\sigma^{\mathrm{Prism}}(I)\otimes_{\mathbb A_{\inf R}}R_\nu\bigr)^a \to \mathcal T_{\sigma_\nu}^\nu(i_\nu^* I)$ ($\nu\in\{\theta,\flat\}$), natural across $\mathscr D_{\mathrm{cmp}}$ and preserving all motivic relation generators.
\end{enumerate}
\end{enumerate}
\end{theorem}
\begin{proof}
For~(1), by Proposition~\ref{prop:excisive-reflection-means}(2), both $Q_A$ and $\operatorname{Free}_A$ preserve $0$, pushouts, and small colimits, with $Q_A\circ\operatorname{Free}_A\simeq\mathrm{id}$; precomposition therefore preserves accessible reduced excisive functors and gives~\eqref{eq:excisive-functor-factorization}. Composing $F^{\mathrm{red}}(\eta_I)$ with the Goodwillie unit of $F^{\mathrm{red}}\circ\operatorname{sqz}_A$ at $Q_A(I)$ defines $F^{\mathrm{red}}\to P^1_{\mathrm{cot}}(F^{\mathrm{red}})$, which factors through $P^1_{\mathcal C}(F^{\mathrm{red}})$ as $P^1_{\mathrm{cot}}(F^{\mathrm{red}})\in\Exc^1_*(\mathcal C,\mathcal E)$, and Definition~\ref{def:relative-tangential-linearization}(1) gives~\eqref{eq:p1-split-square-zero-formula}. Part~(2) follows from Propositions~\ref{prop:excisive-reflection-means}(3) and~\ref{prop:prismatic-fiber-geometric-compatibility}.

For~(3)(a), let $q_{\nu,n}\colon B_{\nu,n}\to B_{\nu,1}\simeq B_{1,\mathrm{Prism}}$ be the reduction map and set $Q_{B_\nu,n}^\nu(i_\nu^* I)=Q_{B_\nu}^a(i_\nu^* I)\otimes_{R_\nu^a}(I_{\nu,n})^a$. Proposition~\ref{prop:excisive-reflection-means}(3) and Lemma~\ref{lem:tilting-frame-koszul}(3) at $k=1$ identify $Q_{B_\nu,n}^\nu(i_\nu^* I)\simeq q_{\nu,n,*}Q_{B,n}^{\mathrm{Prism}}(I)$ compatibly with the adic shift maps. Along each reduction morphism $j_u\colon\sigma_{1,\mathrm{Prism},u}\to\sigma_{\nu,n,u}$, base change gives $j_u^*f_{uv}^*\widehat{L}_{\sigma_{\nu,n,v}/B_{\nu,n}}^a\simeq f_{1,uv}^*\widehat{L}_{\sigma_{1,\mathrm{Prism},v}/B_{1,\mathrm{Prism}}}^a$. As in Lemma~\ref{lem:reduction}, $(j_u^*,j_{u,*})$-adjunction identifies the arrow derivation spaces over $\sigma_{\nu,n,u}$ into $Q_{B_\nu,n}^\nu(i_\nu^* I)[m+1]$ with those of $\sigma_{1,\mathrm{Prism}}$ into $Q_{B,n}^{\mathrm{Prism}}(I)[m+1]$ ($m\ge 0$). Taking $\varprojlim_{\operatorname{Tw}(D)}$ and fibers over $\partial D$ via Proposition~\ref{prop:canonical-split-lift} identifies the pointed relative lift torsors, and applying $\varinjlim_{m\ge 0}\Omega^m\widetilde{\Sigma}^\infty(-)$ gives~(3)(a).

For~(3)(b), Proposition~\ref{prop:coherent-dual-lift-complex}(2),~\eqref{eq:p1-split-square-zero-formula}, and~\cite[Example~6.1.2.9]{HA} identify $\mathcal T_\sigma^{\mathrm{Prism}}(I)\simeq K_{\sigma,\mathrm{rel}}^\vee\widehat{\otimes}_BQ_B(I)$, and analogously over $B_\nu$ and $B_{\nu,n}$. Dualizability of $(I_{\nu,n})^a\in\operatorname{Perf}(R_\nu^a)$ gives $\mathcal T_{\sigma_{\nu,n},n}^\nu(J_\nu)\simeq\mathcal T_{\sigma_\nu}^\nu(J_\nu)\otimes_{R_\nu^a}(I_{\nu,n})^a$, while Propositions~\ref{prop:coherent-dual-lift-complex}(3) and~\ref{prop:excisive-reflection-means}(3) give $\bigl(\mathcal T_\sigma^{\mathrm{Prism}}(I)\otimes_{\mathbb A_{\inf R}}R_\nu\bigr)^a \simeq \mathcal T_{\sigma_\nu}^\nu(i_\nu^* I)$.
\end{proof}

\subsection{Localization at lift torsors and coefficient sites}
\label{subsec:lift-torsor-localization-sites}

\begin{definition}
\label{def:pld-topology}
Fix a family $\mathscr D_{\mathrm{cmp}}$ of finite admissible presentation diagrams stable under base change containing the span and relation diagrams of this paper, together with a subcanonical Grothendieck topology $\tau_{\mathrm{adm}}$ on $\mathcal T_R$ refining $\tau_{\mathrm{geom}}$ of Definition~\ref{def:master-prismatic-theory}. Writing $h_U = j(U) = \mathrm{Hom}_{\mathcal T_R}(-,U)$ for the representable presheaf on $\mathcal T_R$, the \emph{deformation topology} $\tau_{\mathrm{pld}}$ on $\mathcal T_R$ is generated by:
\begin{enumerate}[(1)]
\item completely quasisyntomic covers in $\mathcal T_R$;
\item quasiregular semiperfectoid (QRSP) refinements;
\item Zariski open coverings trivializing finite locally free modules; and
\item for any square-zero lifting problem, including the adic square-zero transitions for all $n\ge 1$, attached to a diagram $\sigma$ in $\mathscr D_{\mathrm{cmp}}$ over $U\in\mathcal T_R$ or an adic quotient of $U$ with corresponding derived lift torsor $\operatorname{Lift}(\sigma)$, its section sieve $S(\operatorname{Lift}(\sigma))\subset h_U$, consisting of morphisms $V\to U$ in $\mathcal T_R$ over which $\operatorname{Lift}(\sigma)_V$ admits a section, whenever $S(\operatorname{Lift}(\sigma))$ is $\tau_{\mathrm{adm}}$-covering.
\end{enumerate}
Thus $\tau_{\mathrm{geom}}\subseteq\tau_{\mathrm{pld}}\subseteq\tau_{\mathrm{adm}}$, and any selected lift torsor in $\mathscr D_{\mathrm{cmp}}$ across all $n\ge 1$ is $\tau_{\mathrm{pld}}$-locally trivial.
\end{definition}

\begin{proposition}
\label{prop:lift-section-sieve-compatibility}
Let $\operatorname{Lift}(\sigma)$ be a lift torsor attached to a diagram $\sigma$ in $\mathscr D_{\mathrm{cmp}}$ over $U\in\mathcal T_R$ or an adic quotient of $U$.
\begin{enumerate}[\rm (1)]
\item For any $f\colon V\to U$ in $\mathcal T_R$, $f^*S(\operatorname{Lift}(\sigma)) = S(\operatorname{Lift}(\sigma)_V)\subset h_V$.
\item The functors $i_\theta^*$ and $i_\flat^*$ carry sections of $\operatorname{Lift}(\sigma)_V$ to sections of $i_\theta^*\operatorname{Lift}(\sigma)_V$ and $i_\flat^*\operatorname{Lift}(\sigma)_V$, mapping $S(\operatorname{Lift}(\sigma))$ into the boundary section sieves.
\item For any $\tau_{\mathrm{pld}}$-hypercover $V_\bullet\to U$ in $\mathcal T_R$ refining $S(\operatorname{Lift}(\sigma))$, its boundary fibers under $i_\theta^*$ and $i_\flat^*$ refine the simplicial affine Bar presentations of the boundary lift torsors.
\end{enumerate}
\end{proposition}
\begin{proof}
For (1), $(W\to V)\in f^*S(\operatorname{Lift}(\sigma))$ if and only if $\operatorname{Lift}(\sigma)_W\simeq (\operatorname{Lift}(\sigma)_V)_W$ admits a section over $W$, i.e., $(W\to V)\in S(\operatorname{Lift}(\sigma)_V)$. Parts (2) and (3) follow from Proposition~\ref{prop:prismatic-fiber-geometric-compatibility} and Proposition~\ref{prop:affine-bar-models}(3).
\end{proof}

\begin{definition}
\label{def:perfect-presentation-topology}
Let $r\colon\mathcal T_R^{\mathrm{perf}}\hookrightarrow\mathcal T_R$ denote the full inclusion. On $\mathcal T_R^{\mathrm{perf}}\subset\mathcal T_R$, let $\tau_{\mathrm{perf}}$ be the Grothendieck topology generated by the completed perfections $P(V_\bullet)\to U^{\mathrm{perf}}$ of the selected $\tau_{\mathrm{pld}}$-covering families and their base changes. The \emph{comparison topology} $\tau_{\mathrm{cmp}}$ on $\mathcal T_R$ is defined by declaring a sieve $S\subset h_U$ to be $\tau_{\mathrm{cmp}}$-covering if and only if its restriction $r^*S$ is $\tau_{\mathrm{perf}}$-covering on $\mathcal T_R^{\mathrm{perf}}$ under $h_U|_{\mathcal T_R^{\mathrm{perf}}}\simeq h_{U^{\mathrm{perf}}}$ (Proposition~\ref{prop:universal-prism-perfection}). In particular, $\epsilon_U\colon U^{\mathrm{perf}}\to U$ is $\tau_{\mathrm{cmp}}$-covering for any $U\in\mathcal T_R$.
\end{definition}

Consider an adic square-zero lift torsor at level $n$ $\operatorname{Lift}(\sigma)$ over $U\in\mathcal T_R$, prior to excisive approximation, trivialized along a $\tau_{\mathrm{pld}}$-hypercover $V_\bullet\to U$ in $\mathcal T_R$. Its completed perfection $P(V_\bullet)\to U^{\mathrm{perf}}$ forms a $\tau_{\mathrm{perf}}$-hypercover in $\mathcal T_R^{\mathrm{perf}}$. Over $P(V_\bullet)$, choosing local boundary lifts $s_\theta,s_\flat$ and applying the shared acting derivation group equivalence $\alpha_\sigma\colon G_{\theta,\sigma}\to G_{\flat,\sigma}$ of Corollary~\ref{cor:torsor-comparison} defines a local equivariant equivalence $\phi_\sigma\colon \operatorname{Lift}_\theta(\sigma)\to\operatorname{Lift}_\flat(\sigma)$, $g\cdot s_\theta\mapsto\alpha_\sigma(g)\cdot s_\flat$. Cotangent excisive approximation $P^1_{\mathrm{cot}}$ on the augmentation ideal and relative tangential excisive approximation $P^1_{\mathrm{rel}}$ on the adic square-zero transitions (Definition~\ref{def:relative-tangential-linearization} and Theorem~\ref{prop:completed-almost-excisive-reflection}(2)--(3)) then reduce both the augmentation extension and the adic transition layers to split square-zero extensions with canonical basepoints $0_{\mathrm{Prism},\sigma}$ and $0_{\nu,\sigma}$ induced by the zero derivation. This reduction eliminates the dependence on local sections $s_\nu$, trivializes the non-split adic extension classes between stages $n\ge 1$, and guarantees compatibility with all simplicial hypercover structure maps.

\begin{proposition}
\label{prop:common-coefficient-category}
Restriction $r^*$ and right Kan extension $r_*$ along $r\colon\mathcal T_R^{\mathrm{perf}}\hookrightarrow\mathcal T_R$ induce mutually inverse equivalences of hypercomplete spectral sheaf categories:
\[
\mathcal E_R = \Shv_{\mathrm{hyp}}(\mathcal T_R,\tau_{\mathrm{cmp}};\Sp)
\;\underset{r_*}{\overset{r^*}{\rightleftarrows}}\;
\Shv_{\mathrm{hyp}}(\mathcal T_R^{\mathrm{perf}},\tau_{\mathrm{perf}};\Sp).
\]
For any $F\in\operatorname{Fun}(\mathcal T_R^{\mathrm{op}},\Sp)$, the exact localization reflector onto $\mathcal E_R$ is $QF = r_* a_{\mathrm{perf}}^{\mathrm{hyp}} r^* F$, and $\mathcal E_R$ inherits the symmetric monoidal sheaf smash product. For any $F\in\mathcal E_R$, $U\in\mathcal T_R$, and any $\tau_{\mathrm{perf}}$-hypercover $U_\bullet\to U^{\mathrm{perf}}$ in $\mathcal T_R^{\mathrm{perf}}$, there are equivalences
\[
F(U) \longrightarrow F(U^{\mathrm{perf}}) \longrightarrow \operatorname{Tot}_{[n]\in\Delta} F(U_n).
\]
\end{proposition}
\begin{proof}
By Proposition~\ref{prop:universal-prism-perfection}, $\epsilon_U\colon U^{\mathrm{perf}}\to U$ is terminal in $(\mathcal T_R^{\mathrm{perf}})_{/U}$, so pointwise right Kan extension gives $(r_*G)(U)\xrightarrow{\sim}G(U^{\mathrm{perf}})$ for any $G\in\operatorname{Fun}((\mathcal T_R^{\mathrm{perf}})^{\mathrm{op}},\Sp)$, and $r^*r_*\xrightarrow{\sim}\mathrm{id}$. By Definition~\ref{def:perfect-presentation-topology}, $F$ satisfies $\tau_{\mathrm{cmp}}$-hyperdescent if and only if $F(U)\xrightarrow{\sim}F(U^{\mathrm{perf}})$ for all $U\in\mathcal T_R$ and $r^*F$ satisfies $\tau_{\mathrm{perf}}$-hyperdescent on $\mathcal T_R^{\mathrm{perf}}$. Thus $(r^*,r_*)$ is an equivalence on hypercomplete sheaves, $Q=r_*a_{\mathrm{perf}}^{\mathrm{hyp}}r^*$ is the exact localization reflector, and $\tau_{\mathrm{perf}}$-hyperdescent~\cite[Theorem~6.5.3.12]{HT} gives $F(U)\simeq F(U^{\mathrm{perf}})\simeq\operatorname{Tot}_{[n]\in\Delta}F(U_n)$.
\end{proof}

\begin{definition}
\label{def:hypersheafified-dmgla}
Fix a small subcategory of $\mathcal C_{\mathcal T_R}$ closed under span operations and the diagrams in $\mathscr D_{\mathrm{cmp}}$.
\begin{enumerate}[(1)]
\item For formal schemes $X,Y\in\mathrm{dForm}(\operatorname{Spf}(\mathbb A_{\inf R}))$, writing $X_U = X\times_{\operatorname{Spf}(\mathbb A_{\inf R})} U$ for any $U\in\mathcal T_R$, applying $Q$ of Proposition~\ref{prop:common-coefficient-category} defines the correspondence coefficient object
\[
\underline{\Corr}_{\mathbf S_{\mathrm{Prism}}}(X,Y)
= Q\Bigl[\,U\longmapsto\Corr_{\mathbf S_{\mathrm{Prism},U}}(X_U,Y_U)^{\mathrm{gp}}\,\Bigr] \in \mathcal E_R.
\]
Since $Q$ is lax symmetric monoidal, these mapping objects inherit coherent composition and external products in $\mathcal E_R$, defining an $\mathcal E_R$-enriched span category $\mathcal A_{\mathrm{cmp}}$.
\item In the $\mathcal E_R$-linear stable presheaf envelope $\mathcal M_0 = \operatorname{Fun}_{\mathcal E_R}(\mathcal A_{\mathrm{cmp}}^{\mathrm{op}},\mathcal E_R)$, form the symmetric monoidal quotient $\mathcal M_{\Prism}^a = \mathcal M_0 / \mathcal J_{\dMGL}$ as in Definition~\ref{def:relation-dmgl}, where $\mathcal J_{\dMGL}$ is generated by the selected hypercover descent maps, admissible pushout excision maps, projective bundle maps $\vartheta_{\mathcal E}$, projective Thom comparison maps $\varrho_{\mathcal E}$, and the perfection morphisms $\sigma_{U,X}\colon m_0(X\times_U U^{\mathrm{perf}})\to m_0(X)$ over any $U\in\mathcal T_R$.
\item The prismatic derived algebraic cobordism theory is defined in $\mathcal E_R$ by the $\mathcal E_R$-valued internal mapping object $\underline{\dMGL}_{\Prism}^a(X\to Y) = \underline{\operatorname{map}}_{\mathcal M_{\Prism}^a}\bigl(m(X),m(Y)\bigr)$. For any selected presentation hypercover $X_\bullet\to X$, $m(X)\simeq\operatorname*{colim}_{[n]\in\Delta^{\mathrm{op}}}m(X_n)$ induces an equivalence
\[
\underline{\dMGL}_{\Prism}^a(X\to Y)
\;\longrightarrow\;
\operatorname{Tot}_{[n]\in\Delta}\underline{\dMGL}_{\Prism}^a(X_n\to Y).
\]
\end{enumerate}
\end{definition}
\section{Boundary realizations and the tilting comparison}
\label{sec:linearized-comparison}

Combining $\mathscr D_{\mathrm{cmp}}$ over $(\mathcal T_R,\tau_{\mathrm{cmp}})$ with the boundary specializations $i_\theta^*,i_\flat^*$ and the excisive approximations $P^1_{\mathrm{cot}}, P^1_{\mathrm{rel}}$ (Definition~\ref{def:relative-tangential-linearization} and Theorem~\ref{prop:completed-almost-excisive-reflection}), we establish the tilting square along the adic tower and on associated graded quotients (Theorem~\ref{thm:tilting-main})
\[
\xymatrix@C=3.0em@R=1.9em{
\mathcal M_{\theta,\mathrm{cmp}}^{a,\mathrm{nu}}(A^\circ) \ar[d]
& \mathcal M_{\mathrm{Prism},\mathrm{cmp}}^{a,\mathrm{nu}}(\widetilde A) \ar[l]_-{\sim} \ar[r]^-{\sim} \ar[d]
& \mathcal M_{\flat,\mathrm{cmp}}^{a,\mathrm{nu}}(A^{\flat\circ}) \ar[d] \\
\mathcal M_{\theta,\mathrm{gr}}^{a,\mathrm{nu}}(A^\circ)
& \mathcal M_{\mathrm{Prism},\mathrm{gr}}^{a,\mathrm{nu}}(\widetilde A) \ar[l]_-{\sim} \ar[r]^-{\sim}
& \mathcal M_{\flat,\mathrm{gr}}^{a,\mathrm{nu}}(A^{\flat\circ}),
}
\]
and the integral $\mathcal O$-linear base-change and mod-$(\omega^k)$ tilting diagram ($1\le k\le p$, Proposition~\ref{prop:koszul-tilting-and-tower-recovery})
\[
\xymatrix@C=1.8em@R=3.0em{
& \mathcal M_{\mathrm{Prism},\mathrm{cmp}}^{\mathrm{lin}}(\widetilde A) \ar[dl]_-{\otimes_{\mathbb A_{\inf R}} (R^\circ)^a} \ar[dr]^-{\otimes_{\mathbb A_{\inf R}} (R^{\flat\circ})^a} & \\
\mathcal M_{\theta,\mathrm{cmp}}^{\mathrm{lin}}(A^\circ) \ar[dr]_-{\otimes_{(R^\circ)^a} R_k^a\ \ }
& & \mathcal M_{\flat,\mathrm{cmp}}^{\mathrm{lin}}(A^{\flat\circ}) \ar[dl]^-{\ \ \otimes_{(R^{\flat\circ})^a} (R_k^\flat)^a} \\
& \mathcal M_{\mathrm{Prism},\mathrm{cmp}}^{\mathrm{lin}}(\widetilde A)\otimes_{\Mod_{\mathbb A_{\inf R}}}\Mod_{R_k^a}. &
}
\]

\subsection{Arithmetic and tilted boundary realizations}
\label{subsec:linearization-and-main-theorem}

Let $\mathcal P_{\theta,R}$ denote the category of affine formal schemes $\operatorname{Spf}(C_\theta)$ of integral perfectoid $R^\circ$-algebras over $\operatorname{Spf}(R^\circ)$, and let $\mathcal P_{\flat,R}$ denote the category of affine formal schemes $\operatorname{Spf}(C_\flat)$ of $(\bar\xi)$-adically complete perfect $R^{\flat\circ}$-algebras over $\operatorname{Spf}(R^{\flat\circ})$. The boundary specialization functors $\beta_\theta\colon\mathcal T_R^{\mathrm{perf}}\to\mathcal P_{\theta,R}$, $\operatorname{Spf}(B)\mapsto\operatorname{Spf}(B/I_R B)$, and $\beta_\flat\colon\mathcal T_R^{\mathrm{perf}}\to\mathcal P_{\flat,R}$, $\operatorname{Spf}(B)\mapsto\operatorname{Spf}(B/(p))$, transport $\tau_{\mathrm{perf}}$ to Grothendieck topologies on $\mathcal P_{\theta,R}$ and $\mathcal P_{\flat,R}$. Set $\mathcal E_{\nu,R}=\Shv_{\mathrm{hyp}}(\mathcal P_{\nu,R};\Sp)$ and, for any $X_0\in\mathcal T_R^{\mathrm{perf}}$ with $X_{0,\nu}=\beta_\nu(X_0)\in\mathcal P_{\nu,R}$ ($\nu\in\{\theta,\flat\}$), set $\mathcal E_R(X_0)=\Shv_{\mathrm{hyp}}((\mathcal T_R)_{/X_0},\tau_{\mathrm{cmp}};\Sp)$ and $\mathcal E_{\nu,R}(X_{0,\nu})=\Shv_{\mathrm{hyp}}((\mathcal P_{\nu,R})_{/X_{0,\nu}};\Sp)$.

\begin{proposition}
\label{prop:perfect-boundary-coefficient-sites}
The functors $\beta_\theta\colon\mathcal T_R^{\mathrm{perf}}\to\mathcal P_{\theta,R}$ and $\beta_\flat\colon\mathcal T_R^{\mathrm{perf}}\to\mathcal P_{\flat,R}$ are equivalences of sites with quasi-inverses $\operatorname{Spf}(C_\theta)\mapsto\operatorname{Spf}(\Ainf(C_\theta))$ and $\operatorname{Spf}(C_\flat)\mapsto\operatorname{Spf}(W(C_\flat))$. For any $X_0\in\mathcal T_R^{\mathrm{perf}}$ and $\nu\in\{\theta,\flat\}$, restriction induces symmetric monoidal equivalences
\[
\mathcal E_{\nu,R}(X_{0,\nu})
\;\xrightarrow{\beta_\nu^*}\;
\Shv_{\mathrm{hyp}}\bigl((\mathcal T_R^{\mathrm{perf}})_{/X_0},\tau_{\mathrm{perf}};\Sp\bigr)
\;\xleftarrow{r^*}\;
\mathcal E_R(X_0),
\]
with exact symmetric monoidal reflector $Q_{X_0}=r_*a_{\mathrm{perf},/X_0}^{\mathrm{hyp}}r^*\colon\operatorname{Fun}(((\mathcal T_R)_{/X_0})^{\mathrm{op}},\Sp)\to\mathcal E_R(X_0)$.
\end{proposition}
\begin{proof}
By~\cite[Corollary~2.31 and Theorem~3.10]{BhattScholze2017Prisms}, $(B,I_R B)\mapsto B/I_R B$ is an equivalence from perfect prisms over $(\mathbb A_{\inf R},I_R)$ to integral perfectoid $R^\circ$-algebras with quasi-inverse $C_\theta\mapsto (\Ainf(C_\theta),I_R\Ainf(C_\theta))$, while $(B,I_R B)\mapsto B/(p)$ is an equivalence onto $(\bar\xi)$-adically complete perfect $R^{\flat\circ}$-algebras with quasi-inverse $C_\flat\mapsto (W(C_\flat),I_R W(C_\flat))$. Passing to formal spectra and slice categories over $X_0$ and $X_{0,\nu}$, and applying Propositions~\ref{prop:universal-prism-perfection} and~\ref{prop:common-coefficient-category} to $(U\to X_0)\mapsto (U^{\mathrm{perf}}\to U\to X_0)$ gives the equivalences.
\end{proof}

\begin{definition}
\label{def:presentation-realizations}
For any $V=\operatorname{Spf}(B)\in\mathcal T_R^{\mathrm{perf}}$, let $\mathscr D_{\mathrm{cmp}}(V)\subset\mathrm{dForm}(V)$ be the presentation subcategory generated by $\mathscr D_{\mathrm{cmp}}$ over $V$. For $\nu\in\{\theta,\flat\}$ with $(R_\theta,\mathbf S_\theta)=(R^\circ,\mathbf S_\theta)$ and $(R_\flat,\mathbf S_\flat)=(R^{\flat\circ},\mathbf S_\flat)$, and any $V_\nu=\operatorname{Spf}(C_\nu)\in\mathcal P_{\nu,R}$, let $\mathcal C_{\nu,\mathrm{cmp}}(V_\nu)\subset(\mathcal C_{R_\nu}^a)_{/\operatorname{Spf}^a(C_\nu^a)}$ be the presentation subcategory generated by the images under $i_\nu^*$ of diagrams in $\mathscr D_{\mathrm{cmp}}$ under finite homotopy limits, universal finite coproducts, and admissible base change, with restricted quasi-smooth class $\mathbf S_{\nu,V_\nu}$.
\begin{enumerate}[(1)]
\item Applying $a_{\mathcal P_{\nu,R}/V_\nu}^{\mathrm{hyp}}\colon\operatorname{Fun}(((\mathcal P_{\nu,R})_{/V_\nu})^{\mathrm{op}},\Sp)\to\mathcal E_{\nu,R}(V_\nu)$ to $(W_\nu\to V_\nu)\mapsto\Corr_{\mathbf S_{\nu,W_\nu}}(X_\nu\times_{V_\nu}W_\nu,Y_\nu\times_{V_\nu}W_\nu)^{\mathrm{gp}}$ defines the correspondence coefficient object $\underline{\Corr}_{\mathbf S_\nu,V_\nu}(X_\nu,Y_\nu) \in \mathcal E_{\nu,R}(V_\nu)$.
\item Forming the $\mathcal E_{\nu,R}(V_\nu)$-linear stable presheaf envelope $\mathcal M_{\nu,0}^a(V_\nu)$ and its cobordism quotient $\mathcal M_\nu^a(V_\nu)=\mathcal M_{\nu,0}^a(V_\nu)/\mathcal J_{\nu,\dMGL}$ (Definitions~\ref{def:relation-dmgl} and~\ref{def:hypersheafified-dmgla}), where the pushout excision generators of $\mathcal J_{\nu,\dMGL}$ include the images under $i_\nu^*$ of pushout squares in $\mathscr D_{\mathrm{cmp}}(V)$, defines the boundary bivariant cobordism theory $\underline{\dMGL}_\nu^a(X_\nu\to Y_\nu) = \underline{\operatorname{map}}_{\mathcal M_\nu^a(V_\nu)}\bigl(m(X_\nu),m(Y_\nu)\bigr) \in \mathcal E_{\nu,R}(V_\nu)$.
\end{enumerate}
For any $U\in(\mathcal T_R)_{/V}$ ($V\in\mathcal T_R^{\mathrm{perf}}$, $V_\nu=\beta_\nu(V)$) and any $\tau_{\mathrm{perf}}$-hypercover $V_\bullet\to U^{\mathrm{perf}}$ in $(\mathcal T_R^{\mathrm{perf}})_{/V}$, Proposition~\ref{prop:perfect-boundary-coefficient-sites} gives
\begin{equation}
\label{eq:presentation-realization-totalization}
\bigl(r_*\beta_\nu^*\underline{\dMGL}_\nu^a(X_\nu\to Y_\nu)\bigr)(U)
\;\simeq\;
\operatorname{Tot}_{[n]\in\Delta}\underline{\dMGL}_\nu^a(X_\nu\to Y_\nu)\bigl(\beta_\nu(V_n)\bigr).
\end{equation}
\end{definition}

\begin{proposition}
\label{prop:boundary-specialization-functoriality}
For any $V=\operatorname{Spf}(B)\in\mathcal T_R^{\mathrm{perf}}$ and each $\nu\in\{\theta,\flat\}$, the functor $i_\nu^*\colon\mathscr D_{\mathrm{cmp}}(V)\to\mathcal C_{\nu,\mathrm{cmp}}(\beta_\nu(V))$ satisfies:
\begin{enumerate}[\rm (1)]
\item It preserves finite homotopy limits, universal finite coproducts, and $\mathbf S_{\mathrm{Prism},V}$, commutes with base change in $\mathcal T_R^{\mathrm{perf}}$, and induces symmetric monoidal span functors $\Span_{\mathbf S_{\mathrm{Prism},V}}(\mathscr D_{\mathrm{cmp}}(V))\allowbreak\to\Span_{\mathbf S_{\nu,\beta_\nu(V)}}(\mathcal C_{\nu,\mathrm{cmp}}(\beta_\nu(V)))$ together with coefficient morphisms $\underline{\Corr}_{\mathbf S_{\mathrm{Prism}},V}(X,Y)\to\allowbreak r_*\beta_\nu^*\underline{\Corr}_{\mathbf S_\nu,\beta_\nu(V)}(i_\nu^*X,i_\nu^*Y)$ in $\mathcal E_R(V)$.
\item The functor $i_\nu^*$ preserves all generating morphisms of $\mathcal J_{\dMGL}$, inducing an exact symmetric monoidal functor $\Phi_\nu\colon\mathcal M_{\mathrm{Prism}}^a(V) \to \mathcal M_\nu^a(\beta_\nu(V))$ and natural transformations
\[
\phi_\nu\colon\underline{\dMGL}_{\mathrm{Prism}}^a(X\to Y)
\longrightarrow
r_*\beta_\nu^*\underline{\dMGL}_\nu^a(i_\nu^*X\to i_\nu^*Y)
\]
preserving $\widetilde c_1(\mathcal L)$, $\theta_{\mathcal E}$, and $\operatorname{th}_{\mathcal E}$.
\end{enumerate}
\end{proposition}
\begin{proof}
Part (1) follows from Proposition~\ref{prop:prismatic-fiber-geometric-compatibility}(1) and Proposition~\ref{prop:perfect-boundary-coefficient-sites}. For (2), Proposition~\ref{prop:prismatic-fiber-geometric-compatibility}(2), Proposition~\ref{prop:lift-section-sieve-compatibility}, and the choice of pushout excision generators in Definition~\ref{def:presentation-realizations}(2) imply that $i_\nu^*$ preserves presentation hypercovers, the selected pushout excision relations, vector and projective bundles, zero sections, derived self-intersections, $\mathcal O(1)$, and Thom objects, and sends $\sigma_{W,X}$ on $W\in(\mathcal T_R^{\mathrm{perf}})_{/V}$ to equivalences. Thus $\mathcal M_0(V)\to\mathcal M_{\nu,0}^a(\beta_\nu(V))$ preserves $\mathcal J_{\dMGL}$, and Proposition~\ref{thm:localized-span-universal-property} induces $\Phi_\nu$ and $\phi_\nu$.
\end{proof}

\subsection{Excisive approximation, assembly along the adic tower, and derived base-change specialization}
\label{subsec:reduction-and-p1}

Fix $X_0=\operatorname{Spf}(B)\in\mathcal T_R^{\mathrm{perf}}$, and set $A_\theta=B/I_R B$, $A_\flat=B/(p)$, $B_\nu=A_\nu^a$, $X_{0,\nu}=\beta_\nu(X_0)=\operatorname{Spf}(A_\nu)\in\mathcal P_{\nu,R}$, and $i_\nu^*(X_0)=\operatorname{Spf}^a(B_\nu)$ for $\nu\in\{\theta,\flat\}$. Write $X_{0,\mathrm{Prism}}=X_0$, $R_\theta^a=(R^\circ)^a$, $R_\flat^a=(R^{\flat\circ})^a$, and $R_{1,\mathrm{Prism}}=\mathbb A_{\inf R}/(\xi,\widetilde\omega)=\mathbb A_{\inf R}/(p,\widetilde\omega)$. When $B=\widetilde A=\Ainf(A^\circ)$, we write $\widetilde A$, $A^\circ$, and $A^{\flat\circ}$ in place of $X_0$, $X_{0,\theta}$, and $X_{0,\flat}$ in the categories and spectra below.

Let $\mathcal A_{\mathrm{cmp}}(X_0)$ be the $\mathcal E_R(X_0)$-enriched presentation span category on $\mathscr D_{\mathrm{cmp}}(X_0)$ with mapping objects $\underline{\Corr}_{\mathbf S_{\mathrm{Prism}},X_0}(X,Y)\in\mathcal E_R(X_0)$. For $\nu\in\{\theta,\flat\}$, the equivalence $\tau_\nu=(\beta_\nu^*)^{-1}r^*\colon\mathcal E_R(X_0)\to\mathcal E_{\nu,R}(X_{0,\nu})$ of Proposition~\ref{prop:perfect-boundary-coefficient-sites} defines the transported $\mathcal E_{\nu,R}(X_{0,\nu})$-enriched presentation span category $\mathcal A_{\nu,\mathrm{cmp}}^{\mathrm{pres}}(X_{0,\nu})=(\tau_\nu)_*\mathcal A_{\mathrm{cmp}}(X_0)$ on $i_\nu^*X$ ($X\in\mathscr D_{\mathrm{cmp}}(X_0)$). Here $X_0=\beta_\nu^{-1}(X_{0,\nu})$ is functorially recovered from $X_{0,\nu}$, and Proposition~\ref{prop:boundary-specialization-functoriality}(1) induces a symmetric monoidal $\mathcal E_{\nu,R}(X_{0,\nu})$-enriched realization functor into the intrinsic boundary span category of Definition~\ref{def:presentation-realizations}(1).

In addition, consider the intrinsic derived first reductions $X_{0,\theta,1}=\Spec^a(B_\theta\otimes_{(R^\circ)^a}R_1^a)$, $X_{0,\mathrm{Prism},1}=\Spec^a((B\otimes_{\mathbb A_{\inf R}}R_{1,\mathrm{Prism}})^a)$, and $X_{0,\flat,1}=\Spec^a(B_\flat\otimes_{(R^{\flat\circ})^a}(R_1^\flat)^a)$ of Subsection~\ref{subsec:excisive-approximation-tilting}. For $\mathcal D\in\{\mathrm{Prism},\theta,\flat\}$, let $\mathcal A_{\mathcal D,1,\mathrm{cmp}}(X_{0,\mathcal D})$ be the intrinsic $\mathcal E_{\mathcal D,R}(X_{0,\mathcal D})$-enriched span category of quasi-compact, quasi-separated quasi-smooth derived almost schemes over $X_{0,\mathcal D,1}$, sheafified over $V_{\mathcal D}\mapsto V_{\mathcal D,1}$ on the slice site of $X_{0,\mathcal D}$.

With respect to the shared almost structure over $\mathfrak m R_1\cong\mathfrak m^\flat R_1^\flat$ of Subsection~\ref{subsec:almost-formal-base}, and writing $\bar\omega$, $\bar{\widetilde\omega}$, and $\bar\omega^\flat$ for the formal degree-$1$ generators of $I_{\theta,1}=I_1$, $I_{\mathrm{Prism},1}$, and $I_{\flat,1}=I_1^\flat$, define the completed associated graded commutative almost rings
\[
\begin{aligned}
\widehat G_\theta &= \Bigl(\prod_{n\ge 0}\bigl(\omega^n R^\circ/\omega^{n+1}R^\circ\bigr)\Bigr)^a \simeq \bigl(R_1[\![\bar\omega]\!]\bigr)^a, \qquad
\widehat G_{\mathrm{Prism}} = \bigl(R_{1,\mathrm{Prism}}[\![\bar{\widetilde\omega}]\!]\bigr)^a, \\
\widehat G_\flat &= \Bigl(\prod_{n\ge 0}\bigl((\omega^\flat)^n R^{\flat\circ}/(\omega^\flat)^{n+1}R^{\flat\circ}\bigr)\Bigr)^a \simeq \bigl(R_1^\flat[\![\bar\omega^\flat]\!]\bigr)^a,
\end{aligned}
\]
linked by the ring isomorphisms $\widehat G_\theta\leftarrow\widehat G_{\mathrm{Prism}}\rightarrow\widehat G_\flat$ sending $\bar{\widetilde\omega}^n$ to $\bar\omega^n$ and $(\bar\omega^\flat)^n$ for each $n\ge 0$ (Lemma~\ref{lem:tilting-frame-koszul}(1)--(2) and Subsection~\ref{subsec:almost-formal-base}). For $\mathcal D\in\{\mathrm{Prism},\theta,\flat\}$, we write $\mathcal E_{\mathrm{Prism},R}(X_0)=\mathcal E_R(X_0)$, $\bar\omega_{\mathrm{Prism}}=\bar{\widetilde\omega}$, $\bar\omega_\theta=\bar\omega$, and $\bar\omega_\flat=\bar\omega^\flat$.
Set $G_{\mathcal D,[1,K]}=\bar\omega_{\mathcal D}\widehat G_{\mathcal D}/\bar\omega_{\mathcal D}^{K+1}\widehat G_{\mathcal D}=\bigoplus_{n=1}^K(I_{\mathcal D,n})^a$ ($K\ge 1$). Let $\Mod_{\widehat G_{\mathcal D}}^\wedge(\mathcal E_{\mathcal D,R}(X_{0,\mathcal D}))\subset\Mod_{\widehat G_{\mathcal D}}(\mathcal E_{\mathcal D,R}(X_{0,\mathcal D}))$ denote the presentably symmetric monoidal stable $\infty$-category of hypercomplete sheaves on the slice site of $X_{0,\mathcal D}$ valued in derived $\bar\omega_{\mathcal D}$-complete $\widehat G_{\mathcal D}$-modules $\Mod_{\widehat G_{\mathcal D}}^\wedge\subset\Mod_{\widehat G_{\mathcal D}}$. Its symmetric monoidal localization $L_{\bar\omega_{\mathcal D}}^\wedge$ is given objectwise by $\varprojlim_{K\ge 1}\cofib(\bar\omega_{\mathcal D}^K\colon -\to -)$~\cite[Tags~091S, 091V, and~0920]{stacks-project},~\cite[Proposition~2.2.1.9]{HA}.

\begin{definition}
\label{def:nonunital-comparison-sources}
Let $X_0=\operatorname{Spf}(B)\in\mathcal T_R^{\mathrm{perf}}$, $\mathcal D\in\{\mathrm{Prism},\theta,\flat\}$, and $\nu\in\{\theta,\flat\}$.
\begin{enumerate}[(1)]
\item Form the enriched comparison and associated graded stable motivic envelopes along the adic tower~\cite{MR3345192}:
\[
\begin{aligned}
\mathcal M_{0,\mathrm{Prism}}^{a,\mathrm{nu}}(X_0)
&= \operatorname{Fun}_{\mathcal E_R(X_0)}\bigl(\mathcal A_{\mathrm{cmp}}(X_0)^{\mathrm{op}},\Mod_{\widehat G_{\mathrm{Prism}}}^\wedge(\mathcal E_R(X_0))\bigr), \\
\mathcal M_{0,\nu}^{a,\mathrm{nu}}(X_{0,\nu})
&= \operatorname{Fun}_{\mathcal E_{\nu,R}(X_{0,\nu})}\bigl(\mathcal A_{\nu,\mathrm{cmp}}^{\mathrm{pres}}(X_{0,\nu})^{\mathrm{op}},\Mod_{\widehat G_\nu}^\wedge(\mathcal E_{\nu,R}(X_{0,\nu}))\bigr), \\
\mathcal M_{0,\mathcal D,\mathrm{gr}}^{a,\mathrm{nu}}(X_{0,\mathcal D})
&= \operatorname{Fun}_{\mathcal E_{\mathcal D,R}(X_{0,\mathcal D})}\bigl(\mathcal A_{\mathcal D,1,\mathrm{cmp}}(X_{0,\mathcal D})^{\mathrm{op}},\Mod_{\widehat G_{\mathcal D}}^\wedge(\mathcal E_{\mathcal D,R}(X_{0,\mathcal D}))\bigr).
\end{aligned}
\]
Forming cobordism quotients defines the localized comparison and associated graded categories
\[
\begin{aligned}
\mathcal M_{\mathrm{Prism},\mathrm{cmp}}^{a,\mathrm{nu}}(X_0) &= \mathcal M_{0,\mathrm{Prism}}^{a,\mathrm{nu}}(X_0)/\mathcal J_{\dMGL}, \\
\mathcal M_{\nu,\mathrm{cmp}}^{a,\mathrm{nu}}(X_{0,\nu}) &= \mathcal M_{0,\nu}^{a,\mathrm{nu}}(X_{0,\nu})/\mathcal J_{\nu,\dMGL}^{\mathrm{pres}}, \\
\mathcal M_{\mathcal D,\mathrm{gr}}^{a,\mathrm{nu}}(X_{0,\mathcal D}) &= \mathcal M_{0,\mathcal D,\mathrm{gr}}^{a,\mathrm{nu}}(X_{0,\mathcal D})/\mathcal J_{\mathcal D,1,\dMGL},
\end{aligned}
\]
where $\mathcal J_{\nu,\dMGL}^{\mathrm{pres}}$ is generated by the images under $i_\nu^*$ of the cobordism relations of Definition~\ref{def:relation-dmgl}(1)--(4) and $\mathcal J_{\mathcal D,1,\dMGL}$ is the intrinsic cobordism ideal of Definition~\ref{def:relation-dmgl}(1)--(4) for Zariski hypercovers and $\tau_{\mathrm{perf}}$-hypercovers on the slice site. Symmetric monoidal associated graded realization functors $\mathcal M_{\mathcal D,\mathrm{cmp}}^{a,\mathrm{nu}}(X_{0,\mathcal D})\to\allowbreak\mathcal M_{\mathcal D,\mathrm{gr}}^{a,\mathrm{nu}}(X_{0,\mathcal D})$ link these quotients. Denote the endomorphism $\mathbb E_\infty$-algebras of their monoidal units $\mathbf 1_{\mathcal D}$ by $\mathbf{dMGL}_{\mathrm{Prism},\mathrm{cmp}}^{a,\mathrm{nu}}(X_0)=\allowbreak\operatorname{End}(\mathbf 1_{\mathrm{Prism}})\in\allowbreak\CAlg(\Mod_{\widehat G_{\mathrm{Prism}}}^\wedge)$ and $\mathbf{dMGL}_{\nu,\mathrm{cmp}}^{a,\mathrm{nu}}(X_{0,\nu})=\allowbreak\operatorname{End}(\mathbf 1_\nu)\in\allowbreak\CAlg(\Mod_{\widehat G_\nu}^\wedge)$, and write $\Mod_E^\wedge=\Mod_E(\Mod_{\widehat G_{\mathcal D}}^\wedge)$ for $E\in\CAlg(\Mod_{\widehat G_{\mathcal D}}^\wedge)$.

\item Fix $I\in\NAlg_{\mathrm{Prism}}(B)$, $J_\nu\in\NAlg_\nu^a(B_\nu)$, and $V=\operatorname{Spf}(C)\in(\mathcal T_R^{\mathrm{perf}})_{/X_0}$, setting $C_\nu=(C\otimes_{\mathbb A_{\inf R}}R_\nu)^a$, $I_C=I\widehat{\otimes}_BC$, $J_{\nu,C_\nu}=J_\nu\widehat{\otimes}_{B_\nu}C_\nu$, and $-_V=(-)\times_{X_0}V$. For any span $\sigma\in\mathscr D_{\mathrm{cmp}}(X_V,Y_V)$ with level-$n$ boundary reduction $i_\nu^*\sigma_n$, the direct sums $\bigoplus_{n=1}^K\mathcal T_{\sigma,n}^{\mathrm{Prism}}(I_C)$ and $\bigoplus_{n=1}^K\mathcal T_{i_\nu^*\sigma_n,n}^\nu(J_{\nu,C_\nu})$ carry the $\widehat G_{\mathcal D}/(\bar\omega_{\mathcal D}^K)$-module structure in which $R_1^a$ acts levelwise while $\bar\omega_{\mathcal D}$ acts via $n\mapsto n+1$ for $1\le n<K$ (Theorem~\ref{prop:completed-almost-excisive-reflection}(3)(a)) and vanishes at $n=K$. Passing to inverse limits over $K\ge 1$ defines $\mathcal T_\sigma^{\mathrm{Prism},\mathrm{adic}}(I_C)=\varprojlim_{K\ge 1}\bigoplus_{n=1}^K\mathcal T_{\sigma,n}^{\mathrm{Prism}}(I_C)\in\Mod_{\widehat G_{\mathrm{Prism}}}^\wedge$ and $\mathcal T_{i_\nu^*\sigma}^{\nu,\mathrm{adic}}(J_{\nu,C_\nu})=\varprojlim_{K\ge 1}\bigoplus_{n=1}^K\mathcal T_{i_\nu^*\sigma_n,n}^\nu(J_{\nu,C_\nu})\in\Mod_{\widehat G_\nu}^\wedge$.

Viewing $\mathscr D_{\mathrm{cmp}}(X_V,Y_V)\subseteq\Span_{\mathbf S_{\mathrm{Prism},V}}(\mathscr D_{\mathrm{cmp}}(V))(X_V,Y_V)$ as the full sub-$\infty$-groupoid of admissible spans, for any connective shifted module $W$ with $V_W=\operatorname{Spf}(C\oplus W)$, base change of spans $s_W\colon\Span_{\mathbf S_{\mathrm{Prism},V}}(\mathscr D_{\mathrm{cmp}}(V))\to\Span_{\mathbf S_{\mathrm{Prism},V_W}}(\mathscr D_{\mathrm{cmp}}(V_W))$ makes $\sigma\mapsto\operatorname{Lift}_{C\oplus W}^{\mathrm{rel}}(\sigma)$, and hence $\mathcal T_\sigma^{\mathrm{Prism},\mathrm{adic}}(I_C)$ and $\mathcal T_{i_\nu^*\sigma}^{\nu,\mathrm{adic}}(J_{\nu,C_\nu})$, coherently contravariant on $\Span_{\mathbf S_{\mathrm{Prism},V}}(\mathscr D_{\mathrm{cmp}}(V))^{\mathrm{op}}$ via pullback composition $s_W(\tau)^*$ (Proposition~\ref{prop:canonical-split-lift}(1)). Applying the finite-product-preserving reflection $(-)^{\mathrm{add}}$ at any $V\in(\mathcal T_R^{\mathrm{perf}})_{/X_0}$ produces a $\operatorname{Fun}(((\mathcal T_R^{\mathrm{perf}})_{/X_0})^{\mathrm{op}},\Sp)$-enriched presheaf on the uncompleted presentation spans; completed slice hypersheafification then defines the $\mathcal E_{\mathcal D,R}(X_{0,\mathcal D})$-enriched presheaves~\cite{MR3345192} in $\mathcal M_{0,\mathrm{Prism}}^{a,\mathrm{nu}}(X_0)$ and $\mathcal M_{0,\nu}^{a,\mathrm{nu}}(X_{0,\nu})$:
\[
\begin{aligned}
\underline{\mathcal T}_{\mathrm{Prism},\mathrm{cmp}}^{\mathrm{adic}}(-,Y;I)
&= L_{\bar{\widetilde\omega}}^\wedge Q_{X_0}\Bigl[\,V\longmapsto\Bigl(\operatorname*{colim}_{\sigma\in\mathscr D_{\mathrm{cmp}}(-_V,Y_V)}\mathcal T_\sigma^{\mathrm{Prism},\mathrm{adic}}(I_C)\Bigr)^{\mathrm{add}}\,\Bigr], \\
\underline{\mathcal T}_{\nu,\mathrm{cmp}}^{\mathrm{adic}}(-,i_\nu^*Y;J_\nu)
&= L_{\bar\omega_\nu}^\wedge a_{\mathcal P_{\nu,R}/X_{0,\nu}}^{\mathrm{hyp}}\Bigl[\,\beta_\nu(V)\longmapsto\Bigl(\operatorname*{colim}_{\sigma\in\mathscr D_{\mathrm{cmp}}(-_V,Y_V)}\mathcal T_{i_\nu^*\sigma}^{\nu,\mathrm{adic}}(J_{\nu,C_\nu})\Bigr)^{\mathrm{add}}\,\Bigr].
\end{aligned}
\]
Applying the cobordism localization reflectors $L_{\dMGL}^{\mathrm{adic}}$ and $L_{\nu,\dMGL}^{\mathrm{adic}}$ and evaluating at $X$ and $i_\nu^*X$ defines the \emph{excisive approximations of the almost cobordism sheaves along the adic tower}
\[
\begin{aligned}
\underline{\dMGL}_{\mathrm{Prism},\mathrm{cmp}}^{a,\mathrm{nu}}(X\to Y;I)
&= \bigl(L_{\dMGL}^{\mathrm{adic}}\underline{\mathcal T}_{\mathrm{Prism},\mathrm{cmp}}^{\mathrm{adic}}(-,Y;I)\bigr)(X) \in \Mod_{\widehat G_{\mathrm{Prism}}}^\wedge\bigl(\mathcal E_R(X_0)\bigr), \\
\underline{\dMGL}_{\nu,\mathrm{cmp}}^{a,\mathrm{nu}}(i_\nu^*X\to i_\nu^*Y;J_\nu)
&= \bigl(L_{\nu,\dMGL}^{\mathrm{adic}}\underline{\mathcal T}_{\nu,\mathrm{cmp}}^{\mathrm{adic}}(-,i_\nu^*Y;J_\nu)\bigr)(i_\nu^*X) \in \Mod_{\widehat G_\nu}^\wedge\bigl(\mathcal E_{\nu,R}(X_{0,\nu})\bigr),
\end{aligned}
\]
with section spectra $\dMGL_{\mathrm{Prism},\mathrm{cmp}}^{a,\mathrm{nu}}(X\to Y;I)\in\Mod_{\widehat G_{\mathrm{Prism}}}^\wedge$ over $X_0$ and $\dMGL_{\nu,\mathrm{cmp}}^{a,\mathrm{nu}}(i_\nu^*X\to i_\nu^*Y;J_\nu)\in\Mod_{\widehat G_\nu}^\wedge$ over $X_{0,\nu}$.
\end{enumerate}
\end{definition}

\begin{definition}
\label{def:integral-linear-comparison-sources}
Let $X_0=\operatorname{Spf}(B)\in\mathcal T_R^{\mathrm{perf}}$, $\nu\in\{\theta,\flat\}$, and $\mathcal O\in\{\mathbb A_{\inf R},R_\theta^a,R_\flat^a\}$.
\begin{enumerate}[(1)]
\item Let $\Mod_{\mathbb A_{\inf R}}(\mathcal E_R(X_0))=\Shv_{\mathrm{hyp}}((\mathcal T_R)_{/X_0},\tau_{\mathrm{cmp}};\Mod_{\mathbb A_{\inf R}})\simeq\mathcal E_R(X_0)\otimes_{\Sp}\Mod_{\mathbb A_{\inf R}}$ and $\Mod_{R_\nu^a}(\mathcal E_{\nu,R}(X_{0,\nu}))=\Shv_{\mathrm{hyp}}((\mathcal P_{\nu,R})_{/X_{0,\nu}};\Mod_{R_\nu^a})\simeq\mathcal E_{\nu,R}(X_{0,\nu})\otimes_{\Sp}\Mod_{R_\nu^a}$. Here $\Mod_{R_\nu^a}$ is dualizable in $\mathcal{P}\mathrm{r}^{\mathrm{L}}_{\mathrm{st}}$ as a retract of $\Mod_{R_\nu}$ via $(-)_!\dashv (-)^a$ (Subsection~\ref{subsec:almost-formal-base}). In the integral $\mathcal O$-linear enriched motivic envelopes
\[
\begin{aligned}
\mathcal M_{0,\mathrm{Prism}}^{\mathrm{lin}}(X_0)
&= \operatorname{Fun}_{\mathcal E_R(X_0)}\bigl(\mathcal A_{\mathrm{cmp}}(X_0)^{\mathrm{op}},\Mod_{\mathbb A_{\inf R}}(\mathcal E_R(X_0))\bigr), \\
\mathcal M_{0,\nu}^{\mathrm{lin}}(X_{0,\nu})
&= \operatorname{Fun}_{\mathcal E_{\nu,R}(X_{0,\nu})}\bigl(\mathcal A_{\nu,\mathrm{cmp}}^{\mathrm{pres}}(X_{0,\nu})^{\mathrm{op}},\Mod_{R_\nu^a}(\mathcal E_{\nu,R}(X_{0,\nu}))\bigr),
\end{aligned}
\]
let $\mathcal M_{\mathrm{Prism},\mathrm{cmp}}^{\mathrm{lin}}(X_0)=\mathcal M_{0,\mathrm{Prism}}^{\mathrm{lin}}(X_0)/\mathcal J_{\dMGL}$ and $\mathcal M_{\nu,\mathrm{cmp}}^{\mathrm{lin}}(X_{0,\nu})=\mathcal M_{0,\nu}^{\mathrm{lin}}(X_{0,\nu})/\mathcal J_{\nu,\dMGL}^{\mathrm{pres}}$ denote the localized $\mathcal O$-linear cobordism $\infty$-categories with unit endomorphism $\mathbb E_\infty$-algebras $\mathbf{dMGL}_{\mathrm{Prism},\mathrm{cmp}}^{\mathrm{lin}}(X_0)\in\CAlg(\Mod_{\mathbb A_{\inf R}})$ and $\mathbf{dMGL}_{\nu,\mathrm{cmp}}^{\mathrm{lin}}(X_{0,\nu})\in\CAlg(\Mod_{R_\nu^a})$.

\item For $I\in\NAlg_{\mathrm{Prism}}(B)$, $J_\nu\in\NAlg_\nu^a(B_\nu)$, and $V=\operatorname{Spf}(C)\in(\mathcal T_R^{\mathrm{perf}})_{/X_0}$ (Definition~\ref{def:nonunital-comparison-sources}(2)), forming the enriched presheaves in $\mathcal M_{0,\mathrm{Prism}}^{\mathrm{lin}}(X_0)$ and $\mathcal M_{0,\nu}^{\mathrm{lin}}(X_{0,\nu})$ given by
\[
\begin{aligned}
\underline{\mathcal T}_{\mathrm{Prism},\mathrm{cmp}}(-,Y;I)
&= Q_{X_0}\Bigl[\,V\longmapsto\Bigl(\operatorname*{colim}_{\sigma\in\mathscr D_{\mathrm{cmp}}(-_V,Y_V)}\mathcal T_\sigma^{\mathrm{Prism}}(I_C)\Bigr)^{\mathrm{add}}\,\Bigr], \\
\underline{\mathcal T}_{\nu,\mathrm{cmp}}(-,i_\nu^*Y;J_\nu)
&= a_{\mathcal P_{\nu,R}/X_{0,\nu}}^{\mathrm{hyp}}\Bigl[\,\beta_\nu(V)\longmapsto\Bigl(\operatorname*{colim}_{\sigma\in\mathscr D_{\mathrm{cmp}}(-_V,Y_V)}\mathcal T_{i_\nu^*\sigma}^\nu(J_{\nu,C_\nu})\Bigr)^{\mathrm{add}}\,\Bigr]
\end{aligned}
\]
from Definition~\ref{def:prismatic-adic-excisive-approximations}(1) and Theorem~\ref{prop:completed-almost-excisive-reflection}(3)(b) and applying the cobordism localization reflectors $L_{\dMGL}$ and $L_{\nu,\dMGL}$ defines the \emph{integral excisive approximations of the cobordism sheaves} $\underline{\dMGL}_{\mathrm{Prism},\mathrm{cmp}}^{\mathrm{lin}}(X\to Y;I)\in\Mod_{\mathbb A_{\inf R}}(\mathcal E_R(X_0))$ and $\underline{\dMGL}_{\nu,\mathrm{cmp}}^{\mathrm{lin}}(i_\nu^*X\to i_\nu^*Y;J_\nu)\in\Mod_{R_\nu^a}(\mathcal E_{\nu,R}(X_{0,\nu}))$, with section spectra $\dMGL_{\mathrm{Prism},\mathrm{cmp}}^{\mathrm{lin}}(X\to Y;I)\in\Mod_{\mathbb A_{\inf R}}$ over $X_0$ and $\dMGL_{\nu,\mathrm{cmp}}^{\mathrm{lin}}(i_\nu^*X\to i_\nu^*Y;J_\nu)\in\Mod_{R_\nu^a}$ over $X_{0,\nu}$.
\end{enumerate}
\end{definition}

For $X_0=\operatorname{Spf}(B)\in\mathcal T_R^{\mathrm{perf}}$ with $X_{0,\nu}=\beta_\nu(X_0)$ ($\nu\in\{\theta,\flat\}$), let $\Xi_\nu\colon\Mod_{\widehat G_{\mathrm{Prism}}}^\wedge(\mathcal E_R(X_0))\to\Mod_{\widehat G_\nu}^\wedge(\mathcal E_{\nu,R}(X_{0,\nu}))$ denote the functor induced by site transport $\tau_\nu=(\beta_\nu^*)^{-1}r^*\colon\mathcal E_R(X_0)\xrightarrow{\sim}\mathcal E_{\nu,R}(X_{0,\nu})$ along $\widehat G_{\mathrm{Prism}}\xrightarrow{\sim}\widehat G_\nu$. For each $R'\in\{R^\circ,R^{\flat\circ}\}\cup\{R_k\}_{1\le k\le p}$ with almostification $(R')^a$, let $\theta_{R'}^{\mathrm{mod}}=\bigl(-\otimes_{\mathbb A_{\inf R}}R'\bigr)^a\colon\Mod_{\mathbb A_{\inf R}}\to\Mod_{(R')^a}$ denote the derived base-change functor, writing $\theta_\nu^{\mathrm{mod}}$ when $R'=R_\nu$ for $\nu\in\{\theta,\flat\}$, set $\Theta_\nu=\tau_\nu\circ\theta_\nu^{\mathrm{mod}}$, and for each $K\ge 1$ view $G_{\nu,[1,K]}=\bigoplus_{n=1}^K(I_{\nu,n})^a\in\operatorname{Perf}(R_\nu^a)$ via restriction of scalars along $R_\nu^a\to R_{\nu,1}^a$.

\begin{lemma}
\label{lem:p1-fiber-comparison}
Fix $X_0=\operatorname{Spf}(B)\in\mathcal T_R^{\mathrm{perf}}$ with $X_{0,\nu}=\beta_\nu(X_0)$, $I\in\NAlg_{\mathrm{Prism}}(B)$, $J_\nu\in\NAlg_\nu^a(B_\nu)$ ($\nu\in\{\theta,\flat\}$), and $X,Y\in\mathscr D_{\mathrm{cmp}}(X_0)$.
\begin{enumerate}[\rm (1)]
\item For each $\nu\in\{\theta,\flat\}$ with $I_\nu=i_\nu^* I$, the functor $\Xi_\nu$ is an exact symmetric monoidal equivalence commuting with $L_{\dMGL}^{\mathrm{adic}}$ and inducing equivalences in $\Mod_{\widehat G_\nu}^\wedge(\mathcal E_{\nu,R}(X_{0,\nu}))$:
\[
\Xi_\nu\Bigl(\underline{\dMGL}_{\mathrm{Prism},\mathrm{cmp}}^{a,\mathrm{nu}}(X\to Y;I)\Bigr)
\;\longrightarrow\;
\underline{\dMGL}_{\nu,\mathrm{cmp}}^{a,\mathrm{nu}}(i_\nu^*X\to i_\nu^*Y;I_\nu).
\]
\item For each $R'\in\{R^\circ,R^{\flat\circ}\}\cup\{R_k\}_{1\le k\le p}$, the functor $\theta_{R'}^{\mathrm{mod}}$ preserves all small limits and all small colimits, preserves $\tau_{\mathrm{perf}}$-hypercomplete sheaves upon objectwise postcomposition, and commutes with $Q_{X_0}$ and $L_{\dMGL}$. For $\nu\in\{\theta,\flat\}$, $\Theta_\nu$ induces specialization equivalences $\Theta_\nu\bigl(\underline{\dMGL}_{\mathrm{Prism},\mathrm{cmp}}^{\mathrm{lin}}(X\to Y;I)\bigr)\to\underline{\dMGL}_{\nu,\mathrm{cmp}}^{\mathrm{lin}}(i_\nu^*X\to i_\nu^*Y;I_\nu)$ in $\Mod_{R_\nu^a}(\mathcal E_{\nu,R}(X_{0,\nu}))$, and there is an equivalence of $\mathbb Z_{\ge 1}^{\mathrm{op}}$-towers along $\widehat G_\nu/(\bar\omega_\nu^{K+1})\twoheadrightarrow\widehat G_\nu/(\bar\omega_\nu^K)$ and $G_{\nu,[1,K+1]}\twoheadrightarrow G_{\nu,[1,K]}$:
\[
\begin{aligned}
\cofib\Bigl(\bar\omega_\nu^K\colon\underline{\dMGL}_{\nu,\mathrm{cmp}}^{a,\mathrm{nu}}(i_\nu^*X\to i_\nu^*Y;J_\nu)&\longrightarrow\underline{\dMGL}_{\nu,\mathrm{cmp}}^{a,\mathrm{nu}}(i_\nu^*X\to i_\nu^*Y;J_\nu)\Bigr) \\
&\longrightarrow
\underline{\dMGL}_{\nu,\mathrm{cmp}}^{\mathrm{lin}}(i_\nu^*X\to i_\nu^*Y;J_\nu)\otimes_{R_\nu^a}G_{\nu,[1,K]}.
\end{aligned}
\]
\end{enumerate}
\end{lemma}
\begin{proof}
For~(1), Proposition~\ref{prop:perfect-boundary-coefficient-sites} and the isomorphism $\widehat G_{\mathrm{Prism}}\xrightarrow{\sim}\widehat G_\nu$ ($\bar{\widetilde\omega}\mapsto\bar\omega_\nu$) make $\Xi_\nu$ an exact symmetric monoidal equivalence preserving $L_{\bar\omega_{\mathcal D}}^\wedge$, small limits, and small colimits. Transport along $(\tau_\nu,\Xi_\nu)$ carries $\mathcal J_{\dMGL}$ onto $\mathcal J_{\nu,\dMGL}^{\mathrm{pres}}$ by Proposition~\ref{prop:boundary-specialization-functoriality}(2) and hence commutes with $L_{\dMGL}^{\mathrm{adic}}$. Theorem~\ref{prop:completed-almost-excisive-reflection}(3)(a) identifies the level-$n$ lift spectra $\mathcal T_{\sigma,n}^{\mathrm{Prism}}(I_C)$ and $\mathcal T_{i_\nu^*\sigma_n,n}^\nu((I_\nu)_{C_\nu})$ compatibly with the shift operators and pullback composition; passing to $\varprojlim_{K\ge 1}\bigoplus_{n=1}^K$, $(\operatorname*{colim}_\sigma)^{\mathrm{add}}$, and $L_{\nu,\dMGL}^{\mathrm{adic}}$ gives~(1).

For~(2), dualizability of $R'\in\operatorname{Perf}(\mathbb A_{\inf R})$ (Lemma~\ref{lem:tilting-frame-koszul}(3)) and the adjoint triple $(-)_!\dashv (-)^a\dashv R(-)_*$ imply that $\theta_{R'}^{\mathrm{mod}}$ preserves all small limits and colimits. Thus postcomposition with $\theta_\nu^{\mathrm{mod}}$ preserves $\Delta$-indexed totalizations and $\tau_{\mathrm{perf}}$-hypersheaves, evaluating $\Theta_\nu$ objectwise and commuting with $Q_{X_0}$ and $(-)^{\mathrm{add}}$. By Proposition~\ref{prop:boundary-specialization-functoriality}(2), transport along $(\tau_\nu,\Theta_\nu)$ also preserves both $\mathcal J_{\nu,\dMGL}^{\mathrm{pres}}$ and $\mathcal J_{\nu,\dMGL}^{\mathrm{pres}}$-local objects, commuting with $L_{\dMGL}$; Proposition~\ref{prop:coherent-dual-lift-complex}(3) and Theorem~\ref{prop:completed-almost-excisive-reflection}(3)(b) then give the specialization equivalence for $\underline{\dMGL}_{\mathrm{Prism},\mathrm{cmp}}^{\mathrm{lin}}$.

Finally, on $\Mod_{\widehat G_\nu}^\wedge(\mathcal E_{\nu,R}(X_{0,\nu}))$, the functor $K\mapsto\cofib(\bar\omega_\nu^K)$ preserves all small limits and, via $(-)_!$ and~\cite[Tags~091V and~0920]{stacks-project}, all small colimits, as does $-\otimes_{R_\nu^a}G_{\nu,[1,K]}$ since $G_{\nu,[1,K]}\in\operatorname{Perf}(R_\nu^a)$ is dualizable. Under $\mathcal T_{i_\nu^*\sigma}^{\nu,\mathrm{adic}}(J_{\nu,C_\nu})\simeq\varprojlim_{M\ge 1}\bigl(\mathcal T_{i_\nu^*\sigma}^\nu(J_{\nu,C_\nu})\otimes_{R_\nu^a}G_{\nu,[1,M]}\bigr)$ (Theorem~\ref{prop:completed-almost-excisive-reflection}(3)(b)), the kernel of $\bar\omega_\nu^K$ on $G_{\nu,[1,M]}\cong\widehat G_\nu/(\bar\omega_\nu^M)$ ($M\ge K$) is pro-zero under $\cdot\,\bar\omega_\nu$ and vanishes in $\varprojlim_{M\ge 1}$, identifying the $\bar\omega_\nu^K$-cofiber of $\mathcal T_{i_\nu^*\sigma}^{\nu,\mathrm{adic}}(J_{\nu,C_\nu})$ with $\mathcal T_{i_\nu^*\sigma}^\nu(J_{\nu,C_\nu})\otimes_{R_\nu^a}G_{\nu,[1,K]}$. Commuting both towers with $(\operatorname*{colim}_\sigma)^{\mathrm{add}}$, $a_{\mathcal P_{\nu,R}/X_{0,\nu}}^{\mathrm{hyp}}$, and $L_{\nu,\dMGL}$ completes~(2).
\end{proof}

\subsection{The tilting equivalence over perfectoid algebras}
\label{subsec:tilting-equivalence-perfectoid}

Throughout this subsection, fix an integral perfectoid $R^\circ$-algebra $A^\circ$ with tilt $A^{\flat\circ}$ over $R^{\flat\circ}$ and perfect prism $\widetilde A=\Ainf(A^\circ)=W(A^{\flat\circ})$. Under $P^1_{\mathrm{cot}}$ and $P^1_{\mathrm{rel}}$ (Definition~\ref{def:relative-tangential-linearization} and Theorem~\ref{prop:completed-almost-excisive-reflection}(3)(a)), the theories along the adic tower pass to the completed associated graded almost rings $\widehat G_\theta\xleftarrow{\sim}\widehat G_{\mathrm{Prism}}\xrightarrow{\sim}\widehat G_\flat$, while retaining the integral base rings $R^\circ$ and $R^{\flat\circ}$ with $P^1_{\mathrm{cot}}$ (Theorem~\ref{prop:completed-almost-excisive-reflection}(3)(b)) gives derived base-change specializations along $\theta$ and $\rho_\flat$.

\begin{theorem}
\label{thm:tilting-main}
Let $X_I=\operatorname{Spf}(\widetilde A\oplus I)\to X_0=\operatorname{Spf}(\widetilde A)$ ($I\in\NAlg_{\mathrm{Prism}}(\widetilde A)$) be a split augmented derived $(p,I_R)$-complete formal scheme with boundary specializations $I_\theta=i_\theta^*I\in\NAlg_\theta^a((A^\circ)^a)$ and $I_\flat=i_\flat^*I\in\NAlg_\flat^a((A^{\flat\circ})^a)$, and let $X,Y\in\mathscr D_{\mathrm{cmp}}(X_0)$.
\begin{enumerate}[\rm (1)]
\item For the transported cobordism categories $\mathcal M_{\mathcal D,\mathrm{cmp}}^{a,\mathrm{nu}}$ along the adic tower and their intrinsic mod-$\omega_{\mathcal D}$ associated graded versions $\mathcal M_{\mathcal D,\mathrm{gr}}^{a,\mathrm{nu}}$ (Definition~\ref{def:nonunital-comparison-sources}), there are symmetric monoidal almost tilting equivalences of $\infty$-categories and of derived $\bar\omega_{\mathcal D}$-complete module $\infty$-categories $\Mod_E^\wedge$ over $E=\mathbf{dMGL}_{\mathcal D,\mathrm{cmp}}^{a,\mathrm{nu}}=\operatorname{End}(\mathbf 1_{\mathcal D})$ ($\mathcal D\in\{\mathrm{Prism},\theta,\flat\}$):
\[
\begin{aligned}
\mathcal M_{\theta,\mathrm{cmp}}^{a,\mathrm{nu}}(A^\circ)
&\;\longleftarrow\;
\mathcal M_{\mathrm{Prism},\mathrm{cmp}}^{a,\mathrm{nu}}(\widetilde A)
\;\longrightarrow\;
\mathcal M_{\flat,\mathrm{cmp}}^{a,\mathrm{nu}}(A^{\flat\circ}), \\
\mathcal M_{\theta,\mathrm{gr}}^{a,\mathrm{nu}}(A^\circ)
&\;\longleftarrow\;
\mathcal M_{\mathrm{Prism},\mathrm{gr}}^{a,\mathrm{nu}}(\widetilde A)
\;\longrightarrow\;
\mathcal M_{\flat,\mathrm{gr}}^{a,\mathrm{nu}}(A^{\flat\circ}), \\
\Mod_{\mathbf{dMGL}_{\theta,\mathrm{cmp}}^{a,\mathrm{nu}}(A^\circ)}^\wedge
&\;\longleftarrow\;
\Mod_{\mathbf{dMGL}_{\mathrm{Prism},\mathrm{cmp}}^{a,\mathrm{nu}}(\widetilde A)}^\wedge
\;\longrightarrow\;
\Mod_{\mathbf{dMGL}_{\flat,\mathrm{cmp}}^{a,\mathrm{nu}}(A^{\flat\circ})}^\wedge,
\end{aligned}
\]
compatible under the associated graded realization functors $\mathcal M_{\mathcal D,\mathrm{cmp}}^{a,\mathrm{nu}}\to\mathcal M_{\mathcal D,\mathrm{gr}}^{a,\mathrm{nu}}$ ($\mathcal D\in\{\mathrm{Prism},\theta,\flat\}$), together with almost tilting equivalences of cobordism sheaves and section spectra across all $n\ge 1$:
\begin{equation}
\label{eq:evaluated-tilting-zigzag}
\dMGL_{\theta,\mathrm{cmp}}^{a,\mathrm{nu}}(i_\theta^*X\to i_\theta^*Y;I_\theta)
\;\longleftarrow\;
\dMGL_{\mathrm{Prism},\mathrm{cmp}}^{a,\mathrm{nu}}(X\to Y;I)
\;\longrightarrow\;
\dMGL_{\flat,\mathrm{cmp}}^{a,\mathrm{nu}}(i_\flat^*X\to i_\flat^*Y;I_\flat),
\end{equation}
compatible with refinements of presentation hypercovers.

\item For the transported integral $\mathcal O$-linear cobordism categories and intrinsic integral cobordism spectra over $\mathcal O\in\{\mathbb A_{\inf R},(R^\circ)^a,(R^{\flat\circ})^a\}$ (Definition~\ref{def:integral-linear-comparison-sources}), derived base change along $\theta\colon\mathbb A_{\inf R}\to R^\circ$ and $\rho_\flat\colon\mathbb A_{\inf R}\to R^{\flat\circ}$ induces symmetric monoidal equivalences of $\infty$-categories and endomorphism $\mathbb E_\infty$-algebras of the monoidal unit:
\[
\begin{aligned}
\mathcal M_{\mathrm{Prism},\mathrm{cmp}}^{\mathrm{lin}}(\widetilde A)\otimes_{\Mod_{\mathbb A_{\inf R}}}\Mod_{(R^\circ)^a}
&\;\longrightarrow\;
\mathcal M_{\theta,\mathrm{cmp}}^{\mathrm{lin}}(A^\circ), \\
\mathcal M_{\mathrm{Prism},\mathrm{cmp}}^{\mathrm{lin}}(\widetilde A)\otimes_{\Mod_{\mathbb A_{\inf R}}}\Mod_{(R^{\flat\circ})^a}
&\;\longrightarrow\;
\mathcal M_{\flat,\mathrm{cmp}}^{\mathrm{lin}}(A^{\flat\circ}), \\
\Bigl(\mathbf{dMGL}_{\mathrm{Prism},\mathrm{cmp}}^{\mathrm{lin}}(\widetilde A)\otimes_{\mathbb A_{\inf R},\theta}R^\circ\Bigr)^a
&\;\longrightarrow\;
\mathbf{dMGL}_{\theta,\mathrm{cmp}}^{\mathrm{lin}}(A^\circ), \\
\Bigl(\mathbf{dMGL}_{\mathrm{Prism},\mathrm{cmp}}^{\mathrm{lin}}(\widetilde A)\otimes_{\mathbb A_{\inf R},\rho_\flat}R^{\flat\circ}\Bigr)^a
&\;\longrightarrow\;
\mathbf{dMGL}_{\flat,\mathrm{cmp}}^{\mathrm{lin}}(A^{\flat\circ}),
\end{aligned}
\]
and equivalences of the integral cobordism spectra:
\begin{equation}
\label{eq:functorial-tilting-zigzag}
\begin{aligned}
\Bigl(\dMGL_{\mathrm{Prism},\mathrm{cmp}}^{\mathrm{lin}}(X\to Y;I)\otimes_{\mathbb A_{\inf R},\theta}R^\circ\Bigr)^a
&\;\longrightarrow\;
\dMGL_{\theta,\mathrm{cmp}}^{\mathrm{lin}}(i_\theta^*X\to i_\theta^*Y;I_\theta), \\
\Bigl(\dMGL_{\mathrm{Prism},\mathrm{cmp}}^{\mathrm{lin}}(X\to Y;I)\otimes_{\mathbb A_{\inf R},\rho_\flat}R^{\flat\circ}\Bigr)^a
&\;\longrightarrow\;
\dMGL_{\flat,\mathrm{cmp}}^{\mathrm{lin}}(i_\flat^*X\to i_\flat^*Y;I_\flat).
\end{aligned}
\end{equation}
\end{enumerate}
\end{theorem}
\begin{proof}
By~\cite[Theorem~3.10 and Lemma~3.8]{BhattScholze2017Prisms} and Proposition~\ref{prop:prismatic-fiber-completion}, $X_0=\operatorname{Spf}(\widetilde A)\in\mathcal T_R^{\mathrm{perf}}$ satisfies $(\widetilde A\otimes_{\mathbb A_{\inf R}}R_\nu)^a\simeq A_\nu^a$ ($A_\theta=A^\circ$, $A_\flat=A^{\flat\circ}$), and Lemma~\ref{lem:tilting-frame-koszul}(2)--(3) identifies the mod-$\omega_{\mathcal D}$ reductions $C_\theta^a\otimes_{(R^\circ)^a}R_1^a\xleftarrow{\sim}(C\otimes_{\mathbb A_{\inf R}}R_{1,\mathrm{Prism}})^a\xrightarrow{\sim}C_\flat^a\otimes_{(R^{\flat\circ})^a}(R_1^\flat)^a$ for any $V=\operatorname{Spf}(C)\in(\mathcal T_R^{\mathrm{perf}})_{/X_0}$.

For~(1), transport along $(\tau_\nu,\Xi_\nu)$ (Proposition~\ref{prop:perfect-boundary-coefficient-sites} and Lemma~\ref{lem:p1-fiber-comparison}(1)) gives exact symmetric monoidal equivalences $\mathcal M_{0,\mathrm{Prism}}^{a,\mathrm{nu}}(X_0)\allowbreak\to\allowbreak\mathcal M_{0,\nu}^{a,\mathrm{nu}}(X_{0,\nu})$ and $\mathcal M_{0,\mathrm{Prism},\mathrm{gr}}^{a,\mathrm{nu}}(X_0)\allowbreak\to\allowbreak\mathcal M_{0,\nu,\mathrm{gr}}^{a,\mathrm{nu}}(X_{0,\nu})$ carrying the cobordism relation ideals onto the boundary relation ideals; Proposition~\ref{thm:localized-span-universal-property} gives the category equivalences on $\mathcal M_{\mathcal D,\mathrm{cmp}}^{a,\mathrm{nu}}$ and $\mathcal M_{\mathcal D,\mathrm{gr}}^{a,\mathrm{nu}}$.

Under enriched Day convolution~\cite[Example~2.2.6.9 and Section~4.8.5]{HA},~\cite{MR3345192}, taking global sections of $\underline{\operatorname{End}}(\mathbf 1_{\mathcal D})$ over $X_{0,\mathcal D}$ identifies $\Xi_\nu$ on endomorphisms of the unit with $\mathbb E_\infty$-algebra equivalences on $\mathbf{dMGL}_{\mathcal D,\mathrm{cmp}}^{a,\mathrm{nu}}$, inducing the equivalences of module categories~\cite[Theorem~4.8.4.6]{HA}. Evaluating Lemma~\ref{lem:p1-fiber-comparison}(1) at $X_0$ gives~\eqref{eq:evaluated-tilting-zigzag}.

For~(2), dualizability of $\Mod_{R_\nu^a}$ over $\Sp$ (Definition~\ref{def:integral-linear-comparison-sources}(1)), Proposition~\ref{prop:perfect-boundary-coefficient-sites},~\cite{MR3345192}, and~\cite[Proposition~4.8.1.17 and Theorem~4.8.4.6]{HA} give $\mathcal M_{0,\mathrm{Prism}}^{\mathrm{lin}}(X_0)\otimes_{\Mod_{\mathbb A_{\inf R}}}\Mod_{R_\nu^a}\xrightarrow{\sim}\mathcal M_{0,\nu}^{\mathrm{lin}}(X_{0,\nu})$. Propositions~\ref{prop:boundary-specialization-functoriality}(2) and~\ref{thm:localized-span-universal-property} pass this equivalence to $\mathcal M_{\mathcal D,\mathrm{cmp}}^{\mathrm{lin}}$, and evaluating Lemma~\ref{lem:p1-fiber-comparison}(2) on $\mathbf 1_{\mathcal M_{0,\mathrm{Prism}}^{\mathrm{lin}}(X_0)}$ and on $(X,Y;I)$ over $X_0$ gives the unit endomorphism equivalences and~\eqref{eq:functorial-tilting-zigzag}.
\end{proof}

\begin{proposition}
\label{prop:koszul-tilting-and-tower-recovery}
Let $I\in\NAlg_{\mathrm{Prism}}(\widetilde A)$ have boundary specializations $I_\theta=i_\theta^*I$ and $I_\flat=i_\flat^*I$, and let $X,Y\in\mathscr D_{\mathrm{cmp}}(\operatorname{Spf}(\widetilde A))$.
\begin{enumerate}[\rm (1)]
\item For each $\nu\in\{\theta,\flat\}$, the cobordism spectra along the adic tower are recovered from the integral $\mathcal O$-linear theory over $G_{\nu,[1,K]}=\bigoplus_{n=1}^K (I_{\nu,n})^a$ ($K\ge 1$) via an equivalence in $\Mod_{\widehat G_\nu}^\wedge$:
\[
\dMGL_{\nu,\mathrm{cmp}}^{a,\mathrm{nu}}(i_\nu^*X\to i_\nu^*Y;I_\nu)
\;\longrightarrow\;
\varprojlim_{K\ge 1}\Bigl(\dMGL_{\nu,\mathrm{cmp}}^{\mathrm{lin}}(i_\nu^*X\to i_\nu^*Y;I_\nu)\otimes_{R_\nu^a}G_{\nu,[1,K]}\Bigr).
\]
\item For each integer $1\le k\le p$, and in particular at $k=p$ over $R_p^a=(R^\circ/(p))^a\simeq (R^{\flat\circ}/(\bar\xi))^a=(R_p^\flat)^a$, derived reduction along $R^\circ\to R_k$ and $R^{\flat\circ}\to R_k^\flat$ induces symmetric monoidal equivalences of $\infty$-categories and module $\infty$-categories:
\[
\begin{aligned}
\mathcal M_{\theta,\mathrm{cmp}}^{\mathrm{lin}}(A^\circ)\otimes_{\Mod_{(R^\circ)^a}}\Mod_{R_k^a}
&\;\longleftarrow\;
\mathcal M_{\mathrm{Prism},\mathrm{cmp}}^{\mathrm{lin}}(\widetilde A)\otimes_{\Mod_{\mathbb A_{\inf R}}}\Mod_{R_k^a} \\
&\;\longrightarrow\;
\mathcal M_{\flat,\mathrm{cmp}}^{\mathrm{lin}}(A^{\flat\circ})\otimes_{\Mod_{(R^{\flat\circ})^a}}\Mod_{(R_k^\flat)^a}, \\
\Mod_{(\mathbf{dMGL}_{\theta,\mathrm{cmp}}^{\mathrm{lin}}(A^\circ)\otimes_{(R^\circ)^a}R_k^a)}
&\;\longleftarrow\;
\Mod_{(\mathbf{dMGL}_{\mathrm{Prism},\mathrm{cmp}}^{\mathrm{lin}}(\widetilde A)\otimes_{\mathbb A_{\inf R}}R_k)^a} \\
&\;\longrightarrow\;
\Mod_{(\mathbf{dMGL}_{\flat,\mathrm{cmp}}^{\mathrm{lin}}(A^{\flat\circ})\otimes_{(R^{\flat\circ})^a}(R_k^\flat)^a)},
\end{aligned}
\]
and equivalences of cobordism spectra in $\Mod_{R_k^a}\simeq\Mod_{(R_k^\flat)^a}$:
\[
\begin{aligned}
\dMGL_{\theta,\mathrm{cmp}}^{\mathrm{lin}}(i_\theta^*X\to i_\theta^*Y;I_\theta)\otimes_{(R^\circ)^a}R_k^a
&\;\longleftarrow\;
\Bigl(\dMGL_{\mathrm{Prism},\mathrm{cmp}}^{\mathrm{lin}}(X\to Y;I)\otimes_{\mathbb A_{\inf R}}R_k\Bigr)^a \\
&\;\longrightarrow\;
\dMGL_{\flat,\mathrm{cmp}}^{\mathrm{lin}}(i_\flat^*X\to i_\flat^*Y;I_\flat)\otimes_{(R^{\flat\circ})^a}(R_k^\flat)^a.
\end{aligned}
\]
\end{enumerate}
\end{proposition}
\begin{proof}
For~(1), derived $\bar\omega_\nu$-completeness of $\dMGL_{\nu,\mathrm{cmp}}^{a,\mathrm{nu}}(i_\nu^*X\to i_\nu^*Y;I_\nu)$ identifies it with the inverse limit of its $\bar\omega_\nu^K$-cofibers; taking global sections over $X_{0,\nu}=\operatorname{Spf}(A_\nu)$ and applying $\varprojlim_{K\ge 1}$ to the tower equivalence of Lemma~\ref{lem:p1-fiber-comparison}(2) gives~(1).

For~(2), fix $1\le k\le p$. Lemma~\ref{lem:tilting-frame-koszul}(2)--(3) and Subsection~\ref{subsec:almost-formal-base} identify $R_k^a\simeq (R_k^\flat)^a$ and show that $R_k\in\operatorname{Perf}(\mathbb A_{\inf R})$ and $R_k^a\in\operatorname{Perf}(R_\nu^a)$ ($\nu\in\{\theta,\flat\}$) are dualizable over $\mathbb A_{\inf R}$ and $R_\nu^a$. This dualizability gives an equivalence $\Mod_{R_k^a}\to\Mod_{(R_k^\flat)^a}$ over $\Mod_{\mathbb A_{\inf R}}$ and $\Mod_{R_\nu^a}$, together with transitivity $(M\otimes_{\mathbb A_{\inf R},\nu}R_\nu)^a\allowbreak\otimes_{R_\nu^a}R_{\nu,k}^a\simeq\allowbreak (M\otimes_{\mathbb A_{\inf R}}R_k)^a$. Combining Theorem~\ref{thm:tilting-main}(2),~\cite[Theorem~4.8.4.6]{HA}, and Lemma~\ref{lem:p1-fiber-comparison}(2) gives the equivalences in~(2).
\end{proof}

For $X,Y\in\mathscr D_{\mathrm{cmp}}(\operatorname{Spf}(\widetilde A))$ and $I\in\NAlg_{\mathrm{Prism}}(\widetilde A)$ with boundary specializations $I_\theta=i_\theta^*I$ and $I_\flat=i_\flat^*I$, write $\widehat{\dMGL}_{\theta,\mathrm{cmp}}^{\mathrm{lin}}(i_\theta^*X\to i_\theta^*Y;I_\theta)\in\Mod_{(R^\circ)^a}^\wedge$ and $\widehat{\dMGL}_{\flat,\mathrm{cmp}}^{\mathrm{lin}}(i_\flat^*X\to i_\flat^*Y;I_\flat)\in\Mod_{(R^{\flat\circ})^a}^\wedge$ for the derived $\omega$- and $\omega^\flat$-completions of the integral boundary almost cobordism spectra.

\begin{corollary}
\label{cor:mixed-to-equal-conservativity}
Let $I\in\NAlg_{\mathrm{Prism}}(\widetilde A)$ have boundary specializations $I_\theta=i_\theta^*I$ and $I_\flat=i_\flat^*I$, and let $X,Y\in\mathscr D_{\mathrm{cmp}}(\operatorname{Spf}(\widetilde A))$. Then $\widehat{\dMGL}_{\theta,\mathrm{cmp}}^{\mathrm{lin}}(i_\theta^*X\to i_\theta^*Y;I_\theta)\simeq 0$ in mixed characteristic over $(R^\circ)^a$ if and only if $\dMGL_{\theta,\mathrm{cmp}}^{a,\mathrm{nu}}(i_\theta^*X\to i_\theta^*Y;I_\theta)\simeq 0$, if and only if $\dMGL_{\flat,\mathrm{cmp}}^{a,\mathrm{nu}}(i_\flat^*X\to i_\flat^*Y;I_\flat)\simeq 0$, if and only if $\widehat{\dMGL}_{\flat,\mathrm{cmp}}^{\mathrm{lin}}(i_\flat^*X\to i_\flat^*Y;I_\flat)\simeq 0$ in characteristic $p$ over $(R^{\flat\circ})^a$. More generally, any morphism of prismatic presentation data, or any morphism in $\mathcal M_{\mathrm{Prism},\mathrm{cmp}}^{\mathrm{lin}}(\widetilde A)$ evaluated objectwise, induces an almost equivalence after derived $\omega$-completion over $(R^\circ)^a$ if and only if its tilted specialization induces an almost equivalence after derived $\omega^\flat$-completion over $(R^{\flat\circ})^a$.
\end{corollary}
\begin{proof}
Let $M\in\Mod_{\mathbb A_{\inf R}}$ have boundary specializations $M_\theta=(M\otimes_{\mathbb A_{\inf R},\theta}R^\circ)^a$ and $M_\flat=(M\otimes_{\mathbb A_{\inf R},\rho_\flat}R^{\flat\circ})^a$, with derived completions $\widehat M_\theta=\varprojlim_{n\ge 1}\cofib(\omega^n\colon M_\theta\to M_\theta)$ and $\widehat M_\flat=\varprojlim_{n\ge 1}\cofib((\omega^\flat)^n\colon M_\flat\to M_\flat)$. In $\Mod_{(R^\circ)^a}$, the fiber of $M_\theta\to\widehat M_\theta$ is $T_\omega(M_\theta)=\varprojlim(\cdots\xrightarrow{\omega}M_\theta\xrightarrow{\omega}M_\theta)$, on which multiplication by $\omega$ is an equivalence; thus $\cofib(\omega\colon T_\omega(M_\theta)\to T_\omega(M_\theta))\simeq 0$ and $M_\theta\otimes_{(R^\circ)^a}R_1^a\xrightarrow{\sim}\widehat M_\theta\otimes_{(R^\circ)^a}R_1^a$.

By derived $\omega$-completeness, $\widehat M_\theta\otimes_{(R^\circ)^a}R_1^a\simeq 0$ holds if and only if multiplication by $\omega$ on $\widehat M_\theta$ is an equivalence, which happens if and only if $\widehat M_\theta\simeq T_\omega(\widehat M_\theta)\simeq 0$. The same characterization applies to $\widehat M_\flat$ with respect to $\omega^\flat$ and to derived $\bar\omega_\nu$-complete modules with respect to $\bar\omega_\nu$~\cite[Tags~091V and~0G1U]{stacks-project}. Combining these equivalences with Theorem~\ref{thm:tilting-main}(1)--(2), Proposition~\ref{prop:koszul-tilting-and-tower-recovery}(1)--(2) at $k=1$, Lemma~\ref{lem:p1-fiber-comparison}(2) at $K=1$, and Lemma~\ref{lem:tilting-frame-koszul}(3) at $k=1$ proves the assertion for $M=\dMGL_{\mathrm{Prism},\mathrm{cmp}}^{\mathrm{lin}}(X\to Y;I)$ and for objectwise cofibers of morphisms in $\mathcal M_{\mathrm{Prism},\mathrm{cmp}}^{\mathrm{lin}}(\widetilde A)$.
\end{proof}

\begin{remark}
\label{rem:generic-fiber-localization}
Over the Tate perfectoid generic fibers $R=R^\circ[\omega^{-1}]=R^\circ[1/p]$ and $R^\flat=R^{\flat\circ}[(\omega^\flat)^{-1}]$, since $\omega\in\mathfrak{m}$ makes the flat inclusion $\mathfrak{m}\to R^\circ$ an isomorphism $\mathfrak{m}\otimes_{R^\circ}R\to R$, any $L\in\ker((-)^a)$ satisfies $L\otimes_{R^\circ}R\simeq (\mathfrak{m}\otimes_{R^\circ}L)\otimes_{R^\circ}R\simeq 0$. Thus $-\otimes_{R^\circ}R$ and $-\otimes_{R^{\flat\circ}}R^\flat$ factor through almostification as exact symmetric monoidal localizations $(-)[\omega^{-1}]\colon\Mod_{(R^\circ)^a}\to\Mod_R$ and $(-)[(\omega^\flat)^{-1}]\colon\Mod_{(R^{\flat\circ})^a}\to\Mod_{R^\flat}$ eliminating almost-zero modules.

Inverting $\widetilde\omega=[\omega^\flat]\in\mathbb A_{\inf R}$ in Theorem~\ref{thm:tilting-main}(2) specializes $\dMGL_{\mathrm{Prism},\mathrm{cmp}}^{\mathrm{lin}}(X\to Y;I)[\widetilde\omega^{-1}]$ and $\mathcal M_{\mathrm{Prism},\mathrm{cmp}}^{\mathrm{lin}}(\widetilde A)\otimes_{\Mod_{\mathbb A_{\inf R}}}\Mod_{\mathbb A_{\inf R}[\widetilde\omega^{-1}]}$ along $\theta[\widetilde\omega^{-1}]$ and $\rho_\flat[\widetilde\omega^{-1}]$ to ordinary cobordism spectra and symmetric monoidal cobordism $\infty$-categories over $R$ and $R^\flat$. Applying $(-)[\omega^{-1}]$ to Corollary~\ref{cor:mixed-to-equal-conservativity} transfers tilted equivalences on $\dMGL_{\flat,\mathrm{cmp}}^{a,\mathrm{nu}}$ or $\widehat{\dMGL}_{\flat,\mathrm{cmp}}^{\mathrm{lin}}$ to equivalences on $\widehat{\dMGL}_{\theta,\mathrm{cmp}}^{\mathrm{lin}}[\omega^{-1}]$ over $R$.
\end{remark}
\section{Hyperdescent and the neutral locus}
\label{sec:hyperdescent-and-neutral-locus}

Iterating loop spaces shifts the connectivity $\conn$ of cotangent fibers linearly, so the excisive approximation of the \v{C}ech--Amitsur totalization tower stabilizes at an explicit finite stage (Theorem~\ref{thm:finite-stage-detection}). At the end of the section, we describe the neutral locus on which $\tau_{\mathrm{geom}}$- and $\tau_{\mathrm{pld}}$-hypercomplete sheaves coincide.

\subsection{Iterated loop spaces and connectivity shifts}
\label{subsec:loop-charts-bar}

\begin{definition}\label{def:connectivity}
For a spectrum $E \in \mathrm{Sp}$, its \emph{connectivity} is defined by
\[
  \conn(E)
  = \sup\bigl\{\, m \in \mathbb{Z} \;\big|\; \pi_j(E) = 0 \text{ for all } j < m \,\bigr\}
  \in \mathbb{Z} \cup \{-\infty, +\infty\},
\]
with the convention $\conn(0) = +\infty$. Thus $\conn(E) \ge m$ if and only if $E$ is $m$-connective, and $\conn(E[n]) = \conn(E) + n$ for each $n \in \mathbb{Z}$. For an $\mathbb{E}_\infty$-ring $A$, an $A$-module $M \in \Mod_A$, or a non-unital $A$-algebra $J \in \NAlg_A$, we define $\conn(A)$, $\conn(M)$, and $\conn(J)$ to be the connectivity of the respective underlying spectrum. Unlike a nonzero unital ring, a non-unital $A$-algebra $J = \fib(\epsilon \colon B \to A) \in \NAlg_A$ can have arbitrarily high connectivity, measuring the connectivity of the augmentation $\epsilon$.
\end{definition}

Let $U = \Spec B$ be a pointed derived affine scheme over $\Spec A$ with unit $\iota \colon A \to B$, augmentation $\epsilon \colon B \to A$, section $s = \Spec(\epsilon) \colon \Spec A \to U$, and augmentation ideal $J = \fib(\epsilon) \in \NAlg_A$ corresponding to $B$ under the unitalization equivalence $\NAlg_A \to \CAlg_{A//A}$, $I \mapsto A \oplus I$. Setting $B^{(0)} = B$, define inductively $B^{(n+1)} = A \otimes_{B^{(n)}} A$ ($n \ge 0$) with augmentation ideal $J^{(n)} = \fib(B^{(n)} \to A) \in \NAlg_A$. By contravariance of $\Spec \colon \CAlg_{A//A}^{\mathrm{op}} \to (\mathrm{dAff}_{/\Spec A})_*$, $\Spec B^{(n)} \simeq \Omega_A^n U$ as pointed derived affine schemes over $\Spec A$. Under $\CAlg_{A//A} \simeq \NAlg_A$, the pushout $B^{(n+1)}$ corresponds to the categorical suspension $\Sigma_{\NAlg_A} J^{(n)} = 0 \amalg_{J^{(n)}} 0$, so $J^{(n)} \simeq \Sigma_{\NAlg_A}^n J$. Since the cotangent fiber functor $Q_A \colon \NAlg_A \to \Mod_A^{\ge 0}$ is a left adjoint and preserves pushouts (Proposition~\ref{prop:excisive-reflection-means}(2) and~\cite[Theorem~7.3.4.13]{HA},~\cite[Sections~2--3]{KatoCalcNUAlg}), we obtain $Q_A(J^{(n)}) \simeq Q_A(J)[n]$ in $\Mod_A^{\ge 0}$, split square-zero excisive approximation $A \oplus \operatorname{sqz}_A\bigl(Q_A(J^{(n)})\bigr) \simeq A \oplus Q_A(J)[n]$, and the connectivity shift
\[
  \conn\bigl(Q_A(J^{(n)})\bigr) = \conn\bigl(Q_A(J)\bigr) + n.
\]

Similarly, for a pointed affine formal scheme $U = \operatorname{Spf}(B)$ over $\operatorname{Spf}(A)$ in $\mathcal{T}_R^{\mathrm{an}}$ with augmentation $\epsilon \colon B \to A$ and $\widehat{J} = \fib(\epsilon) \in \NAlg_{\mathrm{Prism}}(A)$, set $\widehat{B}^{(0)} = B$ and $\widehat{B}^{(n+1)} = (A \otimes_{\widehat{B}^{(n)}} A)^\wedge_{(p,I_R)}$ ($n \ge 0$), with $\widehat{J}^{(n)} = \fib(\widehat{B}^{(n)} \to A) \in \NAlg_{\mathrm{Prism}}(A)$. Proposition~\ref{prop:excisive-reflection-means}(1)--(2) identifies $\widehat{J}^{(n)} \simeq \Sigma_{\NAlg_{\mathrm{Prism}}}^n \widehat{J}$ and $Q_A(\widehat{J}^{(n)}) \simeq Q_A(\widehat{J})[n]$ in $\Mod_A^{\ge 0,\wedge}$.

Whenever $\widehat{B}^{(k)}$ ($0 \le k \le n$) defines a strict bounded animated prism in $\mathcal{T}_R^{\mathrm{an}}$, derived $(p,I_R)$-completed tensor products compute homotopy pullbacks over $\operatorname{Spf}(A)$ by Definition~\ref{def:category-animated-prisms}(3), so $\operatorname{Spf}(\widehat{B}^{(n)})$ represents $\Omega_{\operatorname{Spf}(A)}^n U$ in $\mathcal{T}_R^{\mathrm{an}}$. For general $n$, we view $\Omega_{\operatorname{Spf}(A)}^n U$ as a presheaf on $\mathcal{T}_R$ and evaluate hypercovers via $\tau_{\mathrm{cmp}}$-hyperdescent (Definition~\ref{def:perfect-presentation-topology}). For any $c$-connective $A$-module $M$ ($c \in \mathbb{Z}$), exactness of derived $(p,I_R)$-completion on $\Mod_A$ and preservation of connectivity~\cite[Tag~091V]{stacks-project} imply $\conn(M^\wedge_{(p,I_R)}) \ge c$, preserving the linear connectivity slope in Subsection~\ref{subsec:finite-stage-detection}.
\subsection{The totalization tower and loop hyperdescent}
\label{subsec:totalization-tower}

For a cosimplicial spectrum $C^\bullet$, write $\operatorname{Tot} C^\bullet = \operatorname{Tot}_{[q] \in \Delta} C^q$ for its totalization and $\operatorname{Tot}_{\le n} C^\bullet = \operatorname{Tot}_{[q] \in \Delta_{\le n}} C^q$ for its $n$-th partial totalization along the full subcategory $\Delta_{\le n} \subset \Delta$ spanned by $[0], \dots, [n]$. Its $(q-1)$-th cosimplicial matching object is $M^{q-1} C^\bullet = \varprojlim_{([q] \twoheadrightarrow [m]),\, m < q} C^m$, with $M^{-1} C^\bullet = 0$ and normalized matching fiber $\fib(C^q \to M^{q-1} C^\bullet)$ in degree $q$.

\begin{lemma}\label{lem:tot-tower-connectivity}
Let $C^\bullet$ be a cosimplicial spectrum. If there exist constants $c > 1$, $d \in \mathbb{R}$, and an integer $q_0 \ge 1$ such that $\conn\bigl(\fib(C^q \to M^{q-1} C^\bullet)\bigr) \ge cq - d$ for all $q \ge q_0$, then
\[
\conn\bigl(\fib(\operatorname{Tot}_{\le q} C^\bullet \to \operatorname{Tot}_{\le q-1} C^\bullet)\bigr) \ge (c - 1)q - d \quad \text{for all } q \ge q_0,
\]
and for each $i \in \mathbb{Z}$, the truncation map induces an isomorphism
\[
\pi_i(\operatorname{Tot} C^\bullet) \longrightarrow \pi_i(\operatorname{Tot}_{\le n} C^\bullet) \quad \text{for all integers } n \ge \max\Bigl(q_0 - 1, \frac{i + d + 2}{c - 1}\Bigr).
\]
\end{lemma}
\begin{proof}
Since $\fib(\operatorname{Tot}_{\le q} C^\bullet \to \operatorname{Tot}_{\le q-1} C^\bullet) \simeq \Omega^q \fib(C^q \to M^{q-1} C^\bullet) \simeq \fib(C^q \to M^{q-1} C^\bullet)[-q]$ \cite{MS16}, Definition~\ref{def:connectivity} gives
\[
\conn\bigl(\fib(\operatorname{Tot}_{\le q} C^\bullet \to \operatorname{Tot}_{\le q-1} C^\bullet)\bigr) = \conn\bigl(\fib(C^q \to M^{q-1} C^\bullet)\bigr) - q \ge (c - 1)q - d
\]
for all $q \ge q_0$. For any integer $q > n \ge \max\bigl(q_0 - 1, \frac{i + d + 2}{c - 1}\bigr)$, we have $q \ge q_0$ and $(c - 1)q - d > i + 2$, so $\fib(\operatorname{Tot}_{\le q} C^\bullet \to \operatorname{Tot}_{\le q-1} C^\bullet)$ vanishes on $\pi_{i-1}$ and $\pi_i$; hence $\operatorname{Tot}_{\le q} C^\bullet \to \operatorname{Tot}_{\le q-1} C^\bullet$ is an isomorphism on $\pi_i$ and surjective on $\pi_{i+1}$. The Milnor exact sequence gives $\varprojlim_q^1 \pi_{i+1}(\operatorname{Tot}_{\le q} C^\bullet) = 0$ and $\pi_i(\operatorname{Tot} C^\bullet) \xrightarrow{\sim} \pi_i(\operatorname{Tot}_{\le n} C^\bullet)$.
\end{proof}

For any $E \in \mathcal{E}_R = \Shv_{\mathrm{hyp}}(\mathcal{T}_R, \tau_{\mathrm{cmp}}; \Sp)$, let $\widehat{E} \colon \mathrm{Pre}(\mathcal{T}_R)^{\mathrm{op}} \to \Sp$ denote its right Kan extension along the Yoneda embedding $j \colon \mathcal{T}_R \to \mathrm{Pre}(\mathcal{T}_R)$, given on $F \in \mathrm{Pre}(\mathcal{T}_R)$ by $\widehat{E}(F) = \varprojlim_{(V \to F) \in (\mathcal{T}_R)_{/F}} E(V)$, with $\widehat{E}(j(U)) \simeq E(U)$ for $U \in \mathcal{T}_R$. Recall that a simplicial morphism $P_\bullet \to F$ in $\mathrm{Pre}(\mathcal{T}_R)$ is a $\tau_{\mathrm{cmp}}$-hypercover if, for each $q \ge 0$, the relative matching map $P_q \to M_q(P_\bullet / F)$ becomes an effective epimorphism after $\tau_{\mathrm{cmp}}$-hypersheafification. For a pointed object $U$ over $X$ in $\mathcal{T}_R$ with section $s \colon X \to U$, let $\rho_n \colon \Omega_X^n U \to U$ be $\mathrm{id}_U$ for $n = 0$ and $\Omega_X^n U \to X \xrightarrow{s} U$ for $n \ge 1$ in $\mathrm{Pre}(\mathcal{T}_R)$.

\begin{theorem}\label{thm:derived-p-complete-descent}
For any pointed object $U$ over $X$ in $\mathcal{T}_R$, any $\tau_{\mathrm{cmp}}$-hypercover $P_\bullet \to U$, and each $n \ge 0$, the pullback $P_\bullet^{(n)} = P_\bullet \times_U \Omega_X^n U \to \Omega_X^n U$ along $\rho_n$ is a $\tau_{\mathrm{cmp}}$-hypercover in $\mathrm{Pre}(\mathcal{T}_R)$, and any hypercomplete sheaf $E \in \mathcal{E}_R = \Shv_{\mathrm{hyp}}(\mathcal{T}_R, \tau_{\mathrm{cmp}}; \Sp)$ with right Kan extension $\widehat{E} \colon \mathrm{Pre}(\mathcal{T}_R)^{\mathrm{op}} \to \Sp$ satisfies equivalences
\[
  \widehat{E}(\Omega_X^n U) \longrightarrow \operatorname{Tot}_{[q] \in \Delta} \widehat{E}\bigl(P_q^{(n)}\bigr),
  \qquad
  \bigl(\widehat{E}(\Omega_X^n U)\bigr)^\wedge \longrightarrow \operatorname{Tot}_{[q] \in \Delta} \bigl(\widehat{E}(P_q^{(n)})^\wedge\bigr).
\]
\end{theorem}
\begin{proof}
Pullback along $\rho_n$ commutes with the finite homotopy limits defining relative matching objects $M_q(P_\bullet / U) = \varprojlim_{([m] \hookrightarrow [q]) \in (\partial\Delta^q)^{\mathrm{op}}} P_m$. Since $\tau_{\mathrm{cmp}}$-hypersheafification is left exact and effective epimorphisms are stable under pullback in the $\tau_{\mathrm{cmp}}$-hypercomplete $\infty$-topos, $P_\bullet^{(n)} \to \Omega_X^n U$ is a $\tau_{\mathrm{cmp}}$-hypercover. By~\cite[Theorem~6.5.3.12]{HT}, $\tau_{\mathrm{cmp}}$-hypersheafification sends $\operatorname*{colim}_{[q] \in \Delta^{\mathrm{op}}} P_q^{(n)} \to \Omega_X^n U$ to an equivalence in $\Shv_{\mathrm{hyp}}(\mathcal{T}_R, \tau_{\mathrm{cmp}})$, so mapping into $E \in \mathcal{E}_R$ gives the first equivalence. Since smashing with the dualizable Moore spectrum $\mathbb{S}/p^m$ commutes with small limits, interchanging $\varprojlim_{m \ge 1}$ and $\operatorname{Tot}_{[q] \in \Delta}$ gives the $p$-completed equivalence.
\end{proof}

\subsection{Stabilization at a finite stage after excisive approximation}
\label{subsec:finite-stage-detection}

Let $A$ be a connective animated ring, $D \in \mathscr{D}_{\mathrm{cmp}}$, and $\sigma \in \mathscr{S}_D(A)$ a base diagram of quasi-compact, quasi-separated quasi-smooth derived schemes over $\Spec A$, with associated dual lift complex $K_\sigma^\vee \in \Mod_A$ of Proposition~\ref{prop:coherent-dual-lift-complex}(1). Write $P^1_{\Mod}$ and $P^1_{\mathrm{cot}}$ for the uncompleted excisive approximation functors of Definition~\ref{def:relative-tangential-linearization}(1).

\begin{lemma}\label{lem:lift-excisive-connectivity}
The module $K_\sigma^\vee \in \Mod_A$ is bounded below, and there is an equivalence in $\operatorname{Exc}_*^1(\Mod_A^{\ge 0}, \Sp)$
\[
P^1_{\Mod}\Bigl(\widetilde{\Sigma}^\infty\bigl(\operatorname{Lift}_{A \oplus (-)}(\sigma), 0_\sigma\bigr)\Bigr)(M) \longrightarrow K_\sigma^\vee \otimes_A M,
\]
and for any $M \in \Mod_A$,
\[
\conn(K_\sigma^\vee \otimes_A M) \ge \conn(M) - \max\bigl(0, -\conn(K_\sigma^\vee)\bigr).
\]
\end{lemma}
\begin{proof}
By Proposition~\ref{prop:coherent-dual-lift-complex}(1), $K_\sigma^\vee$ is bounded below and $(\operatorname{Lift}_{A \oplus (-)}(\sigma), 0_\sigma) \simeq \Omega^\infty(K_\sigma^\vee \otimes_A (-))$ in $\operatorname{Fun}_*(\Mod_A^{\ge 0}, \mathcal{S}_*)$, so~\cite[Example~6.1.2.9]{HA} gives the equivalence for $P^1_{\Mod}$. Set $\delta = \max\bigl(0, -\conn(K_\sigma^\vee)\bigr) \ge 0$. For any $M \in \Mod_A$ with $\conn(M) \ge m$ ($m \in \mathbb{Z}$), both $K_\sigma^\vee[\delta]$ and $M[-m]$ are connective over the connective ring $A$, so~\cite[Proposition~7.1.1.13]{HA} implies that $K_\sigma^\vee[\delta] \otimes_A M[-m] \simeq (K_\sigma^\vee \otimes_A M)[\delta - m]$ is connective, giving $\conn(K_\sigma^\vee \otimes_A M) \ge \conn(M) - \max\bigl(0, -\conn(K_\sigma^\vee)\bigr)$.
\end{proof}

Let $U = \Spec(A \oplus J)$ be a pointed derived affine scheme over
$\Spec A$ with connective augmentation ideal $J \in \NAlg_A$, and fix $k
\ge 1$. Under operadic Smith ideal
theory~\cite{WY2024},~\cite[Subsection~4.1]{KatoCalcNUAlg},
the loop scheme $\Omega_A^k U \simeq \Spec B^{(k)}$ with $B^{(k)} = A
\oplus J^{(k)}$ corresponds to the split commutative Smith ideal
$j^{(k)} \colon J^{(k)} \to B^{(k)}$ in
$\bigl(\operatorname{Ar}(\Mod_A^{\ge 0}), \square\bigr)$ with cofiber
$A$. For the canonical section $s_{k-1} \colon \Spec A \to
\Omega_A^{k-1} U$, the \v{C}ech--Amitsur nerve
$\check{C}_\bullet(s_{k-1}) \to \Omega_A^{k-1} U$ has $q$-th term
$\check{C}_q(s_{k-1}) = (\Spec A)^{\times_{\Omega_A^{k-1} U} (q+1)}
\simeq (\Omega_A^k U)^{\times_{\Spec A} q}$, represented contravariantly
by the cosimplicial augmented $A$-algebra
\[
[q] \longmapsto \bigl(B^{(k)}\bigr)^{\otimes_A q} = \bigl(A \oplus J^{(k)}\bigr)^{\otimes_A q}.
\]
In degree $q$, its augmentation corresponds to the $q$-fold pushout-product Smith ideal $(j^{(k)})^{\square q} \colon I^{(k)}_q \to (B^{(k)})^{\otimes_A q}$ with cofiber $A^{\otimes_A q} \simeq A$. The resulting cosimplicial augmentation Smith ideal $I^{(k)}_\bullet$ in $\Mod_A^{\ge 0}$ satisfies $I^{(k)}_0 = 0$ and
\[
I^{(k)}_q = \fib\Bigl(\bigl(A \oplus J^{(k)}\bigr)^{\otimes_A q} \to A\Bigr) \simeq \bigoplus_{\varnothing \neq S \subseteq \{1,\dots,q\}} \bigl(J^{(k)}\bigr)^{\otimes_A |S|} \quad (q \ge 1),
\]
where the summand indexed by $S$ takes $J^{(k)}$ in the tensor positions $j \in S$ and $A$ in the positions $j \notin S$, and each codegeneracy $s^i$ ($0 \le i \le q-1$) applies the augmentation $A \oplus J^{(k)} \to A$ to the $(i+1)$-th tensor factor. Applying the excisive approximation $P^1_{\Mod}\bigl(\widetilde{\Sigma}^\infty(\operatorname{Lift}_{A \oplus (-)}(\sigma), 0_\sigma)\bigr) \simeq K_\sigma^\vee \otimes_A (-)$ of Lemma~\ref{lem:lift-excisive-connectivity} to the underlying cosimplicial Smith ideal $I^{(k)}_\bullet$ defines the cosimplicial spectrum $C^\bullet$ with $q$-th term
\[
C^q = K_\sigma^\vee \otimes_A I^{(k)}_q = K_\sigma^\vee \otimes_A \fib\Bigl(\bigl(A \oplus J^{(k)}\bigr)^{\otimes_A q} \to A\Bigr).
\]

\begin{lemma}\label{lem:normalized-fiber-slope}
For each $q \ge 1$, there is an equivalence
\[
\fib(C^q \to M^{q-1} C^\bullet) \longrightarrow K_\sigma^\vee \otimes_A \bigl(J^{(k)}\bigr)^{\otimes_A q},
\]
with $\conn\bigl(J^{(k)}\bigr) = \conn\bigl(Q_A(J)\bigr) + k$, and the connectivity bound
\[
\conn\bigl(\fib(C^q \to M^{q-1} C^\bullet)\bigr) \ge \bigl(\conn(Q_A(J)) + k\bigr)q - \max\bigl(0, -\conn(K_\sigma^\vee)\bigr).
\]
\end{lemma}
\begin{proof}
Exactness of $K_\sigma^\vee \otimes_A (-)$ and the cosimplicial Dold--Kan correspondence~\cite[Theorem~1.2.4.1]{HA} identify $\fib(C^q \to M^{q-1} C^\bullet)$ with $K_\sigma^\vee \otimes_A \bigcap_{i=0}^{q-1}\fib(s^i)$ on $I^{(k)}_\bullet$. On the summand indexed by $\varnothing \neq S \subseteq \{1,\dots,q\}$, $s^i$ vanishes if $i+1 \in S$ and is an equivalence onto the corresponding summand in degree $q-1$ if $i+1 \notin S$; intersecting the fibers over $0 \le i \le q-1$ leaves only $S = \{1,\dots,q\}$, namely $(J^{(k)})^{\otimes_A q}$.

Next, $\pi_0(B^{(1)}) \cong \pi_0(A) \otimes_{\pi_0(A)\oplus\pi_0(J)} \pi_0(A) \cong \pi_0(A)$ gives $\conn(J^{(1)}) \ge 1$, and the normalized bar filtration on $B^{(n+1)} = A \otimes_{B^{(n)}} A$ has $r$-th graded piece $(J^{(n)}[1])^{\otimes_A r}$ ($r \ge 1$), giving $\conn(J^{(k)}) \ge k \ge 1$ by induction. For $I = J^{(k)} \not\simeq 0$ with $c = \conn(I) \ge 1$ and $\epsilon \colon B = A \oplus I \to A$, $\fib(I \to I \otimes_B A) \simeq I \otimes_B I$ is $2c$-connective, and the Hurewicz estimate~\cite[Theorem~7.4.3.1]{HA} shows that $(I \otimes_B A)[1] \to L_{A/B} \simeq Q_A(I)[1]$ has $(2c+1)$-connective fiber. Thus $\fib(I \to Q_A(I))$ is $2c$-connective with $2c > c$, so $\conn(J^{(k)}) = \conn(Q_A(J^{(k)})) = \conn(Q_A(J)) + k$ by Subsection~\ref{subsec:loop-charts-bar}, and Lemma~\ref{lem:lift-excisive-connectivity} gives the connectivity bound.
\end{proof}

\begin{theorem}\label{thm:finite-stage-detection}
Let $U = \Spec(A \oplus J)$ be a pointed derived affine scheme over $\Spec A$ with connective augmentation ideal $J \in \NAlg_A$, let $\sigma \in \mathscr{S}_D(A)$ ($D \in \mathscr{D}_{\mathrm{cmp}}$) have dual lift complex $K_\sigma^\vee \in \Mod_A$, and for $k \ge 1$ let $C^\bullet$ be the cosimplicial spectrum with $C^q = K_\sigma^\vee \otimes_A \fib\bigl((A \oplus J^{(k)})^{\otimes_A q} \to A\bigr)$ constructed above. If $\conn(Q_A(J)) + k > 1$, then for each $i \in \mathbb{Z}$ and each integer $n$ in the stated range, the truncation maps induce isomorphisms
\[
\pi_i(\operatorname{Tot} C^\bullet) \longrightarrow \pi_i(\operatorname{Tot}_{\le n} C^\bullet)
\quad \text{for } n \ge \max\Biggl(0, \frac{i + \max\bigl(0, -\conn(K_\sigma^\vee)\bigr) + 2}{\conn(Q_A(J)) + k - 1}\Biggr),
\]
\[
\pi_i\bigl((\operatorname{Tot} C^\bullet)^\wedge\bigr) \longrightarrow \pi_i\bigl((\operatorname{Tot}_{\le n} C^\bullet)^\wedge\bigr)
\quad \text{for } n \ge \max\Biggl(0, \frac{i + \max\bigl(0, -\conn(K_\sigma^\vee)\bigr) + 3}{\conn(Q_A(J)) + k - 1}\Biggr).
\]
\end{theorem}
\begin{proof}
Set $c = \conn(Q_A(J)) + k > 1$ and $d = \max\bigl(0, -\conn(K_\sigma^\vee)\bigr) \ge 0$. By Lemma~\ref{lem:normalized-fiber-slope}, $\conn\bigl(\fib(C^q \to M^{q-1} C^\bullet)\bigr) \ge cq - d$ for all $q \ge 1$. Applying Lemma~\ref{lem:tot-tower-connectivity} with $q_0 = 1$ gives
\[
\conn\bigl(\fib(\operatorname{Tot}_{\le q} C^\bullet \to \operatorname{Tot}_{\le q-1} C^\bullet)\bigr) \ge (c - 1)q - d \quad \text{for all } q \ge 1,
\]
and implies the uncompleted isomorphism for $n \ge \max\bigl(0, \frac{i + d + 2}{c - 1}\bigr)$. Since smashing with the dualizable Moore spectrum $\mathbb{S}/p^m$ commutes with small limits, derived $p$-completion $(-)^\wedge = \varprojlim_{m \ge 1} (- \wedge \mathbb{S}/p^m)$ preserves fiber sequences and the totalization tower limit, and the Milnor exact sequence gives
\begin{align*}
\conn\Bigl(\fib\bigl((\operatorname{Tot}_{\le q} C^\bullet)^\wedge \to (\operatorname{Tot}_{\le q-1} C^\bullet)^\wedge\bigr)\Bigr)
&\ge \conn\bigl(\fib(\operatorname{Tot}_{\le q} C^\bullet \to \operatorname{Tot}_{\le q-1} C^\bullet)\bigr) - 1 \\
&\ge (c - 1)q - (d + 1)
\end{align*}
for all $q \ge 1$. The same Milnor exact sequence argument on the tower $(\operatorname{Tot}_{\le q} C^\bullet)^\wedge$ as in Lemma~\ref{lem:tot-tower-connectivity} implies the $p$-completed isomorphism for $n \ge \max\bigl(0, \frac{i + d + 3}{c - 1}\bigr)$.
\end{proof}

\begin{example}\label{ex:slope-loop-order}
In Theorem~\ref{thm:finite-stage-detection}, since $Q_A(J) \in \Mod_A^{\ge 0}$ always satisfies $\conn(Q_A(J)) \ge 0$, the slope condition $\conn(Q_A(J)) + k > 1$ holds unconditionally at loop order $k = 2$, and already at loop order $k = 1$ whenever $Q_A(J)$ is $1$-connective ($\conn(Q_A(J)) \ge 1$). Since derived $(p,I_R)$-completion and almostification are exact and preserve connectivity (Subsection~\ref{subsec:loop-charts-bar}), the same slope bound and stabilization at a finite stage apply to the derived complete and almost dual lift complexes of Proposition~\ref{prop:coherent-dual-lift-complex}(2).
\end{example}

\subsection{The neutral locus}
\label{subsec:neutral-locus}

\begin{definition}\label{def:neutral-locus}
The \emph{neutral locus} $\mathcal{T}_R^{\mathrm{neu}} \subseteq \mathcal{T}_R$ is the full subcategory of objects $U \in \mathcal{T}_R$ such that, for any $V \to U$ in $\mathcal{T}_R$ and any lift torsor $\operatorname{Lift}(\sigma)$ associated with a diagram $\sigma$ in $\mathscr{D}_{\mathrm{cmp}}$ over $V$, $S(\operatorname{Lift}(\sigma))$ contains a $\tau_{\mathrm{geom}}$-covering sieve.
\end{definition}

\begin{proposition}\label{prop:neutral-locus-equivalence}
For any $U \in \mathcal{T}_R^{\mathrm{neu}}$, the topologies $\tau_{\mathrm{geom}}$ and $\tau_{\mathrm{pld}}$ generate identical hypercomplete $\infty$-topoi on $(\mathcal{T}_R)_{/U}$, inducing an equivalence
\[
  \Shv_{\mathrm{hyp}}\bigl((\mathcal{T}_R)_{/U}, \tau_{\mathrm{geom}}; \Sp\bigr)
  \longrightarrow
  \Shv_{\mathrm{hyp}}\bigl((\mathcal{T}_R)_{/U}, \tau_{\mathrm{pld}}; \Sp\bigr).
\]
\end{proposition}
\begin{proof}
By Definition~\ref{def:neutral-locus} and Definition~\ref{def:pld-topology}, for any $U \in \mathcal{T}_R^{\mathrm{neu}}$ and any $V \in (\mathcal{T}_R)_{/U}$, any section sieve $S(\operatorname{Lift}(\sigma))$ over $V$ contains a $\tau_{\mathrm{geom}}$-covering sieve in $(\mathcal{T}_R)_{/U}$, while $\tau_{\mathrm{geom}} \subseteq \tau_{\mathrm{pld}}$. Hence $\tau_{\mathrm{geom}}$ and $\tau_{\mathrm{pld}}$ generate identical covering sieves on $(\mathcal{T}_R)_{/U}$.
\end{proof}

\begin{example}\label{ex:neutral-locus-qrsp}
The neutral locus $\mathcal{T}_R^{\mathrm{neu}}$ includes all objects $U \in \mathcal{T}_R$ over which the selected diagrams in $\mathscr{D}_{\mathrm{cmp}}$ consist of single affine smooth objects, single spans of affine smooth objects, or single affine quasi-smooth objects by Corollary~\ref{cor:honest-chart-lift-existence}, and more generally all $U \in \mathcal{T}_R$ over which completely quasisyntomic prism covers modeled on QRSP refinements~\cite[Section~7]{BhattScholze2017Prisms} split the lift obstructions for diagrams in $\mathscr{D}_{\mathrm{cmp}}$; over any such $U$, Proposition~\ref{prop:neutral-locus-equivalence} and Proposition~\ref{thm:localized-span-universal-property} imply that any $\tau_{\mathrm{geom}}$-hypercomplete additive span theory coincides with its $\tau_{\mathrm{pld}}$-localization.
\end{example}
\section{Universal \texorpdfstring{$\mathbb G_a$}{Ga}- and \texorpdfstring{$\mathbb G_m$}{Gm}-oriented span motives and tilting equivalence}
\label{sec:oriented-span-motives}

Following~\cite{AHIConnerFloyd} and Subsection~\ref{subsec:chern-law-refinements}, we construct universal $\mathbb G_a$- and $\mathbb G_m$-oriented derived span motives from graded spectral correspondence categories, compare the $\mathbb G_a$-oriented theory with classical rational Chow motives on $\pi_0$ (Theorem~\ref{thm:p1-chow-comparison}), and prove tilting equivalences over a perfectoid base (Theorem~\ref{thm:tilting-oriented-motives}).

\subsection{Graded derived spans and universal \texorpdfstring{$\mathbb G_a$}{Ga}-oriented motives}
\label{subsec:graded-spans-realization}

Fix a field $k$, let $\mathbf{SmProj}_k$ denote the category of equidimensional smooth projective schemes over $k$, and let $\mathbf{S}_k$ denote the class of projective quasi-smooth morphisms of derived $k$-schemes. For $X, Y \in \mathbf{SmProj}_k$, a span $[X \xleftarrow{p} W \xrightarrow{g} Y]$ with $p \in \mathbf{S}_k$ has locally constant relative virtual dimension $\operatorname{vdim}(p) = \operatorname{rk}(L_{W/X}) = \operatorname{vdim}(W) - \dim X$, where $\operatorname{rk}$ denotes the Euler characteristic rank of a perfect complex and $\operatorname{vdim}(W)$ denotes the virtual dimension. This relative virtual dimension decomposes the $\mathbb E_\infty$-space of spans as $\Corr_{\mathbf{S}_k}(X, Y) \simeq \coprod_{d \in \mathbb{Z}} \Corr_{\mathbf{S}_k}^d(X, Y)$ compatibly with pullback composition.

\begin{definition}
\label{def:graded-span-category}
Write $H\mathbb Q$ for the rational Eilenberg--MacLane spectrum. For any $X\in\mathbf{SmProj}_k$ and $r\in\mathbb Z$, write $\mathbf M(X)(r)[r]$ for the formal degree-$r$ shifted object associated with $X$, with $\mathbf 1(r)[r]=\mathbf M(\Spec k)(r)[r]$. The \emph{graded spectral span category} $\mathcal A_k$ is the $\mathbb Q$-linear symmetric monoidal spectral category with objects $\mathbf M(X)(r)[r]$ ($X\in\mathbf{SmProj}_k$, $r\in\mathbb Z$), symmetric monoidal product $\mathbf M(X)(r)[r]\otimes\mathbf M(Y)(s)[s] = \mathbf M(X\times_k Y)(r+s)[r+s]$, and mapping spectra $\mathcal A_k^{r,s}(X,Y) = \mathcal A_k\bigl(\mathbf M(X)(r)[r],\mathbf M(Y)(s)[s]\bigr) = \Corr_{\mathbf S_k}^{r-s}(X,Y)^{\mathrm{gp}}\otimes H\mathbb Q$ graded by relative virtual dimension. Its stable presheaf envelope is $\mathcal M_0(k)=\operatorname{Fun}_{\Sp}(\mathcal A_k^{\mathrm{op}},\Mod_{H\mathbb Q})$, equipped with the Yoneda embedding $m_0\colon\mathcal A_k\to\mathcal M_0(k)$.
\end{definition}

\begin{definition}
\label{def:derived-chow-motives}
For $X,Y\in\mathbf{SmProj}_k$ and $r,s\in\mathbb Z$, write $(-)^{\mathrm{gp}}_{\mathbb Q} = (-)^{\mathrm{gp}}\otimes H\mathbb Q$. We introduce three classes of generating relation morphisms in $\mathcal A_k$:
\begin{enumerate}[(1)]
\item Let $\mathscr F_{\mathbb P^1}^{r,s}(X,Y)$ be the maximal $\infty$-groupoid of spans $\Sigma = [X\times\mathbb P^1 \xleftarrow{p} W \xrightarrow{g} Y]$ with $p\in\mathbf S_k$ of relative virtual dimension $r-s$. Derived endpoint pullbacks $e_t$ along $i_t\colon X\hookrightarrow X\times\mathbb P^1$ ($t\in\{0,\infty\}$) define the \emph{$\mathbb P^1$-cobordism boundary morphism} $\partial_{X,Y}^{r,s} = (e_0)^{\mathrm{gp}}_{\mathbb Q}-(e_\infty)^{\mathrm{gp}}_{\mathbb Q} \colon \mathscr F_{\mathbb P^1}^{r,s}(X,Y)^{\mathrm{gp}}_{\mathbb Q} \to \mathcal A_k^{r,s}(X,Y)$.
\item Let $\sigma=[X\xleftarrow{p}Z\xrightarrow{g}Y]$ be a span with $p\in\mathbf S_k$ of relative virtual dimension $r-s+1$, and let $\mathcal L$ be an invertible $\mathcal O_Z$-module with zero section $s_0\colon Z\hookrightarrow\mathbb V_Z(\mathcal L)$ of relative virtual dimension $-1$ and derived self-intersection $\pi_{\mathcal L}\colon Z\times_{\mathbb V_Z(\mathcal L)}Z\to Z$. Set $\widetilde c_1(\mathcal L)\cap\sigma = [X\xleftarrow{p\circ\pi_{\mathcal L}}Z\times_{\mathbb V_Z(\mathcal L)}Z\xrightarrow{g\circ\pi_{\mathcal L}}Y]\in\mathcal A_k^{r,s}(X,Y)$; at $\sigma=\mathrm{id}_X$, this recovers the \emph{first Chern span} $\widetilde c_1(\mathcal L)=\widetilde c_1(\mathcal L)\cap\mathrm{id}_X\in\mathcal A_k^{0,1}(X,X)$. On the maximal $\infty$-groupoids of pairs $(\sigma,\mathcal O_Z)$ and triples $(\sigma,\mathcal L_1,\mathcal L_2)$ with invertible $\mathcal O_Z$-modules $\mathcal L_1,\mathcal L_2$, we define the \emph{additive Chern relation morphisms} into $\mathcal A_k^{r,s}(X,Y)$ with components $\mathcal R_{\mathbb G_a}^0(\sigma) = \widetilde c_1(\mathcal O_Z)\cap\sigma$ and $\mathcal R_{\mathbb G_a}(\sigma;\mathcal L_1,\mathcal L_2) = \widetilde c_1(\mathcal L_1\otimes\mathcal L_2)\cap\sigma - \widetilde c_1(\mathcal L_1)\cap\sigma - \widetilde c_1(\mathcal L_2)\cap\sigma$.
\item Let $\mathscr B_k^{r,s}(X,Y)$ be the maximal $\infty$-groupoid of pairs $(\sigma,\mathcal E)$ consisting of a span $\sigma=[X\xleftarrow{p}Z\xrightarrow{g}Y]$ with $p\in\mathbf S_k$ of relative virtual dimension $r-s-(n-1)$ and a finite locally free module $\mathcal E$ on $Z$ of rank $n\ge 2$ with projective bundle $\pi_{\mathcal E}\colon\mathbb P_Z(\mathcal E)\to Z$. Setting $\mathrm{ev}_{\mathbb P(\mathcal E)}(\sigma,\mathcal E)=\pi_{\mathcal E}^*\sigma=[X\xleftarrow{p\circ\pi_{\mathcal E}}\mathbb P_Z(\mathcal E)\xrightarrow{g\circ\pi_{\mathcal E}}Y]$ defines the \emph{projective bundle vanishing morphism} $\rho_{X,Y}^{r,s} = (\mathrm{ev}_{\mathbb P(\mathcal E)})^{\mathrm{gp}}_{\mathbb Q}\colon\mathscr B_k^{r,s}(X,Y)^{\mathrm{gp}}_{\mathbb Q}\to\mathcal A_k^{r,s}(X,Y)$.
\end{enumerate}
Let $\mathcal I_{\mathbb G_a}$ be the two-sided tensor ideal in $\mathcal A_k$ generated by $\partial_{X,Y}^{r,s}$, $\mathcal R_{\mathbb G_a}^0(\sigma)$, $\mathcal R_{\mathbb G_a}(\sigma;\mathcal L_1,\mathcal L_2)$, and $\rho_{X,Y}^{r,s}$ under span composition and external products, and let $q_{\mathbb G_a}\colon\mathcal A_k \to \mathcal A_{k,\mathbb G_a} = \mathcal A_k / \mathcal I_{\mathbb G_a}$ be the universal symmetric monoidal $\mathbb Q$-linear spectral quotient coning off $\mathcal I_{\mathbb G_a}$ (Definition~\ref{def:chern-law-span-quotient}). Its stable presheaf envelope is the category of \emph{universal $\mathbb G_a$-oriented derived Chow motives} $\mathcal{DM}_{\mathrm{Chow}}(k)_\mathbb Q = \operatorname{Fun}_{\Sp}(\mathcal A_{k,\mathbb G_a}^{\mathrm{op}},\Mod_{H\mathbb Q})$, with representable motives $M_{\mathrm{Chow}}(X)(r)[r]$ spanning $\mathcal A_{k,\mathbb G_a}\subset\mathcal{DM}_{\mathrm{Chow}}(k)_{\mathbb Q}$. Since $\mathcal A_k^{r,s}(X,Y)$ and the source spectra in (1)--(3) are connective, $\pi_0\mathcal A_{k,\mathbb G_a}^{r,s}(X,Y)$ is the algebraic quotient of $\pi_0\mathcal A_k^{r,s}(X,Y)$ by the two-sided ideal $R_{X,Y}^{r,s}$ generated by the classes of $\partial_{X,Y}^{r,s}$, $\mathcal R_{\mathbb G_a}^0(\sigma)$, $\mathcal R_{\mathbb G_a}(\sigma;\mathcal L_1,\mathcal L_2)$, and $\rho_{X,Y}^{r,s}$.
\end{definition}

\begin{proposition}
\label{prop:universal-ga-motives}
For any $H\mathbb Q$-linear stable presentable symmetric monoidal $\infty$-category $\mathcal C$, precomposition with $(q_{\mathbb G_a})_!\colon\mathcal M_0(k)\to\mathcal{DM}_{\mathrm{Chow}}(k)_\mathbb Q$ identifies $\operatorname{Fun}_{\mathrm{Mod}_{H\mathbb Q}}^{\mathrm{L},\otimes}\bigl(\mathcal{DM}_{\mathrm{Chow}}(k)_\mathbb Q,\mathcal C\bigr)$ with the $\infty$-category of $H\mathbb Q$-linear symmetric monoidal functors $E\colon\mathcal M_0(k)\to\mathcal C$ preserving small colimits and equipped with nullhomotopies of the images under $E\circ m_0$ of $\partial_{X,Y}^{r,s}$, $\mathcal R_{\mathbb G_a}^0(\sigma)$, $\mathcal R_{\mathbb G_a}(\sigma;\mathcal L_1,\mathcal L_2)$, and $\rho_{X,Y}^{r,s}$.
\end{proposition}
\begin{proof}
By~\cite[Proposition~4.8.1.10 and Corollary~4.8.1.12]{HA}, restriction along the Yoneda embeddings identifies $\operatorname{Fun}_{\mathrm{Mod}_{H\mathbb Q}}^{\mathrm{L},\otimes}(\mathcal M_0(k),\mathcal C)\simeq\operatorname{Fun}_{\mathrm{Mod}_{H\mathbb Q}}^{\otimes}(\mathcal A_k,\mathcal C)$ and $\operatorname{Fun}_{\mathrm{Mod}_{H\mathbb Q}}^{\mathrm{L},\otimes}(\mathcal{DM}_{\mathrm{Chow}}(k)_{\mathbb Q},\mathcal C)\simeq\operatorname{Fun}_{\mathrm{Mod}_{H\mathbb Q}}^{\otimes}(\mathcal A_{k,\mathbb G_a},\mathcal C)$, where the pushout defining $\mathcal A_{k,\mathbb G_a}=\mathcal A_k/\mathcal I_{\mathbb G_a}$ identifies extensions along $q_{\mathbb G_a}$ with nullhomotopies on the generating morphisms of $\mathcal I_{\mathbb G_a}$.
\end{proof}

\subsection{Virtual cycle class realization and comparison with classical Chow motives}
\label{subsec:classical-recovery-pi0}

Assume that $k$ is a perfect field, and denote by $(-)_{\mathrm{cl}}$ the underlying classical scheme, morphism, or line bundle of a derived object, and by $CH_d(-)_{\mathbb Q}$ and $CH^c(-)_{\mathbb Q}$ the rational homological and cohomological Chow groups. For $X,Y\in\mathbf{SmProj}_k$ and $r,s\in\mathbb Z$, set $c=\dim Y+s-r$ and $d=\dim X+r-s=\dim(X\times Y)-c$. Any span $\sigma=[X\xleftarrow{p}W\xrightarrow{g}Y]$ with $p\in\mathbf S_k$ of relative virtual dimension $r-s$ has projective quasi-smooth structure morphism $\pi_W\colon W\to\Spec k$ of virtual dimension $d$ and projective morphism $f_\sigma=(p,g)\colon W\to X\times Y$. Virtual Gysin pullback $\pi_W^!$ along $\pi_W$ defines the virtual fundamental class $[W]^{\mathrm{vir}}=\pi_W^![\Spec k]\in CH_d(W_{\mathrm{cl}})_\mathbb Q$~\cite{MR4092773,LS16}, and proper pushforward along $f_{\sigma,\mathrm{cl}}$ induces the \emph{virtual cycle class map} $\operatorname{cl}_{X,Y}^{r,s}\colon\pi_0\mathcal A_k^{r,s}(X,Y)\to CH^c(X\times Y)_\mathbb Q$, $[\sigma]\mapsto(f_{\sigma,\mathrm{cl}})_*[W]^{\mathrm{vir}}$. Functoriality and base change for virtual Gysin pullbacks along derived Cartesian squares~\cite[Section~3.2]{MR4092773} make $\operatorname{cl}_{X,Y}^{r,s}$ compatible with span composition and external products.

\begin{proposition}
\label{prop:pi0-cycle-surjectivity}
Let $k$ be a perfect field, let $X,Y\in\mathbf{SmProj}_k$, let $r,s\in\mathbb Z$, and set $c=\dim Y+s-r$. The virtual cycle class realization
\[
\operatorname{cl}_{X,Y}^{r,s}\colon\pi_0\mathcal A_k^{r,s}(X,Y)\longrightarrow CH^c(X\times Y)_\mathbb Q
\]
is surjective.
\end{proposition}
\begin{proof}
It suffices to consider an integral closed subvariety $W\subset X\times Y$ of codimension $c$. Resolution of singularities in characteristic zero and de Jong alterations in positive characteristic~\cite[Theorem~4.1]{dJ96} provide a smooth projective variety $\widetilde W$ over $k$ and a projective generically finite morphism $\mu\colon\widetilde W\to W$ of degree $m>0$. Since $\widetilde W$ and $X\times Y$ are smooth projective over $k$, the composite $f=(p,g)\colon\widetilde W\to W\hookrightarrow X\times Y$ has $p\in\mathbf S_k$ of relative virtual dimension $\dim\widetilde W-\dim X=r-s$, defining a span $\widetilde\sigma=[X\xleftarrow{p}\widetilde W\xrightarrow{g}Y]$ in $\mathcal A_k^{r,s}(X,Y)$ with $\operatorname{cl}_{X,Y}^{r,s}([\widetilde\sigma])=f_*[\widetilde W]=m[W]$ in $CH^c(X\times Y)_\mathbb Q$. Multiplying by $m^{-1}\in\mathbb Q$ gives a preimage of $[W]$.
\end{proof}

For any $X,Y\in\mathbf{SmProj}_k$ and $r,s\in\mathbb Z$, set $c=\dim Y+s-r$, let $K_{X,Y}^{r,s}=\ker(\operatorname{cl}_{X,Y}^{r,s})$, and let $R_{X,Y}^{r,s}\subseteq\pi_0\mathcal A_k^{r,s}(X,Y)$ be the two-sided correspondence ideal generated by the $\mathbb P^1$-cobordism boundary morphisms $\partial_{X,Y}^{r,s}$, additive Chern relations $\mathcal R_{\mathbb G_a}^0(\sigma)$ and $\mathcal R_{\mathbb G_a}(\sigma;\mathcal L_1,\mathcal L_2)$, and projective bundle vanishing morphisms $\rho_{X,Y}^{r,s}$.

\begin{proposition}
\label{prop:lazard-ga-compatibility}
Let $k$ be a field of characteristic $0$, let $X,Y\in\mathbf{SmProj}_k$ and $r,s\in\mathbb Z$, and set $V=X\times Y$, $c=\dim Y+s-r$, and $d=\dim V-c$.
\begin{enumerate}[\rm (1)]
\item The identity on projective quasi-smooth generators annihilates all derived double-point cobordism relations in $\pi_0\mathcal A_k^{r,s}(X,Y)/R_{X,Y}^{r,s}$ and induces a $\mathbb Q$-linear surjection
\[
\Omega_d(V)_{\mathbb Q} \xrightarrow{\ \sim\ } d\Omega_d(V)_{\mathbb Q} \longrightarrow \pi_0\mathcal A_k^{r,s}(X,Y)/R_{X,Y}^{r,s}.
\]
\item At $X=Y=\Spec k$, the induced graded ring homomorphism $\mathbb L_*\otimes\mathbb Q\xrightarrow{\sim}\Omega_*(\Spec k)_{\mathbb Q}\to\bigoplus_{d\in\mathbb Z}\pi_0\mathcal A_k^{d,0}(\Spec k,\Spec k)/R_{\Spec k,\Spec k}^{d,0}$ factors through the additive formal group law specialization $\lambda_{\mathbb G_a}\colon\mathbb L_*\otimes\mathbb Q\twoheadrightarrow\mathbb Q$ with kernel $\mathbb L_{>0}\otimes\mathbb Q$. For any $X,Y\in\mathbf{SmProj}_k$, the surjection of~(1) annihilates $\mathbb L_{>0}\cdot\Omega_*(V)_{\mathbb Q}$ and descends to a $\mathbb Q$-linear surjection
\[
CH^c(V)_{\mathbb Q} \cong \Omega_d(V)_{\mathbb Q}/\mathbb L_{>0}\Omega_*(V)_{\mathbb Q} \longrightarrow \pi_0\mathcal A_k^{r,s}(X,Y)/R_{X,Y}^{r,s}.
\]
\end{enumerate}
\end{proposition}
\begin{proof}
Since $X$ and $V=X\times Y$ are smooth projective over $k$, $\Omega^1_{X/k}$ and $\Omega^1_{V/k}$ are locally free in degree~$0$, so the cotangent transitivity triangles for $W\to X\to\Spec k$ and $W\to V\to\Spec k$ identify spans $[X\xleftarrow{p}W\xrightarrow{g}Y]$ with $p\in\mathbf S_k$ of relative virtual dimension $r-s$ with projective quasi-smooth morphisms $f=(p,g)\colon W\to V$ of virtual dimension $d=\dim V-c$. Let $\sigma_Z=[X\xleftarrow{p_Z}Z\xrightarrow{g_Z}Y]$ be a span with $\operatorname{vdim}(p_Z)=r-s+1$, let $\mathcal L$ be an invertible $\mathcal O_Z$-module carrying a section $s\in\Gamma(Z,\mathcal L)$, and let $t_0\in\Gamma(\mathbb P^1,\mathcal O_{\mathbb P^1}(1))$ be a homogeneous coordinate vanishing at $0\in\mathbb P^1$ and nonzero at $\infty\in\mathbb P^1$. Writing $V^{\mathrm{der}}(s) = Z\times_{s,\mathbb V_Z(\mathcal L),0} Z$ for the derived zero locus of $s$, the boundary relation $\partial_{X,Y}^{r,s}$ applied to $V^{\mathrm{der}}(s\boxtimes t_0)$ has derived fibers $Z\times_{\mathbb V_Z(\mathcal L)}Z$ at $0$ and $V^{\mathrm{der}}(s)$ at $\infty$, identifying $[V^{\mathrm{der}}(s)]=[\widetilde c_1(\mathcal L)\cap\sigma_Z]$ in $\pi_0\mathcal A_k^{r,s}(X,Y)/R_{X,Y}^{r,s}$.

For any projective quasi-smooth degeneration $W\to V\times\mathbb P^1$ with span $\sigma_W$ over $(X,Y)$ whose derived special fiber $W_\infty=V^{\mathrm{der}}(s_1\otimes s_2)$ is a sum of virtual Cartier divisors $A=V^{\mathrm{der}}(s_1)$ and $B=V^{\mathrm{der}}(s_2)$ ($s_i\in\Gamma(W,\mathcal L_i)$) with intersection $D=A\times_W B$, the boundary relation $\partial_{X,Y}^{r,s}$, the additive Chern relation $\mathcal R_{\mathbb G_a}$, and the projective bundle vanishing $\rho_{X,Y}^{r,s}$ give
\begin{align*}
[W_0] = [W_\infty]
&= [\widetilde c_1(\mathcal L_1\otimes\mathcal L_2)\cap\sigma_W]
= [\widetilde c_1(\mathcal L_1)\cap\sigma_W] + [\widetilde c_1(\mathcal L_2)\cap\sigma_W] \\
&= [A] + [B]
= [A] + [B] - [\mathbb P_D(\mathcal O_D(A)\oplus\mathcal O_D)],
\end{align*}
annihilating all derived double-point cobordism relations~\cite{zbMATH05550622,LS16,KR19,Ann21}. Thus the identity on projective quasi-smooth generators descends to a $\mathbb Q$-linear surjection $d\Omega_d(V)_\mathbb Q\to\pi_0\mathcal A_k^{r,s}(X,Y)/R_{X,Y}^{r,s}$ from Lowrey--Sch\"{u}rg's derived algebraic cobordism, and composing with the isomorphism $\Omega_d(V)_\mathbb Q\to d\Omega_d(V)_\mathbb Q$ in characteristic zero~\cite[Theorem~5.12]{LS16} proves~(1).

For~(2), writing $\mathbb L_*$ for the Lazard ring and $\mathbb L_{>0}$ for its ideal of elements of positive degree, $\mathbb L_*\otimes\mathbb Q\cong\Omega_*(\Spec k)_{\mathbb Q}$ is generated in positive degrees by $[\mathbb P^{n-1}]$ for $n\ge 2$~\cite[Section~4.1]{LevineMorel2007}. For any span $[X\xleftarrow{p}Z\xrightarrow{g}Y]$, the relation $\rho_{X,Y}^{r,s}$ kills $[Z\times\mathbb P^{n-1}\to Z\to V]=[\mathbb P^{n-1}]\cdot[Z\to V]$ via $\mathcal E=\mathcal O_Z^{\oplus n}$. At $X=Y=\Spec k$, this annihilates $\mathbb L_{>0}\otimes\mathbb Q$, so the induced ring homomorphism factors through $\lambda_{\mathbb G_a}\colon\mathbb L_*\otimes\mathbb Q\twoheadrightarrow\mathbb Q$. For general $X,Y\in\mathbf{SmProj}_k$, the surjection of~(1) therefore annihilates $\mathbb L_{>0}\cdot\Omega_*(V)_\mathbb Q$ and descends via $CH^c(V)_\mathbb Q\cong\Omega_d(V)_\mathbb Q/\mathbb L_{>0}\Omega_*(V)_\mathbb Q$~\cite[Theorem~4.5.1]{LevineMorel2007} to the stated surjection.
\end{proof}

\begin{theorem}
\label{thm:p1-chow-comparison}
Let $k$ be a perfect field. For any $X,Y\in\mathbf{SmProj}_k$ and $r,s\in\mathbb Z$, we have $R_{X,Y}^{r,s}\subseteq K_{X,Y}^{r,s}$, and $\operatorname{cl}_{X,Y}^{r,s}$ induces a natural short exact sequence
\[
\begin{aligned}
0&\longrightarrow K_{X,Y}^{r,s}/R_{X,Y}^{r,s}
\longrightarrow \pi_0\operatorname{Map}_{\mathcal{DM}_{\mathrm{Chow}}(k)_\mathbb Q}
\bigl(M_{\mathrm{Chow}}(X)(r)[r],M_{\mathrm{Chow}}(Y)(s)[s]\bigr)\\
&\longrightarrow CH^c(X\times Y)_\mathbb Q\longrightarrow0.
\end{aligned}
\]
When $\operatorname{char}(k)=0$, the equality $R_{X,Y}^{r,s}=K_{X,Y}^{r,s}$ holds for all $X,Y\in\mathbf{SmProj}_k$ and $r,s\in\mathbb Z$, the induced virtual cycle class map is an isomorphism, and the additive idempotent completion of the homotopy category $\operatorname{h}(\mathcal A_{k,\mathbb G_a})$ is equivalent to the classical category of homological Grothendieck Chow motives over $k$ with rational coefficients~\cite{Manin1968,Scholl1994}.
\end{theorem}
\begin{proof}
Set $V=X\times Y$ and $d=\dim V-c$. For any $\Sigma=[X\times\mathbb P^1\xleftarrow{p}W\xrightarrow{g}Y]\in\mathscr F_{\mathbb P^1}^{r,s}(X,Y)$ with proper map $F=(p,g)\colon W\to V\times\mathbb P^1$ and derived fibers $j_t\colon W_t\to W$ over $i_t\colon\Spec k\hookrightarrow\mathbb P^1$ ($t\in\{0,\infty\}$), functoriality and base change of virtual Gysin pullbacks~\cite[Section~3.2]{MR4092773} give $[W_t]^{\mathrm{vir}}=j_t^![W]^{\mathrm{vir}}$ and $(f_{\Sigma_t,\mathrm{cl}})_*j_t^!=(\mathrm{id}_V\times i_t)^!F_{\mathrm{cl}*}$. Since $(\mathrm{id}_V\times i_0)^!=(\mathrm{id}_V\times i_\infty)^!$ on $CH_{d+1}(V\times\mathbb P^1)_\mathbb Q$~\cite[Theorem~6.2(a)]{FultonIT}, we obtain $\operatorname{cl}_{X,Y}^{r,s}\circ\partial_{X,Y}^{r,s}=0$.

For any span $\sigma=[X\xleftarrow{p}Z\xrightarrow{g}Y]$ with $\operatorname{vdim}(p)=r-s+1$ and line bundle $\mathcal L$ on $Z$, the self-intersection formula~\cite[Section~3.2]{MR4092773},~\cite[Corollary~6.3]{FultonIT} gives $\operatorname{cl}_{X,Y}^{r,s}([\widetilde c_1(\mathcal L)\cap\sigma])=(f_{\sigma,\mathrm{cl}})_*(c_1(\mathcal L_{\mathrm{cl}})\cap[Z]^{\mathrm{vir}})$, which vanishes on $\mathcal R_{\mathbb G_a}^0(\sigma)$ and $\mathcal R_{\mathbb G_a}(\sigma;\mathcal L_1,\mathcal L_2)$ by additivity of $c_1$~\cite[Theorem~2.5]{FultonIT}. For $(\sigma,\mathcal E)\in\mathscr B_k^{r,s}(X,Y)$ of rank $n\ge 2$, smoothness of $\pi_{\mathcal E}\colon\mathbb P_Z(\mathcal E)\to Z$ gives $[\mathbb P_Z(\mathcal E)]^{\mathrm{vir}}=\pi_{\mathcal E,\mathrm{cl}}^*[Z]^{\mathrm{vir}}$, whose pushforward along $\pi_{\mathcal E,\mathrm{cl}}$ vanishes as $n-1\ge 1$~\cite[Proposition~3.1(a)]{FultonIT}. Thus $R_{X,Y}^{r,s}\subseteq K_{X,Y}^{r,s}$, and Proposition~\ref{prop:pi0-cycle-surjectivity} gives the exact sequence.

When $\operatorname{char}(k)=0$, combining Proposition~\ref{prop:lazard-ga-compatibility} with $\operatorname{cl}_{X,Y}^{r,s}$ gives the commutative comparison diagram
\[
\xymatrix@C=3.6em@R=2.0em{
\Omega_d(V)_{\mathbb Q} \ar[r]^-{\sim}_-{\text{\cite{LS16}}} \ar@{->>}[d]
& d\Omega_d(V)_{\mathbb Q} \ar@{->>}[d] & \\
CH^c(V)_{\mathbb Q} \ar@{->>}[r] \ar@/_2.2pc/[rr]_-{\mathrm{id}}
& \pi_0\mathcal A_k^{r,s}(X,Y)/R_{X,Y}^{r,s} \ar@{->>}[r]^-{\operatorname{cl}_{X,Y}^{r,s}}
& CH^c(V)_{\mathbb Q}.
}
\]
Here the lower composite sends $f_*[W]$ to $f_*[W]$ for any $W\in\mathbf{SmProj}_k$ of dimension $d$, hence equals $\mathrm{id}_{CH^c(V)_\mathbb Q}$. Both lower maps are isomorphisms, $R_{X,Y}^{r,s}=K_{X,Y}^{r,s}$, and compatibility with composition and external products identifies the additive idempotent completion of $\operatorname{h}(\mathcal A_{k,\mathbb G_a})$ with classical homological Grothendieck Chow motives over $k$ with rational coefficients.
\end{proof}

\begin{remark}
\label{rem:chow-motives-variance}
In Theorem~\ref{thm:p1-chow-comparison}, classical homological Grothendieck Chow motives over $k$ with rational coefficients~\cite{Manin1968,Scholl1994} are covariant in $X\in\mathbf{SmProj}_k$ via graph correspondences.
\end{remark}

\subsection{Universal \texorpdfstring{$\mathbb G_m$}{Gm}-oriented derived span motives}
\label{subsec:derived-gm-oriented-motives}

Let $U = \Spec A$ be an affine derived scheme of a connective animated ring $A$. Let $\mathbf{SmProj}_U$ be the category of smooth projective derived schemes over $U$, and let $\mathbf S_U$ be the class of projective quasi-smooth morphisms graded by locally constant relative virtual dimension $\operatorname{vdim}(p)=\operatorname{rk}(L_{W/X})$. The \emph{integral graded spectral span category} $\mathcal A_U$ has objects $\mathbf M(X)(r)[r]$ ($X\in\mathbf{SmProj}_U$, $r\in\mathbb Z$), symmetric monoidal product $\mathbf M(X)(r)[r]\otimes\mathbf M(Y)(s)[s]=\mathbf M(X\times_U Y)(r+s)[r+s]$, mapping spectra $\mathcal A_U^{r,s}(X,Y) = \Corr_{\mathbf S_U}^{r-s}(X,Y)^{\mathrm{gp}}\in\Sp$, and integral stable envelope $\mathcal M_0(U)=\operatorname{Fun}_{\Sp}(\mathcal A_U^{\mathrm{op}},\Sp)$. Coning off the integral analogues of the $\mathbb G_a$-relation morphisms of Definition~\ref{def:derived-chow-motives}(1)--(3) defines $q_{\mathbb G_a}\colon\mathcal A_U\to\mathcal A_{U,\mathbb G_a}=\mathcal A_U/\mathcal I_{\mathbb G_a}$ and $\mathcal{DM}_{\mathrm{Chow}}(U)=\operatorname{Fun}_{\Sp}(\mathcal A_{U,\mathbb G_a}^{\mathrm{op}},\Sp)$.

Define the \emph{Bott morphism} $\beta\colon\mathbf 1(1)[1]\to\mathbf 1(0)[0]$ to be the fundamental span of $\mathbb P^1$
\[
\beta = [\mathbb P^1_U] = \left[U\xleftarrow{\pi_{\mathbb P^1}}\mathbb P^1_U\xrightarrow{\pi_{\mathbb P^1}}U\right] \in \mathcal A_U^{1,0}(U,U),
\]
and let $\mathcal A_U[\beta^{-1}]$ be its symmetric monoidal Bott localization, with degree-zero mapping spectra
\[
\mathcal A_{U,\beta}(X,Y) = \mathcal A_U[\beta^{-1}]\bigl(\mathbf M(X)(0)[0],\mathbf M(Y)(0)[0]\bigr).
\]
For any $X\in\mathbf{SmProj}_U$ and $d\in\mathbb Z$, write $\mathbf M(X)=\mathbf M(X)(0)[0]$, $\beta_X^d\colon\mathbf M(X)(d)[d]\xrightarrow{\sim}\mathbf M(X)(0)[0]$ in $\mathcal A_U[\beta^{-1}]$ for the $d$-fold power of $\beta_X=\mathbf M(X)\otimes\beta$, and $\beta_X^{-d}\colon\mathbf M(X)(0)[0]\xrightarrow{\sim}\mathbf M(X)(d)[d]$ for its inverse; any span $\sigma\in\mathcal A_U^{d,0}(X,Y)$ of relative virtual dimension $d$ then determines $\sigma\circ\beta_X^{-d}\in\mathcal A_{U,\beta}(X,Y)$.

\begin{definition}
\label{def:derived-gm-span-motives}
Using the integral span operations $\partial_{X,Y}^{d,0}$, $\widetilde c_1(\mathcal L)\cap\sigma$, and $\pi_{\mathcal E}^*\sigma$ on $\mathcal A_U$ defined as in Definition~\ref{def:derived-chow-motives}(1)--(3), we introduce three families of relation morphisms into $\mathcal A_{U,\beta}(X,Y)$ for $X,Y\in\mathbf{SmProj}_U$:
\begin{enumerate}[(1)]
\item the periodized \emph{$\mathbb P^1$-cobordism boundary morphism} $\partial_{X,Y}$ on $\bigoplus_{d\in\mathbb Z}\mathscr F_{\mathbb P^1,U}^{d,0}(X,Y)^{\mathrm{gp}}$, given in degree $d$ by $\partial_{X,Y}^{d,0}\circ\beta_X^{-d}$;
\item for spans $\sigma\in\Corr_{\mathbf S_U}^d(X,Y)$ with middle term $Z$ and invertible $\mathcal O_Z$-modules $\mathcal L_1,\mathcal L_2$, the \emph{multiplicative Chern relation morphisms} with components $\mathcal R^0_{\mathbb G_m}(\sigma) = -(\widetilde c_1(\mathcal O_Z)\cap\sigma)\circ\beta_X^{-(d-1)}$ and
\[
\begin{aligned}
\mathcal R_{\mathbb G_m}(\sigma;\mathcal L_1,\mathcal L_2)
&= -\bigl(\widetilde c_1(\mathcal L_1\otimes\mathcal L_2)\cap\sigma - \widetilde c_1(\mathcal L_1)\cap\sigma - \widetilde c_1(\mathcal L_2)\cap\sigma\bigr)\circ\beta_X^{-(d-1)} \\
&\qquad -\bigl(\widetilde c_1(\mathcal L_1)\cap(\widetilde c_1(\mathcal L_2)\cap\sigma)\bigr)\circ\beta_X^{-(d-2)},
\end{aligned}
\]
which at $\sigma=\mathrm{id}_X$ determine $\mathcal R^0_{\mathbb G_m}(X)=u(\mathcal O_X)-\mathbf 1_X$ and $\mathcal R_{\mathbb G_m}(\mathcal L_1,\mathcal L_2)=u(\mathcal L_1\otimes\mathcal L_2)-u(\mathcal L_1)\circ u(\mathcal L_2)$ for $u(\mathcal L)=\mathbf 1_X-\beta_X\circ\widetilde c_1(\mathcal L)\in\mathcal A_{U,\beta}(X,X)$ (Definition~\ref{def:additive-multiplicative-chern-laws}(2)); and
\item for pairs $(\sigma,\mathcal E)$ consisting of a span $\sigma\in\Corr_{\mathbf S_U}^d(X,Y)$ with middle term $W$ and a finite locally free module $\mathcal E$ on $W$ of rank $n\ge 1$, the \emph{projective bundle relation morphism} $\rho_{X,Y}$ on the group completion of the $\infty$-groupoid of such pairs with components $\pi_{\mathcal E}^*\sigma\circ\beta_X^{-(d+n-1)} - \sigma\circ\beta_X^{-d}$~\cite{AHIConnerFloyd}.
\end{enumerate}
Let $\mathcal I_{\mathbb G_m,\beta}$ be the two-sided tensor ideal in $\mathcal A_U[\beta^{-1}]$ generated by $\partial_{X,Y}$, $\mathcal R^0_{\mathbb G_m}(\sigma)$, $\mathcal R_{\mathbb G_m}(\sigma;\mathcal L_1,\mathcal L_2)$, and $\rho_{X,Y}$, and let $q_{\mathbb G_m}\colon\mathcal A_U\to\mathcal A_U[\beta^{-1}]\to\mathcal A_{U,\mathbb G_m} = \mathcal A_U[\beta^{-1}] / \mathcal I_{\mathbb G_m,\beta}$ be the composite symmetric monoidal spectral quotient functor (Definition~\ref{def:chern-law-span-quotient}). The stable $\infty$-category of \emph{universal $\mathbb G_m$-oriented derived span motives} is $\mathcal{DM}_{\mathbb G_m}(U) = \operatorname{Fun}_{\Sp}(\mathcal A_{U,\mathbb G_m}^{\mathrm{op}},\Sp)$. In $\mathcal A_{U,\mathbb G_m}$, the relations $\mathcal R^0_{\mathbb G_m}(X)$ and $\mathcal R_{\mathbb G_m}(\mathcal L,\mathcal L^\vee)$, together with $\mathcal L\otimes\mathcal L^\vee\simeq\mathcal O_X$, imply that $\mathcal L\mapsto u(\mathcal L)$ defines a group homomorphism $\operatorname{Pic}(X)\to\pi_0\operatorname{End}_{\mathcal A_{U,\mathbb G_m}}(\mathbf M(X))^\times$ with $u(\mathcal L)^{-1}=u(\mathcal L^\vee)$.
\end{definition}

\begin{proposition}
\label{prop:universal-gm-motives}
For any stable presentable symmetric monoidal $\infty$-category $\mathcal C$:
\begin{enumerate}[\rm (1)]
\item Precomposition with $(q_{\mathbb G_m})_!\colon\mathcal M_0(U)\to\mathcal{DM}_{\mathbb G_m}(U)$ identifies $\operatorname{Fun}^{\mathrm{L},\otimes}\bigl(\mathcal{DM}_{\mathbb G_m}(U),\mathcal C\bigr)$ with the $\infty$-category of symmetric monoidal functors $E\colon\mathcal M_0(U)\to\mathcal C$ preserving small colimits that invert the Bott morphism $\beta=[\mathbb P^1_U]$ and are equipped with nullhomotopies of the images of $\partial_{X,Y}$, $\mathcal R^0_{\mathbb G_m}(\sigma)$, $\mathcal R_{\mathbb G_m}(\sigma;\mathcal L_1,\mathcal L_2)$, and $\rho_{X,Y}$ in $\mathcal A_U[\beta^{-1}]$.
\item Precomposition with $(q_{\mathbb G_a})_!\colon\mathcal M_0(U)\to\mathcal{DM}_{\mathrm{Chow}}(U)$ identifies $\operatorname{Fun}^{\mathrm{L},\otimes}\bigl(\mathcal{DM}_{\mathrm{Chow}}(U),\mathcal C\bigr)$ with the $\infty$-category of symmetric monoidal functors $E\colon\mathcal M_0(U)\to\mathcal C$ preserving small colimits and equipped with nullhomotopies of the images of the integral $\mathbb G_a$-relation morphisms $\partial_{X,Y}^{r,s}$, $\mathcal R^0_{\mathbb G_a}(\sigma)$, $\mathcal R_{\mathbb G_a}(\sigma;\mathcal L_1,\mathcal L_2)$, and $\rho_{X,Y}^{r,s}$.
\end{enumerate}
\end{proposition}
\begin{proof}
By~\cite[Proposition~4.8.1.10 and Corollary~4.8.1.12]{HA}, symmetric monoidal functors $E\colon\mathcal M_0(U)\to\mathcal C$ preserving small colimits correspond to symmetric monoidal spectral functors $\mathcal A_U\to\mathcal C$. For~(1), those inverting $\beta$ factor uniquely through $\mathcal A_U[\beta^{-1}]$, and extensions across the pushout $\mathcal A_U[\beta^{-1}]\to\mathcal A_{U,\mathbb G_m}=\mathcal A_U[\beta^{-1}]/\mathcal I_{\mathbb G_m,\beta}$ correspond to nullhomotopies on the generating morphisms of $\mathcal I_{\mathbb G_m,\beta}$. Part~(2) follows identically from the pushout $\mathcal A_U\to\mathcal A_{U,\mathbb G_a}=\mathcal A_U/\mathcal I_{\mathbb G_a}$.
\end{proof}

\begin{remark}
\label{rem:formal-almost-oriented-specialization}
The construction of Definition~\ref{def:derived-gm-span-motives} and the universal properties of Proposition~\ref{prop:universal-gm-motives} extend without change to any affine prismatic formal base $U=\operatorname{Spf}(B)\in\mathcal T_R$ using smooth projective formal schemes in $\mathrm{dForm}(U)$ and $\operatorname{vdim}(p)=\operatorname{rk}(\widehat L_{W/X})$, and to any boundary almost formal base $U=\operatorname{Spf}^a(B_\nu)\in\mathcal C_{R_\nu}^a$ ($\nu\in\{\theta,\flat\}$) using smooth projective formal schemes in $\mathcal C_{R_\nu}^a$ and $\operatorname{vdim}(p)=\operatorname{rk}(\widehat L_{W/X}^a)$. By Proposition~\ref{prop:prismatic-axioms-verification} and Definition~\ref{def:almost-formal-base-categories}, the formal spaces $\mathbb P^1_U$, $\mathbb V_Z(\mathcal L)$, and $\mathbb P_W(\mathcal E)$ are represented in $\mathrm{dForm}(U)$ and $\mathcal C_{R_\nu}^a$, and zero sections lie in $\mathbf S_U$, so all relation morphisms of Definition~\ref{def:derived-gm-span-motives} are well-defined over these bases. In Subsection~\ref{subsec:tilting-oriented-motives}, we pass from a fixed base $U$ to the $\mathcal E_{\mathcal D,R}(X_{0,\mathcal D})$-enriched setting over the slice site of $X_0\in\mathcal T_R^{\mathrm{perf}}$ and combine these oriented quotients with $P^1_{\mathrm{cot}}$ and $P^1_{\mathrm{rel}}$.
\end{remark}

\subsection{Tilting equivalence for oriented derived span motives}
\label{subsec:tilting-oriented-motives}

Throughout this subsection, fix an integral perfectoid $R^\circ$-algebra $A^\circ$ with tilt $A^{\flat\circ}$ over $R^{\flat\circ}$ and perfect prism $\widetilde A=\Ainf(A^\circ)=W(A^{\flat\circ})$. Set $X_{0,\mathrm{Prism}}=X_0=\operatorname{Spf}(\widetilde A)\in\mathcal T_R^{\mathrm{perf}}$ with boundary fibers $X_{0,\theta}=\beta_\theta(X_0)=\operatorname{Spf}(A^\circ)$ and $X_{0,\flat}=\beta_\flat(X_0)=\operatorname{Spf}(A^{\flat\circ})$ ($A_\theta=A^\circ$, $A_\flat=A^{\flat\circ}$). Applying $Q_{X_0}$ and $a_{\mathcal P_{\nu,R}/X_{0,\nu}}^{\mathrm{hyp}}$ (Proposition~\ref{prop:perfect-boundary-coefficient-sites}) over $\mathscr D_{\mathrm{cmp}}(X_0)$ defines three graded presentation span categories on smooth projective objects, graded by relative virtual dimension $r-s$:
\begin{itemize}
\item the $\mathcal E_R(X_0)$-enriched prismatic category $\mathcal A_{\mathrm{Prism},\mathrm{cmp}}^{\mathrm{gr}}(X_0)$ of smooth projective formal schemes over $X_0$ in $\mathscr D_{\mathrm{cmp}}(X_0)$;
\item the $\mathcal E_{\nu,R}(X_{0,\nu})$-enriched boundary category $\mathcal A_{\nu,\mathrm{cmp}}^{\mathrm{gr}}(X_{0,\nu})=(\tau_\nu)_*\mathcal A_{\mathrm{Prism},\mathrm{cmp}}^{\mathrm{gr}}(X_0)$, obtained by transport along $\tau_\nu=(\beta_\nu^*)^{-1}r^*\colon\mathcal E_R(X_0)\xrightarrow{\sim}\mathcal E_{\nu,R}(X_{0,\nu})$; and
\item the intrinsic graded category $\mathcal A_{\mathcal D,1,\mathrm{cmp}}^{\mathrm{gr}}(X_{0,\mathcal D})$ ($\mathcal D\in\{\mathrm{Prism},\theta,\flat\}$) of smooth projective derived almost schemes over $X_{0,\mathcal D,1}$ (Subsection~\ref{subsec:reduction-and-p1}).
\end{itemize}

In the symmetric monoidal $\mathcal E_{\mathcal D,R}(X_{0,\mathcal D})$-enriched setting~\cite{MR3345192} (Definitions~\ref{def:chern-law-span-quotient},~\ref{def:derived-chow-motives}, and~\ref{def:derived-gm-span-motives}), coning off the integral $\mathbb G_a$-relations gives the quotient functor $q_{\mathcal D,\mathbb G_a}\colon\mathcal A_{\mathcal D,\mathrm{cmp}}^{\mathrm{gr}}(X_{0,\mathcal D})\to\allowbreak\mathcal A_{\mathcal D,\mathrm{cmp},\mathbb G_a}^{\mathrm{gr}}(X_{0,\mathcal D})$, while inverting $\beta=[\mathbb P^1]$ and coning off the integral $\mathbb G_m$-relations gives $q_{\mathcal D,\mathbb G_m}\colon\mathcal A_{\mathcal D,\mathrm{cmp}}^{\mathrm{gr}}(X_{0,\mathcal D})\to\allowbreak\mathcal A_{\mathcal D,\mathrm{cmp},\mathbb G_m}^{\mathrm{gr}}(X_{0,\mathcal D})$. Write $\mathcal A_{\mathcal D,1,\mathrm{cmp},\mathbb G_a}^{\mathrm{gr}}(X_{0,\mathcal D})$ (resp.\ $\mathcal A_{\mathcal D,1,\mathrm{cmp},\mathbb G_m}^{\mathrm{gr}}(X_{0,\mathcal D})$) for the corresponding mod-$\omega_{\mathcal D}$ quotients and $\mathcal A_{\nu,\mathrm{int}}^{\mathrm{gr}}(X_{0,\nu})\to\allowbreak\mathcal A_{\nu,\mathrm{int},\mathbb G_a}^{\mathrm{gr}}(X_{0,\nu})$ (resp.\ $\mathcal A_{\nu,\mathrm{int}}^{\mathrm{gr}}(X_{0,\nu})\to\allowbreak\mathcal A_{\nu,\mathrm{int},\mathbb G_m}^{\mathrm{gr}}(X_{0,\nu})$) ($\nu\in\{\theta,\flat\}$) for the intrinsic boundary graded span categories and their $\mathbb G_a$-quotients (resp.\ $\mathbb G_m$-quotients) on $\mathcal C_{\nu,\mathrm{cmp}}(X_{0,\nu})$ (Definition~\ref{def:presentation-realizations}).

With coefficient categories $\Mod_{\widehat G_{\mathcal D}}^\wedge(\mathcal E_{\mathcal D,R}(X_{0,\mathcal D}))$, $\Mod_{\mathbb A_{\inf R}}(\mathcal E_R(X_0))$, and $\Mod_{R_\nu^a}(\mathcal E_{\nu,R}(X_{0,\nu}))$ from Subsection~\ref{subsec:reduction-and-p1} and Definition~\ref{def:integral-linear-comparison-sources}(1) ($R_\theta^a=(R^\circ)^a$, $R_\flat^a=(R^{\flat\circ})^a$), we form the enriched presheaf envelopes~\cite{MR3345192} along the adic tower and over the integral $\mathcal O$-linear bases on the $\mathbb G_a$-oriented enriched quotients:
\[
\begin{aligned}
\mathcal M_{0,\mathrm{Chow},\Prism}^{a,\mathrm{nu}}(X_0)
&= \operatorname{Fun}_{\mathcal E_R(X_0)}\bigl(\mathcal A_{\mathrm{Prism},\mathrm{cmp},\mathbb G_a}^{\mathrm{gr}}(X_0)^{\mathrm{op}},\Mod_{\widehat G_{\mathrm{Prism}}}^\wedge(\mathcal E_R(X_0))\bigr), \\
\mathcal M_{0,\mathrm{Chow},\nu}^{a,\mathrm{nu}}(X_{0,\nu})
&= \operatorname{Fun}_{\mathcal E_{\nu,R}(X_{0,\nu})}\bigl(\mathcal A_{\nu,\mathrm{cmp},\mathbb G_a}^{\mathrm{gr}}(X_{0,\nu})^{\mathrm{op}},\Mod_{\widehat G_\nu}^\wedge(\mathcal E_{\nu,R}(X_{0,\nu}))\bigr), \\
\mathcal M_{0,\mathrm{Chow},\mathcal D,\mathrm{gr}}^{a,\mathrm{nu}}(X_{0,\mathcal D})
&= \operatorname{Fun}_{\mathcal E_{\mathcal D,R}(X_{0,\mathcal D})}\bigl(\mathcal A_{\mathcal D,1,\mathrm{cmp},\mathbb G_a}^{\mathrm{gr}}(X_{0,\mathcal D})^{\mathrm{op}},\Mod_{\widehat G_{\mathcal D}}^\wedge(\mathcal E_{\mathcal D,R}(X_{0,\mathcal D}))\bigr), \\
\mathcal M_{0,\mathrm{Chow},\Prism}^{\mathrm{lin}}(X_0)
&= \operatorname{Fun}_{\mathcal E_R(X_0)}\bigl(\mathcal A_{\mathrm{Prism},\mathrm{cmp},\mathbb G_a}^{\mathrm{gr}}(X_0)^{\mathrm{op}},\Mod_{\mathbb A_{\inf R}}(\mathcal E_R(X_0))\bigr), \\
\mathcal M_{0,\mathrm{Chow},\nu}^{\mathrm{lin}}(X_{0,\nu})
&= \operatorname{Fun}_{\mathcal E_{\nu,R}(X_{0,\nu})}\bigl(\mathcal A_{\nu,\mathrm{cmp},\mathbb G_a}^{\mathrm{gr}}(X_{0,\nu})^{\mathrm{op}},\Mod_{R_\nu^a}(\mathcal E_{\nu,R}(X_{0,\nu}))\bigr),
\end{aligned}
\]
and identically define $\mathcal M_{0,\mathbb G_m,\Prism}^{a,\mathrm{nu}}(X_0)$, $\mathcal M_{0,\mathbb G_m,\nu}^{a,\mathrm{nu}}(X_{0,\nu})$, $\mathcal M_{0,\mathbb G_m,\mathcal D,\mathrm{gr}}^{a,\mathrm{nu}}(X_{0,\mathcal D})$, $\mathcal M_{0,\mathbb G_m,\Prism}^{\mathrm{lin}}(X_0)$, and $\mathcal M_{0,\mathbb G_m,\nu}^{\mathrm{lin}}(X_{0,\nu})$ by replacing $\mathbb G_a$ with $\mathbb G_m$. Symmetric monoidal localization $L_{\mathrm{hyp}}$ at the
tensor ideal $\mathcal J_{\mathrm{hyp}}$ generated by presentation
hypercover assembly morphisms
(Proposition~\ref{thm:localized-span-universal-property}) defines the
presentably symmetric monoidal stable $\infty$-categories
$\mathcal{DM}_{\mathrm{Chow},\Prism}^{a,\mathrm{nu}}(\widetilde A)$,
$\mathcal{DM}_{\mathrm{Chow},\theta}^{a,\mathrm{nu}}(A^\circ)$,
$\mathcal{DM}_{\mathrm{Chow},\flat}^{a,\mathrm{nu}}(A^{\flat\circ})$ and
$\mathcal{DM}_{\mathbb G_m,\Prism}^{a,\mathrm{nu}}(\widetilde A)$,
$\mathcal{DM}_{\mathbb G_m,\theta}^{a,\mathrm{nu}}(A^\circ)$,
$\mathcal{DM}_{\mathbb G_m,\flat}^{a,\mathrm{nu}}(A^{\flat\circ})$. It
likewise produces their intrinsic mod-$\omega_{\mathcal D}$
associated graded versions $\mathcal{DM}_{\mathrm{Chow},\mathcal
D,\mathrm{gr}}^{a,\mathrm{nu}}$ and $\mathcal{DM}_{\mathbb G_m,\mathcal
D,\mathrm{gr}}^{a,\mathrm{nu}}$, linked by symmetric monoidal
associated graded realization functors $\mathcal{DM}_{\mathrm{Chow},\mathcal
D}^{a,\mathrm{nu}}\to\mathcal{DM}_{\mathrm{Chow},\mathcal
D,\mathrm{gr}}^{a,\mathrm{nu}}$ (resp.\ $\mathcal{DM}_{\mathbb G_m,\mathcal
D}^{a,\mathrm{nu}}\to\mathcal{DM}_{\mathbb G_m,\mathcal
D,\mathrm{gr}}^{a,\mathrm{nu}}$), together with the integral $\mathcal
O$-linear categories
$\mathcal{DM}_{\mathrm{Chow},\mathcal D}^{\mathrm{lin}}$ and
$\mathcal{DM}_{\mathbb G_m,\mathcal D}^{\mathrm{lin}}$.

Now fix a split augmented formal scheme $X_I=\operatorname{Spf}(\widetilde A\oplus I)\to X_0$ ($I\in\NAlg_{\mathrm{Prism}}(\widetilde A)$) with boundary specializations $I_\theta=i_\theta^*I$ and $I_\flat=i_\flat^*I$, and a target object $\mathbf M(Y)(s)[s]$. Evaluating the enriched lift presheaves of Definitions~\ref{def:nonunital-comparison-sources}(2) and~\ref{def:integral-linear-comparison-sources}(2) on $\mathbf M(X)(r)[r]$ via spans $\sigma\in\mathscr D_{\mathrm{cmp}}^{r-s}(X_V,Y_V)$ defines graded lift presheaves on $\mathcal A_{\mathcal D,\mathrm{cmp}}^{\mathrm{gr}}(X_{0,\mathcal D})^{\mathrm{op}}$. Enriched left Kan extension $(q_{\mathcal D,\mathbb G_a})_!$ (resp.\ $(q_{\mathcal D,\mathbb G_m})_!$) followed by $L_{\mathrm{hyp}}$ then defines the excisive approximations along the adic tower and over the integral bases for the oriented tangent motives of $(I,I_\theta,I_\flat)$.

\begin{proposition}
\label{prop:oriented-boundary-specialization-functoriality}
The boundary specialization functors $i_\theta^*$ and $i_\flat^*$ preserve smooth projective formal schemes, relative virtual dimensions, the Bott morphism $\beta=[\mathbb P^1]$, the $\mathbb P^1$-cobordism boundary morphisms $\partial_{X,Y}$ and $\partial_{X,Y}^{r,s}$, derived first Chern spans $\widetilde c_1(\mathcal L)\cap\sigma$, additive and multiplicative Chern relations, and the projective bundle relation morphisms $\rho_{X,Y}$ and $\rho_{X,Y}^{r,s}$, inducing symmetric monoidal specialization functors
\[
\begin{aligned}
\mathcal A_{\theta,\mathrm{int},\mathbb G_m}^{\mathrm{gr}}(X_{0,\theta})
&\;\xleftarrow{\ \Phi_{\theta,\mathbb G_m}\ }\;
\mathcal A_{\mathrm{Prism},\mathrm{cmp},\mathbb G_m}^{\mathrm{gr}}(X_0)
\;\xrightarrow{\ \Phi_{\flat,\mathbb G_m}\ }\;
\mathcal A_{\flat,\mathrm{int},\mathbb G_m}^{\mathrm{gr}}(X_{0,\flat}), \\
\mathcal A_{\theta,\mathrm{int},\mathbb G_a}^{\mathrm{gr}}(X_{0,\theta})
&\;\xleftarrow{\ \Phi_{\theta,\mathbb G_a}\ }\;
\mathcal A_{\mathrm{Prism},\mathrm{cmp},\mathbb G_a}^{\mathrm{gr}}(X_0)
\;\xrightarrow{\ \Phi_{\flat,\mathbb G_a}\ }\;
\mathcal A_{\flat,\mathrm{int},\mathbb G_a}^{\mathrm{gr}}(X_{0,\flat}),
\end{aligned}
\]
whose enriched left Kan extensions descend across $L_{\mathrm{hyp}}$ to symmetric monoidal realization functors preserving small colimits from $\mathcal{DM}_{\mathrm{Chow},\nu}^{a,\mathrm{nu}}(A_\nu)$ and $\mathcal{DM}_{\mathrm{Chow},\nu}^{\mathrm{lin}}(A_\nu)$ (resp.\ $\mathcal{DM}_{\mathbb G_m,\nu}^{a,\mathrm{nu}}(A_\nu)$ and $\mathcal{DM}_{\mathbb G_m,\nu}^{\mathrm{lin}}(A_\nu)$) into the intrinsic boundary motivic categories obtained by localizing presheaves on $\mathcal A_{\nu,\mathrm{int},\mathbb G_a}^{\mathrm{gr}}(X_{0,\nu})$ (resp.\ $\mathcal A_{\nu,\mathrm{int},\mathbb G_m}^{\mathrm{gr}}(X_{0,\nu})$) at presentation hyperdescent ($\nu\in\{\theta,\flat\}$).
\end{proposition}
\begin{proof}
Fix $\nu \in \{\theta, \flat\}$. By
Proposition~\ref{prop:prismatic-fiber-geometric-compatibility}, the
specialization functor $i_\nu^*$ preserves smooth projective formal
schemes, projective quasi-smooth morphisms,
$\operatorname{vdim}(p)=\operatorname{rk}(\widehat L_{W/X})$, derived
fiber products, finite disjoint unions, and $\beta=[\mathbb P^1]$,
inducing symmetric monoidal functors $\mathcal
A_{\mathrm{Prism},\mathrm{cmp}}^{\mathrm{gr}}(X_0)\to\mathcal
A_{\nu,\mathrm{int}}^{\mathrm{gr}}(X_{0,\nu})$ and $\mathcal
A_{\mathrm{Prism},\mathrm{cmp}}^{\mathrm{gr}}(X_0)[\beta^{-1}] \to
\mathcal A_{\nu,\mathrm{int}}^{\mathrm{gr}}(X_{0,\nu})[\beta^{-1}]$. By
Proposition~\ref{prop:prismatic-fiber-geometric-compatibility} and
Proposition~\ref{prop:lift-section-sieve-compatibility}(2), $i_\nu^*$
also preserves presentation hypercovers, $\mathbb P^1$-cobordism spans,
derived endpoint pullbacks $e_0,e_\infty$, zero sections of line bundles
and their derived self-intersections, and projective bundles $\mathbb
P_W(\mathcal E)$; it therefore preserves $\partial_{X,Y}$,
$\partial_{X,Y}^{r,s}$, $\widetilde c_1(\mathcal L)\cap\sigma$, the
additive and multiplicative Chern relations, $\rho_{X,Y}$, and
$\rho_{X,Y}^{r,s}$.

Evaluating over $X_0=\operatorname{Spf}(\widetilde A)$ via
Proposition~\ref{prop:universal-gm-motives},
Remark~\ref{rem:formal-almost-oriented-specialization}, and
Proposition~\ref{prop:perfect-boundary-coefficient-sites} therefore
induces $\Phi_{\nu,\mathbb G_m}$ and $\Phi_{\nu,\mathbb G_a}$, which
factor through $(\tau_\nu)_*\mathcal
A_{\mathrm{Prism},\mathrm{cmp},\mathbb G_m}^{\mathrm{gr}}(X_0)=\mathcal
A_{\nu,\mathrm{cmp},\mathbb G_m}^{\mathrm{gr}}(X_{0,\nu})$ and
$(\tau_\nu)_*\mathcal A_{\mathrm{Prism},\mathrm{cmp},\mathbb
G_a}^{\mathrm{gr}}(X_0)=\mathcal A_{\nu,\mathrm{cmp},\mathbb
G_a}^{\mathrm{gr}}(X_{0,\nu})$. Their enriched left Kan
extensions~\cite{MR3345192} preserve presentation hyperdescent and
descend by Proposition~\ref{thm:localized-span-universal-property} to
the boundary realization functors.
\end{proof}

Combining Proposition~\ref{prop:oriented-boundary-specialization-functoriality} with $P^1_{\mathrm{cot}}$, $P^1_{\mathrm{rel}}$ (Theorem~\ref{prop:completed-almost-excisive-reflection}(3)), and Lemma~\ref{lem:p1-fiber-comparison} gives the tilting and base-change equivalences for oriented derived span motives below.

\begin{theorem}
\label{thm:tilting-oriented-motives}
Let $I\in\NAlg_{\mathrm{Prism}}(\widetilde A)$ have boundary specializations $I_\theta=i_\theta^*I$ and $I_\flat=i_\flat^*I$.
\begin{enumerate}[\rm (1)]
\item For the transported categories of universal $\mathbb G_m$-oriented derived span motives and $\mathbb G_a$-oriented derived Chow motives along the adic tower, together with their intrinsic associated graded versions, there are symmetric monoidal almost tilting equivalences:
\[
\begin{aligned}
\mathcal{DM}_{\mathbb G_m,\theta}^{a,\mathrm{nu}}(A^\circ)
&\;\longleftarrow\;
\mathcal{DM}_{\mathbb G_m,\Prism}^{a,\mathrm{nu}}(\widetilde A)
\;\longrightarrow\;
\mathcal{DM}_{\mathbb G_m,\flat}^{a,\mathrm{nu}}(A^{\flat\circ}), \\
\mathcal{DM}_{\mathrm{Chow},\theta}^{a,\mathrm{nu}}(A^\circ)
&\;\longleftarrow\;
\mathcal{DM}_{\mathrm{Chow},\Prism}^{a,\mathrm{nu}}(\widetilde A)
\;\longrightarrow\;
\mathcal{DM}_{\mathrm{Chow},\flat}^{a,\mathrm{nu}}(A^{\flat\circ}),
\end{aligned}
\]
compatible with the associated graded equivalences $\mathcal{DM}_{\mathbb G_m,\theta,\mathrm{gr}}^{a,\mathrm{nu}}(A^\circ)\leftarrow\mathcal{DM}_{\mathbb G_m,\Prism,\mathrm{gr}}^{a,\mathrm{nu}}(\widetilde A)\rightarrow\mathcal{DM}_{\mathbb G_m,\flat,\mathrm{gr}}^{a,\mathrm{nu}}(A^{\flat\circ})$ (resp.\ $\mathcal{DM}_{\mathrm{Chow},\theta,\mathrm{gr}}^{a,\mathrm{nu}}(A^\circ)\leftarrow\mathcal{DM}_{\mathrm{Chow},\Prism,\mathrm{gr}}^{a,\mathrm{nu}}(\widetilde A)\rightarrow\mathcal{DM}_{\mathrm{Chow},\flat,\mathrm{gr}}^{a,\mathrm{nu}}(A^{\flat\circ})$), identifying the excisive approximations of the oriented tangent motives of $(I,I_\theta,I_\flat)$ along the adic tower across all $n\ge 1$.
\item For the transported integral $\mathcal O$-linear categories of oriented motives over $\mathbb A_{\inf R}$, $R^\circ$, and $R^{\flat\circ}$, derived base change along $\theta\colon\mathbb A_{\inf R}\to R^\circ$ and $\rho_\flat\colon\mathbb A_{\inf R}\to R^{\flat\circ}$ induces symmetric monoidal equivalences:
\[
\begin{aligned}
\mathcal{DM}_{\mathbb G_m,\Prism}^{\mathrm{lin}}(\widetilde A)\otimes_{\Mod_{\mathbb A_{\inf R}}}\Mod_{(R^\circ)^a}
&\;\longrightarrow\;
\mathcal{DM}_{\mathbb G_m,\theta}^{\mathrm{lin}}(A^\circ), \\
\mathcal{DM}_{\mathbb G_m,\Prism}^{\mathrm{lin}}(\widetilde A)\otimes_{\Mod_{\mathbb A_{\inf R}}}\Mod_{(R^{\flat\circ})^a}
&\;\longrightarrow\;
\mathcal{DM}_{\mathbb G_m,\flat}^{\mathrm{lin}}(A^{\flat\circ}), \\
\mathcal{DM}_{\mathrm{Chow},\Prism}^{\mathrm{lin}}(\widetilde A)\otimes_{\Mod_{\mathbb A_{\inf R}}}\Mod_{(R^\circ)^a}
&\;\longrightarrow\;
\mathcal{DM}_{\mathrm{Chow},\theta}^{\mathrm{lin}}(A^\circ), \\
\mathcal{DM}_{\mathrm{Chow},\Prism}^{\mathrm{lin}}(\widetilde A)\otimes_{\Mod_{\mathbb A_{\inf R}}}\Mod_{(R^{\flat\circ})^a}
&\;\longrightarrow\;
\mathcal{DM}_{\mathrm{Chow},\flat}^{\mathrm{lin}}(A^{\flat\circ}),
\end{aligned}
\]
identifying the integral cotangent excisive approximations of the oriented tangent motives of $(I,I_\theta,I_\flat)$ under derived base change along $\theta$ and $\rho_\flat$.
\end{enumerate}
\end{theorem}
\begin{proof}
Set $X_0=\operatorname{Spf}(\widetilde A)\in\mathcal T_R^{\mathrm{perf}}$ and $X_{0,\nu}=\beta_\nu(X_0)$ ($\nu\in\{\theta,\flat\}$). For~(1), as in the proof of Theorem~\ref{thm:tilting-main}(1), Theorem~\ref{prop:completed-almost-excisive-reflection}(3)(a), Lemma~\ref{lem:tilting-frame-koszul}(2)--(3), and site transport along $(\tau_\nu,\Xi_\nu)$ (Lemma~\ref{lem:p1-fiber-comparison}(1)) give equivalences $(\tau_\nu)_*\mathcal A_{\mathrm{Prism},\mathrm{cmp},\mathbb G_a}^{\mathrm{gr}}(X_0)\to\allowbreak\mathcal A_{\nu,\mathrm{cmp},\mathbb G_a}^{\mathrm{gr}}(X_{0,\nu})$ and $(\tau_\nu)_*\mathcal A_{\mathrm{Prism},1,\mathrm{cmp},\mathbb G_a}^{\mathrm{gr}}(X_0)\to\allowbreak\mathcal A_{\nu,1,\mathrm{cmp},\mathbb G_a}^{\mathrm{gr}}(X_{0,\nu})$, and identically for $\mathbb G_m$. Passing to enriched presheaf categories gives exact symmetric monoidal equivalences on $\mathcal M_{0,\mathrm{Chow},-}^{a,\mathrm{nu}}$ (resp.\ $\mathcal M_{0,\mathbb G_m,-}^{a,\mathrm{nu}}$) that commute with $(q_{\mathcal D,\mathbb G_a})_!$ (resp.\ $(q_{\mathcal D,\mathbb G_m})_!$) and preserve $\mathcal J_{\mathrm{hyp}}$ by Proposition~\ref{prop:oriented-boundary-specialization-functoriality}. Proposition~\ref{thm:localized-span-universal-property} induces the compatible tilting equivalences on $\mathcal{DM}_{\mathrm{Chow},\mathcal D}^{a,\mathrm{nu}}$ and $\mathcal{DM}_{\mathrm{Chow},\mathcal D,\mathrm{gr}}^{a,\mathrm{nu}}$ (resp.\ $\mathcal{DM}_{\mathbb G_m,\mathcal D}^{a,\mathrm{nu}}$ and $\mathcal{DM}_{\mathbb G_m,\mathcal D,\mathrm{gr}}^{a,\mathrm{nu}}$), and Theorem~\ref{prop:completed-almost-excisive-reflection}(3)(a) identifies the oriented tangent motives along the adic tower.

For~(2), as in the proof of Theorem~\ref{thm:tilting-main}(2), dualizability of $\Mod_{R_\nu^a}$ over $\Sp$, Proposition~\ref{prop:perfect-boundary-coefficient-sites},~\cite{MR3345192}, and~\cite[Proposition~4.8.1.17 and Theorem~4.8.4.6]{HA} give an equivalence $\mathcal M_{0,\mathrm{Chow},\Prism}^{\mathrm{lin}}(X_0)\allowbreak\otimes_{\Mod_{\mathbb A_{\inf R}}}\allowbreak\Mod_{R_\nu^a}\to\allowbreak\mathcal M_{0,\mathrm{Chow},\nu}^{\mathrm{lin}}(X_{0,\nu})$, and similarly for $\mathcal M_{0,\mathbb G_m,-}^{\mathrm{lin}}$, compatibly with $(q_{\mathcal D,\mathbb G_a})_!$ (resp.\ $(q_{\mathcal D,\mathbb G_m})_!$). Propositions~\ref{prop:oriented-boundary-specialization-functoriality} and~\ref{thm:localized-span-universal-property} induce the base-change equivalences along $\theta$ and $\rho_\flat$, while Lemma~\ref{lem:p1-fiber-comparison}(2) and Theorem~\ref{prop:completed-almost-excisive-reflection}(3)(b) identify the integral cotangent excisive approximations of the oriented tangent motives.
\end{proof}

\begin{proposition}
\label{prop:oriented-koszul-tilting}
Let $I\in\NAlg_{\mathrm{Prism}}(\widetilde A)$ have boundary specializations $I_\theta=i_\theta^*I$ and $I_\flat=i_\flat^*I$. For each integer $1\le k\le p$, and in particular at $k=p$ over $R_p^a=(R^\circ/(p))^a\simeq (R^{\flat\circ}/(\bar\xi))^a=(R_p^\flat)^a$, derived reduction along $R^\circ\to R_k$ and $R^{\flat\circ}\to R_k^\flat$ induces symmetric monoidal equivalences
\[
\begin{aligned}
\mathcal{DM}_{\mathbb G_m,\theta}^{\mathrm{lin}}(A^\circ)\otimes_{\Mod_{(R^\circ)^a}}\Mod_{R_k^a}
&\;\longleftarrow\;
\mathcal{DM}_{\mathbb G_m,\Prism}^{\mathrm{lin}}(\widetilde A)\otimes_{\Mod_{\mathbb A_{\inf R}}}\Mod_{R_k^a} \\
&\;\longrightarrow\;
\mathcal{DM}_{\mathbb G_m,\flat}^{\mathrm{lin}}(A^{\flat\circ})\otimes_{\Mod_{(R^{\flat\circ})^a}}\Mod_{(R_k^\flat)^a}, \\
\mathcal{DM}_{\mathrm{Chow},\theta}^{\mathrm{lin}}(A^\circ)\otimes_{\Mod_{(R^\circ)^a}}\Mod_{R_k^a}
&\;\longleftarrow\;
\mathcal{DM}_{\mathrm{Chow},\Prism}^{\mathrm{lin}}(\widetilde A)\otimes_{\Mod_{\mathbb A_{\inf R}}}\Mod_{R_k^a} \\
&\;\longrightarrow\;
\mathcal{DM}_{\mathrm{Chow},\flat}^{\mathrm{lin}}(A^{\flat\circ})\otimes_{\Mod_{(R^{\flat\circ})^a}}\Mod_{(R_k^\flat)^a},
\end{aligned}
\]
identifying the integral cotangent excisive approximations of the oriented tangent motives of $(I,I_\theta,I_\flat)$ after derived tensoring with $R_k^a\simeq (R_k^\flat)^a$.
\end{proposition}
\begin{proof}
As in the proof of Proposition~\ref{prop:koszul-tilting-and-tower-recovery}(2), Lemma~\ref{lem:tilting-frame-koszul}(2)--(3) and Subsection~\ref{subsec:almost-formal-base} identify $R_k^a\simeq (R_k^\flat)^a$ as dualizable algebras over $\mathbb A_{\inf R}$ and $R_\nu^a$ ($\nu\in\{\theta,\flat\}$). Tensoring Theorem~\ref{thm:tilting-oriented-motives}(2) with $\Mod_{R_k^a}$ over $\Mod_{R_\nu^a}$ gives the category equivalences, and Lemma~\ref{lem:p1-fiber-comparison}(2), Theorem~\ref{prop:completed-almost-excisive-reflection}(3)(b), and Proposition~\ref{prop:oriented-boundary-specialization-functoriality} identify the oriented tangent motives over $R_k^a\simeq (R_k^\flat)^a$.
\end{proof}

\begin{remark}
\label{rem:oriented-conservativity-and-generic-fiber}
As in Corollary~\ref{cor:mixed-to-equal-conservativity} and Remark~\ref{rem:generic-fiber-localization}, any morphism in $\mathcal{DM}_{\mathbb G_m,\Prism}^{\mathrm{lin}}(\widetilde A)$ or $\mathcal{DM}_{\mathrm{Chow},\Prism}^{\mathrm{lin}}(\widetilde A)$ induces objectwise almost equivalences after derived $\omega$-completion over $(R^\circ)^a$ if and only if its tilted specialization induces objectwise almost equivalences after derived $\omega^\flat$-completion over $(R^{\flat\circ})^a$. This condition is detected modulo $(\xi,\widetilde\omega)$ over $R_1^a$, and equivalently in $\mathcal{DM}_{\mathbb G_m,\nu}^{a,\mathrm{nu}}(A_\nu)$ (resp.\ $\mathcal{DM}_{\mathrm{Chow},\nu}^{a,\mathrm{nu}}(A_\nu)$) for morphisms of oriented tangent motives induced by prismatic augmentation data. Inverting $\widetilde\omega$ specializes both prismatic oriented categories along $\theta[\widetilde\omega^{-1}]$ and $\rho_\flat[\widetilde\omega^{-1}]$ to ordinary symmetric monoidal oriented motivic $\infty$-categories over $R=R^\circ[1/p]$ and $R^\flat=R^{\flat\circ}[(\omega^\flat)^{-1}]$.
\end{remark}
\subsection*{Acknowledgements and declaration of Generative AI usage}
The author utilized Google Gemini 4 Carbon and OpenAI Prism (GPT-5.2) to
assist in exploring proof strategies, checking technical details of
intermediate lemmas, and refining the clarity and logical flow of the
English exposition. All mathematical ideas, definitions, main theorems,
and proof architectures were conceived entirely by the author.  Each
AI-assisted argument and calculation was independently verified line by
line by the author, who assumes sole responsibility for the mathematical
accuracy and content of this publication.

\bibliographystyle{amsalpha}
\bibliography{bibkato}
\end{document}